\documentclass[12pt]{book}
\usepackage{amsmath}
\usepackage{amsthm}
\usepackage[all]{xy}
\usepackage[margin=1.0in]{geometry}
\usepackage{times}
\usepackage{newtxmath}

\newtheorem{proposition}{Proposition}

\theoremstyle{definition}
\newtheorem{definition}{Definition}
\newtheorem{remark}{Remark}
\newtheorem{reference}{Reference}
\newtheorem{assumption}{Assumption}
\newtheorem{conjecture}{Conjecture}

\title{$\infty$-Categorical Generalized Langlands Program I: Mixed-Parity Modules and Sheaves}
\author{Xin Tong}
\date{}

\begin{document}

\maketitle

\newpage
\subsection*{Abstract}
Mixed-parity module emerges for instance when a de Rham Galois representation is being tensored with a square root of cyclotomic character, which produces half odd integers as the corresponding Hodge-Tate weights. We build the whole foundation on the $p$-adic Hodge theory in this setting over small $v$-stacks after Scholze and we also consider certain moduli $v$-stack which parametrizes families of mixed-parity Hodge modules. Examples of the small $v$-stacks in our mind are rigid analytic spaces over $p$-adic fields and moduli $v$-stack of vector bundles over Fargues-Fontaine curves. The preparation implemented at this level will be expected to provide further essential foundationalization for generalized Langlands program after Langlands, Drinfeld, Fargues-Scholze. One side of the generalized Langlands correspondence in the geometric setting is the perverse motivic derived $\infty$-category over $\mathrm{Moduli}_G$ related to smooth representations of reductive groups, while the other side of the generalized Langlands correspondence in the geometric setting is the corresponding derived $\infty$-category over the stack of mixed-parity $L$-parametrizations (i.e. from two-fold covering of the Weil group into $\ell$-adic groups) related to the representations of Weil group in our setting into Langlands dual groups. Although after Scholze and Fargues-Scholze our generalized Langlands program can go along $\ell$-adic cohomologicalization to immediately achieve various solid derived $\infty$-categories $\mathrm{DerCat}_\text{\'et}(\mathrm{Moduli}_G,\square)$, $\mathrm{DerCat}_\mathrm{lisse, \blacksquare}(\mathrm{Moduli}_G,\square)$, $\mathrm{DerCat}_{\blacksquare}(\mathrm{Moduli}_G,\square)$ and so on with well-established formalism regarding 6-functors, we already provide certain $p$-adic cohomologicalization of the story over $\mathrm{Moduli}_G$.

\newpage
\tableofcontents

\newpage

\begin{reference}\mbox{}
\begin{itemize}
\item[$\square$] Chapter 1 Main References: \cite{Sch1}, \cite{KL1}, \cite{KL2}, \cite{BL1}, \cite{BL2}, \cite{BS}, \cite{BHS}, \cite{Fon1}, \cite{CS1}, \cite{CS2}, \cite{BK}, \cite{BBK}, \cite{BBBK}, \cite{KKM}, \cite{KM}, \cite{LZ}, \cite{TT}, \cite{B}, \cite{Shi}, \cite{AI1}, \cite{AI2}, \cite{AI3}, \cite{AB1}, \cite{AB2}, \cite{Fon2}, \cite{Fon3}, \cite{Fa1}, \cite{M}, \cite{Fa2}, \cite{Fa3}, \cite{Fon4}, \cite{Fon5}, \cite{Fon6}; 
\item[$\square$] Chapter 2 Main References: \cite{Sch1}, \cite{Sch2}, \cite{FS}, \cite{KL1}, \cite{KL2}, \cite{BL1}, \cite{BL2}, \cite{BS}, \cite{BHS}, \cite{Fon1}, \cite{CS1}, \cite{CS2}, \cite{BK}, \cite{BBK}, \cite{BBBK}, \cite{KKM}, \cite{KM}, \cite{LZ}, \cite{TT}, \cite{B}, \cite{Shi}, \cite{AI1}, \cite{AI2}, \cite{AI3}, \cite{AB1}, \cite{AB2}, \cite{Fon2}, \cite{Fon3}, \cite{Fa1}, \cite{M}, \cite{Fa2}, \cite{Fa3}, \cite{Fon4}, \cite{Fon5}, \cite{Fon6};
\item[$\square$] Chapter 3 Main References: \cite{Sch1}, \cite{Sch2}, \cite{FS}, \cite{FF}, \cite{KL1}, \cite{KL2}, \cite{BL1}, \cite{BL2}, \cite{BS}, \cite{BHS}, \cite{Fon1}, \cite{CS1}, \cite{CS2}, \cite{BK}, \cite{BBK}, \cite{BBBK}, \cite{KKM}, \cite{KM}, \cite{LBV}, \cite{B}, \cite{SW}; \cite{Sch1}, \cite{Sch2}, \cite{FS}, \cite{KL1}, \cite{KL2}, \cite{BL1}, \cite{BL2}, \cite{BS}, \cite{BHS}, \cite{Fon1}, \cite{CS1}, \cite{CS2}, \cite{BK}, \cite{BBK}, \cite{BBBK}, \cite{KKM}, \cite{KM}, \cite{LZ}, \cite{TT}, \cite{B}, \cite{Shi}, \cite{AI1}, \cite{AI2}, \cite{AI3}, \cite{AB1}, \cite{AB2}, \cite{Fon2}, \cite{Fon3}, \cite{Fa1}, \cite{M}, \cite{Fa2}, \cite{Fa3}, \cite{Fon4}, \cite{Fon5}, \cite{Fon6};
\item[$\square$] Chapter 4 Main References I: \cite{FS}, \cite{FF}, \cite{Sch1}, \cite{Sch2}, \cite{KL1}, \cite{KL2}, \cite{LBV}, \cite{B}, \cite{SW}, \cite{BS}, \cite{Lan1}, \cite{Drin1}, \cite{Drin2}, \cite{Zhu}, \cite{DHKM};
\item[$\square$] Chapter 4 Main References II: \cite{pHodgeT}, \cite{pHodgeF}, \cite{pHodgeS1}, \cite{pHodgeS2}, \cite{pHodgeKL1}, \cite{pHodgeKL2}, \cite{pHodgeBS}, \cite{pHodgeKPX}, \cite{pToAnCS1}, \cite{pToAnCS2}, \cite{pToAnCS3}, \cite{pToAnCS4},  
\cite{pToAnBBBK}, \cite{LPL}, \cite{LPD1}, \cite{LPLL}, \cite{LPVL}, \cite{LPC}, \cite{LPFS}, \cite{LPGL}, \cite{LPEGH}, \cite{LPEG}, \cite{LPZ}, \cite{LPDHKM}, \cite{LPD2}, \cite{LPL}, \cite{LPD1}, \cite{LPLL}, \cite{LPVL}, \cite{LPC}, \cite{LPFS}, \cite{LPGL}, \cite{LPEGH}, \cite{LPEG}, \cite{LPZ}, \cite{LPDHKM}, \cite{LPD2}.
\end{itemize}
\end{reference}

\newpage

\noindent \textbf{Notations}:
\begin{itemize}
\item[$\square$] Chapter 1: The period sheaves in the pro-\'etale topology in this chapter are assumed to be already tensored with a finite extension of $\mathbb{Q}_p$ containing square roots of $p$, although we do write that notation in an explicit way. We assume the corresonding interval $I$ contains 1. 
\begin{align}
&\Gamma_{\text{deRham},X,\text{pro\'et}}\{t^{1/2}\},\Gamma^\mathcal{O}_{\text{deRham},X,\text{pro\'et}}\{t^{1/2}\},\\
&\Gamma_{\text{deRham},X,\text{pro\'et}}\{t^{1/2},\log(t)\},\Gamma^\mathcal{O}_{\text{deRham},X,\text{pro\'et}}\{t^{1/2},\log(t)\},\\
&\Gamma^\mathrm{perfect}_{\text{Robba},X,\text{pro\'et}}\{t^{1/2}\},\Gamma^\mathrm{perfect}_{\text{Robba},X,\text{pro\'et},\infty}\{t^{1/2}\},\Gamma^\mathrm{perfect}_{\text{Robba},X,\text{pro\'et},I}\{t^{1/2}\},\\
&\Gamma^\mathrm{perfect}_{\text{Robba},X,\text{pro\'et}}\{t^{1/2},\log(t)\},\Gamma^\mathrm{perfect}_{\text{Robba},X,\text{pro\'et},\infty}\{t^{1/2},\log(t)\},\Gamma^\mathrm{perfect}_{\text{Robba},X,\text{pro\'et},I}\{t^{1/2},\log(t)\};
\end{align} 
\begin{align}
&\Gamma_{\text{cristalline},X,\text{pro\'et}}\{t^{1/2}\},\Gamma^\mathcal{O}_{\text{cristalline},X,\text{pro\'et}}\{t^{1/2}\},\\
&\Gamma_{\text{cristalline},X,\text{pro\'et}}\{t^{1/2},\log(t)\},\Gamma^\mathcal{O}_{\text{cristalline},X,\text{pro\'et}}\{t^{1/2},\log(t)\},\\
&\Gamma^\mathrm{perfect}_{\text{Robba},X,\text{pro\'et}}\{t^{1/2}\},\Gamma^\mathrm{perfect}_{\text{Robba},X,\text{pro\'et},\infty}\{t^{1/2}\},\Gamma^\mathrm{perfect}_{\text{Robba},X,\text{pro\'et},I}\{t^{1/2}\},\\
&\Gamma^\mathrm{perfect}_{\text{Robba},X,\text{pro\'et}}\{t^{1/2},\log(t)\},\Gamma^\mathrm{perfect}_{\text{Robba},X,\text{pro\'et},\infty}\{t^{1/2},\log(t)\},\Gamma^\mathrm{perfect}_{\text{Robba},X,\text{pro\'et},I}\{t^{1/2},\log(t)\}.
\end{align}
\item[$\square$] Chapter 2: The period sheaves in the $v$-topology in this chapter are assumed to be already tensored with a finite extension of $\mathbb{Q}_p$ containing square roots of $p$, although we do write that notation in an explicit way. We assume the corresonding interval $I$ contains 1.
\begin{align}
&\Gamma_{\text{deRham},X,v}\{t^{1/2}\},\Gamma^\mathcal{O}_{\text{deRham},X,v}\{t^{1/2}\},\\
&\Gamma_{\text{deRham},X,v}\{t^{1/2},\log(t)\},\Gamma^\mathcal{O}_{\text{deRham},X,v}\{t^{1/2},\log(t)\},\\
&\Gamma^\mathrm{perfect}_{\text{Robba},X,v}\{t^{1/2}\},\Gamma^\mathrm{perfect}_{\text{Robba},X,v,\infty}\{t^{1/2}\},\Gamma^\mathrm{perfect}_{\text{Robba},X,v,I}\{t^{1/2}\},\\
&\Gamma^\mathrm{perfect}_{\text{Robba},X,v}\{t^{1/2},\log(t)\},\Gamma^\mathrm{perfect}_{\text{Robba},X,v,\infty}\{t^{1/2},\log(t)\},\Gamma^\mathrm{perfect}_{\text{Robba},X,v,I}\{t^{1/2},\log(t)\};
\end{align}
\begin{align}
&\Gamma_{\text{cristalline},X,v}\{t^{1/2}\},\Gamma^\mathcal{O}_{\text{cristalline},X,v}\{t^{1/2}\},\\
&\Gamma_{\text{cristalline},X,v}\{t^{1/2},\log(t)\},\Gamma^\mathcal{O}_{\text{cristalline},X,v}\{t^{1/2},\log(t)\},\\
&\Gamma^\mathrm{perfect}_{\text{Robba},X,v}\{t^{1/2}\},\Gamma^\mathrm{perfect}_{\text{Robba},X,v,\infty}\{t^{1/2}\},\Gamma^\mathrm{perfect}_{\text{Robba},X,v,I}\{t^{1/2}\},\\
&\Gamma^\mathrm{perfect}_{\text{Robba},X,v}\{t^{1/2},\log(t)\},\Gamma^\mathrm{perfect}_{\text{Robba},X,v,\infty}\{t^{1/2},\log(t)\},\Gamma^\mathrm{perfect}_{\text{Robba},X,v,I}\{t^{1/2},\log(t)\}.
\end{align}
\item[$\square$] Chapter 3: The period rings in 3.1, 3.2 are assumed to be not already tensored with a finite extension of $\mathbb{Q}_p$ containing square roots of $p$,  we do write that notation in an explicit way; Then the period sheaves in the $v$-topology in 3.3, 3.4, 3.5 are assumed to be already tensored with a finite extension of $\mathbb{Q}_p$ containing square roots of $p$, although we do write that notation in an explicit way. We assume the corresonding interval $I$ contains 1.
\begin{align}
&\Gamma_{\text{deRham},X,v}\{t^{1/2}\},\Gamma^\mathcal{O}_{\text{deRham},X,v}\{t^{1/2}\},\\
&\Gamma_{\text{deRham},X,v}\{t^{1/2},\log(t)\},\Gamma^\mathcal{O}_{\text{deRham},X,v}\{t^{1/2},\log(t)\},\\
&\Gamma^\mathrm{perfect}_{\text{Robba},X,v}\{t^{1/2}\},\Gamma^\mathrm{perfect}_{\text{Robba},X,v,\infty}\{t^{1/2}\},\Gamma^\mathrm{perfect}_{\text{Robba},X,v,I}\{t^{1/2}\},\\
&\Gamma^\mathrm{perfect}_{\text{Robba},X,v}\{t^{1/2},\log(t)\},\Gamma^\mathrm{perfect}_{\text{Robba},X,v,\infty}\{t^{1/2},\log(t)\},\Gamma^\mathrm{perfect}_{\text{Robba},X,v,I}\{t^{1/2},\log(t)\};
\end{align}
\begin{align}
&\Gamma_{\text{cristalline},X,v}\{t^{1/2}\},\Gamma^\mathcal{O}_{\text{cristalline},X,v}\{t^{1/2}\},\\
&\Gamma_{\text{cristalline},X,v}\{t^{1/2},\log(t)\},\Gamma^\mathcal{O}_{\text{cristalline},X,v}\{t^{1/2},\log(t)\},\\
&\Gamma^\mathrm{perfect}_{\text{Robba},X,v}\{t^{1/2}\},\Gamma^\mathrm{perfect}_{\text{Robba},X,v,\infty}\{t^{1/2}\},\Gamma^\mathrm{perfect}_{\text{Robba},X,v,I}\{t^{1/2}\},\\
&\Gamma^\mathrm{perfect}_{\text{Robba},X,v}\{t^{1/2},\log(t)\},\Gamma^\mathrm{perfect}_{\text{Robba},X,v,\infty}\{t^{1/2},\log(t)\},\Gamma^\mathrm{perfect}_{\text{Robba},X,v,I}\{t^{1/2},\log(t)\}.
\end{align} 
\item[$\square$] Chapter 4: The period rings in this chapter are assumed to be not already tensored with a finite extension of $\mathbb{Q}_p$ containing square roots of $p$,  we do write that notation in an explicit way.
\end{itemize}

\chapter{Mixed-Parity $p$-adic Hodge Modules over Pro-\'Etale Sites}

\newpage
\section{Geometric Family of Mixed-Parity Hodge Modules I: de Rham Situations}

\subsection{Period Rings and Sheaves}

\begin{reference}
\cite{Sch1}, \cite{KL1}, \cite{KL2}, \cite{BL1}, \cite{BL2}, \cite{BS}, \cite{BHS}, \cite{Fon1}, \cite{CS1}, \cite{CS2}, \cite{BK}, \cite{BBK}, \cite{BBBK}, \cite{KKM}, \cite{KM}, \cite{LZ}, \cite{M}.
\end{reference}

\subsubsection{Rings}

\noindent Let $X$ be a rigid analytic space over $\mathbb{Q}_p$. We have the corresponding \'etale site and the corresponding pro-\'etale site of $X$, which we denote them by $X_{\text{pro\'et}},X_\text{\'et}$. The relationship of the two sites can be reflected by the corresponding morphism $f:X_{\text{pro\'et}}\longrightarrow X_\text{\'et}$. Then we have the corresponding de Rham period rings and sheaves from \cite{Sch1}:
\begin{align}
\Gamma_{\text{deRham},X,\text{pro\'et}}, \Gamma^\mathcal{O}_{\text{deRham},X,\text{pro\'et}}.
\end{align}
Our notations are different from \cite{Sch1}, we use $\Gamma$ to mean $B$ in \cite{Sch1}, while $\Gamma^\mathcal{O}$ will be the corresponding $OB$ ring in \cite{Sch1}.\\

\begin{definition}
\indent Now we assume that $p>2$, following \cite{BS} we join the square root of $t$ element in $\Gamma_{\text{deRham},X,\text{pro\'et}}$ which forms the sheaves:
\begin{align}
\Gamma_{\text{deRham},X,\text{pro\'et}}\{t^{1/2}\},\Gamma^\mathcal{O}_{\text{deRham},X,\text{pro\'et}}\{t^{1/2}\}.
\end{align}
And following \cite{BL1}, \cite{BL2}, \cite{Fon1}, \cite{BHS} we further have the following sheaves of rings:
\begin{align}
\Gamma_{\text{deRham},X,\text{pro\'et}}\{t^{1/2},\log(t)\},\Gamma^\mathcal{O}_{\text{deRham},X,\text{pro\'et}}\{t^{1/2},\log(t)\}.
\end{align}
\end{definition}

\begin{definition}
We use the notations:
\begin{align}
\Gamma^\mathrm{perfect}_{\text{Robba},X,\text{pro\'et}},\Gamma^\mathrm{perfect}_{\text{Robba},X,\text{pro\'et},\infty},\Gamma^\mathrm{perfect}_{\text{Robba},X,\text{pro\'et},I}
\end{align}
to denote the perfect Robba rings from \cite{KL1}, \cite{KL2}, where $I\subset (0,\infty)$. Then we join $t^{1/2}$ to these sheaves we have:
\begin{align}
\Gamma^\mathrm{perfect}_{\text{Robba},X,\text{pro\'et}}\{t^{1/2}\},\Gamma^\mathrm{perfect}_{\text{Robba},X,\text{pro\'et},\infty}\{t^{1/2}\},\Gamma^\mathrm{perfect}_{\text{Robba},X,\text{pro\'et},I}\{t^{1/2}\}.
\end{align}
And following \cite{BL1}, \cite{BL2}, \cite{Fon1}, \cite{BHS} we have the following larger sheaves:
\begin{align}
\Gamma^\mathrm{perfect}_{\text{Robba},X,\text{pro\'et}}\{t^{1/2},\log(t)\},\Gamma^\mathrm{perfect}_{\text{Robba},X,\text{pro\'et},\infty}\{t^{1/2},\log(t)\},\Gamma^\mathrm{perfect}_{\text{Robba},X,\text{pro\'et},I}\{t^{1/2},\log(t)\}.
\end{align} 
\end{definition}

\begin{definition}
From now on, we use the same notation to denote the period rings involved tensored with a finite extension of $\mathbb{Q}_p$ containing square root of $p$ as in \cite{BS}.
\begin{align}
\Gamma_{\text{deRham},X,\text{pro\'et}}\{t^{1/2}\},\Gamma^\mathcal{O}_{\text{deRham},X,\text{pro\'et}}\{t^{1/2}\}.
\end{align}
\begin{align}
\Gamma_{\text{deRham},X,\text{pro\'et}}\{t^{1/2},\log(t)\},\Gamma^\mathcal{O}_{\text{deRham},X,\text{pro\'et}}\{t^{1/2},\log(t)\}.
\end{align}
\begin{align}
\Gamma^\mathrm{perfect}_{\text{Robba},X,\text{pro\'et}}\{t^{1/2}\},\Gamma^\mathrm{perfect}_{\text{Robba},X,\text{pro\'et},\infty}\{t^{1/2}\},\Gamma^\mathrm{perfect}_{\text{Robba},X,\text{pro\'et},I}\{t^{1/2}\}.
\end{align}
\begin{align}
\Gamma^\mathrm{perfect}_{\text{Robba},X,\text{pro\'et}}\{t^{1/2},\log(t)\},\Gamma^\mathrm{perfect}_{\text{Robba},X,\text{pro\'et},\infty}\{t^{1/2},\log(t)\},\Gamma^\mathrm{perfect}_{\text{Robba},X,\text{pro\'et},I}\{t^{1/2},\log(t)\}.
\end{align}
This is necessary since we to extend the action of $\varphi$ to the period rings by $\varphi(t^{1/2}\otimes 1)=\varphi(t)^{1/2}\otimes 1$.
\end{definition}

\subsubsection{Modules}

\noindent We consider quasicoherent presheaves in the following two situation:
\begin{itemize}
\item[$\square$] The solid quasicoherent modules from \cite{CS1}, \cite{CS2};
\item[$\square$] The ind-Banach quasicoherent modules from \cite{BK}, \cite{BBK}, \cite{BBBK}, \cite{KKM}, \cite{KM} with the corresponding monomorphic ind-Banach quasicoherent modules from \cite{BK}, \cite{BBK}, \cite{BBBK}, \cite{KKM}, \cite{KM}.
\end{itemize}

\begin{definition}
We use the notation:
\begin{align}
\mathrm{preModule}^\mathrm{solid,quasicoherent}_{\square,\Gamma^\mathrm{perfect}_{\text{Robba},X,\text{pro\'et}}\{t^{1/2}\}},\mathrm{preModule}^\mathrm{solid,quasicoherent}_{\square,\Gamma^\mathrm{perfect}_{\text{Robba},X,\text{pro\'et},\infty}\{t^{1/2}\}},
\mathrm{preModule}^\mathrm{solid,quasicoherent}_{\square,\Gamma^\mathrm{perfect}_{\text{Robba},X,\text{pro\'et},I}\{t^{1/2}\}} 
\end{align}
to denote the $(\infty,1)$-categories of solid quasicoherent presheaves over the corresonding Robba sheaves. Locally the section is defined by taking the corresponding $(\infty,1)$-categories of solid modules.
\end{definition}

\begin{definition}
We use the notation:
\begin{align}
\mathrm{preModule}^\mathrm{ind-Banach,quasicoherent}_{\Gamma^\mathrm{perfect}_{\text{Robba},X,\text{pro\'et}}\{t^{1/2}\}},\\
\mathrm{preModule}^\mathrm{ind-Banach,quasicoherent}_{\Gamma^\mathrm{perfect}_{\text{Robba},X,\text{pro\'et},\infty}\{t^{1/2}\}},\\
\mathrm{preModule}^\mathrm{ind-Banach,quasicoherent}_{\Gamma^\mathrm{perfect}_{\text{Robba},X,\text{pro\'et},I}\{t^{1/2}\}} 
\end{align}
to denote the $(\infty,1)$-categories of solid quasicoherent presheaves over the corresonding Robba sheaves. Locally the section is defined by taking the corresponding $(\infty,1)$-categories of inductive Banach  modules. 
\end{definition}

\begin{definition}
We use the notation:
\begin{align}
\mathrm{Module}^\mathrm{solid,quasicoherent}_{\square,\Gamma^\mathrm{perfect}_{\text{Robba},X,\text{pro\'et}}\{t^{1/2}\}},\mathrm{Module}^\mathrm{solid,quasicoherent}_{\square,\Gamma^\mathrm{perfect}_{\text{Robba},X,\text{pro\'et},\infty}\{t^{1/2}\}},
\mathrm{Module}^\mathrm{solid,quasicoherent}_{\square,\Gamma^\mathrm{perfect}_{\text{Robba},X,\text{pro\'et},I}\{t^{1/2}\}} 
\end{align}
to denote the $(\infty,1)$-categories of solid quasicoherent sheaves over the corresonding Robba sheaves. Locally the section is defined by taking the corresponding $(\infty,1)$-categories of solid modules.
\end{definition}

\subsubsection{Mixed-Parity Hodge Modules without Frobenius}

\noindent Now we consider the key objects in our study namely those complexes generated by certain mixed-parity Hodge modules. We start from the following definition.

\begin{definition}
For any locally free coherent sheaf $F$ over
\begin{align}
\Gamma^\mathrm{perfect}_{\text{Robba},X,\text{pro\'et},\infty}\{t^{1/2}\},\Gamma^\mathrm{perfect}_{\text{Robba},X,\text{pro\'et},I}\{t^{1/2}\},
\end{align} 
we consider the following functor $\mathrm{dR}$ sending $F$ to the following object:
\begin{align}
f_*(F\otimes_{\Gamma^\mathrm{perfect}_{\text{Robba},X,\text{pro\'et},\infty}\{t^{1/2}\}} \Gamma^\mathcal{O}_{\text{deRham},X,\text{pro\'et}}\{t^{1/2}\})
\end{align}
or 
\begin{align}
f_*(F\otimes_{\Gamma^\mathrm{perfect}_{\text{Robba},X,\text{pro\'et},I}\{t^{1/2}\}} \Gamma^\mathcal{O}_{\text{deRham},X,\text{pro\'et}}\{t^{1/2}\}).
\end{align}
We call $F$ mixed-parity de Rham if we have the following isomorphism\footnote{As in \cite[Definition 10.10]{KL}, when we consider the corresponding de Rham, cristalline, semi-stable functors we will assume 1 is belonging to the interval $I$ in all the following corresponding discussion.}:
\begin{align}
f^*f_*(F\otimes_{\Gamma^\mathrm{perfect}_{\text{Robba},X,\text{pro\'et},\infty}\{t^{1/2}\}} \Gamma^\mathcal{O}_{\text{deRham},X,\text{pro\'et}}\{t^{1/2}\}) \otimes \Gamma^\mathcal{O}_{\text{deRham},X,\text{pro\'et}}\{t^{1/2}\} \overset{\sim}{\longrightarrow} F \otimes \Gamma^\mathcal{O}_{\text{deRham},X,\text{pro\'et}}\{t^{1/2}\} 
\end{align}
or 
\begin{align}
f^*f_*(F\otimes_{\Gamma^\mathrm{perfect}_{\text{Robba},X,\text{pro\'et},I}\{t^{1/2}\}} \Gamma^\mathcal{O}_{\text{deRham},X,\text{pro\'et}}\{t^{1/2}\}) \otimes \Gamma^\mathcal{O}_{\text{deRham},X,\text{pro\'et}}\{t^{1/2}\} \overset{\sim}{\longrightarrow} F \otimes \Gamma^\mathcal{O}_{\text{deRham},X,\text{pro\'et}}\{t^{1/2}\}. 
\end{align}
\end{definition}

\begin{definition}
For any locally free coherent sheaf $F$ over
\begin{align}
\Gamma^\mathrm{perfect}_{\text{Robba},X,\text{pro\'et},\infty}\{t^{1/2}\},\Gamma^\mathrm{perfect}_{\text{Robba},X,\text{pro\'et},I}\{t^{1/2}\},
\end{align} 
we consider the following functor $\mathrm{dR}^\mathrm{almost}$ sending $F$ to the following object:
\begin{align}
f_*(F\otimes_{\Gamma^\mathrm{perfect}_{\text{Robba},X,\text{pro\'et},\infty}\{t^{1/2}\}} \Gamma^\mathcal{O}_{\text{deRham},X,\text{pro\'et}}\{t^{1/2},\log(t)\})
\end{align}
or 
\begin{align}
f_*(F\otimes_{\Gamma^\mathrm{perfect}_{\text{Robba},X,\text{pro\'et},I}\{t^{1/2}\}} \Gamma^\mathcal{O}_{\text{deRham},X,\text{pro\'et}}\{t^{1/2},\log(t)\}).
\end{align}
We call $F$ mixed-parity almost de Rham if we have the following isomorphism:
\begin{align}
f^*f_*(F\otimes_{\Gamma^\mathrm{perfect}_{\text{Robba},X,\text{pro\'et},\infty}\{t^{1/2}\}} \Gamma^\mathcal{O}_{\text{deRham},X,\text{pro\'et}}\{t^{1/2},\log(t)\}) \otimes \Gamma^\mathcal{O}_{\text{deRham},X,\text{pro\'et}}\{t^{1/2},\log(t)\} \\
\overset{\sim}{\longrightarrow} F \otimes \Gamma^\mathcal{O}_{\text{deRham},X,\text{pro\'et}}\{t^{1/2},\log(t)\} 
\end{align}
or 
\begin{align}
f^*f_*(F\otimes_{\Gamma^\mathrm{perfect}_{\text{Robba},X,\text{pro\'et},I}\{t^{1/2}\}} \Gamma^\mathcal{O}_{\text{deRham},X,\text{pro\'et}}\{t^{1/2},\log(t)\}) \otimes \Gamma^\mathcal{O}_{\text{deRham},X,\text{pro\'et}}\{t^{1/2},\log(t)\}\\ \overset{\sim}{\longrightarrow} F \otimes \Gamma^\mathcal{O}_{\text{deRham},X,\text{pro\'et}}\{t^{1/2},\log(t)\}. 
\end{align}
\end{definition}

\noindent We now define the $(\infty,1)$-categories of mixed-parity de Rham modules and he corresponding mixed-parity almost de Rham modules by using the objects involved to generated these categories:

\begin{definition}
Considering all the mixed parity de Rham bundles (locally finite free) as defined above, we consider the sub-$(\infty,1)$ category of 
\begin{align}
\mathrm{preModule}^\mathrm{solid,quasicoherent}_{\square,\Gamma^\mathrm{perfect}_{\text{Robba},X,\text{pro\'et},\infty}\{t^{1/2}\}},
\mathrm{preModule}^\mathrm{solid,quasicoherent}_{\square,\Gamma^\mathrm{perfect}_{\text{Robba},X,\text{pro\'et},I}\{t^{1/2}\}} 
\end{align}
generated by the mixed-parity de Rham bundles (locally finite free ones). These are defined to be the $(\infty,1)$-categories of mixed-parity de Rham complexes:
\begin{align}
\mathrm{preModule}^\mathrm{solid,quasicoherent,mixed-paritydeRham}_{\square,\Gamma^\mathrm{perfect}_{\text{Robba},X,\text{pro\'et},\infty}\{t^{1/2}\}},
\mathrm{preModule}^\mathrm{solid,quasicoherent,mixed-paritydeRham}_{\square,\Gamma^\mathrm{perfect}_{\text{Robba},X,\text{pro\'et},I}\{t^{1/2}\}}. 
\end{align}
\end{definition}

\begin{definition}
Considering all the mixed parity almost de Rham bundles (locally finite free) as defined above, we consider the sub-$(\infty,1)$ category of 
\begin{align}
\mathrm{preModule}^\mathrm{solid,quasicoherent}_{\square,\Gamma^\mathrm{perfect}_{\text{Robba},X,\text{pro\'et},\infty}\{t^{1/2}\}},
\mathrm{preModule}^\mathrm{solid,quasicoherent}_{\square,\Gamma^\mathrm{perfect}_{\text{Robba},X,\text{pro\'et},I}\{t^{1/2}\}} 
\end{align}
generated by the mixed-parity almost de Rham bundles (locally finite free ones). These are defined to be the $(\infty,1)$-categories of mixed-parity de Rham complexes:
\begin{align}
\mathrm{preModule}^\mathrm{solid,quasicoherent,mixed-parityalmostdeRham}_{\square,\Gamma^\mathrm{perfect}_{\text{Robba},X,\text{pro\'et},\infty}\{t^{1/2}\}},
\mathrm{preModule}^\mathrm{solid,quasicoherent,mixed-parityalmostdeRham}_{\square,\Gamma^\mathrm{perfect}_{\text{Robba},X,\text{pro\'et},I}\{t^{1/2}\}}. 
\end{align}
\end{definition}

\indent Then the corresponding mixed-parity de Rham functors can be extended to these categories:
\begin{align}
\mathrm{preModule}^\mathrm{solid,quasicoherent,mixed-paritydeRham}_{\square,\Gamma^\mathrm{perfect}_{\text{Robba},X,\text{pro\'et},\infty}\{t^{1/2}\}},
\mathrm{preModule}^\mathrm{solid,quasicoherent,mixed-paritydeRham}_{\square,\Gamma^\mathrm{perfect}_{\text{Robba},X,\text{pro\'et},I}\{t^{1/2}\}}, 
\end{align}
and
\begin{align}
\mathrm{preModule}^\mathrm{solid,quasicoherent,mixed-parityalmostdeRham}_{\square,\Gamma^\mathrm{perfect}_{\text{Robba},X,\text{pro\'et},\infty}\{t^{1/2}\}},\\
\mathrm{preModule}^\mathrm{solid,quasicoherent,mixed-parityalmostdeRham}_{\square,\Gamma^\mathrm{perfect}_{\text{Robba},X,\text{pro\'et},I}\{t^{1/2}\}}. 
\end{align}

\subsubsection{Mixed-Parity Hodge Modules with Frobenius}

\noindent Now we consider the key objects in our study namely those complexes generated by certain mixed-parity Hodge modules. We start from the following definition.

\begin{remark}
All the coherent sheaves over mixed-parity Robba sheaves in this section will carry the corresponding Frobenius morphism $\varphi: F \overset{\sim}{\longrightarrow} \varphi^*F$.
\end{remark}

\begin{definition}
For any locally free coherent sheaf $F$ over
\begin{align}
\Gamma^\mathrm{perfect}_{\text{Robba},X,\text{pro\'et},\infty}\{t^{1/2}\},\Gamma^\mathrm{perfect}_{\text{Robba},X,\text{pro\'et},I}\{t^{1/2}\},
\end{align} 
we consider the following functor $\mathrm{dR}$ sending $F$ to the following object:
\begin{align}
f_*(F\otimes_{\Gamma^\mathrm{perfect}_{\text{Robba},X,\text{pro\'et},\infty}\{t^{1/2}\}} \Gamma^\mathcal{O}_{\text{deRham},X,\text{pro\'et}}\{t^{1/2}\})
\end{align}
or 
\begin{align}
f_*(F\otimes_{\Gamma^\mathrm{perfect}_{\text{Robba},X,\text{pro\'et},I}\{t^{1/2}\}} \Gamma^\mathcal{O}_{\text{deRham},X,\text{pro\'et}}\{t^{1/2}\}).
\end{align}
We call $F$ mixed-parity de Rham if we have the following isomorphism:
\begin{align}
f^*f_*(F\otimes_{\Gamma^\mathrm{perfect}_{\text{Robba},X,\text{pro\'et},\infty}\{t^{1/2}\}} \Gamma^\mathcal{O}_{\text{deRham},X,\text{pro\'et}}\{t^{1/2}\}) \otimes \Gamma^\mathcal{O}_{\text{deRham},X,\text{pro\'et}}\{t^{1/2}\} \overset{\sim}{\longrightarrow} F \otimes \Gamma^\mathcal{O}_{\text{deRham},X,\text{pro\'et}}\{t^{1/2}\} 
\end{align}
or 
\begin{align}
f^*f_*(F\otimes_{\Gamma^\mathrm{perfect}_{\text{Robba},X,\text{pro\'et},I}\{t^{1/2}\}} \Gamma^\mathcal{O}_{\text{deRham},X,\text{pro\'et}}\{t^{1/2}\}) \otimes \Gamma^\mathcal{O}_{\text{deRham},X,\text{pro\'et}}\{t^{1/2}\} \overset{\sim}{\longrightarrow} F \otimes \Gamma^\mathcal{O}_{\text{deRham},X,\text{pro\'et}}\{t^{1/2}\}. 
\end{align}
\end{definition}

\begin{definition}
For any locally free coherent sheaf $F$ over
\begin{align}
\Gamma^\mathrm{perfect}_{\text{Robba},X,\text{pro\'et},\infty}\{t^{1/2}\},\Gamma^\mathrm{perfect}_{\text{Robba},X,\text{pro\'et},I}\{t^{1/2}\},
\end{align} 
we consider the following functor $\mathrm{dR}^\mathrm{almost}$ sending $F$ to the following object:
\begin{align}
f_*(F\otimes_{\Gamma^\mathrm{perfect}_{\text{Robba},X,\text{pro\'et},\infty}\{t^{1/2}\}} \Gamma^\mathcal{O}_{\text{deRham},X,\text{pro\'et}}\{t^{1/2},\log(t)\})
\end{align}
or 
\begin{align}
f_*(F\otimes_{\Gamma^\mathrm{perfect}_{\text{Robba},X,\text{pro\'et},I}\{t^{1/2}\}} \Gamma^\mathcal{O}_{\text{deRham},X,\text{pro\'et}}\{t^{1/2},\log(t)\}).
\end{align}
We call $F$ mixed-parity almost de Rham if we have the following isomorphism:
\begin{align}
f^*f_*(F\otimes_{\Gamma^\mathrm{perfect}_{\text{Robba},X,\text{pro\'et},\infty}\{t^{1/2}\}} \Gamma^\mathcal{O}_{\text{deRham},X,\text{pro\'et}}\{t^{1/2},\log(t)\}) \otimes \Gamma^\mathcal{O}_{\text{deRham},X,\text{pro\'et}}\{t^{1/2},\log(t)\} \\
\overset{\sim}{\longrightarrow} F \otimes \Gamma^\mathcal{O}_{\text{deRham},X,\text{pro\'et}}\{t^{1/2},\log(t)\} 
\end{align}
or 
\begin{align}
f^*f_*(F\otimes_{\Gamma^\mathrm{perfect}_{\text{Robba},X,\text{pro\'et},I}\{t^{1/2}\}} \Gamma^\mathcal{O}_{\text{deRham},X,\text{pro\'et}}\{t^{1/2},\log(t)\}) \otimes \Gamma^\mathcal{O}_{\text{deRham},X,\text{pro\'et}}\{t^{1/2},\log(t)\}\\ \overset{\sim}{\longrightarrow} F \otimes \Gamma^\mathcal{O}_{\text{deRham},X,\text{pro\'et}}\{t^{1/2},\log(t)\}. 
\end{align}
\end{definition}

\noindent We now define the $(\infty,1)$-categories of mixed-parity de Rham modules and he corresponding mixed-parity almost de Rham modules by using the objects involved to generated these categories:

\begin{definition}
Considering all the mixed parity de Rham bundles (locally finite free) as defined above, we consider the sub-$(\infty,1)$-category of 
\begin{align}
\varphi\mathrm{preModule}^\mathrm{solid,quasicoherent}_{\square,\Gamma^\mathrm{perfect}_{\text{Robba},X,\text{pro\'et},\infty}\{t^{1/2}\}},
\varphi\mathrm{preModule}^\mathrm{solid,quasicoherent}_{\square,\Gamma^\mathrm{perfect}_{\text{Robba},X,\text{pro\'et},I}\{t^{1/2}\}} 
\end{align}
generated by the mixed-parity de Rham bundles (locally finite free ones). These are defined to be the $(\infty,1)$-categories of mixed-parity de Rham complexes:
\begin{align}
\varphi\mathrm{preModule}^\mathrm{solid,quasicoherent,mixed-paritydeRham}_{\square,\Gamma^\mathrm{perfect}_{\text{Robba},X,\text{pro\'et},\infty}\{t^{1/2}\}},
\varphi\mathrm{preModule}^\mathrm{solid,quasicoherent,mixed-paritydeRham}_{\square,\Gamma^\mathrm{perfect}_{\text{Robba},X,\text{pro\'et},I}\{t^{1/2}\}}. 
\end{align}
\end{definition}

\begin{definition}
Considering all the mixed parity almost de Rham bundles (locally finite free) as defined above, we consider the sub-$(\infty,1)$ category of 
\begin{align}
\varphi\mathrm{preModule}^\mathrm{solid,quasicoherent}_{\square,\Gamma^\mathrm{perfect}_{\text{Robba},X,\text{pro\'et},\infty}\{t^{1/2}\}},
\varphi\mathrm{preModule}^\mathrm{solid,quasicoherent}_{\square,\Gamma^\mathrm{perfect}_{\text{Robba},X,\text{pro\'et},I}\{t^{1/2}\}} 
\end{align}
generated by the mixed-parity almost de Rham bundles (locally finite free ones). These are defined to be the $(\infty,1)$-categories of mixed-parity de Rham complexes:
\begin{align}
\varphi\mathrm{preModule}^\mathrm{solid,quasicoherent,mixed-parityalmostdeRham}_{\square,\Gamma^\mathrm{perfect}_{\text{Robba},X,\text{pro\'et},\infty}\{t^{1/2}\}},
\varphi\mathrm{preModule}^\mathrm{solid,quasicoherent,mixed-parityalmostdeRham}_{\square,\Gamma^\mathrm{perfect}_{\text{Robba},X,\text{pro\'et},I}\{t^{1/2}\}}. 
\end{align}
\end{definition}

\indent Then the corresponding mixed-parity de Rham functors can be extended to these categories:
\begin{align}
\varphi\mathrm{preModule}^\mathrm{solid,quasicoherent,mixed-paritydeRham}_{\square,\Gamma^\mathrm{perfect}_{\text{Robba},X,\text{pro\'et},\infty}\{t^{1/2}\}},
\varphi\mathrm{preModule}^\mathrm{solid,quasicoherent,mixed-paritydeRham}_{\square,\Gamma^\mathrm{perfect}_{\text{Robba},X,\text{pro\'et},I}\{t^{1/2}\}}, 
\end{align}
and
\begin{align}
\varphi\mathrm{preModule}^\mathrm{solid,quasicoherent,mixed-parityalmostdeRham}_{\square,\Gamma^\mathrm{perfect}_{\text{Robba},X,\text{pro\'et},\infty}\{t^{1/2}\}},\\
\varphi\mathrm{preModule}^\mathrm{solid,quasicoherent,mixed-parityalmostdeRham}_{\square,\Gamma^\mathrm{perfect}_{\text{Robba},X,\text{pro\'et},I}\{t^{1/2}\}}. 
\end{align}

\subsection{Mixed-Parity de Rham Riemann-Hilbert Correspondence}

\indent This chapter will extend the corresponding Riemann-Hilbert correspondence from \cite{Sch1}, \cite{LZ}, \cite{BL1}, \cite{BL2}, \cite{M} to the mixed-parity setting.

\begin{definition}
We define the following Riemann-Hilbert functor $\text{RH}_\text{mixed-parity}$ from the one of categories:
\begin{align}
\mathrm{preModule}^\mathrm{solid,quasicoherent,mixed-paritydeRham}_{\square,\Gamma^\mathrm{perfect}_{\text{Robba},X,\text{pro\'et},\infty}\{t^{1/2}\}},
\mathrm{preModule}^\mathrm{solid,quasicoherent,mixed-paritydeRham}_{\square,\Gamma^\mathrm{perfect}_{\text{Robba},X,\text{pro\'et},I}\{t^{1/2}\}}, 
\end{align}
and
\begin{align}
\mathrm{preModule}^\mathrm{solid,quasicoherent,mixed-parityalmostdeRham}_{\square,\Gamma^\mathrm{perfect}_{\text{Robba},X,\text{pro\'et},\infty}\{t^{1/2}\}},\\
\mathrm{preModule}^\mathrm{solid,quasicoherent,mixed-parityalmostdeRham}_{\square,\Gamma^\mathrm{perfect}_{\text{Robba},X,\text{pro\'et},I}\{t^{1/2}\}} 
\end{align}
to $(\infty,1)$-categories in image denoted by:
\begin{align}
\mathrm{preModule}_{X,\text{\'et}}
\end{align}
to be the following functors sending each $F$ in the domain to:
\begin{align}
&\text{RH}_\text{mixed-parity}(F):=f_*(F\otimes_{\Gamma^\mathrm{perfect}_{\text{Robba},X,\text{pro\'et},\infty}\{t^{1/2}\}} \Gamma^\mathcal{O}_{\text{deRham},X,\text{pro\'et}}\{t^{1/2}\}),\\
&\text{RH}_\text{mixed-parity}(F):=f_*(F\otimes_{\Gamma^\mathrm{perfect}_{\text{Robba},X,\text{pro\'et},I}\{t^{1/2}\}} \Gamma^\mathcal{O}_{\text{deRham},X,\text{pro\'et}}\{t^{1/2}\}),\\
&\text{RH}_\text{mixed-parity}(F):=f_*(F\otimes_{\Gamma^\mathrm{perfect}_{\text{Robba},X,\text{pro\'et},\infty}\{t^{1/2}\}} \Gamma^\mathcal{O}_{\text{deRham},X,\text{pro\'et}}\{t^{1/2},\log(t)\}),\\
&\text{RH}_\text{mixed-parity}(F):=f_*(F\otimes_{\Gamma^\mathrm{perfect}_{\text{Robba},X,\text{pro\'et},I}\{t^{1/2}\}} \Gamma^\mathcal{O}_{\text{deRham},X,\text{pro\'et}}\{t^{1/2},\log(t)\}),\\
\end{align}
respectively.

\end{definition}

\begin{definition}
In the situation where we have the Frobenius action we consider the follwing. We define the following Riemann-Hilbert functor $\text{RH}_\text{mixed-parity}$ from the one of categories:
\begin{align}
\varphi\mathrm{preModule}^\mathrm{solid,quasicoherent,mixed-paritydeRham}_{\square,\Gamma^\mathrm{perfect}_{\text{Robba},X,\text{pro\'et},\infty}\{t^{1/2}\}},
\varphi\mathrm{preModule}^\mathrm{solid,quasicoherent,mixed-paritydeRham}_{\square,\Gamma^\mathrm{perfect}_{\text{Robba},X,\text{pro\'et},I}\{t^{1/2}\}}, 
\end{align}
and
\begin{align}
\varphi\mathrm{preModule}^\mathrm{solid,quasicoherent,mixed-parityalmostdeRham}_{\square,\Gamma^\mathrm{perfect}_{\text{Robba},X,\text{pro\'et},\infty}\{t^{1/2}\}},\\
\varphi\mathrm{preModule}^\mathrm{solid,quasicoherent,mixed-parityalmostdeRham}_{\square,\Gamma^\mathrm{perfect}_{\text{Robba},X,\text{pro\'et},I}\{t^{1/2}\}} 
\end{align}
to $(\infty,1)$-categories in image denoted by:
\begin{align}
\mathrm{preModule}_{X,\text{\'et}}
\end{align}
to be the following functors sending each $F$ in the domain to:
\begin{align}
&\text{RH}_\text{mixed-parity}(F):=f_*(F\otimes_{\Gamma^\mathrm{perfect}_{\text{Robba},X,\text{pro\'et},\infty}\{t^{1/2}\}} \Gamma^\mathcal{O}_{\text{deRham},X,\text{pro\'et}}\{t^{1/2}\}),\\
&\text{RH}_\text{mixed-parity}(F):=f_*(F\otimes_{\Gamma^\mathrm{perfect}_{\text{Robba},X,\text{pro\'et},I}\{t^{1/2}\}} \Gamma^\mathcal{O}_{\text{deRham},X,\text{pro\'et}}\{t^{1/2}\}),\\
&\text{RH}_\text{mixed-parity}(F):=f_*(F\otimes_{\Gamma^\mathrm{perfect}_{\text{Robba},X,\text{pro\'et},\infty}\{t^{1/2}\}} \Gamma^\mathcal{O}_{\text{deRham},X,\text{pro\'et}}\{t^{1/2},\log(t)\}),\\
&\text{RH}_\text{mixed-parity}(F):=f_*(F\otimes_{\Gamma^\mathrm{perfect}_{\text{Robba},X,\text{pro\'et},I}\{t^{1/2}\}} \Gamma^\mathcal{O}_{\text{deRham},X,\text{pro\'et}}\{t^{1/2},\log(t)\}),\\
\end{align}
respectively.

\end{definition}

\newpage
\section{Geometric Family of Mixed-Parity Hodge Modules II: Cristalline Situations}

\subsection{Period Rings and Sheaves}

\begin{reference}
\cite{Sch1}, \cite{KL1}, \cite{KL2}, \cite{BL1}, \cite{BL2}, \cite{BS}, \cite{BHS}, \cite{Fon1}, \cite{CS1}, \cite{CS2}, \cite{BK}, \cite{BBK}, \cite{BBBK}, \cite{KKM}, \cite{KM}, \cite{LZ}, \cite{TT}, \cite{M}.
\end{reference}

\subsubsection{Rings}

\noindent Let $X$ be a rigid analytic space over $\mathbb{Q}_p$. We have the corresponding \'etale site and the corresponding pro-\'etale site of $X$, which we denote them by $X_{\text{pro\'et}},X_\text{\'et}$. The relationship of the two sites can be reflected by the corresponding morphism $f:X_{\text{pro\'et}}\longrightarrow X_\text{\'et}$. Then we have the corresponding cristalline period rings and sheaves from \cite{TT}:
\begin{align}
\Gamma_{\text{cristalline},X,\text{pro\'et}}, \Gamma^\mathcal{O}_{\text{cristalline},X,\text{pro\'et}}.
\end{align}
Our notations are different from \cite{TT}, we use $\Gamma$ to mean $B$ in \cite{TT}, while $\Gamma^\mathcal{O}$ will be the corresponding $OB$ ring in \cite{TT}.\\

\begin{definition}
\indent Now we assume that $p>2$, following \cite{BS} we join the square root of $t$ element in $\Gamma_{\text{cristalline},X,\text{pro\'et}}$ which forms the sheaves:
\begin{align}
\Gamma_{\text{cristalline},X,\text{pro\'et}}\{t^{1/2}\},\Gamma^\mathcal{O}_{\text{cristalline},X,\text{pro\'et}}\{t^{1/2}\}.
\end{align}
And following \cite{BL1}, \cite{BL2}, \cite{Fon1}, \cite{BHS} we further have the following sheaves of rings:
\begin{align}
\Gamma_{\text{cristalline},X,\text{pro\'et}}\{t^{1/2},\log(t)\},\Gamma^\mathcal{O}_{\text{cristalline},X,\text{pro\'et}}\{t^{1/2},\log(t)\}.
\end{align}
\end{definition}

\begin{definition}
We use the notations:
\begin{align}
\Gamma^\mathrm{perfect}_{\text{Robba},X,\text{pro\'et}},\Gamma^\mathrm{perfect}_{\text{Robba},X,\text{pro\'et},\infty},\Gamma^\mathrm{perfect}_{\text{Robba},X,\text{pro\'et},I}
\end{align}
to denote the perfect Robba rings from \cite{KL1}, \cite{KL2}, where $I\subset (0,\infty)$. Then we join $t^{1/2}$ to these sheaves we have:
\begin{align}
\Gamma^\mathrm{perfect}_{\text{Robba},X,\text{pro\'et}}\{t^{1/2}\},\Gamma^\mathrm{perfect}_{\text{Robba},X,\text{pro\'et},\infty}\{t^{1/2}\},\Gamma^\mathrm{perfect}_{\text{Robba},X,\text{pro\'et},I}\{t^{1/2}\}.
\end{align}
And following \cite{BL1}, \cite{BL2}, \cite{Fon1}, \cite{BHS} we have the following larger sheaves:
\begin{align}
\Gamma^\mathrm{perfect}_{\text{Robba},X,\text{pro\'et}}\{t^{1/2},\log(t)\},\Gamma^\mathrm{perfect}_{\text{Robba},X,\text{pro\'et},\infty}\{t^{1/2},\log(t)\},\Gamma^\mathrm{perfect}_{\text{Robba},X,\text{pro\'et},I}\{t^{1/2},\log(t)\}.
\end{align} 
\end{definition}

\begin{definition}
From now on, we use the same notation to denote the period rings involved tensored with a finite extension of $\mathbb{Q}_p$ containing square root of $p$ as in \cite{BS}.
\begin{align}
\Gamma_{\text{cristalline},X,\text{pro\'et}}\{t^{1/2}\},\Gamma^\mathcal{O}_{\text{cristalline},X,\text{pro\'et}}\{t^{1/2}\}.
\end{align}
\begin{align}
\Gamma_{\text{cristalline},X,\text{pro\'et}}\{t^{1/2},\log(t)\},\Gamma^\mathcal{O}_{\text{cristalline},X,\text{pro\'et}}\{t^{1/2},\log(t)\}.
\end{align}
\begin{align}
\Gamma^\mathrm{perfect}_{\text{Robba},X,\text{pro\'et}}\{t^{1/2}\},\Gamma^\mathrm{perfect}_{\text{Robba},X,\text{pro\'et},\infty}\{t^{1/2}\},\Gamma^\mathrm{perfect}_{\text{Robba},X,\text{pro\'et},I}\{t^{1/2}\}.
\end{align}
\begin{align}
\Gamma^\mathrm{perfect}_{\text{Robba},X,\text{pro\'et}}\{t^{1/2},\log(t)\},\Gamma^\mathrm{perfect}_{\text{Robba},X,\text{pro\'et},\infty}\{t^{1/2},\log(t)\},\Gamma^\mathrm{perfect}_{\text{Robba},X,\text{pro\'et},I}\{t^{1/2},\log(t)\}.
\end{align}
This is necessary since we to extend the action of $\varphi$ to the period rings by $\varphi(t^{1/2}\otimes 1)=\varphi(t)^{1/2}\otimes 1$.
\end{definition}

\subsubsection{Modules}

\noindent We consider quasicoherent presheaves in the following two situation:
\begin{itemize}
\item[$\square$] The solid quasicoherent modules from \cite{CS1}, \cite{CS2};
\item[$\square$] The ind-Banach quasicoherent modules from \cite{BK}, \cite{BBK}, \cite{BBBK}, \cite{KKM}, \cite{KM} with the corresponding monomorphic ind-Banach quasicoherent modules from \cite{BK}, \cite{BBK}, \cite{BBBK}, \cite{KKM}, \cite{KM}.
\end{itemize}

\begin{definition}
We use the notation:
\begin{align}
\mathrm{preModule}^\mathrm{solid,quasicoherent}_{\square,\Gamma^\mathrm{perfect}_{\text{Robba},X,\text{pro\'et}}\{t^{1/2}\}},\mathrm{preModule}^\mathrm{solid,quasicoherent}_{\square,\Gamma^\mathrm{perfect}_{\text{Robba},X,\text{pro\'et},\infty}\{t^{1/2}\}},
\mathrm{preModule}^\mathrm{solid,quasicoherent}_{\square,\Gamma^\mathrm{perfect}_{\text{Robba},X,\text{pro\'et},I}\{t^{1/2}\}} 
\end{align}
to denote the $(\infty,1)$-categories of solid quasicoherent presheaves over the corresonding Robba sheaves. Locally the section is defined by taking the corresponding $(\infty,1)$-categories of solid modules.
\end{definition}

\begin{definition}
We use the notation:
\begin{align}
\mathrm{preModule}^\mathrm{ind-Banach,quasicoherent}_{\Gamma^\mathrm{perfect}_{\text{Robba},X,\text{pro\'et}}\{t^{1/2}\}},\\
\mathrm{preModule}^\mathrm{ind-Banach,quasicoherent}_{\Gamma^\mathrm{perfect}_{\text{Robba},X,\text{pro\'et},\infty}\{t^{1/2}\}},\\
\mathrm{preModule}^\mathrm{ind-Banach,quasicoherent}_{\Gamma^\mathrm{perfect}_{\text{Robba},X,\text{pro\'et},I}\{t^{1/2}\}} 
\end{align}
to denote the $(\infty,1)$-categories of solid quasicoherent presheaves over the corresonding Robba sheaves. Locally the section is defined by taking the corresponding $(\infty,1)$-categories of inductive Banach  modules. 
\end{definition}

\begin{definition}
We use the notation:
\begin{align}
\mathrm{Module}^\mathrm{solid,quasicoherent}_{\square,\Gamma^\mathrm{perfect}_{\text{Robba},X,\text{pro\'et}}\{t^{1/2}\}},\mathrm{Module}^\mathrm{solid,quasicoherent}_{\square,\Gamma^\mathrm{perfect}_{\text{Robba},X,\text{pro\'et},\infty}\{t^{1/2}\}},
\mathrm{Module}^\mathrm{solid,quasicoherent}_{\square,\Gamma^\mathrm{perfect}_{\text{Robba},X,\text{pro\'et},I}\{t^{1/2}\}} 
\end{align}
to denote the $(\infty,1)$-categories of solid quasicoherent sheaves over the corresonding Robba sheaves. Locally the section is defined by taking the corresponding $(\infty,1)$-categories of solid modules.
\end{definition}

\subsubsection{Mixed-Parity Hodge Modules without Frobenius}

\noindent Now we consider the key objects in our study namely those complexes generated by certain mixed-parity Hodge modules. We start from the following definition.

\begin{definition}
For any locally free coherent sheaf $F$ over
\begin{align}
\Gamma^\mathrm{perfect}_{\text{Robba},X,\text{pro\'et},\infty}\{t^{1/2}\},\Gamma^\mathrm{perfect}_{\text{Robba},X,\text{pro\'et},I}\{t^{1/2}\},
\end{align} 
we consider the following functor $\mathrm{cristalline}$ sending $F$ to the following object:
\begin{align}
f_*(F\otimes_{\Gamma^\mathrm{perfect}_{\text{Robba},X,\text{pro\'et},\infty}\{t^{1/2}\}} \Gamma^\mathcal{O}_{\text{cristalline},X,\text{pro\'et}}\{t^{1/2}\})
\end{align}
or 
\begin{align}
f_*(F\otimes_{\Gamma^\mathrm{perfect}_{\text{Robba},X,\text{pro\'et},I}\{t^{1/2}\}} \Gamma^\mathcal{O}_{\text{cristalline},X,\text{pro\'et}}\{t^{1/2}\}).
\end{align}
We call $F$ mixed-parity cristalline if we have the following isomorphism:
\begin{align}
f^*f_*(F\otimes_{\Gamma^\mathrm{perfect}_{\text{Robba},X,\text{pro\'et},\infty}\{t^{1/2}\}} \Gamma^\mathcal{O}_{\text{cristalline},X,\text{pro\'et}}\{t^{1/2}\}) \otimes \Gamma^\mathcal{O}_{\text{cristalline},X,\text{pro\'et}}\{t^{1/2}\}\\
 \overset{\sim}{\longrightarrow} F \otimes \Gamma^\mathcal{O}_{\text{cristalline},X,\text{pro\'et}}\{t^{1/2}\} 
\end{align}
or 
\begin{align}
f^*f_*(F\otimes_{\Gamma^\mathrm{perfect}_{\text{Robba},X,\text{pro\'et},I}\{t^{1/2}\}} \Gamma^\mathcal{O}_{\text{cristalline},X,\text{pro\'et}}\{t^{1/2}\}) \otimes \Gamma^\mathcal{O}_{\text{cristalline},X,\text{pro\'et}}\{t^{1/2}\}\\ \overset{\sim}{\longrightarrow} F \otimes \Gamma^\mathcal{O}_{\text{cristalline},X,\text{pro\'et}}\{t^{1/2}\}. 
\end{align}
\end{definition}

\begin{definition}
For any locally free coherent sheaf $F$ over
\begin{align}
\Gamma^\mathrm{perfect}_{\text{Robba},X,\text{pro\'et},\infty}\{t^{1/2}\},\Gamma^\mathrm{perfect}_{\text{Robba},X,\text{pro\'et},I}\{t^{1/2}\},
\end{align} 
we consider the following functor $\mathrm{cristalline}^\mathrm{almost}$ sending $F$ to the following object:
\begin{align}
f_*(F\otimes_{\Gamma^\mathrm{perfect}_{\text{Robba},X,\text{pro\'et},\infty}\{t^{1/2}\}} \Gamma^\mathcal{O}_{\text{cristalline},X,\text{pro\'et}}\{t^{1/2},\log(t)\})
\end{align}
or 
\begin{align}
f_*(F\otimes_{\Gamma^\mathrm{perfect}_{\text{Robba},X,\text{pro\'et},I}\{t^{1/2}\}} \Gamma^\mathcal{O}_{\text{cristalline},X,\text{pro\'et}}\{t^{1/2},\log(t)\}).
\end{align}
We call $F$ mixed-parity almost cristalline if we have the following isomorphism:
\begin{align}
f^*f_*(F\otimes_{\Gamma^\mathrm{perfect}_{\text{Robba},X,\text{pro\'et},\infty}\{t^{1/2}\}} \Gamma^\mathcal{O}_{\text{cristalline},X,\text{pro\'et}}\{t^{1/2},\log(t)\}) \otimes \Gamma^\mathcal{O}_{\text{cristalline},X,\text{pro\'et}}\{t^{1/2},\log(t)\} \\
\overset{\sim}{\longrightarrow} F \otimes \Gamma^\mathcal{O}_{\text{cristalline},X,\text{pro\'et}}\{t^{1/2},\log(t)\} 
\end{align}
or 
\begin{align}
f^*f_*(F\otimes_{\Gamma^\mathrm{perfect}_{\text{Robba},X,\text{pro\'et},I}\{t^{1/2}\}} \Gamma^\mathcal{O}_{\text{cristalline},X,\text{pro\'et}}\{t^{1/2},\log(t)\}) \otimes \Gamma^\mathcal{O}_{\text{cristalline},X,\text{pro\'et}}\{t^{1/2},\log(t)\}\\ \overset{\sim}{\longrightarrow} F \otimes \Gamma^\mathcal{O}_{\text{cristalline},X,\text{pro\'et}}\{t^{1/2},\log(t)\}. 
\end{align}
\end{definition}

\noindent We now define the $(\infty,1)$-categories of mixed-parity cristalline modules and he corresponding mixed-parity almost cristalline modules by using the objects involved to generated these categories:

\begin{definition}
Considering all the mixed parity cristalline bundles (locally finite free) as defined above, we consider the sub-$(\infty,1)$ category of 
\begin{align}
\mathrm{preModule}^\mathrm{solid,quasicoherent}_{\square,\Gamma^\mathrm{perfect}_{\text{Robba},X,\text{pro\'et},\infty}\{t^{1/2}\}},
\mathrm{preModule}^\mathrm{solid,quasicoherent}_{\square,\Gamma^\mathrm{perfect}_{\text{Robba},X,\text{pro\'et},I}\{t^{1/2}\}} 
\end{align}
generated by the mixed-parity cristalline bundles (locally finite free ones). These are defined to be the $(\infty,1)$-categories of mixed-parity cristalline complexes:
\begin{align}
\mathrm{preModule}^\mathrm{solid,quasicoherent,mixed-paritycristalline}_{\square,\Gamma^\mathrm{perfect}_{\text{Robba},X,\text{pro\'et},\infty}\{t^{1/2}\}},
\mathrm{preModule}^\mathrm{solid,quasicoherent,mixed-paritycristalline}_{\square,\Gamma^\mathrm{perfect}_{\text{Robba},X,\text{pro\'et},I}\{t^{1/2}\}}. 
\end{align}
\end{definition}

\begin{definition}
Considering all the mixed parity almost cristalline bundles (locally finite free) as defined above, we consider the sub-$(\infty,1)$ category of 
\begin{align}
\mathrm{preModule}^\mathrm{solid,quasicoherent}_{\square,\Gamma^\mathrm{perfect}_{\text{Robba},X,\text{pro\'et},\infty}\{t^{1/2}\}},
\mathrm{preModule}^\mathrm{solid,quasicoherent}_{\square,\Gamma^\mathrm{perfect}_{\text{Robba},X,\text{pro\'et},I}\{t^{1/2}\}} 
\end{align}
generated by the mixed-parity almost cristalline bundles (locally finite free ones). These are defined to be the $(\infty,1)$-categories of mixed-parity cristalline complexes:
\begin{align}
\mathrm{preModule}^\mathrm{solid,quasicoherent,mixed-parityalmostcristalline}_{\square,\Gamma^\mathrm{perfect}_{\text{Robba},X,\text{pro\'et},\infty}\{t^{1/2}\}},\\
\mathrm{preModule}^\mathrm{solid,quasicoherent,mixed-parityalmostcristalline}_{\square,\Gamma^\mathrm{perfect}_{\text{Robba},X,\text{pro\'et},I}\{t^{1/2}\}}. 
\end{align}
\end{definition}

\indent Then the corresponding mixed-parity cristalline functors can be extended to these categories:
\begin{align}
\mathrm{preModule}^\mathrm{solid,quasicoherent,mixed-paritycristalline}_{\square,\Gamma^\mathrm{perfect}_{\text{Robba},X,\text{pro\'et},\infty}\{t^{1/2}\}},\\
\mathrm{preModule}^\mathrm{solid,quasicoherent,mixed-paritycristalline}_{\square,\Gamma^\mathrm{perfect}_{\text{Robba},X,\text{pro\'et},I}\{t^{1/2}\}}, 
\end{align}
and
\begin{align}
\mathrm{preModule}^\mathrm{solid,quasicoherent,mixed-parityalmostcristalline}_{\square,\Gamma^\mathrm{perfect}_{\text{Robba},X,\text{pro\'et},\infty}\{t^{1/2}\}},\\
\mathrm{preModule}^\mathrm{solid,quasicoherent,mixed-parityalmostcristalline}_{\square,\Gamma^\mathrm{perfect}_{\text{Robba},X,\text{pro\'et},I}\{t^{1/2}\}}. 
\end{align}

\subsubsection{Mixed-Parity Hodge Modules with Frobenius}

\noindent Now we consider the key objects in our study namely those complexes generated by certain mixed-parity Hodge modules. We start from the following definition.

\begin{remark}
All the coherent sheaves over mixed-parity Robba sheaves in this section will carry the corresponding Frobenius morphism $\varphi: F \overset{\sim}{\longrightarrow} \varphi^*F$.
\end{remark}

\begin{definition}
For any locally free coherent sheaf $F$ over
\begin{align}
\Gamma^\mathrm{perfect}_{\text{Robba},X,\text{pro\'et},\infty}\{t^{1/2}\},\Gamma^\mathrm{perfect}_{\text{Robba},X,\text{pro\'et},I}\{t^{1/2}\},
\end{align} 
we consider the following functor $\mathrm{cristalline}$ sending $F$ to the following object:
\begin{align}
f_*(F\otimes_{\Gamma^\mathrm{perfect}_{\text{Robba},X,\text{pro\'et},\infty}\{t^{1/2}\}} \Gamma^\mathcal{O}_{\text{cristalline},X,\text{pro\'et}}\{t^{1/2}\})
\end{align}
or 
\begin{align}
f_*(F\otimes_{\Gamma^\mathrm{perfect}_{\text{Robba},X,\text{pro\'et},I}\{t^{1/2}\}} \Gamma^\mathcal{O}_{\text{cristalline},X,\text{pro\'et}}\{t^{1/2}\}).
\end{align}
We call $F$ mixed-parity cristalline if we have the following isomorphism:
\begin{align}
f^*f_*(F\otimes_{\Gamma^\mathrm{perfect}_{\text{Robba},X,\text{pro\'et},\infty}\{t^{1/2}\}} \Gamma^\mathcal{O}_{\text{cristalline},X,\text{pro\'et}}\{t^{1/2}\}) \otimes \Gamma^\mathcal{O}_{\text{cristalline},X,\text{pro\'et}}\{t^{1/2}\}\\ \overset{\sim}{\longrightarrow} F \otimes \Gamma^\mathcal{O}_{\text{cristalline},X,\text{pro\'et}}\{t^{1/2}\} 
\end{align}
or 
\begin{align}
f^*f_*(F\otimes_{\Gamma^\mathrm{perfect}_{\text{Robba},X,\text{pro\'et},I}\{t^{1/2}\}} \Gamma^\mathcal{O}_{\text{cristalline},X,\text{pro\'et}}\{t^{1/2}\}) \otimes \Gamma^\mathcal{O}_{\text{cristalline},X,\text{pro\'et}}\{t^{1/2}\} \\\overset{\sim}{\longrightarrow} F \otimes \Gamma^\mathcal{O}_{\text{cristalline},X,\text{pro\'et}}\{t^{1/2}\}. 
\end{align}
\end{definition}

\begin{definition}
For any locally free coherent sheaf $F$ over
\begin{align}
\Gamma^\mathrm{perfect}_{\text{Robba},X,\text{pro\'et},\infty}\{t^{1/2}\},\Gamma^\mathrm{perfect}_{\text{Robba},X,\text{pro\'et},I}\{t^{1/2}\},
\end{align} 
we consider the following functor $\mathrm{cristalline}^\mathrm{almost}$ sending $F$ to the following object:
\begin{align}
f_*(F\otimes_{\Gamma^\mathrm{perfect}_{\text{Robba},X,\text{pro\'et},\infty}\{t^{1/2}\}} \Gamma^\mathcal{O}_{\text{cristalline},X,\text{pro\'et}}\{t^{1/2},\log(t)\})
\end{align}
or 
\begin{align}
f_*(F\otimes_{\Gamma^\mathrm{perfect}_{\text{Robba},X,\text{pro\'et},I}\{t^{1/2}\}} \Gamma^\mathcal{O}_{\text{cristalline},X,\text{pro\'et}}\{t^{1/2},\log(t)\}).
\end{align}
We call $F$ mixed-parity almost cristalline if we have the following isomorphism:
\begin{align}
f^*f_*(F\otimes_{\Gamma^\mathrm{perfect}_{\text{Robba},X,\text{pro\'et},\infty}\{t^{1/2}\}} \Gamma^\mathcal{O}_{\text{cristalline},X,\text{pro\'et}}\{t^{1/2},\log(t)\}) \otimes \Gamma^\mathcal{O}_{\text{cristalline},X,\text{pro\'et}}\{t^{1/2},\log(t)\} \\
\overset{\sim}{\longrightarrow} F \otimes \Gamma^\mathcal{O}_{\text{cristalline},X,\text{pro\'et}}\{t^{1/2},\log(t)\} 
\end{align}
or 
\begin{align}
f^*f_*(F\otimes_{\Gamma^\mathrm{perfect}_{\text{Robba},X,\text{pro\'et},I}\{t^{1/2}\}} \Gamma^\mathcal{O}_{\text{cristalline},X,\text{pro\'et}}\{t^{1/2},\log(t)\}) \otimes \Gamma^\mathcal{O}_{\text{cristalline},X,\text{pro\'et}}\{t^{1/2},\log(t)\}\\ \overset{\sim}{\longrightarrow} F \otimes \Gamma^\mathcal{O}_{\text{cristalline},X,\text{pro\'et}}\{t^{1/2},\log(t)\}. 
\end{align}
\end{definition}

\noindent We now define the $(\infty,1)$-categories of mixed-parity cristalline modules and he corresponding mixed-parity almost cristalline modules by using the objects involved to generated these categories:

\begin{definition}
Considering all the mixed parity cristalline bundles (locally finite free) as defined above, we consider the sub-$(\infty,1)$ category of 
\begin{align}
\varphi\mathrm{preModule}^\mathrm{solid,quasicoherent}_{\square,\Gamma^\mathrm{perfect}_{\text{Robba},X,\text{pro\'et},\infty}\{t^{1/2}\}},
\varphi\mathrm{preModule}^\mathrm{solid,quasicoherent}_{\square,\Gamma^\mathrm{perfect}_{\text{Robba},X,\text{pro\'et},I}\{t^{1/2}\}} 
\end{align}
generated by the mixed-parity cristalline bundles (locally finite free ones). These are defined to be the $(\infty,1)$-categories of mixed-parity cristalline complexes:
\begin{align}
\varphi\mathrm{preModule}^\mathrm{solid,quasicoherent,mixed-paritycristalline}_{\square,\Gamma^\mathrm{perfect}_{\text{Robba},X,\text{pro\'et},\infty}\{t^{1/2}\}},
\varphi\mathrm{preModule}^\mathrm{solid,quasicoherent,mixed-paritycristalline}_{\square,\Gamma^\mathrm{perfect}_{\text{Robba},X,\text{pro\'et},I}\{t^{1/2}\}}. 
\end{align}
\end{definition}

\begin{definition}
Considering all the mixed parity almost cristalline bundles (locally finite free) as defined above, we consider the sub-$(\infty,1)$ category of 
\begin{align}
\varphi\mathrm{preModule}^\mathrm{solid,quasicoherent}_{\square,\Gamma^\mathrm{perfect}_{\text{Robba},X,\text{pro\'et},\infty}\{t^{1/2}\}},
\varphi\mathrm{preModule}^\mathrm{solid,quasicoherent}_{\square,\Gamma^\mathrm{perfect}_{\text{Robba},X,\text{pro\'et},I}\{t^{1/2}\}} 
\end{align}
generated by the mixed-parity almost cristalline bundles (locally finite free ones). These are defined to be the $(\infty,1)$-categories of mixed-parity cristalline complexes:
\begin{align}
\varphi\mathrm{preModule}^\mathrm{solid,quasicoherent,mixed-parityalmostcristalline}_{\square,\Gamma^\mathrm{perfect}_{\text{Robba},X,\text{pro\'et},\infty}\{t^{1/2}\}},\\
\varphi\mathrm{preModule}^\mathrm{solid,quasicoherent,mixed-parityalmostcristalline}_{\square,\Gamma^\mathrm{perfect}_{\text{Robba},X,\text{pro\'et},I}\{t^{1/2}\}}. 
\end{align}
\end{definition}

\indent Then the corresponding mixed-parity cristalline functors can be extended to these categories:
\begin{align}
\varphi\mathrm{preModule}^\mathrm{solid,quasicoherent,mixed-paritycristalline}_{\square,\Gamma^\mathrm{perfect}_{\text{Robba},X,\text{pro\'et},\infty}\{t^{1/2}\}},\\
\varphi\mathrm{preModule}^\mathrm{solid,quasicoherent,mixed-paritycristalline}_{\square,\Gamma^\mathrm{perfect}_{\text{Robba},X,\text{pro\'et},I}\{t^{1/2}\}}, 
\end{align}
and
\begin{align}
\varphi\mathrm{preModule}^\mathrm{solid,quasicoherent,mixed-parityalmostcristalline}_{\square,\Gamma^\mathrm{perfect}_{\text{Robba},X,\text{pro\'et},\infty}\{t^{1/2}\}},\\
\varphi\mathrm{preModule}^\mathrm{solid,quasicoherent,mixed-parityalmostcristalline}_{\square,\Gamma^\mathrm{perfect}_{\text{Robba},X,\text{pro\'et},I}\{t^{1/2}\}}. 
\end{align}

\subsection{Mixed-Parity cristalline Riemann-Hilbert Correspondence}

\indent This chapter will extend the corresponding Riemann-Hilbert correspondence from \cite{Sch1}, \cite{LZ}, \cite{BL1}, \cite{BL2}, \cite{M} to the mixed-parity setting.

\begin{definition}
We define the following Riemann-Hilbert functor $\text{RH}_\text{mixed-parity}$ from the one of categories:
\begin{align}
\mathrm{preModule}^\mathrm{solid,quasicoherent,mixed-paritycristalline}_{\square,\Gamma^\mathrm{perfect}_{\text{Robba},X,\text{pro\'et},\infty}\{t^{1/2}\}},
\mathrm{preModule}^\mathrm{solid,quasicoherent,mixed-paritycristalline}_{\square,\Gamma^\mathrm{perfect}_{\text{Robba},X,\text{pro\'et},I}\{t^{1/2}\}}, 
\end{align}
and
\begin{align}
\mathrm{preModule}^\mathrm{solid,quasicoherent,mixed-parityalmostcristalline}_{\square,\Gamma^\mathrm{perfect}_{\text{Robba},X,\text{pro\'et},\infty}\{t^{1/2}\}},\\
\mathrm{preModule}^\mathrm{solid,quasicoherent,mixed-parityalmostcristalline}_{\square,\Gamma^\mathrm{perfect}_{\text{Robba},X,\text{pro\'et},I}\{t^{1/2}\}} 
\end{align}
to $(\infty,1)$-categories in image denoted by:
\begin{align}
\mathrm{preModule}_{X,\text{\'et}}
\end{align}
to be the following functors sending each $F$ in the domain to:
\begin{align}
&\text{RH}_\text{mixed-parity}(F):=f_*(F\otimes_{\Gamma^\mathrm{perfect}_{\text{Robba},X,\text{pro\'et},\infty}\{t^{1/2}\}} \Gamma^\mathcal{O}_{\text{cristalline},X,\text{pro\'et}}\{t^{1/2}\}),\\
&\text{RH}_\text{mixed-parity}(F):=f_*(F\otimes_{\Gamma^\mathrm{perfect}_{\text{Robba},X,\text{pro\'et},I}\{t^{1/2}\}} \Gamma^\mathcal{O}_{\text{cristalline},X,\text{pro\'et}}\{t^{1/2}\}),\\
&\text{RH}_\text{mixed-parity}(F):=f_*(F\otimes_{\Gamma^\mathrm{perfect}_{\text{Robba},X,\text{pro\'et},\infty}\{t^{1/2}\}} \Gamma^\mathcal{O}_{\text{cristalline},X,\text{pro\'et}}\{t^{1/2},\log(t)\}),\\
&\text{RH}_\text{mixed-parity}(F):=f_*(F\otimes_{\Gamma^\mathrm{perfect}_{\text{Robba},X,\text{pro\'et},I}\{t^{1/2}\}} \Gamma^\mathcal{O}_{\text{cristalline},X,\text{pro\'et}}\{t^{1/2},\log(t)\}),\\
\end{align}
respectively.

\end{definition}

\begin{definition}
In the situation where we have the Frobenius action we consider the follwing. We define the following Riemann-Hilbert functor $\text{RH}_\text{mixed-parity}$ from the one of categories:
\begin{align}
\varphi\mathrm{preModule}^\mathrm{solid,quasicoherent,mixed-paritycristalline}_{\square,\Gamma^\mathrm{perfect}_{\text{Robba},X,\text{pro\'et},\infty}\{t^{1/2}\}},
\varphi\mathrm{preModule}^\mathrm{solid,quasicoherent,mixed-paritycristalline}_{\square,\Gamma^\mathrm{perfect}_{\text{Robba},X,\text{pro\'et},I}\{t^{1/2}\}}, 
\end{align}
and
\begin{align}
\varphi\mathrm{preModule}^\mathrm{solid,quasicoherent,mixed-parityalmostcristalline}_{\square,\Gamma^\mathrm{perfect}_{\text{Robba},X,\text{pro\'et},\infty}\{t^{1/2}\}},\\
\varphi\mathrm{preModule}^\mathrm{solid,quasicoherent,mixed-parityalmostcristalline}_{\square,\Gamma^\mathrm{perfect}_{\text{Robba},X,\text{pro\'et},I}\{t^{1/2}\}} 
\end{align}
to $(\infty,1)$-categories in image denoted by:
\begin{align}
\mathrm{preModule}_{X,\text{\'et}}
\end{align}
to be the following functors sending each $F$ in the domain to:
\begin{align}
&\text{RH}_\text{mixed-parity}(F):=f_*(F\otimes_{\Gamma^\mathrm{perfect}_{\text{Robba},X,\text{pro\'et},\infty}\{t^{1/2}\}} \Gamma^\mathcal{O}_{\text{cristalline},X,\text{pro\'et}}\{t^{1/2}\}),\\
&\text{RH}_\text{mixed-parity}(F):=f_*(F\otimes_{\Gamma^\mathrm{perfect}_{\text{Robba},X,\text{pro\'et},I}\{t^{1/2}\}} \Gamma^\mathcal{O}_{\text{cristalline},X,\text{pro\'et}}\{t^{1/2}\}),\\
&\text{RH}_\text{mixed-parity}(F):=f_*(F\otimes_{\Gamma^\mathrm{perfect}_{\text{Robba},X,\text{pro\'et},\infty}\{t^{1/2}\}} \Gamma^\mathcal{O}_{\text{cristalline},X,\text{pro\'et}}\{t^{1/2},\log(t)\}),\\
&\text{RH}_\text{mixed-parity}(F):=f_*(F\otimes_{\Gamma^\mathrm{perfect}_{\text{Robba},X,\text{pro\'et},I}\{t^{1/2}\}} \Gamma^\mathcal{O}_{\text{cristalline},X,\text{pro\'et}}\{t^{1/2},\log(t)\}),\\
\end{align}
respectively.

\end{definition}

\newpage
\section{Geometric Family of Mixed-Parity Hodge Modules III: Semi-Stable Situations}

\begin{reference}
\cite{Sch1}, \cite{KL1}, \cite{KL2}, \cite{BL1}, \cite{BL2}, \cite{BS}, \cite{BHS}, \cite{Fon1}, \cite{CS1}, \cite{CS2}, \cite{BK}, \cite{BBK}, \cite{BBBK}, \cite{KKM}, \cite{KM}, \cite{LZ}, \cite{Shi}, \cite{M}.
\end{reference}

\subsection{Period Rings and Sheaves}

\subsubsection{Rings}

\noindent Let $X$ be a rigid analytic space over $\mathbb{Q}_p$. We have the corresponding \'etale site and the corresponding pro-\'etale site of $X$, which we denote them by $X_{\text{pro\'et}},X_\text{\'et}$. The relationship of the two sites can be reflected by the corresponding morphism $f:X_{\text{pro\'et}}\longrightarrow X_\text{\'et}$. Then we have the corresponding semi-stable period rings and sheaves from \cite{Shi}:
\begin{align}
\Gamma_{\text{semistable},X,\text{pro\'et}}, \Gamma^\mathcal{O}_{\text{semistable},X,\text{pro\'et}}.
\end{align}
Our notations are different from \cite{Shi}, we use $\Gamma$ to mean $B$ in \cite{Shi}, while $\Gamma^\mathcal{O}$ will be the corresponding $OB$ ring in \cite{Shi}.\\

\begin{definition}
\indent Now we assume that $p>2$, following \cite{BS} we join the square root of $t$ element in $\Gamma_{\text{semistable},X,\text{pro\'et}}$ which forms the sheaves:
\begin{align}
\Gamma_{\text{semistable},X,\text{pro\'et}}\{t^{1/2}\},\Gamma^\mathcal{O}_{\text{semistable},X,\text{pro\'et}}\{t^{1/2}\}.
\end{align}
And following \cite{BL1}, \cite{BL2}, \cite{Fon1}, \cite{BHS} we further have the following sheaves of rings:
\begin{align}
\Gamma_{\text{semistable},X,\text{pro\'et}}\{t^{1/2},\log(t)\},\Gamma^\mathcal{O}_{\text{semistable},X,\text{pro\'et}}\{t^{1/2},\log(t)\}.
\end{align}
\end{definition}

\begin{definition}
We use the notations:
\begin{align}
\Gamma^\mathrm{perfect}_{\text{Robba},X,\text{pro\'et}},\Gamma^\mathrm{perfect}_{\text{Robba},X,\text{pro\'et},\infty},\Gamma^\mathrm{perfect}_{\text{Robba},X,\text{pro\'et},I}
\end{align}
to denote the perfect Robba rings from \cite{KL1}, \cite{KL2}, where $I\subset (0,\infty)$. Then we join $t^{1/2}$ to these sheaves we have:
\begin{align}
\Gamma^\mathrm{perfect}_{\text{Robba},X,\text{pro\'et}}\{t^{1/2}\},\Gamma^\mathrm{perfect}_{\text{Robba},X,\text{pro\'et},\infty}\{t^{1/2}\},\Gamma^\mathrm{perfect}_{\text{Robba},X,\text{pro\'et},I}\{t^{1/2}\}.
\end{align}
And following \cite{BL1}, \cite{BL2}, \cite{Fon1}, \cite{BHS} we have the following larger sheaves:
\begin{align}
\Gamma^\mathrm{perfect}_{\text{Robba},X,\text{pro\'et}}\{t^{1/2},\log(t)\},\Gamma^\mathrm{perfect}_{\text{Robba},X,\text{pro\'et},\infty}\{t^{1/2},\log(t)\},\Gamma^\mathrm{perfect}_{\text{Robba},X,\text{pro\'et},I}\{t^{1/2},\log(t)\}.
\end{align} 
\end{definition}

\begin{definition}
From now on, we use the same notation to denote the period rings involved tensored with a finite extension of $\mathbb{Q}_p$ containing square root of $p$ as in \cite{BS}.
\begin{align}
\Gamma_{\text{semistable},X,\text{pro\'et}}\{t^{1/2}\},\Gamma^\mathcal{O}_{\text{semistable},X,\text{pro\'et}}\{t^{1/2}\}.
\end{align}
\begin{align}
\Gamma_{\text{semistable},X,\text{pro\'et}}\{t^{1/2},\log(t)\},\Gamma^\mathcal{O}_{\text{semistable},X,\text{pro\'et}}\{t^{1/2},\log(t)\}.
\end{align}
\begin{align}
\Gamma^\mathrm{perfect}_{\text{Robba},X,\text{pro\'et}}\{t^{1/2}\},\Gamma^\mathrm{perfect}_{\text{Robba},X,\text{pro\'et},\infty}\{t^{1/2}\},\Gamma^\mathrm{perfect}_{\text{Robba},X,\text{pro\'et},I}\{t^{1/2}\}.
\end{align}
\begin{align}
\Gamma^\mathrm{perfect}_{\text{Robba},X,\text{pro\'et}}\{t^{1/2},\log(t)\},\Gamma^\mathrm{perfect}_{\text{Robba},X,\text{pro\'et},\infty}\{t^{1/2},\log(t)\},\Gamma^\mathrm{perfect}_{\text{Robba},X,\text{pro\'et},I}\{t^{1/2},\log(t)\}.
\end{align}
This is necessary since we to extend the action of $\varphi$ to the period rings by $\varphi(t^{1/2}\otimes 1)=\varphi(t)^{1/2}\otimes 1$.
\end{definition}

\subsubsection{Modules}

\noindent We consider quasicoherent presheaves in the following two situation:
\begin{itemize}
\item[$\square$] The solid quasicoherent modules from \cite{CS1}, \cite{CS2};
\item[$\square$] The ind-Banach quasicoherent modules from \cite{BK}, \cite{BBK}, \cite{BBBK}, \cite{KKM}, \cite{KM} with the corresponding monomorphic ind-Banach quasicoherent modules from \cite{BK}, \cite{BBK}, \cite{BBBK}, \cite{KKM}, \cite{KM}.
\end{itemize}

\begin{definition}
We use the notation:
\begin{align}
\mathrm{preModule}^\mathrm{solid,quasicoherent}_{\square,\Gamma^\mathrm{perfect}_{\text{Robba},X,\text{pro\'et}}\{t^{1/2}\}},\mathrm{preModule}^\mathrm{solid,quasicoherent}_{\square,\Gamma^\mathrm{perfect}_{\text{Robba},X,\text{pro\'et},\infty}\{t^{1/2}\}},
\mathrm{preModule}^\mathrm{solid,quasicoherent}_{\square,\Gamma^\mathrm{perfect}_{\text{Robba},X,\text{pro\'et},I}\{t^{1/2}\}} 
\end{align}
to denote the $(\infty,1)$-categories of solid quasicoherent presheaves over the corresonding Robba sheaves. Locally the section is defined by taking the corresponding $(\infty,1)$-categories of solid modules.
\end{definition}

\begin{definition}
We use the notation:
\begin{align}
\mathrm{preModule}^\mathrm{ind-Banach,quasicoherent}_{\Gamma^\mathrm{perfect}_{\text{Robba},X,\text{pro\'et}}\{t^{1/2}\}},\\
\mathrm{preModule}^\mathrm{ind-Banach,quasicoherent}_{\Gamma^\mathrm{perfect}_{\text{Robba},X,\text{pro\'et},\infty}\{t^{1/2}\}},\\
\mathrm{preModule}^\mathrm{ind-Banach,quasicoherent}_{\Gamma^\mathrm{perfect}_{\text{Robba},X,\text{pro\'et},I}\{t^{1/2}\}} 
\end{align}
to denote the $(\infty,1)$-categories of solid quasicoherent presheaves over the corresonding Robba sheaves. Locally the section is defined by taking the corresponding $(\infty,1)$-categories of inductive Banach  modules. 
\end{definition}

\begin{definition}
We use the notation:
\begin{align}
\mathrm{Module}^\mathrm{solid,quasicoherent}_{\square,\Gamma^\mathrm{perfect}_{\text{Robba},X,\text{pro\'et}}\{t^{1/2}\}},\mathrm{Module}^\mathrm{solid,quasicoherent}_{\square,\Gamma^\mathrm{perfect}_{\text{Robba},X,\text{pro\'et},\infty}\{t^{1/2}\}},
\mathrm{Module}^\mathrm{solid,quasicoherent}_{\square,\Gamma^\mathrm{perfect}_{\text{Robba},X,\text{pro\'et},I}\{t^{1/2}\}} 
\end{align}
to denote the $(\infty,1)$-categories of solid quasicoherent sheaves over the corresonding Robba sheaves. Locally the section is defined by taking the corresponding $(\infty,1)$-categories of solid modules.
\end{definition}

\subsubsection{Mixed-Parity Hodge Modules without Frobenius}

\noindent Now we consider the key objects in our study namely those complexes generated by certain mixed-parity Hodge modules. We start from the following definition.

\begin{definition}
For any locally free coherent sheaf $F$ over
\begin{align}
\Gamma^\mathrm{perfect}_{\text{Robba},X,\text{pro\'et},\infty}\{t^{1/2}\},\Gamma^\mathrm{perfect}_{\text{Robba},X,\text{pro\'et},I}\{t^{1/2}\},
\end{align} 
we consider the following functor $\mathrm{semistable}$ sending $F$ to the following object:
\begin{align}
f_*(F\otimes_{\Gamma^\mathrm{perfect}_{\text{Robba},X,\text{pro\'et},\infty}\{t^{1/2}\}} \Gamma^\mathcal{O}_{\text{semistable},X,\text{pro\'et}}\{t^{1/2}\})
\end{align}
or 
\begin{align}
f_*(F\otimes_{\Gamma^\mathrm{perfect}_{\text{Robba},X,\text{pro\'et},I}\{t^{1/2}\}} \Gamma^\mathcal{O}_{\text{semistable},X,\text{pro\'et}}\{t^{1/2}\}).
\end{align}
We call $F$ mixed-parity semi-stable if we have the following isomorphism:
\begin{align}
f^*f_*(F\otimes_{\Gamma^\mathrm{perfect}_{\text{Robba},X,\text{pro\'et},\infty}\{t^{1/2}\}} \Gamma^\mathcal{O}_{\text{semistable},X,\text{pro\'et}}\{t^{1/2}\}) \otimes \Gamma^\mathcal{O}_{\text{semistable},X,\text{pro\'et}}\{t^{1/2}\}\\
 \overset{\sim}{\longrightarrow} F \otimes \Gamma^\mathcal{O}_{\text{semistable},X,\text{pro\'et}}\{t^{1/2}\} 
\end{align}
or 
\begin{align}
f^*f_*(F\otimes_{\Gamma^\mathrm{perfect}_{\text{Robba},X,\text{pro\'et},I}\{t^{1/2}\}} \Gamma^\mathcal{O}_{\text{semistable},X,\text{pro\'et}}\{t^{1/2}\}) \otimes \Gamma^\mathcal{O}_{\text{semistable},X,\text{pro\'et}}\{t^{1/2}\} \\\overset{\sim}{\longrightarrow} F \otimes \Gamma^\mathcal{O}_{\text{semistable},X,\text{pro\'et}}\{t^{1/2}\}. 
\end{align}
\end{definition}

\begin{definition}
For any locally free coherent sheaf $F$ over
\begin{align}
\Gamma^\mathrm{perfect}_{\text{Robba},X,\text{pro\'et},\infty}\{t^{1/2}\},\Gamma^\mathrm{perfect}_{\text{Robba},X,\text{pro\'et},I}\{t^{1/2}\},
\end{align} 
we consider the following functor $\mathrm{semistable}^\mathrm{almost}$ sending $F$ to the following object:
\begin{align}
f_*(F\otimes_{\Gamma^\mathrm{perfect}_{\text{Robba},X,\text{pro\'et},\infty}\{t^{1/2}\}} \Gamma^\mathcal{O}_{\text{semistable},X,\text{pro\'et}}\{t^{1/2},\log(t)\})
\end{align}
or 
\begin{align}
f_*(F\otimes_{\Gamma^\mathrm{perfect}_{\text{Robba},X,\text{pro\'et},I}\{t^{1/2}\}} \Gamma^\mathcal{O}_{\text{semistable},X,\text{pro\'et}}\{t^{1/2},\log(t)\}).
\end{align}
We call $F$ mixed-parity almost semi-stable if we have the following isomorphism:
\begin{align}
f^*f_*(F\otimes_{\Gamma^\mathrm{perfect}_{\text{Robba},X,\text{pro\'et},\infty}\{t^{1/2}\}} \Gamma^\mathcal{O}_{\text{semistable},X,\text{pro\'et}}\{t^{1/2},\log(t)\}) \otimes \Gamma^\mathcal{O}_{\text{semistable},X,\text{pro\'et}}\{t^{1/2},\log(t)\} \\
\overset{\sim}{\longrightarrow} F \otimes \Gamma^\mathcal{O}_{\text{semistable},X,\text{pro\'et}}\{t^{1/2},\log(t)\} 
\end{align}
or 
\begin{align}
f^*f_*(F\otimes_{\Gamma^\mathrm{perfect}_{\text{Robba},X,\text{pro\'et},I}\{t^{1/2}\}} \Gamma^\mathcal{O}_{\text{semistable},X,\text{pro\'et}}\{t^{1/2},\log(t)\}) \otimes \Gamma^\mathcal{O}_{\text{semistable},X,\text{pro\'et}}\{t^{1/2},\log(t)\}\\ \overset{\sim}{\longrightarrow} F \otimes \Gamma^\mathcal{O}_{\text{semistable},X,\text{pro\'et}}\{t^{1/2},\log(t)\}. 
\end{align}
\end{definition}

\noindent We now define the $(\infty,1)$-categories of mixed-parity semi-stable modules and he corresponding mixed-parity almost semi-stable modules by using the objects involved to generated these categories:

\begin{definition}
Considering all the mixed parity semi-stable bundles (locally finite free) as defined above, we consider the sub-$(\infty,1)$ category of 
\begin{align}
\mathrm{preModule}^\mathrm{solid,quasicoherent}_{\square,\Gamma^\mathrm{perfect}_{\text{Robba},X,\text{pro\'et},\infty}\{t^{1/2}\}},
\mathrm{preModule}^\mathrm{solid,quasicoherent}_{\square,\Gamma^\mathrm{perfect}_{\text{Robba},X,\text{pro\'et},I}\{t^{1/2}\}} 
\end{align}
generated by the mixed-parity semi-stable bundles (locally finite free ones). These are defined to be the $(\infty,1)$-categories of mixed-parity semi-stable complexes:
\begin{align}
\mathrm{preModule}^\mathrm{solid,quasicoherent,mixed-paritysemistable}_{\square,\Gamma^\mathrm{perfect}_{\text{Robba},X,\text{pro\'et},\infty}\{t^{1/2}\}},
\mathrm{preModule}^\mathrm{solid,quasicoherent,mixed-paritysemistable}_{\square,\Gamma^\mathrm{perfect}_{\text{Robba},X,\text{pro\'et},I}\{t^{1/2}\}}. 
\end{align}
\end{definition}

\begin{definition}
Considering all the mixed parity almost semi-stable bundles (locally finite free) as defined above, we consider the sub-$(\infty,1)$ category of 
\begin{align}
\mathrm{preModule}^\mathrm{solid,quasicoherent}_{\square,\Gamma^\mathrm{perfect}_{\text{Robba},X,\text{pro\'et},\infty}\{t^{1/2}\}},
\mathrm{preModule}^\mathrm{solid,quasicoherent}_{\square,\Gamma^\mathrm{perfect}_{\text{Robba},X,\text{pro\'et},I}\{t^{1/2}\}} 
\end{align}
generated by the mixed-parity almost semi-stable bundles (locally finite free ones). These are defined to be the $(\infty,1)$-categories of mixed-parity semi-stable complexes:
\begin{align}
\mathrm{preModule}^\mathrm{solid,quasicoherent,mixed-parityalmostsemistable}_{\square,\Gamma^\mathrm{perfect}_{\text{Robba},X,\text{pro\'et},\infty}\{t^{1/2}\}},\\
\mathrm{preModule}^\mathrm{solid,quasicoherent,mixed-parityalmostsemistable}_{\square,\Gamma^\mathrm{perfect}_{\text{Robba},X,\text{pro\'et},I}\{t^{1/2}\}}. 
\end{align}
\end{definition}

\indent Then the corresponding mixed-parity semi-stable functors can be extended to these categories:
\begin{align}
\mathrm{preModule}^\mathrm{solid,quasicoherent,mixed-paritysemistable}_{\square,\Gamma^\mathrm{perfect}_{\text{Robba},X,\text{pro\'et},\infty}\{t^{1/2}\}},\\
\mathrm{preModule}^\mathrm{solid,quasicoherent,mixed-paritysemistable}_{\square,\Gamma^\mathrm{perfect}_{\text{Robba},X,\text{pro\'et},I}\{t^{1/2}\}}, 
\end{align}
and
\begin{align}
\mathrm{preModule}^\mathrm{solid,quasicoherent,mixed-parityalmostsemistable}_{\square,\Gamma^\mathrm{perfect}_{\text{Robba},X,\text{pro\'et},\infty}\{t^{1/2}\}},\\
\mathrm{preModule}^\mathrm{solid,quasicoherent,mixed-parityalmostsemistable}_{\square,\Gamma^\mathrm{perfect}_{\text{Robba},X,\text{pro\'et},I}\{t^{1/2}\}}. 
\end{align}

\subsubsection{Mixed-Parity Hodge Modules with Frobenius}

\noindent Now we consider the key objects in our study namely those complexes generated by certain mixed-parity Hodge modules. We start from the following definition.

\begin{remark}
All the coherent sheaves over mixed-parity Robba sheaves in this section will carry the corresponding Frobenius morphism $\varphi: F \overset{\sim}{\longrightarrow} \varphi^*F$.
\end{remark}

\begin{definition}
For any locally free coherent sheaf $F$ over
\begin{align}
\Gamma^\mathrm{perfect}_{\text{Robba},X,\text{pro\'et},\infty}\{t^{1/2}\},\Gamma^\mathrm{perfect}_{\text{Robba},X,\text{pro\'et},I}\{t^{1/2}\},
\end{align} 
we consider the following functor $\mathrm{semistable}$ sending $F$ to the following object:
\begin{align}
f_*(F\otimes_{\Gamma^\mathrm{perfect}_{\text{Robba},X,\text{pro\'et},\infty}\{t^{1/2}\}} \Gamma^\mathcal{O}_{\text{semistable},X,\text{pro\'et}}\{t^{1/2}\})
\end{align}
or 
\begin{align}
f_*(F\otimes_{\Gamma^\mathrm{perfect}_{\text{Robba},X,\text{pro\'et},I}\{t^{1/2}\}} \Gamma^\mathcal{O}_{\text{semistable},X,\text{pro\'et}}\{t^{1/2}\}).
\end{align}
We call $F$ mixed-parity semi-stable if we have the following isomorphism:
\begin{align}
f^*f_*(F\otimes_{\Gamma^\mathrm{perfect}_{\text{Robba},X,\text{pro\'et},\infty}\{t^{1/2}\}} \Gamma^\mathcal{O}_{\text{semistable},X,\text{pro\'et}}\{t^{1/2}\}) \otimes \Gamma^\mathcal{O}_{\text{semistable},X,\text{pro\'et}}\{t^{1/2}\}\\ \overset{\sim}{\longrightarrow} F \otimes \Gamma^\mathcal{O}_{\text{semistable},X,\text{pro\'et}}\{t^{1/2}\} 
\end{align}
or 
\begin{align}
f^*f_*(F\otimes_{\Gamma^\mathrm{perfect}_{\text{Robba},X,\text{pro\'et},I}\{t^{1/2}\}} \Gamma^\mathcal{O}_{\text{semistable},X,\text{pro\'et}}\{t^{1/2}\}) \otimes \Gamma^\mathcal{O}_{\text{semistable},X,\text{pro\'et}}\{t^{1/2}\}\\ \overset{\sim}{\longrightarrow} F \otimes \Gamma^\mathcal{O}_{\text{semistable},X,\text{pro\'et}}\{t^{1/2}\}. 
\end{align}
\end{definition}

\begin{definition}
For any locally free coherent sheaf $F$ over
\begin{align}
\Gamma^\mathrm{perfect}_{\text{Robba},X,\text{pro\'et},\infty}\{t^{1/2}\},\Gamma^\mathrm{perfect}_{\text{Robba},X,\text{pro\'et},I}\{t^{1/2}\},
\end{align} 
we consider the following functor $\mathrm{semistable}^\mathrm{almost}$ sending $F$ to the following object:
\begin{align}
f_*(F\otimes_{\Gamma^\mathrm{perfect}_{\text{Robba},X,\text{pro\'et},\infty}\{t^{1/2}\}} \Gamma^\mathcal{O}_{\text{semistable},X,\text{pro\'et}}\{t^{1/2},\log(t)\})
\end{align}
or 
\begin{align}
f_*(F\otimes_{\Gamma^\mathrm{perfect}_{\text{Robba},X,\text{pro\'et},I}\{t^{1/2}\}} \Gamma^\mathcal{O}_{\text{semistable},X,\text{pro\'et}}\{t^{1/2},\log(t)\}).
\end{align}
We call $F$ mixed-parity almost semi-stable if we have the following isomorphism:
\begin{align}
f^*f_*(F\otimes_{\Gamma^\mathrm{perfect}_{\text{Robba},X,\text{pro\'et},\infty}\{t^{1/2}\}} \Gamma^\mathcal{O}_{\text{semistable},X,\text{pro\'et}}\{t^{1/2},\log(t)\}) \otimes \Gamma^\mathcal{O}_{\text{semistable},X,\text{pro\'et}}\{t^{1/2},\log(t)\} \\
\overset{\sim}{\longrightarrow} F \otimes \Gamma^\mathcal{O}_{\text{semistable},X,\text{pro\'et}}\{t^{1/2},\log(t)\} 
\end{align}
or 
\begin{align}
f^*f_*(F\otimes_{\Gamma^\mathrm{perfect}_{\text{Robba},X,\text{pro\'et},I}\{t^{1/2}\}} \Gamma^\mathcal{O}_{\text{semistable},X,\text{pro\'et}}\{t^{1/2},\log(t)\}) \otimes \Gamma^\mathcal{O}_{\text{semistable},X,\text{pro\'et}}\{t^{1/2},\log(t)\}\\ \overset{\sim}{\longrightarrow} F \otimes \Gamma^\mathcal{O}_{\text{semistable},X,\text{pro\'et}}\{t^{1/2},\log(t)\}. 
\end{align}
\end{definition}

\noindent We now define the $(\infty,1)$-categories of mixed-parity semi-stable modules and he corresponding mixed-parity almost semi-stable modules by using the objects involved to generated these categories:

\begin{definition}
Considering all the mixed parity semi-stable bundles (locally finite free) as defined above, we consider the sub-$(\infty,1)$ category of 
\begin{align}
\varphi\mathrm{preModule}^\mathrm{solid,quasicoherent}_{\square,\Gamma^\mathrm{perfect}_{\text{Robba},X,\text{pro\'et},\infty}\{t^{1/2}\}},
\varphi\mathrm{preModule}^\mathrm{solid,quasicoherent}_{\square,\Gamma^\mathrm{perfect}_{\text{Robba},X,\text{pro\'et},I}\{t^{1/2}\}} 
\end{align}
generated by the mixed-parity semi-stable bundles (locally finite free ones). These are defined to be the $(\infty,1)$-categories of mixed-parity semi-stable complexes:
\begin{align}
\varphi\mathrm{preModule}^\mathrm{solid,quasicoherent,mixed-paritysemistable}_{\square,\Gamma^\mathrm{perfect}_{\text{Robba},X,\text{pro\'et},\infty}\{t^{1/2}\}},
\varphi\mathrm{preModule}^\mathrm{solid,quasicoherent,mixed-paritysemistable}_{\square,\Gamma^\mathrm{perfect}_{\text{Robba},X,\text{pro\'et},I}\{t^{1/2}\}}. 
\end{align}
\end{definition}

\begin{definition}
Considering all the mixed parity almost semi-stable bundles (locally finite free) as defined above, we consider the sub-$(\infty,1)$ category of 
\begin{align}
\varphi\mathrm{preModule}^\mathrm{solid,quasicoherent}_{\square,\Gamma^\mathrm{perfect}_{\text{Robba},X,\text{pro\'et},\infty}\{t^{1/2}\}},
\varphi\mathrm{preModule}^\mathrm{solid,quasicoherent}_{\square,\Gamma^\mathrm{perfect}_{\text{Robba},X,\text{pro\'et},I}\{t^{1/2}\}} 
\end{align}
generated by the mixed-parity almost semi-stable bundles (locally finite free ones). These are defined to be the $(\infty,1)$-categories of mixed-parity semi-stable complexes:
\begin{align}
\varphi\mathrm{preModule}^\mathrm{solid,quasicoherent,mixed-parityalmostsemistable}_{\square,\Gamma^\mathrm{perfect}_{\text{Robba},X,\text{pro\'et},\infty}\{t^{1/2}\}},\\
\varphi\mathrm{preModule}^\mathrm{solid,quasicoherent,mixed-parityalmostsemistable}_{\square,\Gamma^\mathrm{perfect}_{\text{Robba},X,\text{pro\'et},I}\{t^{1/2}\}}. 
\end{align}
\end{definition}

\indent Then the corresponding mixed-parity semi-stable functors can be extended to these categories:
\begin{align}
\varphi\mathrm{preModule}^\mathrm{solid,quasicoherent,mixed-paritysemistable}_{\square,\Gamma^\mathrm{perfect}_{\text{Robba},X,\text{pro\'et},\infty}\{t^{1/2}\}},
\varphi\mathrm{preModule}^\mathrm{solid,quasicoherent,mixed-paritysemistable}_{\square,\Gamma^\mathrm{perfect}_{\text{Robba},X,\text{pro\'et},I}\{t^{1/2}\}}, 
\end{align}
and
\begin{align}
\varphi\mathrm{preModule}^\mathrm{solid,quasicoherent,mixed-parityalmostsemistable}_{\square,\Gamma^\mathrm{perfect}_{\text{Robba},X,\text{pro\'et},\infty}\{t^{1/2}\}},\\
\varphi\mathrm{preModule}^\mathrm{solid,quasicoherent,mixed-parityalmostsemistable}_{\square,\Gamma^\mathrm{perfect}_{\text{Robba},X,\text{pro\'et},I}\{t^{1/2}\}}. 
\end{align}

\subsection{Mixed-Parity semi-stable Riemann-Hilbert Correspondence}

\indent This chapter will extend the corresponding Riemann-Hilbert correspondence from \cite{Sch1}, \cite{LZ}, \cite{BL1}, \cite{BL2}, \cite{M} to the mixed-parity setting.

\begin{definition}
We define the following Riemann-Hilbert functor $\text{RH}_\text{mixed-parity}$ from the one of categories:
\begin{align}
\mathrm{preModule}^\mathrm{solid,quasicoherent,mixed-paritysemistable}_{\square,\Gamma^\mathrm{perfect}_{\text{Robba},X,\text{pro\'et},\infty}\{t^{1/2}\}},
\mathrm{preModule}^\mathrm{solid,quasicoherent,mixed-paritysemistable}_{\square,\Gamma^\mathrm{perfect}_{\text{Robba},X,\text{pro\'et},I}\{t^{1/2}\}}, 
\end{align}
and
\begin{align}
\mathrm{preModule}^\mathrm{solid,quasicoherent,mixed-parityalmostsemistable}_{\square,\Gamma^\mathrm{perfect}_{\text{Robba},X,\text{pro\'et},\infty}\{t^{1/2}\}},\\
\mathrm{preModule}^\mathrm{solid,quasicoherent,mixed-parityalmostsemistable}_{\square,\Gamma^\mathrm{perfect}_{\text{Robba},X,\text{pro\'et},I}\{t^{1/2}\}} 
\end{align}
to $(\infty,1)$-categories in image denoted by:
\begin{align}
\mathrm{preModule}_{X,\text{\'et}}
\end{align}
to be the following functors sending each $F$ in the domain to:
\begin{align}
&\text{RH}_\text{mixed-parity}(F):=f_*(F\otimes_{\Gamma^\mathrm{perfect}_{\text{Robba},X,\text{pro\'et},\infty}\{t^{1/2}\}} \Gamma^\mathcal{O}_{\text{semistable},X,\text{pro\'et}}\{t^{1/2}\}),\\
&\text{RH}_\text{mixed-parity}(F):=f_*(F\otimes_{\Gamma^\mathrm{perfect}_{\text{Robba},X,\text{pro\'et},I}\{t^{1/2}\}} \Gamma^\mathcal{O}_{\text{semistable},X,\text{pro\'et}}\{t^{1/2}\}),\\
&\text{RH}_\text{mixed-parity}(F):=f_*(F\otimes_{\Gamma^\mathrm{perfect}_{\text{Robba},X,\text{pro\'et},\infty}\{t^{1/2}\}} \Gamma^\mathcal{O}_{\text{semistable},X,\text{pro\'et}}\{t^{1/2},\log(t)\}),\\
&\text{RH}_\text{mixed-parity}(F):=f_*(F\otimes_{\Gamma^\mathrm{perfect}_{\text{Robba},X,\text{pro\'et},I}\{t^{1/2}\}} \Gamma^\mathcal{O}_{\text{semistable},X,\text{pro\'et}}\{t^{1/2},\log(t)\}),\\
\end{align}
respectively.

\end{definition}

\begin{definition}
In the situation where we have the Frobenius action we consider the follwing. We define the following Riemann-Hilbert functor $\text{RH}_\text{mixed-parity}$ from the one of categories:
\begin{align}
\varphi\mathrm{preModule}^\mathrm{solid,quasicoherent,mixed-paritysemistable}_{\square,\Gamma^\mathrm{perfect}_{\text{Robba},X,\text{pro\'et},\infty}\{t^{1/2}\}},
\varphi\mathrm{preModule}^\mathrm{solid,quasicoherent,mixed-paritysemistable}_{\square,\Gamma^\mathrm{perfect}_{\text{Robba},X,\text{pro\'et},I}\{t^{1/2}\}}, 
\end{align}
and
\begin{align}
\varphi\mathrm{preModule}^\mathrm{solid,quasicoherent,mixed-parityalmostsemistable}_{\square,\Gamma^\mathrm{perfect}_{\text{Robba},X,\text{pro\'et},\infty}\{t^{1/2}\}},\\
\varphi\mathrm{preModule}^\mathrm{solid,quasicoherent,mixed-parityalmostsemistable}_{\square,\Gamma^\mathrm{perfect}_{\text{Robba},X,\text{pro\'et},I}\{t^{1/2}\}} 
\end{align}
to $(\infty,1)$-categories in image denoted by:
\begin{align}
\mathrm{preModule}_{X,\text{\'et}}
\end{align}
to be the following functors sending each $F$ in the domain to:
\begin{align}
&\text{RH}_\text{mixed-parity}(F):=f_*(F\otimes_{\Gamma^\mathrm{perfect}_{\text{Robba},X,\text{pro\'et},\infty}\{t^{1/2}\}} \Gamma^\mathcal{O}_{\text{semistable},X,\text{pro\'et}}\{t^{1/2}\}),\\
&\text{RH}_\text{mixed-parity}(F):=f_*(F\otimes_{\Gamma^\mathrm{perfect}_{\text{Robba},X,\text{pro\'et},I}\{t^{1/2}\}} \Gamma^\mathcal{O}_{\text{semistable},X,\text{pro\'et}}\{t^{1/2}\}),\\
&\text{RH}_\text{mixed-parity}(F):=f_*(F\otimes_{\Gamma^\mathrm{perfect}_{\text{Robba},X,\text{pro\'et},\infty}\{t^{1/2}\}} \Gamma^\mathcal{O}_{\text{semistable},X,\text{pro\'et}}\{t^{1/2},\log(t)\}),\\
&\text{RH}_\text{mixed-parity}(F):=f_*(F\otimes_{\Gamma^\mathrm{perfect}_{\text{Robba},X,\text{pro\'et},I}\{t^{1/2}\}} \Gamma^\mathcal{O}_{\text{semistable},X,\text{pro\'et}}\{t^{1/2},\log(t)\}),\\
\end{align}
respectively.

\end{definition}

\newpage
\section{Localizations}

\begin{reference}
\cite{AI1}, \cite{AI2}, \cite{AB1}, \cite{AB2}, \cite{Fon2}, \cite{Fon3}, \cite{Fa1}.
\end{reference}

\subsection{Extension of Fundamental Groups}

\indent In the local setting setting in fact we can have more thorough understanding of more structures. Locally we can have the Galois group of $\mathbb{Q}_p\left<T_1,...,T_n\right>$ for some $n>0$ in the smooth situation for instance. Our current discussion will be in the following situation:

\begin{definition}
We define the corresponding two fold covering of the Galois group:
\begin{align}
\mathrm{Gal}(\overline{\mathbb{Q}_p\left<T_1,...,T_n\right>}^\wedge/\mathbb{Q}_p\left<T_1,...,T_n\right>)_2
\end{align}
by taking the product of
\begin{align}
\mathrm{Gal}(\overline{\mathbb{Q}_p\left<T_1,...,T_n\right>}^\wedge/\mathbb{Q}_p\left<T_1,...,T_n\right>),\mathrm{Gal}(\overline{Q}_p/\mathbb{Q}_p)_2
\end{align}
where the latter is the group defined in \cite[Just before Lemma 7.5]{BS}. This group admits an action on the element $t^{1/2}$ through the action of the group $\mathrm{Gal}(\overline{\mathbb{Q}}_p/\mathbb{Q}_p)_2$.
\end{definition}

\subsection{Modules}

\noindent We consider the following definition of modules with $(\varphi,\mathrm{Gal}(\overline{\mathbb{Q}_p\left<T_1,...,T_n\right>}^\wedge/\mathbb{Q}_p\left<T_1,...,T_n\right>)_2)$-structure.

\begin{definition}
Let $R:=\mathbb{Q}_p\left<T_1,...,T_n\right>$. We use the notation:
\begin{align}
\mathrm{Module}^\mathrm{solid,quasicoherent}_{\square,\Gamma^\mathrm{perfect}_{\text{Robba},R,\text{pro\'et}}\{t^{1/2}\}},\mathrm{Module}^\mathrm{solid,quasicoherent}_{\square,\Gamma^\mathrm{perfect}_{\text{Robba},R,\text{pro\'et},\infty}\{t^{1/2}\}},
\mathrm{Module}^\mathrm{solid,quasicoherent}_{\square,\Gamma^\mathrm{perfect}_{\text{Robba},R,\text{pro\'et},I}\{t^{1/2}\}} 
\end{align}
to denote the $(\infty,1)$-categories of solid modules over the corresonding Robba rings in the local setting namey associated to:
\begin{align}
R^\mathrm{perf\flat}:=\mathbb{Q}_p(p^{1/p^\infty})\left<T^{1/p^\infty}_1,...,T^{1/p^\infty}_n\right>^{\wedge\flat}.
\end{align}
Then we consider all the modules as such carrying commuting operations from $\varphi$ and 
\begin{align}
\Sigma:=\mathrm{Gal}(\overline{\mathbb{Q}_p\left<T_1,...,T_n\right>}^\wedge/\mathbb{Q}_p\left<T_1,...,T_n\right>)_2,
\end{align}
which is assumed to be semilinear. We use the notation 
\begin{align}
\overset{\mathrm{Module}}{\varphi,\Sigma}^\mathrm{solid,quasicoherent}_{\square,\Gamma^\mathrm{perfect}_{\text{Robba},R,\text{pro\'et}}\{t^{1/2}\}},\overset{\mathrm{Module}}{\varphi,\Sigma}^\mathrm{solid,quasicoherent}_{\square,\Gamma^\mathrm{perfect}_{\text{Robba},R,\text{pro\'et},\infty}\{t^{1/2}\}},
\overset{\mathrm{Module}}{\varphi,\Sigma}^\mathrm{solid,quasicoherent}_{\square,\Gamma^\mathrm{perfect}_{\text{Robba},R,\text{pro\'et},I}\{t^{1/2}\}} 
\end{align}
to denote the categories.
\end{definition}

\begin{definition}
For any module $F$ over 
\begin{align}
\Gamma^\mathrm{perfect}_{\text{Robba},R,\text{pro\'et},\infty}\{t^{1/2}\},\Gamma^\mathrm{perfect}_{\text{Robba},R,\text{pro\'et},I}\{t^{1/2}\},
\end{align} 
carrying the structure of $(\varphi,\Sigma)$-action, we consider the following functor $\mathrm{dR}$ sending $F$ to the following object:
\begin{align}
(F\otimes_{\Gamma^\mathrm{perfect}_{\text{Robba},R,\text{pro\'et},\infty}\{t^{1/2}\}} \Gamma^\mathcal{O}_{\text{deRham},R,\text{pro\'et}}\{t^{1/2}\})^\Sigma
\end{align}
or 
\begin{align}
f_*(F\otimes_{\Gamma^\mathrm{perfect}_{\text{Robba},R,\text{pro\'et},I}\{t^{1/2}\}} \Gamma^\mathcal{O}_{\text{deRham},R,\text{pro\'et}}\{t^{1/2}\}).
\end{align}
We call $F$ mixed-parity de Rham if we have the following isomorphism:
\begin{align}
(F\otimes_{\Gamma^\mathrm{perfect}_{\text{Robba},R,\text{pro\'et},\infty}\{t^{1/2}\}} \Gamma^\mathcal{O}_{\text{deRham},R,\text{pro\'et}}\{t^{1/2}\})^\Sigma \otimes \Gamma^\mathcal{O}_{\text{deRham},R,\text{pro\'et}}\{t^{1/2}\} \overset{\sim}{\longrightarrow} F \otimes \Gamma^\mathcal{O}_{\text{deRham},R,\text{pro\'et}}\{t^{1/2}\} 
\end{align}
or 
\begin{align}
(F\otimes_{\Gamma^\mathrm{perfect}_{\text{Robba},R,\text{pro\'et},I}\{t^{1/2}\}} \Gamma^\mathcal{O}_{\text{deRham},R,\text{pro\'et}}\{t^{1/2}\})^\Sigma \otimes \Gamma^\mathcal{O}_{\text{deRham},R,\text{pro\'et}}\{t^{1/2}\} \overset{\sim}{\longrightarrow} F \otimes \Gamma^\mathcal{O}_{\text{deRham},R,\text{pro\'et}}\{t^{1/2}\}. 
\end{align}
\end{definition}

\begin{definition}
For any module $F$ over 
\begin{align}
\Gamma^\mathrm{perfect}_{\text{Robba},R,\text{pro\'et},\infty}\{t^{1/2}\},\Gamma^\mathrm{perfect}_{\text{Robba},R,\text{pro\'et},I}\{t^{1/2}\},
\end{align} 
carrying the structure of $(\varphi,\Sigma)$-action, we consider the following functor $\mathrm{dR}$ sending $F$ to the following object:
\begin{align}
(F\otimes_{\Gamma^\mathrm{perfect}_{\text{Robba},R,\text{pro\'et},\infty}\{t^{1/2}\}} \Gamma^\mathcal{O}_{\text{deRham},R,\text{pro\'et}}\{t^{1/2},\log(t)\})^\Sigma
\end{align}
or 
\begin{align}
f_*(F\otimes_{\Gamma^\mathrm{perfect}_{\text{Robba},R,\text{pro\'et},I}\{t^{1/2}\}} \Gamma^\mathcal{O}_{\text{deRham},R,\text{pro\'et}}\{t^{1/2}\}).
\end{align}
We call $F$ mixed-parity almost de Rham if we have the following isomorphism:
\begin{align}
(F\otimes_{\Gamma^\mathrm{perfect}_{\text{Robba},R,\text{pro\'et},\infty}\{t^{1/2}\}} \Gamma^\mathcal{O}_{\text{deRham},R,\text{pro\'et}}\{t^{1/2},\log(t)\})^\Sigma \otimes \Gamma^\mathcal{O}_{\text{deRham},R,\text{pro\'et}}\{t^{1/2},\log(t)\} \\
\overset{\sim}{\longrightarrow} F \otimes \Gamma^\mathcal{O}_{\text{deRham},R,\text{pro\'et}}\{t^{1/2},\log(t)\} 
\end{align}
or 
\begin{align}
(F\otimes_{\Gamma^\mathrm{perfect}_{\text{Robba},R,\text{pro\'et},I}\{t^{1/2}\}} \Gamma^\mathcal{O}_{\text{deRham},R,\text{pro\'et}}\{t^{1/2},\log(t)\})^\Sigma \otimes \Gamma^\mathcal{O}_{\text{deRham},R,\text{pro\'et}}\{t^{1/2},\log(t)\}\\ \overset{\sim}{\longrightarrow} F \otimes \Gamma^\mathcal{O}_{\text{deRham},R,\text{pro\'et}}\{t^{1/2},\log(t)\}. 
\end{align}
\end{definition}

\begin{definition}
For any module $F$ over 
\begin{align}
\Gamma^\mathrm{perfect}_{\text{Robba},R,\text{pro\'et},\infty}\{t^{1/2}\},\Gamma^\mathrm{perfect}_{\text{Robba},R,\text{pro\'et},I}\{t^{1/2}\},
\end{align} 
carrying the structure of $(\varphi,\Sigma)$-action, we consider the following functor $\mathrm{cristalline}$ sending $F$ to the following object:
\begin{align}
(F\otimes_{\Gamma^\mathrm{perfect}_{\text{Robba},R,\text{pro\'et},\infty}\{t^{1/2}\}} \Gamma^\mathcal{O}_{\text{cristalline},R,\text{pro\'et}}\{t^{1/2}\})^\Sigma
\end{align}
or 
\begin{align}
f_*(F\otimes_{\Gamma^\mathrm{perfect}_{\text{Robba},R,\text{pro\'et},I}\{t^{1/2}\}} \Gamma^\mathcal{O}_{\text{cristalline},R,\text{pro\'et}}\{t^{1/2}\}).
\end{align}
We call $F$ mixed-parity cristalline if we have the following isomorphism:
\begin{align}
(F\otimes_{\Gamma^\mathrm{perfect}_{\text{Robba},R,\text{pro\'et},\infty}\{t^{1/2}\}} \Gamma^\mathcal{O}_{\text{cristalline},R,\text{pro\'et}}\{t^{1/2}\})^\Sigma \otimes \Gamma^\mathcal{O}_{\text{cristalline},R,\text{pro\'et}}\{t^{1/2}\} \overset{\sim}{\longrightarrow} F \otimes \Gamma^\mathcal{O}_{\text{cristalline},R,\text{pro\'et}}\{t^{1/2}\} 
\end{align}
or 
\begin{align}
(F\otimes_{\Gamma^\mathrm{perfect}_{\text{Robba},R,\text{pro\'et},I}\{t^{1/2}\}} \Gamma^\mathcal{O}_{\text{cristalline},R,\text{pro\'et}}\{t^{1/2}\})^\Sigma \otimes \Gamma^\mathcal{O}_{\text{cristalline},R,\text{pro\'et}}\{t^{1/2}\} \overset{\sim}{\longrightarrow} F \otimes \Gamma^\mathcal{O}_{\text{cristalline},R,\text{pro\'et}}\{t^{1/2}\}. 
\end{align}
\end{definition}

\begin{definition}
For any module $F$ over 
\begin{align}
\Gamma^\mathrm{perfect}_{\text{Robba},R,\text{pro\'et},\infty}\{t^{1/2}\},\Gamma^\mathrm{perfect}_{\text{Robba},R,\text{pro\'et},I}\{t^{1/2}\},
\end{align} 
carrying the structure of $(\varphi,\Sigma)$-action, we consider the following functor $\mathrm{cristalline}^\mathrm{almost}$ sending $F$ to the following object:
\begin{align}
(F\otimes_{\Gamma^\mathrm{perfect}_{\text{Robba},R,\text{pro\'et},\infty}\{t^{1/2}\}} \Gamma^\mathcal{O}_{\text{cristalline},R,\text{pro\'et}}\{t^{1/2},\log(t)\})^\Sigma
\end{align}
or 
\begin{align}
f_*(F\otimes_{\Gamma^\mathrm{perfect}_{\text{Robba},R,\text{pro\'et},I}\{t^{1/2}\}} \Gamma^\mathcal{O}_{\text{cristalline},R,\text{pro\'et}}\{t^{1/2}\}).
\end{align}
We call $F$ mixed-parity almost cristalline if we have the following isomorphism:
\begin{align}
(F\otimes_{\Gamma^\mathrm{perfect}_{\text{Robba},R,\text{pro\'et},\infty}\{t^{1/2}\}} \Gamma^\mathcal{O}_{\text{cristalline},R,\text{pro\'et}}\{t^{1/2},\log(t)\})^\Sigma \otimes \Gamma^\mathcal{O}_{\text{cristalline},R,\text{pro\'et}}\{t^{1/2},\log(t)\} \\
\overset{\sim}{\longrightarrow} F \otimes \Gamma^\mathcal{O}_{\text{cristalline},R,\text{pro\'et}}\{t^{1/2},\log(t)\} 
\end{align}
or 
\begin{align}
(F\otimes_{\Gamma^\mathrm{perfect}_{\text{Robba},R,\text{pro\'et},I}\{t^{1/2}\}} \Gamma^\mathcal{O}_{\text{cristalline},R,\text{pro\'et}}\{t^{1/2},\log(t)\})^\Sigma \otimes \Gamma^\mathcal{O}_{\text{cristalline},R,\text{pro\'et}}\{t^{1/2},\log(t)\}\\ \overset{\sim}{\longrightarrow} F \otimes \Gamma^\mathcal{O}_{\text{cristalline},R,\text{pro\'et}}\{t^{1/2},\log(t)\}. 
\end{align}
\end{definition}

\begin{definition}
For any module $F$ over 
\begin{align}
\Gamma^\mathrm{perfect}_{\text{Robba},R,\text{pro\'et},\infty}\{t^{1/2}\},\Gamma^\mathrm{perfect}_{\text{Robba},R,\text{pro\'et},I}\{t^{1/2}\},
\end{align} 
carrying the structure of $(\varphi,\Sigma)$-action, we consider the following functor $\mathrm{semistable}$ sending $F$ to the following object:
\begin{align}
(F\otimes_{\Gamma^\mathrm{perfect}_{\text{Robba},R,\text{pro\'et},\infty}\{t^{1/2}\}} \Gamma^\mathcal{O}_{\text{semistable},R,\text{pro\'et}}\{t^{1/2}\})^\Sigma
\end{align}
or 
\begin{align}
f_*(F\otimes_{\Gamma^\mathrm{perfect}_{\text{Robba},R,\text{pro\'et},I}\{t^{1/2}\}} \Gamma^\mathcal{O}_{\text{deRham},R,\text{pro\'et}}\{t^{1/2}\}).
\end{align}
We call $F$ mixed-parity semi-stable if we have the following isomorphism:
\begin{align}
(F\otimes_{\Gamma^\mathrm{perfect}_{\text{Robba},R,\text{pro\'et},\infty}\{t^{1/2}\}} \Gamma^\mathcal{O}_{\text{semistable},R,\text{pro\'et}}\{t^{1/2}\})^\Sigma \otimes \Gamma^\mathcal{O}_{\text{semistable},R,\text{pro\'et}}\{t^{1/2}\} \overset{\sim}{\longrightarrow} F \otimes \Gamma^\mathcal{O}_{\text{semistable},R,\text{pro\'et}}\{t^{1/2}\} 
\end{align}
or 
\begin{align}
(F\otimes_{\Gamma^\mathrm{perfect}_{\text{Robba},R,\text{pro\'et},I}\{t^{1/2}\}} \Gamma^\mathcal{O}_{\text{semistable},R,\text{pro\'et}}\{t^{1/2}\})^\Sigma \otimes \Gamma^\mathcal{O}_{\text{semistable},R,\text{pro\'et}}\{t^{1/2}\} \overset{\sim}{\longrightarrow} F \otimes \Gamma^\mathcal{O}_{\text{semistable},R,\text{pro\'et}}\{t^{1/2}\}. 
\end{align}
\end{definition}

\begin{definition}
For any module $F$ over 
\begin{align}
\Gamma^\mathrm{perfect}_{\text{Robba},R,\text{pro\'et},\infty}\{t^{1/2}\},\Gamma^\mathrm{perfect}_{\text{Robba},R,\text{pro\'et},I}\{t^{1/2}\},
\end{align} 
carrying the structure of $(\varphi,\Sigma)$-action, we consider the following functor $\mathrm{semistable}^\mathrm{almost}$ sending $F$ to the following object:
\begin{align}
(F\otimes_{\Gamma^\mathrm{perfect}_{\text{Robba},R,\text{pro\'et},\infty}\{t^{1/2}\}} \Gamma^\mathcal{O}_{\text{semistable},R,\text{pro\'et}}\{t^{1/2},\log(t)\})^\Sigma
\end{align}
or 
\begin{align}
f_*(F\otimes_{\Gamma^\mathrm{perfect}_{\text{Robba},R,\text{pro\'et},I}\{t^{1/2}\}} \Gamma^\mathcal{O}_{\text{semistable},R,\text{pro\'et}}\{t^{1/2}\}).
\end{align}
We call $F$ mixed-parity almost de Rham if we have the following isomorphism:
\begin{align}
(F\otimes_{\Gamma^\mathrm{perfect}_{\text{Robba},R,\text{pro\'et},\infty}\{t^{1/2}\}} \Gamma^\mathcal{O}_{\text{semistable},R,\text{pro\'et}}\{t^{1/2},\log(t)\})^\Sigma \otimes \Gamma^\mathcal{O}_{\text{semistable},R,\text{pro\'et}}\{t^{1/2},\log(t)\} \\
\overset{\sim}{\longrightarrow} F \otimes \Gamma^\mathcal{O}_{\text{semistable},R,\text{pro\'et}}\{t^{1/2},\log(t)\} 
\end{align}
or 
\begin{align}
(F\otimes_{\Gamma^\mathrm{perfect}_{\text{Robba},R,\text{pro\'et},I}\{t^{1/2}\}} \Gamma^\mathcal{O}_{\text{semistable},R,\text{pro\'et}}\{t^{1/2},\log(t)\})^\Sigma \otimes \Gamma^\mathcal{O}_{\text{semistable},R,\text{pro\'et}}\{t^{1/2},\log(t)\}\\ \overset{\sim}{\longrightarrow} F \otimes \Gamma^\mathcal{O}_{\text{semistable},R,\text{pro\'et}}\{t^{1/2},\log(t)\}. 
\end{align}
\end{definition}

\chapter{Mixed-Parity $p$-adic Hodge Modules in $v$-Topology}

\newpage
\section{Geometric Family of Mixed-Parity Hodge Modules I: de Rham Situations}

\begin{reference}
\cite{Sch1}, \cite{Sch2}, \cite{FS}, \cite{KL1}, \cite{KL2}, \cite{BL1}, \cite{BL2}, \cite{BS}, \cite{BHS}, \cite{Fon1}, \cite{CS1}, \cite{CS2}, \cite{BK}, \cite{BBK}, \cite{BBBK}, \cite{KKM}, \cite{KM}, \cite{LZ}, \cite{M}.
\end{reference}

\subsection{Period Rings and Sheaves}

\subsubsection{Rings}

\noindent Let $X$ be a rigid analytic space over $\mathbb{Q}_p$. We have the corresponding \'etale site and the corresponding pro-\'etale site of $X$, which we denote them by $X_{v},X_\text{\'et}$. The relationship of the two sites can be reflected by the corresponding morphism $f:X_{v}\longrightarrow X_\text{\'et}$. Then we have the corresponding de Rham period rings and sheaves from \cite{Sch1}:
\begin{align}
\Gamma_{\text{deRham},X,v}, \Gamma^\mathcal{O}_{\text{deRham},X,v}.
\end{align}
Our notations are different from \cite{Sch1}, we use $\Gamma$ to mean $B$ in \cite{Sch1}, while $\Gamma^\mathcal{O}$ will be the corresponding $OB$ ring in \cite{Sch1}.\\

\begin{definition}
\indent Now we assume that $p>2$, following \cite{BS} we join the square root of $t$ element in $\Gamma_{\text{deRham},X,v}$ which forms the sheaves:
\begin{align}
\Gamma_{\text{deRham},X,v}\{t^{1/2}\},\Gamma^\mathcal{O}_{\text{deRham},X,v}\{t^{1/2}\}.
\end{align}
And following \cite{BL1}, \cite{BL2}, \cite{Fon1}, \cite{BHS} we further have the following sheaves of rings:
\begin{align}
\Gamma_{\text{deRham},X,v}\{t^{1/2},\log(t)\},\Gamma^\mathcal{O}_{\text{deRham},X,v}\{t^{1/2},\log(t)\}.
\end{align}
\end{definition}

\begin{definition}
We use the notations:
\begin{align}
\Gamma^\mathrm{perfect}_{\text{Robba},X,v},\Gamma^\mathrm{perfect}_{\text{Robba},X,v,\infty},\Gamma^\mathrm{perfect}_{\text{Robba},X,v,I}
\end{align}
to denote the perfect Robba rings from \cite{KL1}, \cite{KL2}, where $I\subset (0,\infty)$. Then we join $t^{1/2}$ to these sheaves we have:
\begin{align}
\Gamma^\mathrm{perfect}_{\text{Robba},X,v}\{t^{1/2}\},\Gamma^\mathrm{perfect}_{\text{Robba},X,v,\infty}\{t^{1/2}\},\Gamma^\mathrm{perfect}_{\text{Robba},X,v,I}\{t^{1/2}\}.
\end{align}
And following \cite{BL1}, \cite{BL2}, \cite{Fon1}, \cite{BHS} we have the following larger sheaves:
\begin{align}
\Gamma^\mathrm{perfect}_{\text{Robba},X,v}\{t^{1/2},\log(t)\},\Gamma^\mathrm{perfect}_{\text{Robba},X,v,\infty}\{t^{1/2},\log(t)\},\Gamma^\mathrm{perfect}_{\text{Robba},X,v,I}\{t^{1/2},\log(t)\}.
\end{align} 
\end{definition}

\begin{definition}
From now on, we use the same notation to denote the period rings involved tensored with a finite extension of $\mathbb{Q}_p$ containing square root of $p$ as in \cite{BS}.
\begin{align}
\Gamma_{\text{deRham},X,v}\{t^{1/2}\},\Gamma^\mathcal{O}_{\text{deRham},X,v}\{t^{1/2}\}.
\end{align}
\begin{align}
\Gamma_{\text{deRham},X,v}\{t^{1/2},\log(t)\},\Gamma^\mathcal{O}_{\text{deRham},X,v}\{t^{1/2},\log(t)\}.
\end{align}
\begin{align}
\Gamma^\mathrm{perfect}_{\text{Robba},X,v}\{t^{1/2}\},\Gamma^\mathrm{perfect}_{\text{Robba},X,v,\infty}\{t^{1/2}\},\Gamma^\mathrm{perfect}_{\text{Robba},X,v,I}\{t^{1/2}\}.
\end{align}
\begin{align}
\Gamma^\mathrm{perfect}_{\text{Robba},X,v}\{t^{1/2},\log(t)\},\Gamma^\mathrm{perfect}_{\text{Robba},X,v,\infty}\{t^{1/2},\log(t)\},\Gamma^\mathrm{perfect}_{\text{Robba},X,v,I}\{t^{1/2},\log(t)\}.
\end{align}
This is necessary since we to extend the action of $\varphi$ to the period rings by $\varphi(t^{1/2}\otimes 1)=\varphi(t)^{1/2}\otimes 1$.
\end{definition}

\subsubsection{Modules}

\noindent We consider quasicoherent presheaves in the following two situation:
\begin{itemize}
\item[$\square$] The solid quasicoherent modules from \cite{CS1}, \cite{CS2};
\item[$\square$] The ind-Banach quasicoherent modules from \cite{BK}, \cite{BBK}, \cite{BBBK}, \cite{KKM}, \cite{KM} with the corresponding monomorphic ind-Banach quasicoherent modules from \cite{BK}, \cite{BBK}, \cite{BBBK}, \cite{KKM}, \cite{KM}.
\end{itemize}

\begin{definition}
We use the notation:
\begin{align}
\mathrm{preModule}^\mathrm{solid,quasicoherent}_{\square,\Gamma^\mathrm{perfect}_{\text{Robba},X,v}\{t^{1/2}\}},\mathrm{preModule}^\mathrm{solid,quasicoherent}_{\square,\Gamma^\mathrm{perfect}_{\text{Robba},X,v,\infty}\{t^{1/2}\}},
\mathrm{preModule}^\mathrm{solid,quasicoherent}_{\square,\Gamma^\mathrm{perfect}_{\text{Robba},X,v,I}\{t^{1/2}\}} 
\end{align}
to denote the $(\infty,1)$-categories of solid quasicoherent presheaves over the corresonding Robba sheaves. Locally the section is defined by taking the corresponding $(\infty,1)$-categories of solid modules.
\end{definition}

\begin{definition}
We use the notation:
\begin{align}
\mathrm{preModule}^\mathrm{ind-Banach,quasicoherent}_{\Gamma^\mathrm{perfect}_{\text{Robba},X,v}\{t^{1/2}\}},\\\mathrm{preModule}^\mathrm{ind-Banach,quasicoherent}_{\Gamma^\mathrm{perfect}_{\text{Robba},X,v,\infty}\{t^{1/2}\}},\\
\mathrm{preModule}^\mathrm{ind-Banach,quasicoherent}_{\Gamma^\mathrm{perfect}_{\text{Robba},X,v,I}\{t^{1/2}\}} 
\end{align}
to denote the $(\infty,1)$-categories of solid quasicoherent presheaves over the corresonding Robba sheaves. Locally the section is defined by taking the corresponding $(\infty,1)$-categories of inductive Banach  modules. 
\end{definition}

\begin{definition}
We use the notation:
\begin{align}
\mathrm{Module}^\mathrm{solid,quasicoherent}_{\square,\Gamma^\mathrm{perfect}_{\text{Robba},X,v}\{t^{1/2}\}},\mathrm{Module}^\mathrm{solid,quasicoherent}_{\square,\Gamma^\mathrm{perfect}_{\text{Robba},X,v,\infty}\{t^{1/2}\}},
\mathrm{Module}^\mathrm{solid,quasicoherent}_{\square,\Gamma^\mathrm{perfect}_{\text{Robba},X,v,I}\{t^{1/2}\}} 
\end{align}
to denote the $(\infty,1)$-categories of solid quasicoherent sheaves over the corresonding Robba sheaves. Locally the section is defined by taking the corresponding $(\infty,1)$-categories of solid modules.
\end{definition}

\subsubsection{Mixed-Parity Hodge Modules without Frobenius}

\noindent Now we consider the key objects in our study namely those complexes generated by certain mixed-parity Hodge modules. We start from the following definition.

\begin{definition}
For any locally free coherent sheaf $F$ over
\begin{align}
\Gamma^\mathrm{perfect}_{\text{Robba},X,v,\infty}\{t^{1/2}\},\Gamma^\mathrm{perfect}_{\text{Robba},X,v,I}\{t^{1/2}\},
\end{align} 
we consider the following functor $\mathrm{dR}$ sending $F$ to the following object:
\begin{align}
f_*(F\otimes_{\Gamma^\mathrm{perfect}_{\text{Robba},X,v,\infty}\{t^{1/2}\}} \Gamma^\mathcal{O}_{\text{deRham},X,v}\{t^{1/2}\})
\end{align}
or 
\begin{align}
f_*(F\otimes_{\Gamma^\mathrm{perfect}_{\text{Robba},X,v,I}\{t^{1/2}\}} \Gamma^\mathcal{O}_{\text{deRham},X,v}\{t^{1/2}\}).
\end{align}
We call $F$ mixed-parity de Rham if we have the following isomorphism:
\begin{align}
f^*f_*(F\otimes_{\Gamma^\mathrm{perfect}_{\text{Robba},X,v,\infty}\{t^{1/2}\}} \Gamma^\mathcal{O}_{\text{deRham},X,v}\{t^{1/2}\}) \otimes \Gamma^\mathcal{O}_{\text{deRham},X,v}\{t^{1/2}\} \overset{\sim}{\longrightarrow} F \otimes \Gamma^\mathcal{O}_{\text{deRham},X,v}\{t^{1/2}\} 
\end{align}
or 
\begin{align}
f^*f_*(F\otimes_{\Gamma^\mathrm{perfect}_{\text{Robba},X,v,I}\{t^{1/2}\}} \Gamma^\mathcal{O}_{\text{deRham},X,v}\{t^{1/2}\}) \otimes \Gamma^\mathcal{O}_{\text{deRham},X,v}\{t^{1/2}\} \overset{\sim}{\longrightarrow} F \otimes \Gamma^\mathcal{O}_{\text{deRham},X,v}\{t^{1/2}\}. 
\end{align}
\end{definition}

\begin{definition}
For any locally free coherent sheaf $F$ over
\begin{align}
\Gamma^\mathrm{perfect}_{\text{Robba},X,v,\infty}\{t^{1/2}\},\Gamma^\mathrm{perfect}_{\text{Robba},X,v,I}\{t^{1/2}\},
\end{align} 
we consider the following functor $\mathrm{dR}^\mathrm{almost}$ sending $F$ to the following object:
\begin{align}
f_*(F\otimes_{\Gamma^\mathrm{perfect}_{\text{Robba},X,v,\infty}\{t^{1/2}\}} \Gamma^\mathcal{O}_{\text{deRham},X,v}\{t^{1/2},\log(t)\})
\end{align}
or 
\begin{align}
f_*(F\otimes_{\Gamma^\mathrm{perfect}_{\text{Robba},X,v,I}\{t^{1/2}\}} \Gamma^\mathcal{O}_{\text{deRham},X,v}\{t^{1/2},\log(t)\}).
\end{align}
We call $F$ mixed-parity almost de Rham if we have the following isomorphism:
\begin{align}
f^*f_*(F\otimes_{\Gamma^\mathrm{perfect}_{\text{Robba},X,v,\infty}\{t^{1/2}\}} \Gamma^\mathcal{O}_{\text{deRham},X,v}\{t^{1/2},\log(t)\}) \otimes \Gamma^\mathcal{O}_{\text{deRham},X,v}\{t^{1/2},\log(t)\} \\
\overset{\sim}{\longrightarrow} F \otimes \Gamma^\mathcal{O}_{\text{deRham},X,v}\{t^{1/2},\log(t)\} 
\end{align}
or 
\begin{align}
f^*f_*(F\otimes_{\Gamma^\mathrm{perfect}_{\text{Robba},X,v,I}\{t^{1/2}\}} \Gamma^\mathcal{O}_{\text{deRham},X,v}\{t^{1/2},\log(t)\}) \otimes \Gamma^\mathcal{O}_{\text{deRham},X,v}\{t^{1/2},\log(t)\}\\ \overset{\sim}{\longrightarrow} F \otimes \Gamma^\mathcal{O}_{\text{deRham},X,v}\{t^{1/2},\log(t)\}. 
\end{align}
\end{definition}

\noindent We now define the $(\infty,1)$-categories of mixed-parity de Rham modules and he corresponding mixed-parity almost de Rham modules by using the objects involved to generated these categories:

\begin{definition}
Considering all the mixed parity de Rham bundles (locally finite free) as defined above, we consider the sub-$(\infty,1)$ category of 
\begin{align}
\mathrm{preModule}^\mathrm{solid,quasicoherent}_{\square,\Gamma^\mathrm{perfect}_{\text{Robba},X,v,\infty}\{t^{1/2}\}},
\mathrm{preModule}^\mathrm{solid,quasicoherent}_{\square,\Gamma^\mathrm{perfect}_{\text{Robba},X,v,I}\{t^{1/2}\}} 
\end{align}
generated by the mixed-parity de Rham bundles (locally finite free ones). These are defined to be the $(\infty,1)$-categories of mixed-parity de Rham complexes:
\begin{align}
\mathrm{preModule}^\mathrm{solid,quasicoherent,mixed-paritydeRham}_{\square,\Gamma^\mathrm{perfect}_{\text{Robba},X,v,\infty}\{t^{1/2}\}},
\mathrm{preModule}^\mathrm{solid,quasicoherent,mixed-paritydeRham}_{\square,\Gamma^\mathrm{perfect}_{\text{Robba},X,v,I}\{t^{1/2}\}}. 
\end{align}
\end{definition}

\begin{definition}
Considering all the mixed parity almost de Rham bundles (locally finite free) as defined above, we consider the sub-$(\infty,1)$ category of 
\begin{align}
\mathrm{preModule}^\mathrm{solid,quasicoherent}_{\square,\Gamma^\mathrm{perfect}_{\text{Robba},X,v,\infty}\{t^{1/2}\}},
\mathrm{preModule}^\mathrm{solid,quasicoherent}_{\square,\Gamma^\mathrm{perfect}_{\text{Robba},X,v,I}\{t^{1/2}\}} 
\end{align}
generated by the mixed-parity almost de Rham bundles (locally finite free ones). These are defined to be the $(\infty,1)$-categories of mixed-parity de Rham complexes:
\begin{align}
\mathrm{preModule}^\mathrm{solid,quasicoherent,mixed-parityalmostdeRham}_{\square,\Gamma^\mathrm{perfect}_{\text{Robba},X,v,\infty}\{t^{1/2}\}},\\
\mathrm{preModule}^\mathrm{solid,quasicoherent,mixed-parityalmostdeRham}_{\square,\Gamma^\mathrm{perfect}_{\text{Robba},X,v,I}\{t^{1/2}\}}. 
\end{align}
\end{definition}

\indent Then the corresponding mixed-parity de Rham functors can be extended to these categories:
\begin{align}
\mathrm{preModule}^\mathrm{solid,quasicoherent,mixed-paritydeRham}_{\square,\Gamma^\mathrm{perfect}_{\text{Robba},X,v,\infty}\{t^{1/2}\}},\\
\mathrm{preModule}^\mathrm{solid,quasicoherent,mixed-paritydeRham}_{\square,\Gamma^\mathrm{perfect}_{\text{Robba},X,v,I}\{t^{1/2}\}}, 
\end{align}
and
\begin{align}
\mathrm{preModule}^\mathrm{solid,quasicoherent,mixed-parityalmostdeRham}_{\square,\Gamma^\mathrm{perfect}_{\text{Robba},X,v,\infty}\{t^{1/2}\}},\\
\mathrm{preModule}^\mathrm{solid,quasicoherent,mixed-parityalmostdeRham}_{\square,\Gamma^\mathrm{perfect}_{\text{Robba},X,v,I}\{t^{1/2}\}}. 
\end{align}

\subsubsection{Mixed-Parity Hodge Modules with Frobenius}

\noindent Now we consider the key objects in our study namely those complexes generated by certain mixed-parity Hodge modules. We start from the following definition.

\begin{remark}
All the coherent sheaves over mixed-parity Robba sheaves in this section will carry the corresponding Frobenius morphism $\varphi: F \overset{\sim}{\longrightarrow} \varphi^*F$.
\end{remark}

\begin{definition}
For any locally free coherent sheaf $F$ over
\begin{align}
\Gamma^\mathrm{perfect}_{\text{Robba},X,v,\infty}\{t^{1/2}\},\Gamma^\mathrm{perfect}_{\text{Robba},X,v,I}\{t^{1/2}\},
\end{align} 
we consider the following functor $\mathrm{dR}$ sending $F$ to the following object:
\begin{align}
f_*(F\otimes_{\Gamma^\mathrm{perfect}_{\text{Robba},X,v,\infty}\{t^{1/2}\}} \Gamma^\mathcal{O}_{\text{deRham},X,v}\{t^{1/2}\})
\end{align}
or 
\begin{align}
f_*(F\otimes_{\Gamma^\mathrm{perfect}_{\text{Robba},X,v,I}\{t^{1/2}\}} \Gamma^\mathcal{O}_{\text{deRham},X,v}\{t^{1/2}\}).
\end{align}
We call $F$ mixed-parity de Rham if we have the following isomorphism:
\begin{align}
f^*f_*(F\otimes_{\Gamma^\mathrm{perfect}_{\text{Robba},X,v,\infty}\{t^{1/2}\}} \Gamma^\mathcal{O}_{\text{deRham},X,v}\{t^{1/2}\}) \otimes \Gamma^\mathcal{O}_{\text{deRham},X,v}\{t^{1/2}\} \overset{\sim}{\longrightarrow} F \otimes \Gamma^\mathcal{O}_{\text{deRham},X,v}\{t^{1/2}\} 
\end{align}
or 
\begin{align}
f^*f_*(F\otimes_{\Gamma^\mathrm{perfect}_{\text{Robba},X,v,I}\{t^{1/2}\}} \Gamma^\mathcal{O}_{\text{deRham},X,v}\{t^{1/2}\}) \otimes \Gamma^\mathcal{O}_{\text{deRham},X,v}\{t^{1/2}\} \overset{\sim}{\longrightarrow} F \otimes \Gamma^\mathcal{O}_{\text{deRham},X,v}\{t^{1/2}\}. 
\end{align}
\end{definition}

\begin{definition}
For any locally free coherent sheaf $F$ over
\begin{align}
\Gamma^\mathrm{perfect}_{\text{Robba},X,v,\infty}\{t^{1/2}\},\Gamma^\mathrm{perfect}_{\text{Robba},X,v,I}\{t^{1/2}\},
\end{align} 
we consider the following functor $\mathrm{dR}^\mathrm{almost}$ sending $F$ to the following object:
\begin{align}
f_*(F\otimes_{\Gamma^\mathrm{perfect}_{\text{Robba},X,v,\infty}\{t^{1/2}\}} \Gamma^\mathcal{O}_{\text{deRham},X,v}\{t^{1/2},\log(t)\})
\end{align}
or 
\begin{align}
f_*(F\otimes_{\Gamma^\mathrm{perfect}_{\text{Robba},X,v,I}\{t^{1/2}\}} \Gamma^\mathcal{O}_{\text{deRham},X,v}\{t^{1/2},\log(t)\}).
\end{align}
We call $F$ mixed-parity almost de Rham if we have the following isomorphism:
\begin{align}
f^*f_*(F\otimes_{\Gamma^\mathrm{perfect}_{\text{Robba},X,v,\infty}\{t^{1/2}\}} \Gamma^\mathcal{O}_{\text{deRham},X,v}\{t^{1/2},\log(t)\}) \otimes \Gamma^\mathcal{O}_{\text{deRham},X,v}\{t^{1/2},\log(t)\} \\
\overset{\sim}{\longrightarrow} F \otimes \Gamma^\mathcal{O}_{\text{deRham},X,v}\{t^{1/2},\log(t)\} 
\end{align}
or 
\begin{align}
f^*f_*(F\otimes_{\Gamma^\mathrm{perfect}_{\text{Robba},X,v,I}\{t^{1/2}\}} \Gamma^\mathcal{O}_{\text{deRham},X,v}\{t^{1/2},\log(t)\}) \otimes \Gamma^\mathcal{O}_{\text{deRham},X,v}\{t^{1/2},\log(t)\}\\ \overset{\sim}{\longrightarrow} F \otimes \Gamma^\mathcal{O}_{\text{deRham},X,v}\{t^{1/2},\log(t)\}. 
\end{align}
\end{definition}

\noindent We now define the $(\infty,1)$-categories of mixed-parity de Rham modules and he corresponding mixed-parity almost de Rham modules by using the objects involved to generated these categories:

\begin{definition}
Considering all the mixed parity de Rham bundles (locally finite free) as defined above, we consider the sub-$(\infty,1)$ category of 
\begin{align}
\varphi\mathrm{preModule}^\mathrm{solid,quasicoherent}_{\square,\Gamma^\mathrm{perfect}_{\text{Robba},X,v,\infty}\{t^{1/2}\}},
\varphi\mathrm{preModule}^\mathrm{solid,quasicoherent}_{\square,\Gamma^\mathrm{perfect}_{\text{Robba},X,v,I}\{t^{1/2}\}} 
\end{align}
generated by the mixed-parity de Rham bundles (locally finite free ones). These are defined to be the $(\infty,1)$-categories of mixed-parity de Rham complexes:
\begin{align}
\varphi\mathrm{preModule}^\mathrm{solid,quasicoherent,mixed-paritydeRham}_{\square,\Gamma^\mathrm{perfect}_{\text{Robba},X,v,\infty}\{t^{1/2}\}},
\varphi\mathrm{preModule}^\mathrm{solid,quasicoherent,mixed-paritydeRham}_{\square,\Gamma^\mathrm{perfect}_{\text{Robba},X,v,I}\{t^{1/2}\}}. 
\end{align}
\end{definition}

\begin{definition}
Considering all the mixed parity almost de Rham bundles (locally finite free) as defined above, we consider the sub-$(\infty,1)$ category of 
\begin{align}
\varphi\mathrm{preModule}^\mathrm{solid,quasicoherent}_{\square,\Gamma^\mathrm{perfect}_{\text{Robba},X,v,\infty}\{t^{1/2}\}},
\varphi\mathrm{preModule}^\mathrm{solid,quasicoherent}_{\square,\Gamma^\mathrm{perfect}_{\text{Robba},X,v,I}\{t^{1/2}\}} 
\end{align}
generated by the mixed-parity almost de Rham bundles (locally finite free ones). These are defined to be the $(\infty,1)$-categories of mixed-parity de Rham complexes:
\begin{align}
\varphi\mathrm{preModule}^\mathrm{solid,quasicoherent,mixed-parityalmostdeRham}_{\square,\Gamma^\mathrm{perfect}_{\text{Robba},X,v,\infty}\{t^{1/2}\}},\\
\varphi\mathrm{preModule}^\mathrm{solid,quasicoherent,mixed-parityalmostdeRham}_{\square,\Gamma^\mathrm{perfect}_{\text{Robba},X,v,I}\{t^{1/2}\}}. 
\end{align}
\end{definition}

\indent Then the corresponding mixed-parity de Rham functors can be extended to these categories:
\begin{align}
\varphi\mathrm{preModule}^\mathrm{solid,quasicoherent,mixed-paritydeRham}_{\square,\Gamma^\mathrm{perfect}_{\text{Robba},X,v,\infty}\{t^{1/2}\}},
\varphi\mathrm{preModule}^\mathrm{solid,quasicoherent,mixed-paritydeRham}_{\square,\Gamma^\mathrm{perfect}_{\text{Robba},X,v,I}\{t^{1/2}\}}, 
\end{align}
and
\begin{align}
\varphi\mathrm{preModule}^\mathrm{solid,quasicoherent,mixed-parityalmostdeRham}_{\square,\Gamma^\mathrm{perfect}_{\text{Robba},X,v,\infty}\{t^{1/2}\}},\\
\varphi\mathrm{preModule}^\mathrm{solid,quasicoherent,mixed-parityalmostdeRham}_{\square,\Gamma^\mathrm{perfect}_{\text{Robba},X,v,I}\{t^{1/2}\}}. 
\end{align}

\subsection{Mixed-Parity de Rham Riemann-Hilbert Correspondence}

\indent This chapter will extend the corresponding Riemann-Hilbert correspondence from \cite{Sch1}, \cite{LZ}, \cite{BL1}, \cite{BL2}, \cite{M} to the mixed-parity setting.

\begin{definition}
We define the following Riemann-Hilbert functor $\text{RH}_\text{mixed-parity}$ from the one of categories:
\begin{align}
\mathrm{preModule}^\mathrm{solid,quasicoherent,mixed-paritydeRham}_{\square,\Gamma^\mathrm{perfect}_{\text{Robba},X,v,\infty}\{t^{1/2}\}},
\mathrm{preModule}^\mathrm{solid,quasicoherent,mixed-paritydeRham}_{\square,\Gamma^\mathrm{perfect}_{\text{Robba},X,v,I}\{t^{1/2}\}}, 
\end{align}
and
\begin{align}
\mathrm{preModule}^\mathrm{solid,quasicoherent,mixed-parityalmostdeRham}_{\square,\Gamma^\mathrm{perfect}_{\text{Robba},X,v,\infty}\{t^{1/2}\}},\\
\mathrm{preModule}^\mathrm{solid,quasicoherent,mixed-parityalmostdeRham}_{\square,\Gamma^\mathrm{perfect}_{\text{Robba},X,v,I}\{t^{1/2}\}} 
\end{align}
to $(\infty,1)$-categories in image denoted by:
\begin{align}
\mathrm{preModule}_{X,\text{\'et}}
\end{align}
to be the following functors sending each $F$ in the domain to:
\begin{align}
&\text{RH}_\text{mixed-parity}(F):=f_*(F\otimes_{\Gamma^\mathrm{perfect}_{\text{Robba},X,v,\infty}\{t^{1/2}\}} \Gamma^\mathcal{O}_{\text{deRham},X,v}\{t^{1/2}\}),\\
&\text{RH}_\text{mixed-parity}(F):=f_*(F\otimes_{\Gamma^\mathrm{perfect}_{\text{Robba},X,v,I}\{t^{1/2}\}} \Gamma^\mathcal{O}_{\text{deRham},X,v}\{t^{1/2}\}),\\
&\text{RH}_\text{mixed-parity}(F):=f_*(F\otimes_{\Gamma^\mathrm{perfect}_{\text{Robba},X,v,\infty}\{t^{1/2}\}} \Gamma^\mathcal{O}_{\text{deRham},X,v}\{t^{1/2},\log(t)\}),\\
&\text{RH}_\text{mixed-parity}(F):=f_*(F\otimes_{\Gamma^\mathrm{perfect}_{\text{Robba},X,v,I}\{t^{1/2}\}} \Gamma^\mathcal{O}_{\text{deRham},X,v}\{t^{1/2},\log(t)\}),\\
\end{align}
respectively.

\end{definition}

\begin{definition}
In the situation where we have the Frobenius action we consider the follwing. We define the following Riemann-Hilbert functor $\text{RH}_\text{mixed-parity}$ from the one of categories:
\begin{align}
\varphi\mathrm{preModule}^\mathrm{solid,quasicoherent,mixed-paritydeRham}_{\square,\Gamma^\mathrm{perfect}_{\text{Robba},X,v,\infty}\{t^{1/2}\}},
\varphi\mathrm{preModule}^\mathrm{solid,quasicoherent,mixed-paritydeRham}_{\square,\Gamma^\mathrm{perfect}_{\text{Robba},X,v,I}\{t^{1/2}\}}, 
\end{align}
and
\begin{align}
\varphi\mathrm{preModule}^\mathrm{solid,quasicoherent,mixed-parityalmostdeRham}_{\square,\Gamma^\mathrm{perfect}_{\text{Robba},X,v,\infty}\{t^{1/2}\}},\\
\varphi\mathrm{preModule}^\mathrm{solid,quasicoherent,mixed-parityalmostdeRham}_{\square,\Gamma^\mathrm{perfect}_{\text{Robba},X,v,I}\{t^{1/2}\}} 
\end{align}
to $(\infty,1)$-categories in image denoted by:
\begin{align}
\mathrm{preModule}_{X,\text{\'et}}
\end{align}
to be the following functors sending each $F$ in the domain to:
\begin{align}
&\text{RH}_\text{mixed-parity}(F):=f_*(F\otimes_{\Gamma^\mathrm{perfect}_{\text{Robba},X,v,\infty}\{t^{1/2}\}} \Gamma^\mathcal{O}_{\text{deRham},X,v}\{t^{1/2}\}),\\
&\text{RH}_\text{mixed-parity}(F):=f_*(F\otimes_{\Gamma^\mathrm{perfect}_{\text{Robba},X,v,I}\{t^{1/2}\}} \Gamma^\mathcal{O}_{\text{deRham},X,v}\{t^{1/2}\}),\\
&\text{RH}_\text{mixed-parity}(F):=f_*(F\otimes_{\Gamma^\mathrm{perfect}_{\text{Robba},X,v,\infty}\{t^{1/2}\}} \Gamma^\mathcal{O}_{\text{deRham},X,v}\{t^{1/2},\log(t)\}),\\
&\text{RH}_\text{mixed-parity}(F):=f_*(F\otimes_{\Gamma^\mathrm{perfect}_{\text{Robba},X,v,I}\{t^{1/2}\}} \Gamma^\mathcal{O}_{\text{deRham},X,v}\{t^{1/2},\log(t)\}),\\
\end{align}
respectively.

\end{definition}

\newpage
\section{Geometric Family of Mixed-Parity Hodge Modules II: Cristalline Situations}

\noindent References: \cite{Sch1}, \cite{Sch2}, \cite{FS}, \cite{KL1}, \cite{KL2}, \cite{BL1}, \cite{BL2}, \cite{BS}, \cite{BHS}, \cite{Fon1}, \cite{CS1}, \cite{CS2}, \cite{BK}, \cite{BBK}, \cite{BBBK}, \cite{KKM}, \cite{KM}, \cite{LZ}, \cite{TT}, \cite{M}.

\subsection{Period Rings and Sheaves}

\subsubsection{Rings}

\noindent Let $X$ be a rigid analytic space over $\mathbb{Q}_p$. We have the corresponding \'etale site and the corresponding pro-\'etale site of $X$, which we denote them by $X_{v},X_\text{\'et}$. The relationship of the two sites can be reflected by the corresponding morphism $f:X_{v}\longrightarrow X_\text{\'et}$. Then we have the corresponding cristalline period rings and sheaves from \cite{TT}:
\begin{align}
\Gamma_{\text{cristalline},X,v}, \Gamma^\mathcal{O}_{\text{cristalline},X,v}.
\end{align}
Our notations are different from \cite{TT}, we use $\Gamma$ to mean $B$ in \cite{TT}, while $\Gamma^\mathcal{O}$ will be the corresponding $OB$ ring in \cite{TT}.\\

\begin{definition}
\indent Now we assume that $p>2$, following \cite{BS} we join the square root of $t$ element in $\Gamma_{\text{cristalline},X,v}$ which forms the sheaves:
\begin{align}
\Gamma_{\text{cristalline},X,v}\{t^{1/2}\},\Gamma^\mathcal{O}_{\text{cristalline},X,v}\{t^{1/2}\}.
\end{align}
And following \cite{BL1}, \cite{BL2}, \cite{Fon1}, \cite{BHS} we further have the following sheaves of rings:
\begin{align}
\Gamma_{\text{cristalline},X,v}\{t^{1/2},\log(t)\},\Gamma^\mathcal{O}_{\text{cristalline},X,v}\{t^{1/2},\log(t)\}.
\end{align}
\end{definition}

\begin{definition}
We use the notations:
\begin{align}
\Gamma^\mathrm{perfect}_{\text{Robba},X,v},\Gamma^\mathrm{perfect}_{\text{Robba},X,v,\infty},\Gamma^\mathrm{perfect}_{\text{Robba},X,v,I}
\end{align}
to denote the perfect Robba rings from \cite{KL1}, \cite{KL2}, where $I\subset (0,\infty)$. Then we join $t^{1/2}$ to these sheaves we have:
\begin{align}
\Gamma^\mathrm{perfect}_{\text{Robba},X,v}\{t^{1/2}\},\Gamma^\mathrm{perfect}_{\text{Robba},X,v,\infty}\{t^{1/2}\},\Gamma^\mathrm{perfect}_{\text{Robba},X,v,I}\{t^{1/2}\}.
\end{align}
And following \cite{BL1}, \cite{BL2}, \cite{Fon1}, \cite{BHS} we have the following larger sheaves:
\begin{align}
\Gamma^\mathrm{perfect}_{\text{Robba},X,v}\{t^{1/2},\log(t)\},\Gamma^\mathrm{perfect}_{\text{Robba},X,v,\infty}\{t^{1/2},\log(t)\},\Gamma^\mathrm{perfect}_{\text{Robba},X,v,I}\{t^{1/2},\log(t)\}.
\end{align} 
\end{definition}

\begin{definition}
From now on, we use the same notation to denote the period rings involved tensored with a finite extension of $\mathbb{Q}_p$ containing square root of $p$ as in \cite{BS}.
\begin{align}
\Gamma_{\text{cristalline},X,v}\{t^{1/2}\},\Gamma^\mathcal{O}_{\text{cristalline},X,v}\{t^{1/2}\}.
\end{align}
\begin{align}
\Gamma_{\text{cristalline},X,v}\{t^{1/2},\log(t)\},\Gamma^\mathcal{O}_{\text{cristalline},X,v}\{t^{1/2},\log(t)\}.
\end{align}
\begin{align}
\Gamma^\mathrm{perfect}_{\text{Robba},X,v}\{t^{1/2}\},\Gamma^\mathrm{perfect}_{\text{Robba},X,v,\infty}\{t^{1/2}\},\Gamma^\mathrm{perfect}_{\text{Robba},X,v,I}\{t^{1/2}\}.
\end{align}
\begin{align}
\Gamma^\mathrm{perfect}_{\text{Robba},X,v}\{t^{1/2},\log(t)\},\Gamma^\mathrm{perfect}_{\text{Robba},X,v,\infty}\{t^{1/2},\log(t)\},\Gamma^\mathrm{perfect}_{\text{Robba},X,v,I}\{t^{1/2},\log(t)\}.
\end{align}
This is necessary since we to extend the action of $\varphi$ to the period rings by $\varphi(t^{1/2}\otimes 1)=\varphi(t)^{1/2}\otimes 1$.
\end{definition}

\subsubsection{Modules}

\noindent We consider quasicoherent presheaves in the following two situation:
\begin{itemize}
\item[$\square$] The solid quasicoherent modules from \cite{CS1}, \cite{CS2};
\item[$\square$] The ind-Banach quasicoherent modules from \cite{BK}, \cite{BBK}, \cite{BBBK}, \cite{KKM}, \cite{KM} with the corresponding monomorphic ind-Banach quasicoherent modules from \cite{BK}, \cite{BBK}, \cite{BBBK}, \cite{KKM}, \cite{KM}.
\end{itemize}

\begin{definition}
We use the notation:
\begin{align}
\mathrm{preModule}^\mathrm{solid,quasicoherent}_{\square,\Gamma^\mathrm{perfect}_{\text{Robba},X,v}\{t^{1/2}\}},\mathrm{preModule}^\mathrm{solid,quasicoherent}_{\square,\Gamma^\mathrm{perfect}_{\text{Robba},X,v,\infty}\{t^{1/2}\}},
\mathrm{preModule}^\mathrm{solid,quasicoherent}_{\square,\Gamma^\mathrm{perfect}_{\text{Robba},X,v,I}\{t^{1/2}\}} 
\end{align}
to denote the $(\infty,1)$-categories of solid quasicoherent presheaves over the corresonding Robba sheaves. Locally the section is defined by taking the corresponding $(\infty,1)$-categories of solid modules.
\end{definition}

\begin{definition}
We use the notation:
\begin{align}
\mathrm{preModule}^\mathrm{ind-Banach,quasicoherent}_{\Gamma^\mathrm{perfect}_{\text{Robba},X,v}\{t^{1/2}\}},\\\mathrm{preModule}^\mathrm{ind-Banach,quasicoherent}_{\Gamma^\mathrm{perfect}_{\text{Robba},X,v,\infty}\{t^{1/2}\}},\\
\mathrm{preModule}^\mathrm{ind-Banach,quasicoherent}_{\Gamma^\mathrm{perfect}_{\text{Robba},X,v,I}\{t^{1/2}\}} 
\end{align}
to denote the $(\infty,1)$-categories of solid quasicoherent presheaves over the corresonding Robba sheaves. Locally the section is defined by taking the corresponding $(\infty,1)$-categories of inductive Banach  modules. 
\end{definition}

\begin{definition}
We use the notation:
\begin{align}
\mathrm{Module}^\mathrm{solid,quasicoherent}_{\square,\Gamma^\mathrm{perfect}_{\text{Robba},X,v}\{t^{1/2}\}},\mathrm{Module}^\mathrm{solid,quasicoherent}_{\square,\Gamma^\mathrm{perfect}_{\text{Robba},X,v,\infty}\{t^{1/2}\}},
\mathrm{Module}^\mathrm{solid,quasicoherent}_{\square,\Gamma^\mathrm{perfect}_{\text{Robba},X,v,I}\{t^{1/2}\}} 
\end{align}
to denote the $(\infty,1)$-categories of solid quasicoherent sheaves over the corresonding Robba sheaves. Locally the section is defined by taking the corresponding $(\infty,1)$-categories of solid modules.
\end{definition}

\subsubsection{Mixed-Parity Hodge Modules without Frobenius}

\noindent Now we consider the key objects in our study namely those complexes generated by certain mixed-parity Hodge modules. We start from the following definition.

\begin{definition}
For any locally free coherent sheaf $F$ over
\begin{align}
\Gamma^\mathrm{perfect}_{\text{Robba},X,v,\infty}\{t^{1/2}\},\Gamma^\mathrm{perfect}_{\text{Robba},X,v,I}\{t^{1/2}\},
\end{align} 
we consider the following functor $\mathrm{dR}$ sending $F$ to the following object:
\begin{align}
f_*(F\otimes_{\Gamma^\mathrm{perfect}_{\text{Robba},X,v,\infty}\{t^{1/2}\}} \Gamma^\mathcal{O}_{\text{cristalline},X,v}\{t^{1/2}\})
\end{align}
or 
\begin{align}
f_*(F\otimes_{\Gamma^\mathrm{perfect}_{\text{Robba},X,v,I}\{t^{1/2}\}} \Gamma^\mathcal{O}_{\text{cristalline},X,v}\{t^{1/2}\}).
\end{align}
We call $F$ mixed-parity cristalline if we have the following isomorphism:
\begin{align}
f^*f_*(F\otimes_{\Gamma^\mathrm{perfect}_{\text{Robba},X,v,\infty}\{t^{1/2}\}} \Gamma^\mathcal{O}_{\text{cristalline},X,v}\{t^{1/2}\}) \otimes \Gamma^\mathcal{O}_{\text{cristalline},X,v}\{t^{1/2}\} \overset{\sim}{\longrightarrow} F \otimes \Gamma^\mathcal{O}_{\text{cristalline},X,v}\{t^{1/2}\} 
\end{align}
or 
\begin{align}
f^*f_*(F\otimes_{\Gamma^\mathrm{perfect}_{\text{Robba},X,v,I}\{t^{1/2}\}} \Gamma^\mathcal{O}_{\text{cristalline},X,v}\{t^{1/2}\}) \otimes \Gamma^\mathcal{O}_{\text{cristalline},X,v}\{t^{1/2}\} \overset{\sim}{\longrightarrow} F \otimes \Gamma^\mathcal{O}_{\text{cristalline},X,v}\{t^{1/2}\}. 
\end{align}
\end{definition}

\begin{definition}
For any locally free coherent sheaf $F$ over
\begin{align}
\Gamma^\mathrm{perfect}_{\text{Robba},X,v,\infty}\{t^{1/2}\},\Gamma^\mathrm{perfect}_{\text{Robba},X,v,I}\{t^{1/2}\},
\end{align} 
we consider the following functor $\mathrm{dR}^\mathrm{almost}$ sending $F$ to the following object:
\begin{align}
f_*(F\otimes_{\Gamma^\mathrm{perfect}_{\text{Robba},X,v,\infty}\{t^{1/2}\}} \Gamma^\mathcal{O}_{\text{cristalline},X,v}\{t^{1/2},\log(t)\})
\end{align}
or 
\begin{align}
f_*(F\otimes_{\Gamma^\mathrm{perfect}_{\text{Robba},X,v,I}\{t^{1/2}\}} \Gamma^\mathcal{O}_{\text{cristalline},X,v}\{t^{1/2},\log(t)\}).
\end{align}
We call $F$ mixed-parity almost cristalline if we have the following isomorphism:
\begin{align}
f^*f_*(F\otimes_{\Gamma^\mathrm{perfect}_{\text{Robba},X,v,\infty}\{t^{1/2}\}} \Gamma^\mathcal{O}_{\text{cristalline},X,v}\{t^{1/2},\log(t)\}) \otimes \Gamma^\mathcal{O}_{\text{cristalline},X,v}\{t^{1/2},\log(t)\} \\
\overset{\sim}{\longrightarrow} F \otimes \Gamma^\mathcal{O}_{\text{cristalline},X,v}\{t^{1/2},\log(t)\} 
\end{align}
or 
\begin{align}
f^*f_*(F\otimes_{\Gamma^\mathrm{perfect}_{\text{Robba},X,v,I}\{t^{1/2}\}} \Gamma^\mathcal{O}_{\text{cristalline},X,v}\{t^{1/2},\log(t)\}) \otimes \Gamma^\mathcal{O}_{\text{cristalline},X,v}\{t^{1/2},\log(t)\}\\ \overset{\sim}{\longrightarrow} F \otimes \Gamma^\mathcal{O}_{\text{cristalline},X,v}\{t^{1/2},\log(t)\}. 
\end{align}
\end{definition}

\noindent We now define the $(\infty,1)$-categories of mixed-parity cristalline modules and he corresponding mixed-parity almost cristalline modules by using the objects involved to generated these categories:

\begin{definition}
Considering all the mixed parity cristalline bundles (locally finite free) as defined above, we consider the sub-$(\infty,1)$ category of 
\begin{align}
\mathrm{preModule}^\mathrm{solid,quasicoherent}_{\square,\Gamma^\mathrm{perfect}_{\text{Robba},X,v,\infty}\{t^{1/2}\}},
\mathrm{preModule}^\mathrm{solid,quasicoherent}_{\square,\Gamma^\mathrm{perfect}_{\text{Robba},X,v,I}\{t^{1/2}\}} 
\end{align}
generated by the mixed-parity cristalline bundles (locally finite free ones). These are defined to be the $(\infty,1)$-categories of mixed-parity cristalline complexes:
\begin{align}
\mathrm{preModule}^\mathrm{solid,quasicoherent,mixed-paritycristalline}_{\square,\Gamma^\mathrm{perfect}_{\text{Robba},X,v,\infty}\{t^{1/2}\}},
\mathrm{preModule}^\mathrm{solid,quasicoherent,mixed-paritycristalline}_{\square,\Gamma^\mathrm{perfect}_{\text{Robba},X,v,I}\{t^{1/2}\}}. 
\end{align}
\end{definition}

\begin{definition}
Considering all the mixed parity almost cristalline bundles (locally finite free) as defined above, we consider the sub-$(\infty,1)$ category of 
\begin{align}
\mathrm{preModule}^\mathrm{solid,quasicoherent}_{\square,\Gamma^\mathrm{perfect}_{\text{Robba},X,v,\infty}\{t^{1/2}\}},
\mathrm{preModule}^\mathrm{solid,quasicoherent}_{\square,\Gamma^\mathrm{perfect}_{\text{Robba},X,v,I}\{t^{1/2}\}} 
\end{align}
generated by the mixed-parity almost cristalline bundles (locally finite free ones). These are defined to be the $(\infty,1)$-categories of mixed-parity cristalline complexes:
\begin{align}
\mathrm{preModule}^\mathrm{solid,quasicoherent,mixed-parityalmostcristalline}_{\square,\Gamma^\mathrm{perfect}_{\text{Robba},X,v,\infty}\{t^{1/2}\}},\\
\mathrm{preModule}^\mathrm{solid,quasicoherent,mixed-parityalmostcristalline}_{\square,\Gamma^\mathrm{perfect}_{\text{Robba},X,v,I}\{t^{1/2}\}}. 
\end{align}
\end{definition}

\indent Then the corresponding mixed-parity cristalline functors can be extended to these categories:
\begin{align}
\mathrm{preModule}^\mathrm{solid,quasicoherent,mixed-paritycristalline}_{\square,\Gamma^\mathrm{perfect}_{\text{Robba},X,v,\infty}\{t^{1/2}\}},\\
\mathrm{preModule}^\mathrm{solid,quasicoherent,mixed-paritycristalline}_{\square,\Gamma^\mathrm{perfect}_{\text{Robba},X,v,I}\{t^{1/2}\}}, 
\end{align}
and
\begin{align}
\mathrm{preModule}^\mathrm{solid,quasicoherent,mixed-parityalmostcristalline}_{\square,\Gamma^\mathrm{perfect}_{\text{Robba},X,v,\infty}\{t^{1/2}\}},\\
\mathrm{preModule}^\mathrm{solid,quasicoherent,mixed-parityalmostcristalline}_{\square,\Gamma^\mathrm{perfect}_{\text{Robba},X,v,I}\{t^{1/2}\}}. 
\end{align}

\subsubsection{Mixed-Parity Hodge Modules with Frobenius}

\noindent Now we consider the key objects in our study namely those complexes generated by certain mixed-parity Hodge modules. We start from the following definition.

\begin{remark}
All the coherent sheaves over mixed-parity Robba sheaves in this section will carry the corresponding Frobenius morphism $\varphi: F \overset{\sim}{\longrightarrow} \varphi^*F$.
\end{remark}

\begin{definition}
For any locally free coherent sheaf $F$ over
\begin{align}
\Gamma^\mathrm{perfect}_{\text{Robba},X,v,\infty}\{t^{1/2}\},\Gamma^\mathrm{perfect}_{\text{Robba},X,v,I}\{t^{1/2}\},
\end{align} 
we consider the following functor $\mathrm{dR}$ sending $F$ to the following object:
\begin{align}
f_*(F\otimes_{\Gamma^\mathrm{perfect}_{\text{Robba},X,v,\infty}\{t^{1/2}\}} \Gamma^\mathcal{O}_{\text{cristalline},X,v}\{t^{1/2}\})
\end{align}
or 
\begin{align}
f_*(F\otimes_{\Gamma^\mathrm{perfect}_{\text{Robba},X,v,I}\{t^{1/2}\}} \Gamma^\mathcal{O}_{\text{cristalline},X,v}\{t^{1/2}\}).
\end{align}
We call $F$ mixed-parity cristalline if we have the following isomorphism:
\begin{align}
f^*f_*(F\otimes_{\Gamma^\mathrm{perfect}_{\text{Robba},X,v,\infty}\{t^{1/2}\}} \Gamma^\mathcal{O}_{\text{cristalline},X,v}\{t^{1/2}\}) \otimes \Gamma^\mathcal{O}_{\text{cristalline},X,v}\{t^{1/2}\} \overset{\sim}{\longrightarrow} F \otimes \Gamma^\mathcal{O}_{\text{cristalline},X,v}\{t^{1/2}\} 
\end{align}
or 
\begin{align}
f^*f_*(F\otimes_{\Gamma^\mathrm{perfect}_{\text{Robba},X,v,\infty}\{t^{1/2}\}} \Gamma^\mathcal{O}_{\text{cristalline},X,v}\{t^{1/2}\}) \otimes \Gamma^\mathcal{O}_{\text{cristalline},X,v}\{t^{1/2}\} \overset{\sim}{\longrightarrow} F \otimes \Gamma^\mathcal{O}_{\text{cristalline},X,v}\{t^{1/2}\}. 
\end{align}
\end{definition}

\begin{definition}
For any locally free coherent sheaf $F$ over
\begin{align}
\Gamma^\mathrm{perfect}_{\text{Robba},X,v,\infty}\{t^{1/2}\},\Gamma^\mathrm{perfect}_{\text{Robba},X,v,I}\{t^{1/2}\},
\end{align} 
we consider the following functor $\mathrm{dR}^\mathrm{almost}$ sending $F$ to the following object:
\begin{align}
f_*(F\otimes_{\Gamma^\mathrm{perfect}_{\text{Robba},X,v,\infty}\{t^{1/2}\}} \Gamma^\mathcal{O}_{\text{cristalline},X,v}\{t^{1/2},\log(t)\})
\end{align}
or 
\begin{align}
f_*(F\otimes_{\Gamma^\mathrm{perfect}_{\text{Robba},X,v,I}\{t^{1/2}\}} \Gamma^\mathcal{O}_{\text{cristalline},X,v}\{t^{1/2},\log(t)\}).
\end{align}
We call $F$ mixed-parity almost cristalline if we have the following isomorphism:
\begin{align}
f^*f_*(F\otimes_{\Gamma^\mathrm{perfect}_{\text{Robba},X,v,\infty}\{t^{1/2}\}} \Gamma^\mathcal{O}_{\text{cristalline},X,v}\{t^{1/2},\log(t)\}) \otimes \Gamma^\mathcal{O}_{\text{cristalline},X,v}\{t^{1/2},\log(t)\} \\
\overset{\sim}{\longrightarrow} F \otimes \Gamma^\mathcal{O}_{\text{cristalline},X,v}\{t^{1/2},\log(t)\} 
\end{align}
or 
\begin{align}
f^*f_*(F\otimes_{\Gamma^\mathrm{perfect}_{\text{Robba},X,v,I}\{t^{1/2}\}} \Gamma^\mathcal{O}_{\text{cristalline},X,v}\{t^{1/2},\log(t)\}) \otimes \Gamma^\mathcal{O}_{\text{cristalline},X,v}\{t^{1/2},\log(t)\}\\ \overset{\sim}{\longrightarrow} F \otimes \Gamma^\mathcal{O}_{\text{cristalline},X,v}\{t^{1/2},\log(t)\}. 
\end{align}
\end{definition}

\noindent We now define the $(\infty,1)$-categories of mixed-parity cristalline modules and he corresponding mixed-parity almost cristalline modules by using the objects involved to generated these categories:

\begin{definition}
Considering all the mixed parity cristalline bundles (locally finite free) as defined above, we consider the sub-$(\infty,1)$ category of 
\begin{align}
\varphi\mathrm{preModule}^\mathrm{solid,quasicoherent}_{\square,\Gamma^\mathrm{perfect}_{\text{Robba},X,v,\infty}\{t^{1/2}\}},
\varphi\mathrm{preModule}^\mathrm{solid,quasicoherent}_{\square,\Gamma^\mathrm{perfect}_{\text{Robba},X,v,I}\{t^{1/2}\}} 
\end{align}
generated by the mixed-parity cristalline bundles (locally finite free ones). These are defined to be the $(\infty,1)$-categories of mixed-parity cristalline complexes:
\begin{align}
\varphi\mathrm{preModule}^\mathrm{solid,quasicoherent,mixed-paritycristalline}_{\square,\Gamma^\mathrm{perfect}_{\text{Robba},X,v,\infty}\{t^{1/2}\}},
\varphi\mathrm{preModule}^\mathrm{solid,quasicoherent,mixed-paritycristalline}_{\square,\Gamma^\mathrm{perfect}_{\text{Robba},X,v,I}\{t^{1/2}\}}. 
\end{align}
\end{definition}

\begin{definition}
Considering all the mixed parity almost cristalline bundles (locally finite free) as defined above, we consider the sub-$(\infty,1)$ category of 
\begin{align}
\varphi\mathrm{preModule}^\mathrm{solid,quasicoherent}_{\square,\Gamma^\mathrm{perfect}_{\text{Robba},X,v,\infty}\{t^{1/2}\}},
\varphi\mathrm{preModule}^\mathrm{solid,quasicoherent}_{\square,\Gamma^\mathrm{perfect}_{\text{Robba},X,v,I}\{t^{1/2}\}} 
\end{align}
generated by the mixed-parity almost cristalline bundles (locally finite free ones). These are defined to be the $(\infty,1)$-categories of mixed-parity cristalline complexes:
\begin{align}
\varphi\mathrm{preModule}^\mathrm{solid,quasicoherent,mixed-parityalmostcristalline}_{\square,\Gamma^\mathrm{perfect}_{\text{Robba},X,v,\infty}\{t^{1/2}\}},\\
\varphi\mathrm{preModule}^\mathrm{solid,quasicoherent,mixed-parityalmostcristalline}_{\square,\Gamma^\mathrm{perfect}_{\text{Robba},X,v,I}\{t^{1/2}\}}. 
\end{align}
\end{definition}

\indent Then the corresponding mixed-parity cristalline functors can be extended to these categories:
\begin{align}
\varphi\mathrm{preModule}^\mathrm{solid,quasicoherent,mixed-paritycristalline}_{\square,\Gamma^\mathrm{perfect}_{\text{Robba},X,v,\infty}\{t^{1/2}\}},
\varphi\mathrm{preModule}^\mathrm{solid,quasicoherent,mixed-paritycristalline}_{\square,\Gamma^\mathrm{perfect}_{\text{Robba},X,v,I}\{t^{1/2}\}}, 
\end{align}
and
\begin{align}
\varphi\mathrm{preModule}^\mathrm{solid,quasicoherent,mixed-parityalmostcristalline}_{\square,\Gamma^\mathrm{perfect}_{\text{Robba},X,v,\infty}\{t^{1/2}\}},\\
\varphi\mathrm{preModule}^\mathrm{solid,quasicoherent,mixed-parityalmostcristalline}_{\square,\Gamma^\mathrm{perfect}_{\text{Robba},X,v,I}\{t^{1/2}\}}. 
\end{align}

\subsection{Mixed-Parity Cristalline Riemann-Hilbert Correspondence}

\indent This chapter will extend the corresponding Riemann-Hilbert correspondence from \cite{Sch1}, \cite{LZ}, \cite{BL1}, \cite{BL2}, \cite{M} to the mixed-parity setting.

\begin{definition}
We define the following Riemann-Hilbert functor $\text{RH}_\text{mixed-parity}$ from the one of categories:
\begin{align}
\mathrm{preModule}^\mathrm{solid,quasicoherent,mixed-paritycristalline}_{\square,\Gamma^\mathrm{perfect}_{\text{Robba},X,v,\infty}\{t^{1/2}\}},
\mathrm{preModule}^\mathrm{solid,quasicoherent,mixed-paritycristalline}_{\square,\Gamma^\mathrm{perfect}_{\text{Robba},X,v,I}\{t^{1/2}\}}, 
\end{align}
and
\begin{align}
\mathrm{preModule}^\mathrm{solid,quasicoherent,mixed-parityalmostcristalline}_{\square,\Gamma^\mathrm{perfect}_{\text{Robba},X,v,\infty}\{t^{1/2}\}},\\
\mathrm{preModule}^\mathrm{solid,quasicoherent,mixed-parityalmostcristalline}_{\square,\Gamma^\mathrm{perfect}_{\text{Robba},X,v,I}\{t^{1/2}\}} 
\end{align}
to $(\infty,1)$-categories in image denoted by:
\begin{align}
\mathrm{preModule}_{X,\text{\'et}}
\end{align}
to be the following functors sending each $F$ in the domain to:
\begin{align}
&\text{RH}_\text{mixed-parity}(F):=f_*(F\otimes_{\Gamma^\mathrm{perfect}_{\text{Robba},X,v,\infty}\{t^{1/2}\}} \Gamma^\mathcal{O}_{\text{cristalline},X,v}\{t^{1/2}\}),\\
&\text{RH}_\text{mixed-parity}(F):=f_*(F\otimes_{\Gamma^\mathrm{perfect}_{\text{Robba},X,v,I}\{t^{1/2}\}} \Gamma^\mathcal{O}_{\text{cristalline},X,v}\{t^{1/2}\}),\\
&\text{RH}_\text{mixed-parity}(F):=f_*(F\otimes_{\Gamma^\mathrm{perfect}_{\text{Robba},X,v,\infty}\{t^{1/2}\}} \Gamma^\mathcal{O}_{\text{cristalline},X,v}\{t^{1/2},\log(t)\}),\\
&\text{RH}_\text{mixed-parity}(F):=f_*(F\otimes_{\Gamma^\mathrm{perfect}_{\text{Robba},X,v,I}\{t^{1/2}\}} \Gamma^\mathcal{O}_{\text{cristalline},X,v}\{t^{1/2},\log(t)\}),\\
\end{align}
respectively.

\end{definition}

\begin{definition}
In the situation where we have the Frobenius action we consider the follwing. We define the following Riemann-Hilbert functor $\text{RH}_\text{mixed-parity}$ from the one of categories:
\begin{align}
\varphi\mathrm{preModule}^\mathrm{solid,quasicoherent,mixed-paritycristalline}_{\square,\Gamma^\mathrm{perfect}_{\text{Robba},X,v,\infty}\{t^{1/2}\}},
\varphi\mathrm{preModule}^\mathrm{solid,quasicoherent,mixed-paritycristalline}_{\square,\Gamma^\mathrm{perfect}_{\text{Robba},X,v,I}\{t^{1/2}\}}, 
\end{align}
and
\begin{align}
\varphi\mathrm{preModule}^\mathrm{solid,quasicoherent,mixed-parityalmostcristalline}_{\square,\Gamma^\mathrm{perfect}_{\text{Robba},X,v,\infty}\{t^{1/2}\}},\\
\varphi\mathrm{preModule}^\mathrm{solid,quasicoherent,mixed-parityalmostcristalline}_{\square,\Gamma^\mathrm{perfect}_{\text{Robba},X,v,I}\{t^{1/2}\}} 
\end{align}
to $(\infty,1)$-categories in image denoted by:
\begin{align}
\mathrm{preModule}_{X,\text{\'et}}
\end{align}
to be the following functors sending each $F$ in the domain to:
\begin{align}
&\text{RH}_\text{mixed-parity}(F):=f_*(F\otimes_{\Gamma^\mathrm{perfect}_{\text{Robba},X,v,\infty}\{t^{1/2}\}} \Gamma^\mathcal{O}_{\text{cristalline},X,v}\{t^{1/2}\}),\\
&\text{RH}_\text{mixed-parity}(F):=f_*(F\otimes_{\Gamma^\mathrm{perfect}_{\text{Robba},X,v,I}\{t^{1/2}\}} \Gamma^\mathcal{O}_{\text{cristalline},X,v}\{t^{1/2}\}),\\
&\text{RH}_\text{mixed-parity}(F):=f_*(F\otimes_{\Gamma^\mathrm{perfect}_{\text{Robba},X,v,\infty}\{t^{1/2}\}} \Gamma^\mathcal{O}_{\text{cristalline},X,v}\{t^{1/2},\log(t)\}),\\
&\text{RH}_\text{mixed-parity}(F):=f_*(F\otimes_{\Gamma^\mathrm{perfect}_{\text{Robba},X,v,I}\{t^{1/2}\}} \Gamma^\mathcal{O}_{\text{cristalline},X,v}\{t^{1/2},\log(t)\}),\\
\end{align}
respectively.

\end{definition}

\newpage
\section{Geometric Family of Mixed-Parity Hodge Modules III: Semi-Stable Situations}

\noindent References: \cite{Sch1}, \cite{Sch2}, \cite{FS}, \cite{KL1}, \cite{KL2}, \cite{BL1}, \cite{BL2}, \cite{BS}, \cite{BHS}, \cite{Fon1}, \cite{CS1}, \cite{CS2}, \cite{BK}, \cite{BBK}, \cite{BBBK}, \cite{KKM}, \cite{KM}, \cite{LZ}, \cite{Shi}, \cite{M}.

\subsection{Period Rings and Sheaves}

\subsubsection{Rings}

\noindent Let $X$ be a rigid analytic space over $\mathbb{Q}_p$. We have the corresponding \'etale site and the corresponding pro-\'etale site of $X$, which we denote them by $X_{v},X_\text{\'et}$. The relationship of the two sites can be reflected by the corresponding morphism $f:X_{v}\longrightarrow X_\text{\'et}$. Then we have the corresponding semi-stable period rings and sheaves from \cite{Shi}:
\begin{align}
\Gamma_{\text{semistable},X,v}, \Gamma^\mathcal{O}_{\text{semistable},X,v}.
\end{align}
Our notations are different from \cite{Shi}, we use $\Gamma$ to mean $B$ in \cite{Shi}, while $\Gamma^\mathcal{O}$ will be the corresponding $OB$ ring in \cite{Shi}.\\

\begin{definition}
\indent Now we assume that $p>2$, following \cite{BS} we join the square root of $t$ element in $\Gamma_{\text{semistable},X,v}$ which forms the sheaves:
\begin{align}
\Gamma_{\text{semistable},X,v}\{t^{1/2}\},\Gamma^\mathcal{O}_{\text{semistable},X,v}\{t^{1/2}\}.
\end{align}
And following \cite{BL1}, \cite{BL2}, \cite{Fon1}, \cite{BHS} we further have the following sheaves of rings:
\begin{align}
\Gamma_{\text{semistable},X,v}\{t^{1/2},\log(t)\},\Gamma^\mathcal{O}_{\text{semistable},X,v}\{t^{1/2},\log(t)\}.
\end{align}
\end{definition}

\begin{definition}
We use the notations:
\begin{align}
\Gamma^\mathrm{perfect}_{\text{Robba},X,v},\Gamma^\mathrm{perfect}_{\text{Robba},X,v,\infty},\Gamma^\mathrm{perfect}_{\text{Robba},X,v,I}
\end{align}
to denote the perfect Robba rings from \cite{KL1}, \cite{KL2}, where $I\subset (0,\infty)$. Then we join $t^{1/2}$ to these sheaves we have:
\begin{align}
\Gamma^\mathrm{perfect}_{\text{Robba},X,v}\{t^{1/2}\},\Gamma^\mathrm{perfect}_{\text{Robba},X,v,\infty}\{t^{1/2}\},\Gamma^\mathrm{perfect}_{\text{Robba},X,v,I}\{t^{1/2}\}.
\end{align}
And following \cite{BL1}, \cite{BL2}, \cite{Fon1}, \cite{BHS} we have the following larger sheaves:
\begin{align}
\Gamma^\mathrm{perfect}_{\text{Robba},X,v}\{t^{1/2},\log(t)\},\Gamma^\mathrm{perfect}_{\text{Robba},X,v,\infty}\{t^{1/2},\log(t)\},\Gamma^\mathrm{perfect}_{\text{Robba},X,v,I}\{t^{1/2},\log(t)\}.
\end{align} 
\end{definition}

\begin{definition}
From now on, we use the same notation to denote the period rings involved tensored with a finite extension of $\mathbb{Q}_p$ containing square root of $p$ as in \cite{BS}.
\begin{align}
\Gamma_{\text{semistable},X,v}\{t^{1/2}\},\Gamma^\mathcal{O}_{\text{semistable},X,v}\{t^{1/2}\}.
\end{align}
\begin{align}
\Gamma_{\text{semistable},X,v}\{t^{1/2},\log(t)\},\Gamma^\mathcal{O}_{\text{semistable},X,v}\{t^{1/2},\log(t)\}.
\end{align}
\begin{align}
\Gamma^\mathrm{perfect}_{\text{Robba},X,v}\{t^{1/2}\},\Gamma^\mathrm{perfect}_{\text{Robba},X,v,\infty}\{t^{1/2}\},\Gamma^\mathrm{perfect}_{\text{Robba},X,v,I}\{t^{1/2}\}.
\end{align}
\begin{align}
\Gamma^\mathrm{perfect}_{\text{Robba},X,v}\{t^{1/2},\log(t)\},\Gamma^\mathrm{perfect}_{\text{Robba},X,v,\infty}\{t^{1/2},\log(t)\},\Gamma^\mathrm{perfect}_{\text{Robba},X,v,I}\{t^{1/2},\log(t)\}.
\end{align}
This is necessary since we to extend the action of $\varphi$ to the period rings by $\varphi(t^{1/2}\otimes 1)=\varphi(t)^{1/2}\otimes 1$.
\end{definition}

\subsubsection{Modules}

\noindent We consider quasicoherent presheaves in the following two situation:
\begin{itemize}
\item[$\square$] The solid quasicoherent modules from \cite{CS1}, \cite{CS2};
\item[$\square$] The ind-Banach quasicoherent modules from \cite{BK}, \cite{BBK}, \cite{BBBK}, \cite{KKM}, \cite{KM} with the corresponding monomorphic ind-Banach quasicoherent modules from \cite{BK}, \cite{BBK}, \cite{BBBK}, \cite{KKM}, \cite{KM}.
\end{itemize}

\begin{definition}
We use the notation:
\begin{align}
\mathrm{preModule}^\mathrm{solid,quasicoherent}_{\square,\Gamma^\mathrm{perfect}_{\text{Robba},X,v}\{t^{1/2}\}},\mathrm{preModule}^\mathrm{solid,quasicoherent}_{\square,\Gamma^\mathrm{perfect}_{\text{Robba},X,v,\infty}\{t^{1/2}\}},
\mathrm{preModule}^\mathrm{solid,quasicoherent}_{\square,\Gamma^\mathrm{perfect}_{\text{Robba},X,v,I}\{t^{1/2}\}} 
\end{align}
to denote the $(\infty,1)$-categories of solid quasicoherent presheaves over the corresonding Robba sheaves. Locally the section is defined by taking the corresponding $(\infty,1)$-categories of solid modules.
\end{definition}

\begin{definition}
We use the notation:
\begin{align}
\mathrm{preModule}^\mathrm{ind-Banach,quasicoherent}_{\Gamma^\mathrm{perfect}_{\text{Robba},X,v}\{t^{1/2}\}},\\\mathrm{preModule}^\mathrm{ind-Banach,quasicoherent}_{\Gamma^\mathrm{perfect}_{\text{Robba},X,v,\infty}\{t^{1/2}\}},\\
\mathrm{preModule}^\mathrm{ind-Banach,quasicoherent}_{\Gamma^\mathrm{perfect}_{\text{Robba},X,v,I}\{t^{1/2}\}} 
\end{align}
to denote the $(\infty,1)$-categories of solid quasicoherent presheaves over the corresonding Robba sheaves. Locally the section is defined by taking the corresponding $(\infty,1)$-categories of inductive Banach  modules. 
\end{definition}

\begin{definition}
We use the notation:
\begin{align}
\mathrm{Module}^\mathrm{solid,quasicoherent}_{\square,\Gamma^\mathrm{perfect}_{\text{Robba},X,v}\{t^{1/2}\}},\mathrm{Module}^\mathrm{solid,quasicoherent}_{\square,\Gamma^\mathrm{perfect}_{\text{Robba},X,v,\infty}\{t^{1/2}\}},
\mathrm{Module}^\mathrm{solid,quasicoherent}_{\square,\Gamma^\mathrm{perfect}_{\text{Robba},X,v,I}\{t^{1/2}\}} 
\end{align}
to denote the $(\infty,1)$-categories of solid quasicoherent sheaves over the corresonding Robba sheaves. Locally the section is defined by taking the corresponding $(\infty,1)$-categories of solid modules.
\end{definition}

\subsubsection{Mixed-Parity Hodge Modules without Frobenius}

\noindent Now we consider the key objects in our study namely those complexes generated by certain mixed-parity Hodge modules. We start from the following definition.

\begin{definition}
For any locally free coherent sheaf $F$ over
\begin{align}
\Gamma^\mathrm{perfect}_{\text{Robba},X,v,\infty}\{t^{1/2}\},\Gamma^\mathrm{perfect}_{\text{Robba},X,v,I}\{t^{1/2}\},
\end{align} 
we consider the following functor $\mathrm{dR}$ sending $F$ to the following object:
\begin{align}
f_*(F\otimes_{\Gamma^\mathrm{perfect}_{\text{Robba},X,v,\infty}\{t^{1/2}\}} \Gamma^\mathcal{O}_{\text{semistable},X,v}\{t^{1/2}\})
\end{align}
or 
\begin{align}
f_*(F\otimes_{\Gamma^\mathrm{perfect}_{\text{Robba},X,v,I}\{t^{1/2}\}} \Gamma^\mathcal{O}_{\text{semistable},X,v}\{t^{1/2}\}).
\end{align}
We call $F$ mixed-parity semi-stable if we have the following isomorphism:
\begin{align}
f^*f_*(F\otimes_{\Gamma^\mathrm{perfect}_{\text{Robba},X,v,\infty}\{t^{1/2}\}} \Gamma^\mathcal{O}_{\text{semistable},X,v}\{t^{1/2}\}) \otimes \Gamma^\mathcal{O}_{\text{semistable},X,v}\{t^{1/2}\} \overset{\sim}{\longrightarrow} F \otimes \Gamma^\mathcal{O}_{\text{semistable},X,v}\{t^{1/2}\} 
\end{align}
or 
\begin{align}
f^*f_*(F\otimes_{\Gamma^\mathrm{perfect}_{\text{Robba},X,v,I}\{t^{1/2}\}} \Gamma^\mathcal{O}_{\text{semistable},X,v}\{t^{1/2}\}) \otimes \Gamma^\mathcal{O}_{\text{semistable},X,v}\{t^{1/2}\} \overset{\sim}{\longrightarrow} F \otimes \Gamma^\mathcal{O}_{\text{semistable},X,v}\{t^{1/2}\}. 
\end{align}
\end{definition}

\begin{definition}
For any locally free coherent sheaf $F$ over
\begin{align}
\Gamma^\mathrm{perfect}_{\text{Robba},X,v,\infty}\{t^{1/2}\},\Gamma^\mathrm{perfect}_{\text{Robba},X,v,I}\{t^{1/2}\},
\end{align} 
we consider the following functor $\mathrm{dR}^\mathrm{almost}$ sending $F$ to the following object:
\begin{align}
f_*(F\otimes_{\Gamma^\mathrm{perfect}_{\text{Robba},X,v,\infty}\{t^{1/2}\}} \Gamma^\mathcal{O}_{\text{semistable},X,v}\{t^{1/2},\log(t)\})
\end{align}
or 
\begin{align}
f_*(F\otimes_{\Gamma^\mathrm{perfect}_{\text{Robba},X,v,I}\{t^{1/2}\}} \Gamma^\mathcal{O}_{\text{semistable},X,v}\{t^{1/2},\log(t)\}).
\end{align}
We call $F$ mixed-parity almost semi-stable if we have the following isomorphism:
\begin{align}
f^*f_*(F\otimes_{\Gamma^\mathrm{perfect}_{\text{Robba},X,v,\infty}\{t^{1/2}\}} \Gamma^\mathcal{O}_{\text{semistable},X,v}\{t^{1/2},\log(t)\}) \otimes \Gamma^\mathcal{O}_{\text{semistable},X,v}\{t^{1/2},\log(t)\} \\
\overset{\sim}{\longrightarrow} F \otimes \Gamma^\mathcal{O}_{\text{semistable},X,v}\{t^{1/2},\log(t)\} 
\end{align}
or 
\begin{align}
f^*f_*(F\otimes_{\Gamma^\mathrm{perfect}_{\text{Robba},X,v,I}\{t^{1/2}\}} \Gamma^\mathcal{O}_{\text{semistable},X,v}\{t^{1/2},\log(t)\}) \otimes \Gamma^\mathcal{O}_{\text{semistable},X,v}\{t^{1/2},\log(t)\}\\ \overset{\sim}{\longrightarrow} F \otimes \Gamma^\mathcal{O}_{\text{semistable},X,v}\{t^{1/2},\log(t)\}. 
\end{align}
\end{definition}

\noindent We now define the $(\infty,1)$-categories of mixed-parity semi-stable modules and he corresponding mixed-parity almost semi-stable modules by using the objects involved to generated these categories:

\begin{definition}
Considering all the mixed parity semi-stable bundles (locally finite free) as defined above, we consider the sub-$(\infty,1)$ category of 
\begin{align}
\mathrm{preModule}^\mathrm{solid,quasicoherent}_{\square,\Gamma^\mathrm{perfect}_{\text{Robba},X,v,\infty}\{t^{1/2}\}},
\mathrm{preModule}^\mathrm{solid,quasicoherent}_{\square,\Gamma^\mathrm{perfect}_{\text{Robba},X,v,I}\{t^{1/2}\}} 
\end{align}
generated by the mixed-parity semi-stable bundles (locally finite free ones). These are defined to be the $(\infty,1)$-categories of mixed-parity semi-stable complexes:
\begin{align}
\mathrm{preModule}^\mathrm{solid,quasicoherent,mixed-paritysemistable}_{\square,\Gamma^\mathrm{perfect}_{\text{Robba},X,v,\infty}\{t^{1/2}\}},
\mathrm{preModule}^\mathrm{solid,quasicoherent,mixed-paritysemistable}_{\square,\Gamma^\mathrm{perfect}_{\text{Robba},X,v,I}\{t^{1/2}\}}. 
\end{align}
\end{definition}

\begin{definition}
Considering all the mixed parity almost semi-stable bundles (locally finite free) as defined above, we consider the sub-$(\infty,1)$ category of 
\begin{align}
\mathrm{preModule}^\mathrm{solid,quasicoherent}_{\square,\Gamma^\mathrm{perfect}_{\text{Robba},X,v,\infty}\{t^{1/2}\}},
\mathrm{preModule}^\mathrm{solid,quasicoherent}_{\square,\Gamma^\mathrm{perfect}_{\text{Robba},X,v,I}\{t^{1/2}\}} 
\end{align}
generated by the mixed-parity almost semi-stable bundles (locally finite free ones). These are defined to be the $(\infty,1)$-categories of mixed-parity semi-stable complexes:
\begin{align}
\mathrm{preModule}^\mathrm{solid,quasicoherent,mixed-parityalmostsemistable}_{\square,\Gamma^\mathrm{perfect}_{\text{Robba},X,v,\infty}\{t^{1/2}\}},\\
\mathrm{preModule}^\mathrm{solid,quasicoherent,mixed-parityalmostsemistable}_{\square,\Gamma^\mathrm{perfect}_{\text{Robba},X,v,I}\{t^{1/2}\}}. 
\end{align}
\end{definition}

\indent Then the corresponding mixed-parity semi-stable functors can be extended to these categories:
\begin{align}
\mathrm{preModule}^\mathrm{solid,quasicoherent,mixed-paritysemistable}_{\square,\Gamma^\mathrm{perfect}_{\text{Robba},X,v,\infty}\{t^{1/2}\}},
\mathrm{preModule}^\mathrm{solid,quasicoherent,mixed-paritysemistable}_{\square,\Gamma^\mathrm{perfect}_{\text{Robba},X,v,I}\{t^{1/2}\}}, 
\end{align}
and
\begin{align}
\mathrm{preModule}^\mathrm{solid,quasicoherent,mixed-parityalmostsemistable}_{\square,\Gamma^\mathrm{perfect}_{\text{Robba},X,v,\infty}\{t^{1/2}\}},\\
\mathrm{preModule}^\mathrm{solid,quasicoherent,mixed-parityalmostsemistable}_{\square,\Gamma^\mathrm{perfect}_{\text{Robba},X,v,I}\{t^{1/2}\}}. 
\end{align}

\subsubsection{Mixed-Parity Hodge Modules with Frobenius}

\noindent Now we consider the key objects in our study namely those complexes generated by certain mixed-parity Hodge modules. We start from the following definition.

\begin{remark}
All the coherent sheaves over mixed-parity Robba sheaves in this section will carry the corresponding Frobenius morphism $\varphi: F \overset{\sim}{\longrightarrow} \varphi^*F$.
\end{remark}

\begin{definition}
For any locally free coherent sheaf $F$ over
\begin{align}
\Gamma^\mathrm{perfect}_{\text{Robba},X,v,\infty}\{t^{1/2}\},\Gamma^\mathrm{perfect}_{\text{Robba},X,v,I}\{t^{1/2}\},
\end{align} 
we consider the following functor $\mathrm{dR}$ sending $F$ to the following object:
\begin{align}
f_*(F\otimes_{\Gamma^\mathrm{perfect}_{\text{Robba},X,v,\infty}\{t^{1/2}\}} \Gamma^\mathcal{O}_{\text{semistable},X,v}\{t^{1/2}\})
\end{align}
or 
\begin{align}
f_*(F\otimes_{\Gamma^\mathrm{perfect}_{\text{Robba},X,v,I}\{t^{1/2}\}} \Gamma^\mathcal{O}_{\text{semistable},X,v}\{t^{1/2}\}).
\end{align}
We call $F$ mixed-parity semi-stable if we have the following isomorphism:
\begin{align}
f^*f_*(F\otimes_{\Gamma^\mathrm{perfect}_{\text{Robba},X,v,\infty}\{t^{1/2}\}} \Gamma^\mathcal{O}_{\text{semistable},X,v}\{t^{1/2}\}) \otimes \Gamma^\mathcal{O}_{\text{semistable},X,v}\{t^{1/2}\} \overset{\sim}{\longrightarrow} F \otimes \Gamma^\mathcal{O}_{\text{semistable},X,v}\{t^{1/2}\} 
\end{align}
or 
\begin{align}
f^*f_*(F\otimes_{\Gamma^\mathrm{perfect}_{\text{Robba},X,v,I}\{t^{1/2}\}} \Gamma^\mathcal{O}_{\text{semistable},X,v}\{t^{1/2}\}) \otimes \Gamma^\mathcal{O}_{\text{semistable},X,v}\{t^{1/2}\} \overset{\sim}{\longrightarrow} F \otimes \Gamma^\mathcal{O}_{\text{semistable},X,v}\{t^{1/2}\}. 
\end{align}
\end{definition}

\begin{definition}
For any locally free coherent sheaf $F$ over
\begin{align}
\Gamma^\mathrm{perfect}_{\text{Robba},X,v,\infty}\{t^{1/2}\},\Gamma^\mathrm{perfect}_{\text{Robba},X,v,I}\{t^{1/2}\},
\end{align} 
we consider the following functor $\mathrm{dR}^\mathrm{almost}$ sending $F$ to the following object:
\begin{align}
f_*(F\otimes_{\Gamma^\mathrm{perfect}_{\text{Robba},X,v,\infty}\{t^{1/2}\}} \Gamma^\mathcal{O}_{\text{semistable},X,v}\{t^{1/2},\log(t)\})
\end{align}
or 
\begin{align}
f_*(F\otimes_{\Gamma^\mathrm{perfect}_{\text{Robba},X,v,I}\{t^{1/2}\}} \Gamma^\mathcal{O}_{\text{semistable},X,v}\{t^{1/2},\log(t)\}).
\end{align}
We call $F$ mixed-parity almost semi-stable if we have the following isomorphism:
\begin{align}
f^*f_*(F\otimes_{\Gamma^\mathrm{perfect}_{\text{Robba},X,v,\infty}\{t^{1/2}\}} \Gamma^\mathcal{O}_{\text{semistable},X,v}\{t^{1/2},\log(t)\}) \otimes \Gamma^\mathcal{O}_{\text{semistable},X,v}\{t^{1/2},\log(t)\} \\
\overset{\sim}{\longrightarrow} F \otimes \Gamma^\mathcal{O}_{\text{semistable},X,v}\{t^{1/2},\log(t)\} 
\end{align}
or 
\begin{align}
f^*f_*(F\otimes_{\Gamma^\mathrm{perfect}_{\text{Robba},X,v,I}\{t^{1/2}\}} \Gamma^\mathcal{O}_{\text{semistable},X,v}\{t^{1/2},\log(t)\}) \otimes \Gamma^\mathcal{O}_{\text{semistable},X,v}\{t^{1/2},\log(t)\}\\ \overset{\sim}{\longrightarrow} F \otimes \Gamma^\mathcal{O}_{\text{semistable},X,v}\{t^{1/2},\log(t)\}. 
\end{align}
\end{definition}

\noindent We now define the $(\infty,1)$-categories of mixed-parity semi-stable modules and he corresponding mixed-parity almost semi-stable modules by using the objects involved to generated these categories:

\begin{definition}
Considering all the mixed parity semi-stable bundles (locally finite free) as defined above, we consider the sub-$(\infty,1)$ category of 
\begin{align}
\varphi\mathrm{preModule}^\mathrm{solid,quasicoherent}_{\square,\Gamma^\mathrm{perfect}_{\text{Robba},X,v,\infty}\{t^{1/2}\}},
\varphi\mathrm{preModule}^\mathrm{solid,quasicoherent}_{\square,\Gamma^\mathrm{perfect}_{\text{Robba},X,v,I}\{t^{1/2}\}} 
\end{align}
generated by the mixed-parity semi-stable bundles (locally finite free ones). These are defined to be the $(\infty,1)$-categories of mixed-parity semi-stable complexes:
\begin{align}
\varphi\mathrm{preModule}^\mathrm{solid,quasicoherent,mixed-paritysemistable}_{\square,\Gamma^\mathrm{perfect}_{\text{Robba},X,v,\infty}\{t^{1/2}\}},
\varphi\mathrm{preModule}^\mathrm{solid,quasicoherent,mixed-paritysemistable}_{\square,\Gamma^\mathrm{perfect}_{\text{Robba},X,v,I}\{t^{1/2}\}}. 
\end{align}
\end{definition}

\begin{definition}
Considering all the mixed parity almost semi-stable bundles (locally finite free) as defined above, we consider the sub-$(\infty,1)$ category of 
\begin{align}
\varphi\mathrm{preModule}^\mathrm{solid,quasicoherent}_{\square,\Gamma^\mathrm{perfect}_{\text{Robba},X,v,\infty}\{t^{1/2}\}},
\varphi\mathrm{preModule}^\mathrm{solid,quasicoherent}_{\square,\Gamma^\mathrm{perfect}_{\text{Robba},X,v,I}\{t^{1/2}\}} 
\end{align}
generated by the mixed-parity almost semi-stable bundles (locally finite free ones). These are defined to be the $(\infty,1)$-categories of mixed-parity semi-stable complexes:
\begin{align}
\varphi\mathrm{preModule}^\mathrm{solid,quasicoherent,mixed-parityalmostsemistable}_{\square,\Gamma^\mathrm{perfect}_{\text{Robba},X,v,\infty}\{t^{1/2}\}},\\
\varphi\mathrm{preModule}^\mathrm{solid,quasicoherent,mixed-parityalmostsemistable}_{\square,\Gamma^\mathrm{perfect}_{\text{Robba},X,v,I}\{t^{1/2}\}}. 
\end{align}
\end{definition}

\indent Then the corresponding mixed-parity semi-stable functors can be extended to these categories:
\begin{align}
\varphi\mathrm{preModule}^\mathrm{solid,quasicoherent,mixed-paritysemistable}_{\square,\Gamma^\mathrm{perfect}_{\text{Robba},X,v,\infty}\{t^{1/2}\}},
\varphi\mathrm{preModule}^\mathrm{solid,quasicoherent,mixed-paritysemistable}_{\square,\Gamma^\mathrm{perfect}_{\text{Robba},X,v,I}\{t^{1/2}\}}, 
\end{align}
and
\begin{align}
\varphi\mathrm{preModule}^\mathrm{solid,quasicoherent,mixed-parityalmostsemistable}_{\square,\Gamma^\mathrm{perfect}_{\text{Robba},X,v,\infty}\{t^{1/2}\}},\\
\varphi\mathrm{preModule}^\mathrm{solid,quasicoherent,mixed-parityalmostsemistable}_{\square,\Gamma^\mathrm{perfect}_{\text{Robba},X,v,I}\{t^{1/2}\}}. 
\end{align}

\subsection{Mixed-Parity semi-stable Riemann-Hilbert Correspondence}

\indent This chapter will extend the corresponding Riemann-Hilbert correspondence from \cite{Sch1}, \cite{LZ}, \cite{BL1}, \cite{BL2}, \cite{M} to the mixed-parity setting.

\begin{definition}
We define the following Riemann-Hilbert functor $\text{RH}_\text{mixed-parity}$ from the one of categories:
\begin{align}
\mathrm{preModule}^\mathrm{solid,quasicoherent,mixed-paritysemistable}_{\square,\Gamma^\mathrm{perfect}_{\text{Robba},X,v,\infty}\{t^{1/2}\}},
\mathrm{preModule}^\mathrm{solid,quasicoherent,mixed-paritysemistable}_{\square,\Gamma^\mathrm{perfect}_{\text{Robba},X,v,I}\{t^{1/2}\}}, 
\end{align}
and
\begin{align}
\mathrm{preModule}^\mathrm{solid,quasicoherent,mixed-parityalmostsemistable}_{\square,\Gamma^\mathrm{perfect}_{\text{Robba},X,v,\infty}\{t^{1/2}\}},\\
\mathrm{preModule}^\mathrm{solid,quasicoherent,mixed-parityalmostsemistable}_{\square,\Gamma^\mathrm{perfect}_{\text{Robba},X,v,I}\{t^{1/2}\}} 
\end{align}
to $(\infty,1)$-categories in image denoted by:
\begin{align}
\mathrm{preModule}_{X,\text{\'et}}
\end{align}
to be the following functors sending each $F$ in the domain to:
\begin{align}
&\text{RH}_\text{mixed-parity}(F):=f_*(F\otimes_{\Gamma^\mathrm{perfect}_{\text{Robba},X,v,\infty}\{t^{1/2}\}} \Gamma^\mathcal{O}_{\text{semistable},X,v}\{t^{1/2}\}),\\
&\text{RH}_\text{mixed-parity}(F):=f_*(F\otimes_{\Gamma^\mathrm{perfect}_{\text{Robba},X,v,I}\{t^{1/2}\}} \Gamma^\mathcal{O}_{\text{semistable},X,v}\{t^{1/2}\}),\\
&\text{RH}_\text{mixed-parity}(F):=f_*(F\otimes_{\Gamma^\mathrm{perfect}_{\text{Robba},X,v,\infty}\{t^{1/2}\}} \Gamma^\mathcal{O}_{\text{semistable},X,v}\{t^{1/2},\log(t)\}),\\
&\text{RH}_\text{mixed-parity}(F):=f_*(F\otimes_{\Gamma^\mathrm{perfect}_{\text{Robba},X,v,I}\{t^{1/2}\}} \Gamma^\mathcal{O}_{\text{semistable},X,v}\{t^{1/2},\log(t)\}),\\
\end{align}
respectively.

\end{definition}

\begin{definition}
In the situation where we have the Frobenius action we consider the follwing. We define the following Riemann-Hilbert functor $\text{RH}_\text{mixed-parity}$ from the one of categories:
\begin{align}
\varphi\mathrm{preModule}^\mathrm{solid,quasicoherent,mixed-paritysemistable}_{\square,\Gamma^\mathrm{perfect}_{\text{Robba},X,v,\infty}\{t^{1/2}\}},
\varphi\mathrm{preModule}^\mathrm{solid,quasicoherent,mixed-paritysemistable}_{\square,\Gamma^\mathrm{perfect}_{\text{Robba},X,v,I}\{t^{1/2}\}}, 
\end{align}
and
\begin{align}
\varphi\mathrm{preModule}^\mathrm{solid,quasicoherent,mixed-parityalmostsemistable}_{\square,\Gamma^\mathrm{perfect}_{\text{Robba},X,v,\infty}\{t^{1/2}\}},\\
\varphi\mathrm{preModule}^\mathrm{solid,quasicoherent,mixed-parityalmostsemistable}_{\square,\Gamma^\mathrm{perfect}_{\text{Robba},X,v,I}\{t^{1/2}\}} 
\end{align}
to $(\infty,1)$-categories in image denoted by:
\begin{align}
\mathrm{preModule}_{X,\text{\'et}}
\end{align}
to be the following functors sending each $F$ in the domain to:
\begin{align}
&\text{RH}_\text{mixed-parity}(F):=f_*(F\otimes_{\Gamma^\mathrm{perfect}_{\text{Robba},X,v,\infty}\{t^{1/2}\}} \Gamma^\mathcal{O}_{\text{semistable},X,v}\{t^{1/2}\}),\\
&\text{RH}_\text{mixed-parity}(F):=f_*(F\otimes_{\Gamma^\mathrm{perfect}_{\text{Robba},X,v,I}\{t^{1/2}\}} \Gamma^\mathcal{O}_{\text{semistable},X,v}\{t^{1/2}\}),\\
&\text{RH}_\text{mixed-parity}(F):=f_*(F\otimes_{\Gamma^\mathrm{perfect}_{\text{Robba},X,v,\infty}\{t^{1/2}\}} \Gamma^\mathcal{O}_{\text{semistable},X,v}\{t^{1/2},\log(t)\}),\\
&\text{RH}_\text{mixed-parity}(F):=f_*(F\otimes_{\Gamma^\mathrm{perfect}_{\text{Robba},X,v,I}\{t^{1/2}\}} \Gamma^\mathcal{O}_{\text{semistable},X,v}\{t^{1/2},\log(t)\}),\\
\end{align}
respectively.

\end{definition}

\begin{remark}
We now have discussed the corresponding two different morphisms:
\begin{align}
f: X_\text{pro\'et}\longrightarrow X_\text{\'et};\\
f': X_\text{v}\longrightarrow X_\text{\'et}.
\end{align}
One can consider the following relation among the sites:
\begin{align}
X_\text{v}\longrightarrow X_\text{pro\'et}\longrightarrow X_\text{\'et}
\end{align}
which produces $f'$. The map:
\begin{align}
g: X_\text{v}\longrightarrow X_\text{pro\'et}
\end{align}
can help us relate the corresponding constructions above as in \cite[Proposition 2.37]{B}. Namely we have:
\begin{align}
&\mathrm{dR}_{v}=\mathrm{dR}_{\text{pro\'et}}g_*;\\
&\mathrm{dR}_{v,\text{almost}}=\mathrm{dR}_{\text{pro\'et},\text{almost}}g_*;\\
&\mathrm{cristalline}_{v}=\mathrm{cristalline}_{\text{pro\'et}}g_*;\\
&\mathrm{cristalline}_{v,\text{almost}}=\mathrm{cristalline}_{\text{pro\'et},\text{almost}}g_*;\\
&\mathrm{semistable}_{v}=\mathrm{semistable}_{\text{pro\'et}}g_*;\\
&\mathrm{semistable}_{v,\text{almost}}=\mathrm{semistable}_{\text{pro\'et},\text{almost}}g_*.
\end{align} 
\end{remark}

\chapter{Mixed-Parity Hodge Modules over $v$-Stacks}

\newpage
\section{$(\infty,1)$-Quasicoherent Sheaves over Extended Fargues-Fontaine Curves I}

\noindent We now consider the sheaves over extended Fargues-Fontain stacks:

\begin{remark}
Let $X$ be a general small $v$-stack over $\mathbb{Q}_p$ (as a $v$-stack\footnote{All $v$-stacks in this chapter are assumed to be over a $v$-stack associated to $\mathbb{Q}_p$ like this.}).  
\end{remark}

\begin{definition}
\begin{align}
\mathrm{FF}_X:=\bigcup_{I\subset (0,\infty)}\mathrm{Spa}(\Gamma^\text{perfect}_{\text{Robba},X,I\subset (0,\infty)}\{t^{1/2}\}\otimes_{\mathbb{Q}_p}E,\Gamma^{\text{perfect},+}_{\text{Robba},X,I\subset (0,\infty)}\{t^{1/2}\}\otimes_{\mathbb{Q}_p}E)/\varphi^\mathbb{Z},
\end{align}
which has the corresonding structure map as in the following:
\[\displayindent=-0.4in
\xymatrix@R+1pc{
&\mathrm{FF}_X  \ar[d] \\
&\mathrm{FF}_{\mathrm{Spd}^\diamond(\mathbb{Q}_p)}.  
}
\]
\end{definition}

\begin{definition}
We use the notation
\begin{align}
\mathrm{Quasicoherent}^{\mathrm{solid}}_{\mathrm{FF}_X,\mathcal{O}_{\mathrm{FF}_X}}
\end{align}
to denote $(\infty,1)$-category of all the solid quasicoherent sheaves over the stack $\mathrm{FF}_X$. For any local perfectoid $Y\in X_v$ we define the corresponding $(\infty,1)$-category in the local sense.\\
We use the notation
\begin{align}
\mathrm{Quasicoherent}^{\mathrm{solid,perfectcomplexes}}_{\mathrm{FF}_X,\mathcal{O}_{\mathrm{FF}_X}}
\end{align}
to denote $(\infty,1)$-category of all the solid quasicoherent sheaves over the stack $\mathrm{FF}_X$ which are perfect complexes. For any local perfectoid $Y\in X_v$ we define the corresponding $(\infty,1)$-category in the local sense.
\end{definition}

\begin{definition}
We use the notation
\begin{align}
\mathrm{Quasicoherent}^{\mathrm{indBanach}}_{\mathrm{FF}_X,\mathcal{O}_{\mathrm{FF}_X}}
\end{align}
to denote $(\infty,1)$-category of all the ind-Banach quasicoherent sheaves over the stack $\mathrm{FF}_X$. For any local perfectoid $Y\in X_v$ we define the corresponding $(\infty,1)$-category in the local sense.\\
We use the notation
\begin{align}
\mathrm{Quasicoherent}^{\mathrm{indBanach,perfectcomplexes}}_{\mathrm{FF}_X,\mathcal{O}_{\mathrm{FF}_X}}
\end{align}
to denote $(\infty,1)$-category of all the indBanach quasicoherent sheaves over the stack $\mathrm{FF}_X$ which are perfect complexes. For any local perfectoid $Y\in X_v$ we define the corresponding $(\infty,1)$-category in the local sense.
\end{definition}

\begin{definition}
We use the notation
\begin{align}
\{{\varphi\mathrm{Module}^{\mathrm{solid}}}({\Gamma^\text{perfect}_{\text{Robba},X,I}\{t^{1/2}\}\otimes_{\mathbb{Q}_p}E})\}_{I\subset (0,\infty)}
\end{align}
to denote $(\infty,1)$-category of all the solid $\varphi$-modules over the extended Robba ring. The modules satisfy the Frobenius pullback condition and glueing condition for overlapped intervals $I\subset J\subset K$. For any local perfectoid $Y\in X_v$ we define the corresponding $(\infty,1)$-category in the local sense.\\
We use the notation
\begin{align}
\{\underset{\mathrm{solid,perfectcomplexes}}{\varphi\mathrm{Module}}(\Gamma^\text{perfect}_{\text{Robba},X,I}\{t^{1/2}\}\otimes_{\mathbb{Q}_p}E)\}_{I\subset (0,\infty)}
\end{align}
to denote $(\infty,1)$-category of all the solid $\varphi$-modules over the extended Robba ring which are perfect complexes. The modules satisfy the Frobenius pullback condition and glueing condition for overlapped intervals $I\subset J\subset K$. For any local perfectoid $Y\in X_v$ we define the corresponding $(\infty,1)$-category in the local sense. 

\end{definition}

\begin{definition}
We use the notation
\begin{align}
\{{\varphi\mathrm{Module}}^{\mathrm{indBanach}}(\Gamma^\text{perfect}_{\text{Robba},X,I}\{t^{1/2}\}\otimes_{\mathbb{Q}_p}E)\}_{I\subset (0,\infty)}
\end{align}
to denote $(\infty,1)$-category of all the ind-Banach $\varphi$-modules over the extended Robba ring. The modules satisfy the Frobenius pullback condition and glueing condition for overlapped intervals $I\subset J\subset K$. For any local perfectoid $Y\in X_v$ we define the corresponding $(\infty,1)$-category in the local sense.\\
We use the notation
\begin{align}
\left\{\underset{\mathrm{indBanach,perfectcomplexes},\Gamma^\text{perfect}_{\text{Robba},X,I}\{t^{1/2}\}\otimes_{\mathbb{Q}_p}E}{\varphi\mathrm{Module}}\right\}_{I\subset (0,\infty)}
\end{align}
to denote $(\infty,1)$-category of all the ind-Banach $\varphi$-modules over the extended Robba ring which are perfect complexes. The modules satisfy the Frobenius pullback condition and glueing condition for overlapped intervals $I\subset J\subset K$. For any local perfectoid $Y\in X_v$ we define the corresponding $(\infty,1)$-category in the local sense. 

\end{definition}

\begin{proposition}
We have the following commutative diagram by taking the global section functor in the horizontal rows:\\
\[\displayindent=+0in
\xymatrix@R+7pc{
\mathrm{Quasicoherent}^{\mathrm{solid}}_{\mathrm{FF}_{\mathrm{Spd}(\mathbb{Q}_p)^\diamond},\mathcal{O}_{\mathrm{FF}_{\mathrm{Spd}(\mathbb{Q}_p)^\diamond}}} \ar[r] \ar[d] &\{\varphi\mathrm{Module}^{\mathrm{solid}}(\Gamma^\text{perfect}_{\text{Robba},{\mathrm{Spd}(\mathbb{Q}_p)^\diamond},I}\{t^{1/2}\}\otimes_{\mathbb{Q}_p}E)\}_{I\subset (0,\infty)} \ar[d]\\
\mathrm{Quasicoherent}^{\mathrm{solid}}_{\mathrm{FF}_X,\mathcal{O}_{\mathrm{FF}_X}}  \ar[r] \ar[r] \ar[r] &\{\varphi\mathrm{Module}^{\mathrm{solid}}(\Gamma^\text{perfect}_{\text{Robba},X,I}\{t^{1/2}\}\otimes_{\mathbb{Q}_p}E)\}_{I\subset (0,\infty)}.  \\  
}
\]
\end{proposition}

\begin{proposition}
We have the following commutative diagram by taking the global section functor in the horizontal rows:\\
\[\displayindent=+0in
\xymatrix@R+7pc{
\mathrm{Quasicoherent}^{\mathrm{solid,perfectcomplexes}}_{\mathrm{FF}_{\mathrm{Spd}(\mathbb{Q}_p)^\diamond},\mathcal{O}_{\mathrm{FF}_{\mathrm{Spd}(\mathbb{Q}_p)^\diamond}}} \ar[r] \ar[d] &\{\underset{\mathrm{solid,perfectcomplexes}}{\varphi\mathrm{Module}}(\Gamma^\text{perfect}_{\text{Robba},{\mathrm{Spd}(\mathbb{Q}_p)^\diamond},I}\{t^{1/2}\}\otimes_{\mathbb{Q}_p}E)\}_{I\subset (0,\infty)} \ar[d]  \\
\mathrm{Quasicoherent}^{\mathrm{solid,perfectcomplexes}}_{\mathrm{FF}_X,\mathcal{O}_{\mathrm{FF}_X}}  \ar[r] \ar[r] \ar[r] &\{\underset{\mathrm{solid,perfectcomplexes}}{\varphi\mathrm{Module}}(\Gamma^\text{perfect}_{\text{Robba},X,I}\{t^{1/2}\}\otimes_{\mathbb{Q}_p}E)\}_{I\subset (0,\infty)}.   
}
\]
\end{proposition}

\begin{proposition}
We have the following commutative diagram by taking the global section functor in the horizontal rows:\\
\[\displayindent=+0in
\xymatrix@R+7pc{
\mathrm{Quasicoherent}^{\mathrm{indBanach}}_{\mathrm{FF}_{\mathrm{Spd}(\mathbb{Q}_p)^\diamond},\mathcal{O}_{\mathrm{FF}_{\mathrm{Spd}(\mathbb{Q}_p)^\diamond}}} \ar[r] \ar[d] &\{\varphi\mathrm{Module}^\text{indBanach}(\Gamma^\text{perfect}_{\text{Robba},{\mathrm{Spd}(\mathbb{Q}_p)^\diamond},I}\{t^{1/2}\}\otimes_{\mathbb{Q}_p}E)\}_{I\subset (0,\infty)}\ar[d]\\
\mathrm{Quasicoherent}^{\mathrm{indBanach}}_{\mathrm{FF}_X,\mathcal{O}_{\mathrm{FF}_X}}  \ar[r] \ar[r] \ar[r] &\{\varphi\mathrm{Module}^{\mathrm{indBanach}}(\Gamma^\text{perfect}_{\text{Robba},X,I}\{t^{1/2}\}\otimes_{\mathbb{Q}_p}E)\}_{I\subset (0,\infty)}.    
}
\]
\end{proposition}

\begin{proposition}
We have the following commutative diagram by taking the global section functor in the horizontal rows:
\[\displayindent=+0in
\xymatrix@C-0.2in@R+7pc{
\mathrm{Quasicoherent}^{\mathrm{indBanach,perfectcomplexes}}_{\mathrm{FF}_{\mathrm{Spd}(\mathbb{Q}_p)^\diamond},\mathcal{O}_{\mathrm{FF}_{\mathrm{Spd}(\mathbb{Q}_p)^\diamond}}} \ar[d]\ar[r] &\left\{\underset{\mathrm{indBanach,perfectcomplexes},\Gamma^\text{perfect}_{\text{Robba},\mathrm{Spd}(\mathbb{Q}_p)^\diamond,I}\{t^{1/2}\}\otimes_{\mathbb{Q}_p}E}{\varphi\mathrm{Module}}\right\}_{I\subset (0,\infty)}\ar[d]\\
\mathrm{Quasicoherent}^{\mathrm{indBanach,perfectcomplexes}}_{\mathrm{FF}_X,\mathcal{O}_{\mathrm{FF}_X}}  \ar[r] \ar[r] \ar[r] &\left\{\underset{\mathrm{indBanach,perfectcomplexes},\Gamma^\text{perfect}_{\text{Robba},X,I}\{t^{1/2}\}\otimes_{\mathbb{Q}_p}E}{\varphi\mathrm{Module}}\right\}_{I\subset (0,\infty)}.    
}
\]

\end{proposition}

\indent Taking the corresponding simplicial commutative object we have the following propositions:

\begin{proposition}
We have the following commutative diagram by taking the global section functor in the horizontal rows:
\[\displayindent=+0in
\xymatrix@R+7pc{
\underset{\mathrm{Quasicoherent}^{\mathrm{solid}}_{\mathrm{FF}_{\mathrm{Spd}(\mathbb{Q}_p)^\diamond},\mathcal{O}_{\mathrm{FF}_{\mathrm{Spd}(\mathbb{Q}_p)^\diamond}}}}{\mathrm{SimplicialRings}} \ar[d] \ar[r] &\underset{\{\varphi\mathrm{Module}^{\mathrm{solid}}(\Gamma^\text{perfect}_{\text{Robba},{\mathrm{Spd}(\mathbb{Q}_p)^\diamond},I}\{t^{1/2}\}\otimes_{\mathbb{Q}_p}E)\}_{I\subset (0,\infty)}}{\mathrm{SimplicialRings}} \ar[d]\\
\underset{\mathrm{Quasicoherent}^{\mathrm{solid}}_{\mathrm{FF}_X,\mathcal{O}_{\mathrm{FF}_X}}}{\mathrm{SimplicialRings}}  \ar[r] \ar[r] \ar[r] &\underset{\{\varphi\mathrm{Module}^{\mathrm{solid}}(\Gamma^\text{perfect}_{\text{Robba},X,I}\{t^{1/2}\}\otimes_{\mathbb{Q}_p}E)\}_{I\subset (0,\infty)}}{\mathrm{SimplicialRings}}.
}
\]
\end{proposition}

\begin{proposition}
We have the following commutative diagram by taking the global section functor in the horizontal rows:
\[\displayindent=+0in
\xymatrix@C+0in@R+7pc{
\underset{\mathrm{Quasicoherent}^{\mathrm{solid,perfectcomplexes}}_{\mathrm{FF}_{\mathrm{Spd}(\mathbb{Q}_p)^\diamond},\mathcal{O}_{\mathrm{FF}_{\mathrm{Spd}(\mathbb{Q}_p)^\diamond}}}}{\mathrm{SimplicialRings}}\ar[d] \ar[r] &\underset{\{\varphi\mathrm{Module}^{\mathrm{solid,perfectcomplexes}}(\Gamma^\text{perfect}_{\text{Robba},{\mathrm{Spd}(\mathbb{Q}_p)^\diamond},I}\{t^{1/2}\}\otimes_{\mathbb{Q}_p}E)\}_{I\subset (0,\infty)}}{\mathrm{SimplicialRings}}\ar[d]\\
\underset{\mathrm{Quasicoherent}^{\mathrm{solid,perfectcomplexes}}_{\mathrm{FF}_X,\mathcal{O}_{\mathrm{FF}_X}}}{\mathrm{SimplicialRings}}  \ar[r] \ar[r] \ar[r] &\underset{\{\varphi\mathrm{Module}^{\mathrm{solid,perfectcomplexes}}(\Gamma^\text{perfect}_{\text{Robba},X,I}\{t^{1/2}\}\otimes_{\mathbb{Q}_p}E)\}_{I\subset (0,\infty)}}{\mathrm{SimplicialRings}}.  
}
\]
\end{proposition}

\begin{proposition}
We have the following commutative diagram by taking the global section functor in the horizontal rows:\\
\[\displayindent=+0in
\xymatrix@R+7pc{
\underset{\mathrm{Quasicoherent}^{\mathrm{indBanach}}_{\mathrm{FF}_X,\mathcal{O}_{\mathrm{FF}_X}}}{\mathrm{SimplicialRings}}  \ar[r] \ar[r] \ar[r] &\underset{\{\varphi\mathrm{Module}^{\mathrm{indBanach}}(\Gamma^\text{perfect}_{\text{Robba},X,I}\{t^{1/2}\}\otimes_{\mathbb{Q}_p}E)\}_{I\subset (0,\infty)}}{\mathrm{SimplicialRings}}   \\
\underset{\mathrm{Quasicoherent}^{\mathrm{indBanach}}_{\mathrm{FF}_{\mathrm{Spd}(\mathbb{Q}_p)^\diamond},\mathcal{O}_{\mathrm{FF}_{\mathrm{Spd}(\mathbb{Q}_p)^\diamond}}}}{\mathrm{SimplicialRings}} \ar[u]\ar[r] &\underset{\{\varphi\mathrm{Module}^\text{indBanach}(\Gamma^\text{perfect}_{\text{Robba},{\mathrm{Spd}(\mathbb{Q}_p)^\diamond},I}\{t^{1/2}\}\otimes_{\mathbb{Q}_p}E)\}_{I\subset (0,\infty)}}{\mathrm{SimplicialRings}}.\ar[u]  
}
\]
\end{proposition}

\begin{proposition}
We have the following commutative diagram by taking the global section functor in the horizontal rows:\\
\[\displayindent=+0in
\xymatrix@C+0in@R+7pc{
\underset{\mathrm{Quasicoherent}^{\mathrm{indBanach,perfectcomplexes}}_{\mathrm{FF}_X,\mathcal{O}_{\mathrm{FF}_X}}}{\mathrm{SimplicialRings}}  \ar[r] \ar[r] \ar[r] &\underset{\{\varphi\mathrm{Module}^{\mathrm{indBanach,perfectcomplexes}}(\Gamma^\text{perfect}_{\text{Robba},X,I}\{t^{1/2}\}\otimes_{\mathbb{Q}_p}E)\}_{I\subset (0,\infty)}}{\mathrm{SimplicialRings}}   \\
\underset{\mathrm{Quasicoherent}^{\mathrm{indBanach,perfectcomplexes}}_{\mathrm{FF}_{\mathrm{Spd}(\mathbb{Q}_p)^\diamond},\mathcal{O}_{\mathrm{FF}_{\mathrm{Spd}(\mathbb{Q}_p)^\diamond}}}}{\mathrm{SimplicialRings}} \ar[u]\ar[r] &\underset{\{\varphi\mathrm{Module}^{\mathrm{indBanach,perfectcomplexes}}(\Gamma^\text{perfect}_{\text{Robba},{\mathrm{Spd}(\mathbb{Q}_p)^\diamond},I}\{t^{1/2}\}\otimes_{\mathbb{Q}_p}E)\}_{I\subset (0,\infty)}}{\mathrm{SimplicialRings}}.\ar[u]  
}
\]

\end{proposition}

\newpage
\section{$(\infty,1)$-Quasicoherent Sheaves over Extended Fargues-Fontaine Curves II}

\noindent We now consider the sheaves over extended Fargues-Fontain stacks:

\begin{remark}
Let $X$ be a general small $v$-stack over $\mathbb{Q}_p$ (as a $v$-stack\footnote{All $v$-stacks in this chapter are assumed to be over a $v$-stack associated to $\mathbb{Q}_p$ like this.}). $\mathrm{Spa}$ will denote Clausen-Scholze analytic space in \cite{CS2}. 
\end{remark}

\begin{definition}
\begin{align}
\mathrm{FF}_X:=\bigcup_{I\subset (0,\infty)}\mathrm{Spa}(\Gamma^\text{perfect}_{\text{Robba},X,I\subset (0,\infty)}\{t^{1/2},\log(t)\}\otimes_{\mathbb{Q}_p}E,\Gamma^{\text{perfect},+}_{\text{Robba},X,I\subset (0,\infty)}\{t^{1/2},\log(t)\}\otimes_{\mathbb{Q}_p}E)/\varphi^\mathbb{Z},
\end{align}
\footnote{Here the ring $\Gamma^\text{perfect}_{\text{Robba},X,I\subset (0,\infty)}\{t^{1/2},\log(t)\}$ is defined to be just:
\begin{align}
\Gamma^\text{perfect}_{\text{Robba},X,I\subset (0,\infty)}\{t^{1/2}\}[\log(t)]
\end{align}
which carries the corresponding adic topology from the corresponding Banach ring 
$\Gamma^\text{perfect}_{\text{Robba},X,I\subset (0,\infty)}\{t^{1/2}\}$, which induces a topological adic ring structure (therefore a corresponding condensed animated ring structure in \cite{CS2}). Then the corresponding spectrum will be defined to be the correponding analytic spectrum from Clausen-Scholze.  
}which has the corresonding structure map as in the following:
\[\displayindent=-0.4in
\xymatrix@R+1pc{
&\mathrm{FF}_X  \ar[d] \\
&\mathrm{FF}_{\mathrm{Spd}^\diamond(\mathbb{Q}_p)}.  
}
\]
\end{definition}

\begin{definition}
We use the notation
\begin{align}
\mathrm{Quasicoherent}^{\mathrm{solid}}_{\mathrm{FF}_X,\mathcal{O}_{\mathrm{FF}_X}}
\end{align}
to denote $(\infty,1)$-category of all the solid quasicoherent sheaves over the stack $\mathrm{FF}_X$. For any local perfectoid $Y\in X_v$ we define the corresponding $(\infty,1)$-category in the local sense.\\
We use the notation
\begin{align}
\mathrm{Quasicoherent}^{\mathrm{solid,perfectcomplexes}}_{\mathrm{FF}_X,\mathcal{O}_{\mathrm{FF}_X}}
\end{align}
to denote $(\infty,1)$-category of all the solid quasicoherent sheaves over the stack $\mathrm{FF}_X$ which are perfect complexes. For any local perfectoid $Y\in X_v$ we define the corresponding $(\infty,1)$-category in the local sense.
\end{definition}

\begin{definition}
We use the notation
\begin{align}
\mathrm{Quasicoherent}^{\mathrm{indBanach}}_{\mathrm{FF}_X,\mathcal{O}_{\mathrm{FF}_X}}
\end{align}
to denote $(\infty,1)$-category of all the ind-Banach quasicoherent sheaves over the stack $\mathrm{FF}_X$. For any local perfectoid $Y\in X_v$ we define the corresponding $(\infty,1)$-category in the local sense.\\
We use the notation
\begin{align}
\mathrm{Quasicoherent}^{\mathrm{indBanach,perfectcomplexes}}_{\mathrm{FF}_X,\mathcal{O}_{\mathrm{FF}_X}}
\end{align}
to denote $(\infty,1)$-category of all the indBanach quasicoherent sheaves over the stack $\mathrm{FF}_X$ which are perfect complexes. For any local perfectoid $Y\in X_v$ we define the corresponding $(\infty,1)$-category in the local sense.
\end{definition}

\begin{definition}
We use the notation
\begin{align}
\{{\varphi\mathrm{Module}^{\mathrm{solid}}}({\Gamma^\text{perfect}_{\text{Robba},X,I}\{t^{1/2},\log(t)\}\otimes_{\mathbb{Q}_p}E})\}_{I\subset (0,\infty)}
\end{align}
to denote $(\infty,1)$-category of all the solid $\varphi$-modules over the extended Robba ring. The modules satisfy the Frobenius pullback condition and glueing condition for overlapped intervals $I\subset J\subset K$. For any local perfectoid $Y\in X_v$ we define the corresponding $(\infty,1)$-category in the local sense.\\
We use the notation
\begin{align}
\left\{\underset{\mathrm{solid,perfectcomplexes},(\Gamma^\text{perfect}_{\text{Robba},X,I}\{t^{1/2},\log(t)\}\otimes_{\mathbb{Q}_p}E)}{\varphi\mathrm{Module}}\right\}_{I\subset (0,\infty)}
\end{align}
to denote $(\infty,1)$-category of all the solid $\varphi$-modules over the extended Robba ring which are perfect complexes. The modules satisfy the Frobenius pullback condition and glueing condition for overlapped intervals $I\subset J\subset K$. For any local perfectoid $Y\in X_v$ we define the corresponding $(\infty,1)$-category in the local sense. 

\end{definition}

\begin{definition}
We use the notation
\begin{align}
\{\underset{\mathrm{indBanach}}{\varphi\mathrm{Module}}(\Gamma^\text{perfect}_{\text{Robba},X,I}\{t^{1/2},\log(t)\}\otimes_{\mathbb{Q}_p}E)\}_{I\subset (0,\infty)}
\end{align}
to denote $(\infty,1)$-category of all the ind-Banach $\varphi$-modules over the extended Robba ring. The modules satisfy the Frobenius pullback condition and glueing condition for overlapped intervals $I\subset J\subset K$. For any local perfectoid $Y\in X_v$ we define the corresponding $(\infty,1)$-category in the local sense.\\
We use the notation
\begin{align}
\left\{\underset{\mathrm{indBanach,perfectcomplexes},\Gamma^\text{perfect}_{\text{Robba},X,I}\{t^{1/2},\log(t)\}\otimes_{\mathbb{Q}_p}E}{\varphi\mathrm{Module}}\right\}_{I\subset (0,\infty)}
\end{align}
to denote $(\infty,1)$-category of all the ind-Banach $\varphi$-modules over the extended Robba ring which are perfect complexes. The modules satisfy the Frobenius pullback condition and glueing condition for overlapped intervals $I\subset J\subset K$. For any local perfectoid $Y\in X_v$ we define the corresponding $(\infty,1)$-category in the local sense. 

\end{definition}

\begin{proposition}
We have the following commutative diagram by taking the global section functor in the horizontal rows:
\[\displayindent=+0in
\xymatrix@R+7pc{
\mathrm{Quasicoherent}^{\mathrm{solid}}_{\mathrm{FF}_{\mathrm{Spd}(\mathbb{Q}_p)^\diamond},\mathcal{O}_{\mathrm{FF}_{\mathrm{Spd}(\mathbb{Q}_p)^\diamond}}} \ar[r] \ar[d] &\{\varphi\mathrm{Module}^{\mathrm{solid}}(\Gamma^\text{perfect}_{\text{Robba},{\mathrm{Spd}(\mathbb{Q}_p)^\diamond},I}\{t^{1/2},\log(t)\}\otimes_{\mathbb{Q}_p}E)\}_{I\subset (0,\infty)} \ar[d]\\
\mathrm{Quasicoherent}^{\mathrm{solid}}_{\mathrm{FF}_X,\mathcal{O}_{\mathrm{FF}_X}}  \ar[r] \ar[r] \ar[r] &\{\varphi\mathrm{Module}^{\mathrm{solid}}(\Gamma^\text{perfect}_{\text{Robba},X,I}\{t^{1/2},\log(t)\}\otimes_{\mathbb{Q}_p}E)\}_{I\subset (0,\infty)}.   
}
\]
\end{proposition}

\begin{proposition}
We have the following commutative diagram by taking the global section functor in the horizontal rows:
\[\displayindent=+0in
\xymatrix@R+7pc{
\mathrm{Quasicoherent}^{\mathrm{solid,perfectcomplexes}}_{\mathrm{FF}_{\mathrm{Spd}(\mathbb{Q}_p)^\diamond},\mathcal{O}_{\mathrm{FF}_{\mathrm{Spd}(\mathbb{Q}_p)^\diamond}}} \ar[r] \ar[d] &\left\{\underset{\mathrm{solid,perfectcomplexes},(\Gamma^\text{perfect}_{\text{Robba},\mathrm{Spd}(\mathbb{Q}_p)^\diamond,I}\{t^{1/2},\log(t)\}\otimes_{\mathbb{Q}_p}E)}{\varphi\mathrm{Module}}\right\}_{I\subset (0,\infty)} \ar[d]  \\
\mathrm{Quasicoherent}^{\mathrm{solid,perfectcomplexes}}_{\mathrm{FF}_X,\mathcal{O}_{\mathrm{FF}_X}}  \ar[r] \ar[r] \ar[r] &\left\{\underset{\mathrm{solid,perfectcomplexes},(\Gamma^\text{perfect}_{\text{Robba},X,I}\{t^{1/2},\log(t)\}\otimes_{\mathbb{Q}_p}E)}{\varphi\mathrm{Module}}\right\}_{I\subset (0,\infty)}.   
}
\]
\end{proposition}

\begin{proposition}
We have the following commutative diagram by taking the global section functor in the horizontal rows:
\[\displayindent=+0in
\xymatrix@C-0.12in@R+7pc{
\mathrm{Quasicoherent}^{\mathrm{indBanach}}_{\mathrm{FF}_{\mathrm{Spd}(\mathbb{Q}_p)^\diamond},\mathcal{O}_{\mathrm{FF}_{\mathrm{Spd}(\mathbb{Q}_p)^\diamond}}} \ar[r] \ar[d] &\{\underset{\mathrm{indBanach}}{\varphi\mathrm{Module}}(\Gamma^\text{perfect}_{\text{Robba},\mathrm{Spd}(\mathbb{Q}_p)^\diamond,I}\{t^{1/2},\log(t)\}\otimes_{\mathbb{Q}_p}E)\}_{I\subset (0,\infty)}\ar[d]\\
\mathrm{Quasicoherent}^{\mathrm{indBanach}}_{\mathrm{FF}_X,\mathcal{O}_{\mathrm{FF}_X}}  \ar[r] \ar[r] \ar[r] &\{\underset{\mathrm{indBanach}}{\varphi\mathrm{Module}}(\Gamma^\text{perfect}_{\text{Robba},X,I}\{t^{1/2},\log(t)\}\otimes_{\mathbb{Q}_p}E)\}_{I\subset (0,\infty)}.    
}
\]
\end{proposition}

\begin{proposition}
We have the following commutative diagram by taking the global section functor in the horizontal rows:
\[\displayindent=+0in
\xymatrix@C-0.2in@R+7pc{
\mathrm{Quasicoherent}^{\mathrm{indBanach,perfectcomplexes}}_{\mathrm{FF}_{\mathrm{Spd}(\mathbb{Q}_p)^\diamond},\mathcal{O}_{\mathrm{FF}_{\mathrm{Spd}(\mathbb{Q}_p)^\diamond}}} \ar[d]\ar[r] &\left\{\underset{\mathrm{indBanach,perfectcomplexes},\Gamma^\text{perfect}_{\text{Robba},\mathrm{Spd}(\mathbb{Q}_p)^\diamond,I}\{t^{1/2},\log(t)\}\otimes_{\mathbb{Q}_p}E}{\varphi\mathrm{Module}}\right\}_{I\subset (0,\infty)}\ar[d]\\
\mathrm{Quasicoherent}^{\mathrm{indBanach,perfectcomplexes}}_{\mathrm{FF}_X,\mathcal{O}_{\mathrm{FF}_X}}  \ar[r] \ar[r] \ar[r] &\left\{\underset{\mathrm{indBanach,perfectcomplexes},\Gamma^\text{perfect}_{\text{Robba},X,I}\{t^{1/2},\log(t)\}\otimes_{\mathbb{Q}_p}E}{\varphi\mathrm{Module}}\right\}_{I\subset (0,\infty)}.    
}
\]

\end{proposition}

\indent Taking the corresponding simplicial commutative object we have the following propositions:

\begin{proposition}
We have the following commutative diagram by taking the global section functor in the horizontal rows:
\[\displayindent=+0in
\xymatrix@R+7pc{
\underset{\mathrm{Quasicoherent}^{\mathrm{solid}}_{\mathrm{FF}_{\mathrm{Spd}(\mathbb{Q}_p)^\diamond},\mathcal{O}_{\mathrm{FF}_{\mathrm{Spd}(\mathbb{Q}_p)^\diamond}}}}{\mathrm{SimplicialRings}} \ar[d] \ar[r] &\underset{\{\varphi\mathrm{Module}^{\mathrm{solid}}(\Gamma^\text{perfect}_{\text{Robba},{\mathrm{Spd}(\mathbb{Q}_p)^\diamond},I}\{t^{1/2},\log(t)\}\otimes_{\mathbb{Q}_p}E)\}_{I\subset (0,\infty)}}{\mathrm{SimplicialRings}} \ar[d]\\
\underset{\mathrm{Quasicoherent}^{\mathrm{solid}}_{\mathrm{FF}_X,\mathcal{O}_{\mathrm{FF}_X}}}{\mathrm{SimplicialRings}}  \ar[r] \ar[r] \ar[r] &\underset{\{\varphi\mathrm{Module}^{\mathrm{solid}}(\Gamma^\text{perfect}_{\text{Robba},X,I}\{t^{1/2},\log(t)\}\otimes_{\mathbb{Q}_p}E)\}_{I\subset (0,\infty)}}{\mathrm{SimplicialRings}}.
}
\]
\end{proposition}

\begin{proposition}
We have the following commutative diagram by taking the global section functor in the horizontal rows:
\[\displayindent=+0in
\xymatrix@C+0in@R+7pc{
\underset{\mathrm{Quasicoherent}^{\mathrm{solid,perfectcomplexes}}_{\mathrm{FF}_{\mathrm{Spd}(\mathbb{Q}_p)^\diamond},\mathcal{O}_{\mathrm{FF}_{\mathrm{Spd}(\mathbb{Q}_p)^\diamond}}}}{\mathrm{SimplicialRings}}\ar[d] \ar[r] &\underset{\{\varphi\mathrm{Module}^{\mathrm{solid,perfectcomplexes}}(\Gamma^\text{perfect}_{\text{Robba},{\mathrm{Spd}(\mathbb{Q}_p)^\diamond},I}\{t^{1/2},\log(t)\}\otimes_{\mathbb{Q}_p}E)\}_{I\subset (0,\infty)}}{\mathrm{SimplicialRings}}\ar[d]\\
\underset{\mathrm{Quasicoherent}^{\mathrm{solid,perfectcomplexes}}_{\mathrm{FF}_X,\mathcal{O}_{\mathrm{FF}_X}}}{\mathrm{SimplicialRings}}  \ar[r] \ar[r] \ar[r] &\underset{\{\varphi\mathrm{Module}^{\mathrm{solid,perfectcomplexes}}(\Gamma^\text{perfect}_{\text{Robba},X,I}\{t^{1/2},\log(t)\}\otimes_{\mathbb{Q}_p}E)\}_{I\subset (0,\infty)}}{\mathrm{SimplicialRings}}.  
}
\]
\end{proposition}

\begin{proposition}
We have the following commutative diagram by taking the global section functor in the horizontal rows:
\[\displayindent=+0in
\xymatrix@R+7pc{
\underset{\mathrm{Quasicoherent}^{\mathrm{indBanach}}_{\mathrm{FF}_X,\mathcal{O}_{\mathrm{FF}_X}}}{\mathrm{SimplicialRings}}  \ar[r] \ar[r] \ar[r] &\underset{\{\varphi\mathrm{Module}^{\mathrm{indBanach}}(\Gamma^\text{perfect}_{\text{Robba},X,I}\{t^{1/2},\log(t)\}\otimes_{\mathbb{Q}_p}E)\}_{I\subset (0,\infty)}}{\mathrm{SimplicialRings}}   \\
\underset{\mathrm{Quasicoherent}^{\mathrm{indBanach}}_{\mathrm{FF}_{\mathrm{Spd}(\mathbb{Q}_p)^\diamond},\mathcal{O}_{\mathrm{FF}_{\mathrm{Spd}(\mathbb{Q}_p)^\diamond}}}}{\mathrm{SimplicialRings}} \ar[u]\ar[r] &\underset{\{\varphi\mathrm{Module}^\text{indBanach}(\Gamma^\text{perfect}_{\text{Robba},{\mathrm{Spd}(\mathbb{Q}_p)^\diamond},I}\{t^{1/2},\log(t)\}\otimes_{\mathbb{Q}_p}E)\}_{I\subset (0,\infty)}}{\mathrm{SimplicialRings}}.\ar[u]  
}
\]
\end{proposition}

\begin{proposition}
We have the following commutative diagram by taking the global section functor in the horizontal rows:
\[\displayindent=+0in
\xymatrix@C+0in@R+7pc{
\underset{\mathrm{Quasicoherent}^{\mathrm{indBanach,perfectcomplexes}}_{\mathrm{FF}_X,\mathcal{O}_{\mathrm{FF}_X}}}{\mathrm{SimplicialRings}}  \ar[r] \ar[r] \ar[r] &\underset{\{\varphi\mathrm{Module}^{\mathrm{indBanach,perfectcomplexes}}(\Gamma^\text{perfect}_{\text{Robba},X,I}\{t^{1/2},\log(t)\}\otimes_{\mathbb{Q}_p}E)\}_{I\subset (0,\infty)}}{\mathrm{SimplicialRings}}   \\
\underset{\mathrm{Quasicoherent}^{\mathrm{indBanach,perfectcomplexes}}_{\mathrm{FF}_{\mathrm{Spd}(\mathbb{Q}_p)^\diamond},\mathcal{O}_{\mathrm{FF}_{\mathrm{Spd}(\mathbb{Q}_p)^\diamond}}}}{\mathrm{SimplicialRings}} \ar[u]\ar[r] &\underset{\{\varphi\mathrm{Module}^{\mathrm{indBanach,perfectcomplexes}}(\Gamma^\text{perfect}_{\text{Robba},{\mathrm{Spd}(\mathbb{Q}_p)^\diamond},I}\{t^{1/2},\log(t)\}\otimes_{\mathbb{Q}_p}E)\}_{I\subset (0,\infty)}}{\mathrm{SimplicialRings}}.\ar[u]  
}
\]

\end{proposition}

\newpage
\section{Geometric Family of Mixed-Parity Hodge Modules I: de Rham Situations}

\begin{reference}
\cite{Sch1}, \cite{Sch2}, \cite{FS}, \cite{KL1}, \cite{KL2}, \cite{BL1}, \cite{BL2}, \cite{BS}, \cite{BHS}, \cite{Fon1}, \cite{CS1}, \cite{CS2}, \cite{BK}, \cite{BBK}, \cite{BBBK}, \cite{KKM}, \cite{KM}, \cite{LZ}, \cite{M}.
\end{reference}

\subsection{Period Rings and Sheaves}

\subsubsection{Rings}

\noindent Let $X$ be a $v$-stack over $\mathrm{Spd}\mathbb{Q}_p$, which is required to be restricted to be a diamond which is further assumed to be spacial in the local setting. We have the corresponding \'etale site and the corresponding pro-\'etale site of $X$, which we denote them by $X_{v},X_\text{\'et}$. The relationship of the two sites can be reflected by the corresponding morphism $f:X_{v}\longrightarrow X_\text{\'et}$. Then we have the corresponding de Rham period rings and sheaves from \cite{Sch1}:
\begin{align}
\Gamma_{\text{deRham},X,v}, \Gamma^\mathcal{O}_{\text{deRham},X,v}.
\end{align}
Our notations are different from \cite{Sch1}, we use $\Gamma$ to mean $B$ in \cite{Sch1}, while $\Gamma^\mathcal{O}$ will be the corresponding $OB$ ring in \cite{Sch1}.\\

\begin{definition}
\indent Now we assume that $p>2$, following \cite{BS} we join the square root of $t$ element in $\Gamma_{\text{deRham},X,v}$ which forms the sheaves:
\begin{align}
\Gamma_{\text{deRham},X,v}\{t^{1/2}\},\Gamma^\mathcal{O}_{\text{deRham},X,v}\{t^{1/2}\}.
\end{align}
And following \cite{BL1}, \cite{BL2}, \cite{Fon1}, \cite{BHS} we further have the following sheaves of rings:
\begin{align}
\Gamma_{\text{deRham},X,v}\{t^{1/2},\log(t)\},\Gamma^\mathcal{O}_{\text{deRham},X,v}\{t^{1/2},\log(t)\}.
\end{align}
\end{definition}

\begin{definition}
We use the notations:
\begin{align}
\Gamma^\mathrm{perfect}_{\text{Robba},X,v},\Gamma^\mathrm{perfect}_{\text{Robba},X,v,\infty},\Gamma^\mathrm{perfect}_{\text{Robba},X,v,I}
\end{align}
to denote the perfect Robba rings from \cite{KL1}, \cite{KL2}, where $I\subset (0,\infty)$. Then we join $t^{1/2}$ to these sheaves we have:
\begin{align}
\Gamma^\mathrm{perfect}_{\text{Robba},X,v}\{t^{1/2}\},\Gamma^\mathrm{perfect}_{\text{Robba},X,v,\infty}\{t^{1/2}\},\Gamma^\mathrm{perfect}_{\text{Robba},X,v,I}\{t^{1/2}\}.
\end{align}
And following \cite{BL1}, \cite{BL2}, \cite{Fon1}, \cite{BHS} we have the following larger sheaves:
\begin{align}
\Gamma^\mathrm{perfect}_{\text{Robba},X,v}\{t^{1/2},\log(t)\},\Gamma^\mathrm{perfect}_{\text{Robba},X,v,\infty}\{t^{1/2},\log(t)\},\Gamma^\mathrm{perfect}_{\text{Robba},X,v,I}\{t^{1/2},\log(t)\}.
\end{align} 
\end{definition}

\begin{definition}
From now on, we use the same notation to denote the period rings involved tensored with a finite extension of $\mathbb{Q}_p$ containing square root of $p$ as in \cite{BS}.
\begin{align}
\Gamma_{\text{deRham},X,v}\{t^{1/2}\},\Gamma^\mathcal{O}_{\text{deRham},X,v}\{t^{1/2}\}.
\end{align}
\begin{align}
\Gamma_{\text{deRham},X,v}\{t^{1/2},\log(t)\},\Gamma^\mathcal{O}_{\text{deRham},X,v}\{t^{1/2},\log(t)\}.
\end{align}
\begin{align}
\Gamma^\mathrm{perfect}_{\text{Robba},X,v}\{t^{1/2}\},\Gamma^\mathrm{perfect}_{\text{Robba},X,v,\infty}\{t^{1/2}\},\Gamma^\mathrm{perfect}_{\text{Robba},X,v,I}\{t^{1/2}\}.
\end{align}
\begin{align}
\Gamma^\mathrm{perfect}_{\text{Robba},X,v}\{t^{1/2},\log(t)\},\Gamma^\mathrm{perfect}_{\text{Robba},X,v,\infty}\{t^{1/2},\log(t)\},\Gamma^\mathrm{perfect}_{\text{Robba},X,v,I}\{t^{1/2},\log(t)\}.
\end{align}
This is necessary since we to extend the action of $\varphi$ to the period rings by $\varphi(t^{1/2}\otimes 1)=\varphi(t)^{1/2}\otimes 1$.
\end{definition}

\subsubsection{Modules}

\noindent We consider quasicoherent presheaves in the following two situation:
\begin{itemize}
\item[$\square$] The solid quasicoherent modules from \cite{CS1}, \cite{CS2};
\item[$\square$] The ind-Banach quasicoherent modules from \cite{BK}, \cite{BBK}, \cite{BBBK}, \cite{KKM}, \cite{KM} with the corresponding monomorphic ind-Banach quasicoherent modules from \cite{BK}, \cite{BBK}, \cite{BBBK}, \cite{KKM}, \cite{KM}.
\end{itemize}

\begin{definition}
We use the notation:
\begin{align}
\mathrm{preModule}^\mathrm{solid,quasicoherent}_{\square,\Gamma^\mathrm{perfect}_{\text{Robba},X,v}\{t^{1/2}\}},\mathrm{preModule}^\mathrm{solid,quasicoherent}_{\square,\Gamma^\mathrm{perfect}_{\text{Robba},X,v,\infty}\{t^{1/2}\}},
\mathrm{preModule}^\mathrm{solid,quasicoherent}_{\square,\Gamma^\mathrm{perfect}_{\text{Robba},X,v,I}\{t^{1/2}\}} 
\end{align}
to denote the $(\infty,1)$-categories of solid quasicoherent presheaves over the corresonding Robba sheaves. Locally the section is defined by taking the corresponding $(\infty,1)$-categories of solid modules.
\end{definition}

\begin{definition}
We use the notation:
\begin{align}
\mathrm{preModule}^\mathrm{ind-Banach,quasicoherent}_{\Gamma^\mathrm{perfect}_{\text{Robba},X,v}\{t^{1/2}\}},\\\mathrm{preModule}^\mathrm{ind-Banach,quasicoherent}_{\Gamma^\mathrm{perfect}_{\text{Robba},X,v,\infty}\{t^{1/2}\}},\\
\mathrm{preModule}^\mathrm{ind-Banach,quasicoherent}_{\Gamma^\mathrm{perfect}_{\text{Robba},X,v,I}\{t^{1/2}\}} 
\end{align}
to denote the $(\infty,1)$-categories of solid quasicoherent presheaves over the corresonding Robba sheaves. Locally the section is defined by taking the corresponding $(\infty,1)$-categories of inductive Banach  modules. 
\end{definition}

\begin{definition}
We use the notation:
\begin{align}
\mathrm{Module}^\mathrm{solid,quasicoherent}_{\square,\Gamma^\mathrm{perfect}_{\text{Robba},X,v}\{t^{1/2}\}},\mathrm{Module}^\mathrm{solid,quasicoherent}_{\square,\Gamma^\mathrm{perfect}_{\text{Robba},X,v,\infty}\{t^{1/2}\}},
\mathrm{Module}^\mathrm{solid,quasicoherent}_{\square,\Gamma^\mathrm{perfect}_{\text{Robba},X,v,I}\{t^{1/2}\}} 
\end{align}
to denote the $(\infty,1)$-categories of solid quasicoherent sheaves over the corresonding Robba sheaves. Locally the section is defined by taking the corresponding $(\infty,1)$-categories of solid modules.
\end{definition}

\subsubsection{Mixed-Parity Hodge Modules without Frobenius}

\noindent Now we consider the key objects in our study namely those complexes generated by certain mixed-parity Hodge modules. We start from the following definition.

\begin{definition}
For any locally free coherent sheaf $F$ over
\begin{align}
\Gamma^\mathrm{perfect}_{\text{Robba},X,v,\infty}\{t^{1/2}\},\Gamma^\mathrm{perfect}_{\text{Robba},X,v,I}\{t^{1/2}\},
\end{align} 
we consider the following functor $\mathrm{dR}$ sending $F$ to the following object:
\begin{align}
f_*(F\otimes_{\Gamma^\mathrm{perfect}_{\text{Robba},X,v,\infty}\{t^{1/2}\}} \Gamma^\mathcal{O}_{\text{deRham},X,v}\{t^{1/2}\})
\end{align}
or 
\begin{align}
f_*(F\otimes_{\Gamma^\mathrm{perfect}_{\text{Robba},X,v,I}\{t^{1/2}\}} \Gamma^\mathcal{O}_{\text{deRham},X,v}\{t^{1/2}\}).
\end{align}
We call $F$ mixed-parity de Rham if we have the following isomorphism:
\begin{align}
f^*f_*(F\otimes_{\Gamma^\mathrm{perfect}_{\text{Robba},X,v,\infty}\{t^{1/2}\}} \Gamma^\mathcal{O}_{\text{deRham},X,v}\{t^{1/2}\}) \otimes \Gamma^\mathcal{O}_{\text{deRham},X,v}\{t^{1/2}\} \overset{\sim}{\longrightarrow} F \otimes \Gamma^\mathcal{O}_{\text{deRham},X,v}\{t^{1/2}\} 
\end{align}
or 
\begin{align}
f^*f_*(F\otimes_{\Gamma^\mathrm{perfect}_{\text{Robba},X,v,I}\{t^{1/2}\}} \Gamma^\mathcal{O}_{\text{deRham},X,v}\{t^{1/2}\}) \otimes \Gamma^\mathcal{O}_{\text{deRham},X,v}\{t^{1/2}\} \overset{\sim}{\longrightarrow} F \otimes \Gamma^\mathcal{O}_{\text{deRham},X,v}\{t^{1/2}\}. 
\end{align}
\end{definition}

\begin{definition}
For any locally free coherent sheaf $F$ over
\begin{align}
\Gamma^\mathrm{perfect}_{\text{Robba},X,v,\infty}\{t^{1/2}\},\Gamma^\mathrm{perfect}_{\text{Robba},X,v,I}\{t^{1/2}\},
\end{align} 
we consider the following functor $\mathrm{dR}^\mathrm{almost}$ sending $F$ to the following object:
\begin{align}
f_*(F\otimes_{\Gamma^\mathrm{perfect}_{\text{Robba},X,v,\infty}\{t^{1/2}\}} \Gamma^\mathcal{O}_{\text{deRham},X,v}\{t^{1/2},\log(t)\})
\end{align}
or 
\begin{align}
f_*(F\otimes_{\Gamma^\mathrm{perfect}_{\text{Robba},X,v,I}\{t^{1/2}\}} \Gamma^\mathcal{O}_{\text{deRham},X,v}\{t^{1/2},\log(t)\}).
\end{align}
We call $F$ mixed-parity almost de Rham if we have the following isomorphism:
\begin{align}
f^*f_*(F\otimes_{\Gamma^\mathrm{perfect}_{\text{Robba},X,v,\infty}\{t^{1/2}\}} \Gamma^\mathcal{O}_{\text{deRham},X,v}\{t^{1/2},\log(t)\}) \otimes \Gamma^\mathcal{O}_{\text{deRham},X,v}\{t^{1/2},\log(t)\} \\
\overset{\sim}{\longrightarrow} F \otimes \Gamma^\mathcal{O}_{\text{deRham},X,v}\{t^{1/2},\log(t)\} 
\end{align}
or 
\begin{align}
f^*f_*(F\otimes_{\Gamma^\mathrm{perfect}_{\text{Robba},X,v,I}\{t^{1/2}\}} \Gamma^\mathcal{O}_{\text{deRham},X,v}\{t^{1/2},\log(t)\}) \otimes \Gamma^\mathcal{O}_{\text{deRham},X,v}\{t^{1/2},\log(t)\}\\ \overset{\sim}{\longrightarrow} F \otimes \Gamma^\mathcal{O}_{\text{deRham},X,v}\{t^{1/2},\log(t)\}. 
\end{align}
\end{definition}

\noindent We now define the $(\infty,1)$-categories of mixed-parity de Rham modules and he corresponding mixed-parity almost de Rham modules by using the objects involved to generated these categories:

\begin{definition}
Considering all the mixed parity de Rham bundles (locally finite free) as defined above, we consider the sub-$(\infty,1)$ category of 
\begin{align}
\mathrm{preModule}^\mathrm{solid,quasicoherent}_{\square,\Gamma^\mathrm{perfect}_{\text{Robba},X,v,\infty}\{t^{1/2}\}},
\mathrm{preModule}^\mathrm{solid,quasicoherent}_{\square,\Gamma^\mathrm{perfect}_{\text{Robba},X,v,I}\{t^{1/2}\}} 
\end{align}
generated by the mixed-parity de Rham bundles (locally finite free ones). These are defined to be the $(\infty,1)$-categories of mixed-parity de Rham complexes:
\begin{align}
\mathrm{preModule}^\mathrm{solid,quasicoherent,mixed-paritydeRham}_{\square,\Gamma^\mathrm{perfect}_{\text{Robba},X,v,\infty}\{t^{1/2}\}},
\mathrm{preModule}^\mathrm{solid,quasicoherent,mixed-paritydeRham}_{\square,\Gamma^\mathrm{perfect}_{\text{Robba},X,v,I}\{t^{1/2}\}}. 
\end{align}
\end{definition}

\begin{definition}
Considering all the mixed parity almost de Rham bundles (locally finite free) as defined above, we consider the sub-$(\infty,1)$ category of 
\begin{align}
\mathrm{preModule}^\mathrm{solid,quasicoherent}_{\square,\Gamma^\mathrm{perfect}_{\text{Robba},X,v,\infty}\{t^{1/2}\}},
\mathrm{preModule}^\mathrm{solid,quasicoherent}_{\square,\Gamma^\mathrm{perfect}_{\text{Robba},X,v,I}\{t^{1/2}\}} 
\end{align}
generated by the mixed-parity almost de Rham bundles (locally finite free ones). These are defined to be the $(\infty,1)$-categories of mixed-parity de Rham complexes:
\begin{align}
\mathrm{preModule}^\mathrm{solid,quasicoherent,mixed-parityalmostdeRham}_{\square,\Gamma^\mathrm{perfect}_{\text{Robba},X,v,\infty}\{t^{1/2}\}},\\
\mathrm{preModule}^\mathrm{solid,quasicoherent,mixed-parityalmostdeRham}_{\square,\Gamma^\mathrm{perfect}_{\text{Robba},X,v,I}\{t^{1/2}\}}. 
\end{align}
\end{definition}

\indent Then the corresponding mixed-parity de Rham functors can be extended to these categories:
\begin{align}
\mathrm{preModule}^\mathrm{solid,quasicoherent,mixed-paritydeRham}_{\square,\Gamma^\mathrm{perfect}_{\text{Robba},X,v,\infty}\{t^{1/2}\}},\\
\mathrm{preModule}^\mathrm{solid,quasicoherent,mixed-paritydeRham}_{\square,\Gamma^\mathrm{perfect}_{\text{Robba},X,v,I}\{t^{1/2}\}}, 
\end{align}
and
\begin{align}
\mathrm{preModule}^\mathrm{solid,quasicoherent,mixed-parityalmostdeRham}_{\square,\Gamma^\mathrm{perfect}_{\text{Robba},X,v,\infty}\{t^{1/2}\}},\\
\mathrm{preModule}^\mathrm{solid,quasicoherent,mixed-parityalmostdeRham}_{\square,\Gamma^\mathrm{perfect}_{\text{Robba},X,v,I}\{t^{1/2}\}}. 
\end{align}

\subsubsection{Mixed-Parity Hodge Modules with Frobenius}

\noindent Now we consider the key objects in our study namely those complexes generated by certain mixed-parity Hodge modules. We start from the following definition.

\begin{remark}
All the coherent sheaves over mixed-parity Robba sheaves in this section will carry the corresponding Frobenius morphism $\varphi: F \overset{\sim}{\longrightarrow} \varphi^*F$.
\end{remark}

\begin{definition}
For any locally free coherent sheaf $F$ over
\begin{align}
\Gamma^\mathrm{perfect}_{\text{Robba},X,v,\infty}\{t^{1/2}\},\Gamma^\mathrm{perfect}_{\text{Robba},X,v,I}\{t^{1/2}\},
\end{align} 
we consider the following functor $\mathrm{dR}$ sending $F$ to the following object:
\begin{align}
f_*(F\otimes_{\Gamma^\mathrm{perfect}_{\text{Robba},X,v,\infty}\{t^{1/2}\}} \Gamma^\mathcal{O}_{\text{deRham},X,v}\{t^{1/2}\})
\end{align}
or 
\begin{align}
f_*(F\otimes_{\Gamma^\mathrm{perfect}_{\text{Robba},X,v,I}\{t^{1/2}\}} \Gamma^\mathcal{O}_{\text{deRham},X,v}\{t^{1/2}\}).
\end{align}
We call $F$ mixed-parity de Rham if we have the following isomorphism:
\begin{align}
f^*f_*(F\otimes_{\Gamma^\mathrm{perfect}_{\text{Robba},X,v,\infty}\{t^{1/2}\}} \Gamma^\mathcal{O}_{\text{deRham},X,v}\{t^{1/2}\}) \otimes \Gamma^\mathcal{O}_{\text{deRham},X,v}\{t^{1/2}\} \overset{\sim}{\longrightarrow} F \otimes \Gamma^\mathcal{O}_{\text{deRham},X,v}\{t^{1/2}\} 
\end{align}
or 
\begin{align}
f^*f_*(F\otimes_{\Gamma^\mathrm{perfect}_{\text{Robba},X,v,I}\{t^{1/2}\}} \Gamma^\mathcal{O}_{\text{deRham},X,v}\{t^{1/2}\}) \otimes \Gamma^\mathcal{O}_{\text{deRham},X,v}\{t^{1/2}\} \overset{\sim}{\longrightarrow} F \otimes \Gamma^\mathcal{O}_{\text{deRham},X,v}\{t^{1/2}\}. 
\end{align}
\end{definition}

\begin{definition}
For any locally free coherent sheaf $F$ over
\begin{align}
\Gamma^\mathrm{perfect}_{\text{Robba},X,v,\infty}\{t^{1/2}\},\Gamma^\mathrm{perfect}_{\text{Robba},X,v,I}\{t^{1/2}\},
\end{align} 
we consider the following functor $\mathrm{dR}^\mathrm{almost}$ sending $F$ to the following object:
\begin{align}
f_*(F\otimes_{\Gamma^\mathrm{perfect}_{\text{Robba},X,v,\infty}\{t^{1/2}\}} \Gamma^\mathcal{O}_{\text{deRham},X,v}\{t^{1/2},\log(t)\})
\end{align}
or 
\begin{align}
f_*(F\otimes_{\Gamma^\mathrm{perfect}_{\text{Robba},X,v,I}\{t^{1/2}\}} \Gamma^\mathcal{O}_{\text{deRham},X,v}\{t^{1/2},\log(t)\}).
\end{align}
We call $F$ mixed-parity almost de Rham if we have the following isomorphism:
\begin{align}
f^*f_*(F\otimes_{\Gamma^\mathrm{perfect}_{\text{Robba},X,v,\infty}\{t^{1/2}\}} \Gamma^\mathcal{O}_{\text{deRham},X,v}\{t^{1/2},\log(t)\}) \otimes \Gamma^\mathcal{O}_{\text{deRham},X,v}\{t^{1/2},\log(t)\} \\
\overset{\sim}{\longrightarrow} F \otimes \Gamma^\mathcal{O}_{\text{deRham},X,v}\{t^{1/2},\log(t)\} 
\end{align}
or 
\begin{align}
f^*f_*(F\otimes_{\Gamma^\mathrm{perfect}_{\text{Robba},X,v,I}\{t^{1/2}\}} \Gamma^\mathcal{O}_{\text{deRham},X,v}\{t^{1/2},\log(t)\}) \otimes \Gamma^\mathcal{O}_{\text{deRham},X,v}\{t^{1/2},\log(t)\}\\ \overset{\sim}{\longrightarrow} F \otimes \Gamma^\mathcal{O}_{\text{deRham},X,v}\{t^{1/2},\log(t)\}. 
\end{align}
\end{definition}

\noindent We now define the $(\infty,1)$-categories of mixed-parity de Rham modules and he corresponding mixed-parity almost de Rham modules by using the objects involved to generated these categories:

\begin{definition}
Considering all the mixed parity de Rham bundles (locally finite free) as defined above, we consider the sub-$(\infty,1)$ category of 
\begin{align}
\varphi\mathrm{preModule}^\mathrm{solid,quasicoherent}_{\square,\Gamma^\mathrm{perfect}_{\text{Robba},X,v,\infty}\{t^{1/2}\}},
\varphi\mathrm{preModule}^\mathrm{solid,quasicoherent}_{\square,\Gamma^\mathrm{perfect}_{\text{Robba},X,v,I}\{t^{1/2}\}} 
\end{align}
generated by the mixed-parity de Rham bundles (locally finite free ones). These are defined to be the $(\infty,1)$-categories of mixed-parity de Rham complexes:
\begin{align}
\varphi\mathrm{preModule}^\mathrm{solid,quasicoherent,mixed-paritydeRham}_{\square,\Gamma^\mathrm{perfect}_{\text{Robba},X,v,\infty}\{t^{1/2}\}},
\varphi\mathrm{preModule}^\mathrm{solid,quasicoherent,mixed-paritydeRham}_{\square,\Gamma^\mathrm{perfect}_{\text{Robba},X,v,I}\{t^{1/2}\}}. 
\end{align}
\end{definition}

\begin{definition}
Considering all the mixed parity almost de Rham bundles (locally finite free) as defined above, we consider the sub-$(\infty,1)$ category of 
\begin{align}
\varphi\mathrm{preModule}^\mathrm{solid,quasicoherent}_{\square,\Gamma^\mathrm{perfect}_{\text{Robba},X,v,\infty}\{t^{1/2}\}},
\varphi\mathrm{preModule}^\mathrm{solid,quasicoherent}_{\square,\Gamma^\mathrm{perfect}_{\text{Robba},X,v,I}\{t^{1/2}\}} 
\end{align}
generated by the mixed-parity almost de Rham bundles (locally finite free ones). These are defined to be the $(\infty,1)$-categories of mixed-parity de Rham complexes:
\begin{align}
\varphi\mathrm{preModule}^\mathrm{solid,quasicoherent,mixed-parityalmostdeRham}_{\square,\Gamma^\mathrm{perfect}_{\text{Robba},X,v,\infty}\{t^{1/2}\}},\\
\varphi\mathrm{preModule}^\mathrm{solid,quasicoherent,mixed-parityalmostdeRham}_{\square,\Gamma^\mathrm{perfect}_{\text{Robba},X,v,I}\{t^{1/2}\}}. 
\end{align}
\end{definition}

\indent Then the corresponding mixed-parity de Rham functors can be extended to these categories:
\begin{align}
\varphi\mathrm{preModule}^\mathrm{solid,quasicoherent,mixed-paritydeRham}_{\square,\Gamma^\mathrm{perfect}_{\text{Robba},X,v,\infty}\{t^{1/2}\}},
\varphi\mathrm{preModule}^\mathrm{solid,quasicoherent,mixed-paritydeRham}_{\square,\Gamma^\mathrm{perfect}_{\text{Robba},X,v,I}\{t^{1/2}\}}, 
\end{align}
and
\begin{align}
\varphi\mathrm{preModule}^\mathrm{solid,quasicoherent,mixed-parityalmostdeRham}_{\square,\Gamma^\mathrm{perfect}_{\text{Robba},X,v,\infty}\{t^{1/2}\}},\\
\varphi\mathrm{preModule}^\mathrm{solid,quasicoherent,mixed-parityalmostdeRham}_{\square,\Gamma^\mathrm{perfect}_{\text{Robba},X,v,I}\{t^{1/2}\}}. 
\end{align}

\subsection{Mixed-Parity de Rham Riemann-Hilbert Correspondence}

\indent This chapter will extend the corresponding Riemann-Hilbert correspondence from \cite{Sch1}, \cite{LZ}, \cite{BL1}, \cite{BL2}, \cite{M} to the mixed-parity setting.

\begin{definition}
We define the following Riemann-Hilbert functor $\text{RH}_\text{mixed-parity}$ from the one of categories:
\begin{align}
\mathrm{preModule}^\mathrm{solid,quasicoherent,mixed-paritydeRham}_{\square,\Gamma^\mathrm{perfect}_{\text{Robba},X,v,\infty}\{t^{1/2}\}},
\mathrm{preModule}^\mathrm{solid,quasicoherent,mixed-paritydeRham}_{\square,\Gamma^\mathrm{perfect}_{\text{Robba},X,v,I}\{t^{1/2}\}}, 
\end{align}
and
\begin{align}
\mathrm{preModule}^\mathrm{solid,quasicoherent,mixed-parityalmostdeRham}_{\square,\Gamma^\mathrm{perfect}_{\text{Robba},X,v,\infty}\{t^{1/2}\}},\\
\mathrm{preModule}^\mathrm{solid,quasicoherent,mixed-parityalmostdeRham}_{\square,\Gamma^\mathrm{perfect}_{\text{Robba},X,v,I}\{t^{1/2}\}} 
\end{align}
to $(\infty,1)$-categories in image denoted by:
\begin{align}
\mathrm{preModule}_{X,\text{\'et}}
\end{align}
to be the following functors sending each $F$ in the domain to:
\begin{align}
&\text{RH}_\text{mixed-parity}(F):=f_*(F\otimes_{\Gamma^\mathrm{perfect}_{\text{Robba},X,v,\infty}\{t^{1/2}\}} \Gamma^\mathcal{O}_{\text{deRham},X,v}\{t^{1/2}\}),\\
&\text{RH}_\text{mixed-parity}(F):=f_*(F\otimes_{\Gamma^\mathrm{perfect}_{\text{Robba},X,v,I}\{t^{1/2}\}} \Gamma^\mathcal{O}_{\text{deRham},X,v}\{t^{1/2}\}),\\
&\text{RH}_\text{mixed-parity}(F):=f_*(F\otimes_{\Gamma^\mathrm{perfect}_{\text{Robba},X,v,\infty}\{t^{1/2}\}} \Gamma^\mathcal{O}_{\text{deRham},X,v}\{t^{1/2},\log(t)\}),\\
&\text{RH}_\text{mixed-parity}(F):=f_*(F\otimes_{\Gamma^\mathrm{perfect}_{\text{Robba},X,v,I}\{t^{1/2}\}} \Gamma^\mathcal{O}_{\text{deRham},X,v}\{t^{1/2},\log(t)\}),\\
\end{align}
respectively.

\end{definition}

\begin{definition}
In the situation where we have the Frobenius action we consider the follwing. We define the following Riemann-Hilbert functor $\text{RH}_\text{mixed-parity}$ from the one of categories:
\begin{align}
\varphi\mathrm{preModule}^\mathrm{solid,quasicoherent,mixed-paritydeRham}_{\square,\Gamma^\mathrm{perfect}_{\text{Robba},X,v,\infty}\{t^{1/2}\}},
\varphi\mathrm{preModule}^\mathrm{solid,quasicoherent,mixed-paritydeRham}_{\square,\Gamma^\mathrm{perfect}_{\text{Robba},X,v,I}\{t^{1/2}\}}, 
\end{align}
and
\begin{align}
\varphi\mathrm{preModule}^\mathrm{solid,quasicoherent,mixed-parityalmostdeRham}_{\square,\Gamma^\mathrm{perfect}_{\text{Robba},X,v,\infty}\{t^{1/2}\}},\\
\varphi\mathrm{preModule}^\mathrm{solid,quasicoherent,mixed-parityalmostdeRham}_{\square,\Gamma^\mathrm{perfect}_{\text{Robba},X,v,I}\{t^{1/2}\}} 
\end{align}
to $(\infty,1)$-categories in image denoted by:
\begin{align}
\mathrm{preModule}_{X,\text{\'et}}
\end{align}
to be the following functors sending each $F$ in the domain to:
\begin{align}
&\text{RH}_\text{mixed-parity}(F):=f_*(F\otimes_{\Gamma^\mathrm{perfect}_{\text{Robba},X,v,\infty}\{t^{1/2}\}} \Gamma^\mathcal{O}_{\text{deRham},X,v}\{t^{1/2}\}),\\
&\text{RH}_\text{mixed-parity}(F):=f_*(F\otimes_{\Gamma^\mathrm{perfect}_{\text{Robba},X,v,I}\{t^{1/2}\}} \Gamma^\mathcal{O}_{\text{deRham},X,v}\{t^{1/2}\}),\\
&\text{RH}_\text{mixed-parity}(F):=f_*(F\otimes_{\Gamma^\mathrm{perfect}_{\text{Robba},X,v,\infty}\{t^{1/2}\}} \Gamma^\mathcal{O}_{\text{deRham},X,v}\{t^{1/2},\log(t)\}),\\
&\text{RH}_\text{mixed-parity}(F):=f_*(F\otimes_{\Gamma^\mathrm{perfect}_{\text{Robba},X,v,I}\{t^{1/2}\}} \Gamma^\mathcal{O}_{\text{deRham},X,v}\{t^{1/2},\log(t)\}),\\
\end{align}
respectively.

\end{definition}

\newpage
\section{Geometric Family of Mixed-Parity Hodge Modules II: Cristalline Situations}

\noindent References: \cite{Sch1}, \cite{Sch2}, \cite{FS}, \cite{KL1}, \cite{KL2}, \cite{BL1}, \cite{BL2}, \cite{BS}, \cite{BHS}, \cite{Fon1}, \cite{CS1}, \cite{CS2}, \cite{BK}, \cite{BBK}, \cite{BBBK}, \cite{KKM}, \cite{KM}, \cite{LZ}, \cite{TT}, \cite{M}.

\subsection{Period Rings and Sheaves}

\subsubsection{Rings}

\noindent Let $X$ be a $v$-stack over $\mathrm{Spd}\mathbb{Q}_p$, which is required to be restricted to be a diamond which is further assumed to be spacial in the local setting. We have the corresponding \'etale site and the corresponding pro-\'etale site of $X$, which we denote them by $X_{v},X_\text{\'et}$. The relationship of the two sites can be reflected by the corresponding morphism $f:X_{v}\longrightarrow X_\text{\'et}$. Then we have the corresponding cristalline period rings and sheaves from \cite{TT}:
\begin{align}
\Gamma_{\text{cristalline},X,v}, \Gamma^\mathcal{O}_{\text{cristalline},X,v}.
\end{align}
Our notations are different from \cite{TT}, we use $\Gamma$ to mean $B$ in \cite{TT}, while $\Gamma^\mathcal{O}$ will be the corresponding $OB$ ring in \cite{TT}.\\

\begin{definition}
\indent Now we assume that $p>2$, following \cite{BS} we join the square root of $t$ element in $\Gamma_{\text{cristalline},X,v}$ which forms the sheaves:
\begin{align}
\Gamma_{\text{cristalline},X,v}\{t^{1/2}\},\Gamma^\mathcal{O}_{\text{cristalline},X,v}\{t^{1/2}\}.
\end{align}
And following \cite{BL1}, \cite{BL2}, \cite{Fon1}, \cite{BHS} we further have the following sheaves of rings:
\begin{align}
\Gamma_{\text{cristalline},X,v}\{t^{1/2},\log(t)\},\Gamma^\mathcal{O}_{\text{cristalline},X,v}\{t^{1/2},\log(t)\}.
\end{align}
\end{definition}

\begin{definition}
We use the notations:
\begin{align}
\Gamma^\mathrm{perfect}_{\text{Robba},X,v},\Gamma^\mathrm{perfect}_{\text{Robba},X,v,\infty},\Gamma^\mathrm{perfect}_{\text{Robba},X,v,I}
\end{align}
to denote the perfect Robba rings from \cite{KL1}, \cite{KL2}, where $I\subset (0,\infty)$. Then we join $t^{1/2}$ to these sheaves we have:
\begin{align}
\Gamma^\mathrm{perfect}_{\text{Robba},X,v}\{t^{1/2}\},\Gamma^\mathrm{perfect}_{\text{Robba},X,v,\infty}\{t^{1/2}\},\Gamma^\mathrm{perfect}_{\text{Robba},X,v,I}\{t^{1/2}\}.
\end{align}
And following \cite{BL1}, \cite{BL2}, \cite{Fon1}, \cite{BHS} we have the following larger sheaves:
\begin{align}
\Gamma^\mathrm{perfect}_{\text{Robba},X,v}\{t^{1/2},\log(t)\},\Gamma^\mathrm{perfect}_{\text{Robba},X,v,\infty}\{t^{1/2},\log(t)\},\Gamma^\mathrm{perfect}_{\text{Robba},X,v,I}\{t^{1/2},\log(t)\}.
\end{align} 
\end{definition}

\begin{definition}
From now on, we use the same notation to denote the period rings involved tensored with a finite extension of $\mathbb{Q}_p$ containing square root of $p$ as in \cite{BS}.
\begin{align}
\Gamma_{\text{cristalline},X,v}\{t^{1/2}\},\Gamma^\mathcal{O}_{\text{cristalline},X,v}\{t^{1/2}\}.
\end{align}
\begin{align}
\Gamma_{\text{cristalline},X,v}\{t^{1/2},\log(t)\},\Gamma^\mathcal{O}_{\text{cristalline},X,v}\{t^{1/2},\log(t)\}.
\end{align}
\begin{align}
\Gamma^\mathrm{perfect}_{\text{Robba},X,v}\{t^{1/2}\},\Gamma^\mathrm{perfect}_{\text{Robba},X,v,\infty}\{t^{1/2}\},\Gamma^\mathrm{perfect}_{\text{Robba},X,v,I}\{t^{1/2}\}.
\end{align}
\begin{align}
\Gamma^\mathrm{perfect}_{\text{Robba},X,v}\{t^{1/2},\log(t)\},\Gamma^\mathrm{perfect}_{\text{Robba},X,v,\infty}\{t^{1/2},\log(t)\},\Gamma^\mathrm{perfect}_{\text{Robba},X,v,I}\{t^{1/2},\log(t)\}.
\end{align}
This is necessary since we to extend the action of $\varphi$ to the period rings by $\varphi(t^{1/2}\otimes 1)=\varphi(t)^{1/2}\otimes 1$.
\end{definition}

\subsubsection{Modules}

\noindent We consider quasicoherent presheaves in the following two situation:
\begin{itemize}
\item[$\square$] The solid quasicoherent modules from \cite{CS1}, \cite{CS2};
\item[$\square$] The ind-Banach quasicoherent modules from \cite{BK}, \cite{BBK}, \cite{BBBK}, \cite{KKM}, \cite{KM} with the corresponding monomorphic ind-Banach quasicoherent modules from \cite{BK}, \cite{BBK}, \cite{BBBK}, \cite{KKM}, \cite{KM}.
\end{itemize}

\begin{definition}
We use the notation:
\begin{align}
\mathrm{preModule}^\mathrm{solid,quasicoherent}_{\square,\Gamma^\mathrm{perfect}_{\text{Robba},X,v}\{t^{1/2}\}},\mathrm{preModule}^\mathrm{solid,quasicoherent}_{\square,\Gamma^\mathrm{perfect}_{\text{Robba},X,v,\infty}\{t^{1/2}\}},
\mathrm{preModule}^\mathrm{solid,quasicoherent}_{\square,\Gamma^\mathrm{perfect}_{\text{Robba},X,v,I}\{t^{1/2}\}} 
\end{align}
to denote the $(\infty,1)$-categories of solid quasicoherent presheaves over the corresonding Robba sheaves. Locally the section is defined by taking the corresponding $(\infty,1)$-categories of solid modules.
\end{definition}

\begin{definition}
We use the notation:
\begin{align}
\mathrm{preModule}^\mathrm{ind-Banach,quasicoherent}_{\Gamma^\mathrm{perfect}_{\text{Robba},X,v}\{t^{1/2}\}},\\\mathrm{preModule}^\mathrm{ind-Banach,quasicoherent}_{\Gamma^\mathrm{perfect}_{\text{Robba},X,v,\infty}\{t^{1/2}\}},\\
\mathrm{preModule}^\mathrm{ind-Banach,quasicoherent}_{\Gamma^\mathrm{perfect}_{\text{Robba},X,v,I}\{t^{1/2}\}} 
\end{align}
to denote the $(\infty,1)$-categories of solid quasicoherent presheaves over the corresonding Robba sheaves. Locally the section is defined by taking the corresponding $(\infty,1)$-categories of inductive Banach  modules. 
\end{definition}

\begin{definition}
We use the notation:
\begin{align}
\mathrm{Module}^\mathrm{solid,quasicoherent}_{\square,\Gamma^\mathrm{perfect}_{\text{Robba},X,v}\{t^{1/2}\}},\mathrm{Module}^\mathrm{solid,quasicoherent}_{\square,\Gamma^\mathrm{perfect}_{\text{Robba},X,v,\infty}\{t^{1/2}\}},
\mathrm{Module}^\mathrm{solid,quasicoherent}_{\square,\Gamma^\mathrm{perfect}_{\text{Robba},X,v,I}\{t^{1/2}\}} 
\end{align}
to denote the $(\infty,1)$-categories of solid quasicoherent sheaves over the corresonding Robba sheaves. Locally the section is defined by taking the corresponding $(\infty,1)$-categories of solid modules.
\end{definition}

\subsubsection{Mixed-Parity Hodge Modules without Frobenius}

\noindent Now we consider the key objects in our study namely those complexes generated by certain mixed-parity Hodge modules. We start from the following definition.

\begin{definition}
For any locally free coherent sheaf $F$ over
\begin{align}
\Gamma^\mathrm{perfect}_{\text{Robba},X,v,\infty}\{t^{1/2}\},\Gamma^\mathrm{perfect}_{\text{Robba},X,v,I}\{t^{1/2}\},
\end{align} 
we consider the following functor $\mathrm{dR}$ sending $F$ to the following object:
\begin{align}
f_*(F\otimes_{\Gamma^\mathrm{perfect}_{\text{Robba},X,v,\infty}\{t^{1/2}\}} \Gamma^\mathcal{O}_{\text{cristalline},X,v}\{t^{1/2}\})
\end{align}
or 
\begin{align}
f_*(F\otimes_{\Gamma^\mathrm{perfect}_{\text{Robba},X,v,I}\{t^{1/2}\}} \Gamma^\mathcal{O}_{\text{cristalline},X,v}\{t^{1/2}\}).
\end{align}
We call $F$ mixed-parity cristalline if we have the following isomorphism:
\begin{align}
f^*f_*(F\otimes_{\Gamma^\mathrm{perfect}_{\text{Robba},X,v,\infty}\{t^{1/2}\}} \Gamma^\mathcal{O}_{\text{cristalline},X,v}\{t^{1/2}\}) \otimes \Gamma^\mathcal{O}_{\text{cristalline},X,v}\{t^{1/2}\} \overset{\sim}{\longrightarrow} F \otimes \Gamma^\mathcal{O}_{\text{cristalline},X,v}\{t^{1/2}\} 
\end{align}
or 
\begin{align}
f^*f_*(F\otimes_{\Gamma^\mathrm{perfect}_{\text{Robba},X,v,I}\{t^{1/2}\}} \Gamma^\mathcal{O}_{\text{cristalline},X,v}\{t^{1/2}\}) \otimes \Gamma^\mathcal{O}_{\text{cristalline},X,v}\{t^{1/2}\} \overset{\sim}{\longrightarrow} F \otimes \Gamma^\mathcal{O}_{\text{cristalline},X,v}\{t^{1/2}\}. 
\end{align}
\end{definition}

\begin{definition}
For any locally free coherent sheaf $F$ over
\begin{align}
\Gamma^\mathrm{perfect}_{\text{Robba},X,v,\infty}\{t^{1/2}\},\Gamma^\mathrm{perfect}_{\text{Robba},X,v,I}\{t^{1/2}\},
\end{align} 
we consider the following functor $\mathrm{dR}^\mathrm{almost}$ sending $F$ to the following object:
\begin{align}
f_*(F\otimes_{\Gamma^\mathrm{perfect}_{\text{Robba},X,v,\infty}\{t^{1/2}\}} \Gamma^\mathcal{O}_{\text{cristalline},X,v}\{t^{1/2},\log(t)\})
\end{align}
or 
\begin{align}
f_*(F\otimes_{\Gamma^\mathrm{perfect}_{\text{Robba},X,v,I}\{t^{1/2}\}} \Gamma^\mathcal{O}_{\text{cristalline},X,v}\{t^{1/2},\log(t)\}).
\end{align}
We call $F$ mixed-parity almost cristalline if we have the following isomorphism:
\begin{align}
f^*f_*(F\otimes_{\Gamma^\mathrm{perfect}_{\text{Robba},X,v,\infty}\{t^{1/2}\}} \Gamma^\mathcal{O}_{\text{cristalline},X,v}\{t^{1/2},\log(t)\}) \otimes \Gamma^\mathcal{O}_{\text{cristalline},X,v}\{t^{1/2},\log(t)\} \\
\overset{\sim}{\longrightarrow} F \otimes \Gamma^\mathcal{O}_{\text{cristalline},X,v}\{t^{1/2},\log(t)\} 
\end{align}
or 
\begin{align}
f^*f_*(F\otimes_{\Gamma^\mathrm{perfect}_{\text{Robba},X,v,I}\{t^{1/2}\}} \Gamma^\mathcal{O}_{\text{cristalline},X,v}\{t^{1/2},\log(t)\}) \otimes \Gamma^\mathcal{O}_{\text{cristalline},X,v}\{t^{1/2},\log(t)\}\\ \overset{\sim}{\longrightarrow} F \otimes \Gamma^\mathcal{O}_{\text{cristalline},X,v}\{t^{1/2},\log(t)\}. 
\end{align}
\end{definition}

\noindent We now define the $(\infty,1)$-categories of mixed-parity cristalline modules and he corresponding mixed-parity almost cristalline modules by using the objects involved to generated these categories:

\begin{definition}
Considering all the mixed parity cristalline bundles (locally finite free) as defined above, we consider the sub-$(\infty,1)$ category of 
\begin{align}
\mathrm{preModule}^\mathrm{solid,quasicoherent}_{\square,\Gamma^\mathrm{perfect}_{\text{Robba},X,v,\infty}\{t^{1/2}\}},
\mathrm{preModule}^\mathrm{solid,quasicoherent}_{\square,\Gamma^\mathrm{perfect}_{\text{Robba},X,v,I}\{t^{1/2}\}} 
\end{align}
generated by the mixed-parity cristalline bundles (locally finite free ones). These are defined to be the $(\infty,1)$-categories of mixed-parity cristalline complexes:
\begin{align}
\mathrm{preModule}^\mathrm{solid,quasicoherent,mixed-paritycristalline}_{\square,\Gamma^\mathrm{perfect}_{\text{Robba},X,v,\infty}\{t^{1/2}\}},
\mathrm{preModule}^\mathrm{solid,quasicoherent,mixed-paritycristalline}_{\square,\Gamma^\mathrm{perfect}_{\text{Robba},X,v,I}\{t^{1/2}\}}. 
\end{align}
\end{definition}

\begin{definition}
Considering all the mixed parity almost cristalline bundles (locally finite free) as defined above, we consider the sub-$(\infty,1)$ category of 
\begin{align}
\mathrm{preModule}^\mathrm{solid,quasicoherent}_{\square,\Gamma^\mathrm{perfect}_{\text{Robba},X,v,\infty}\{t^{1/2}\}},
\mathrm{preModule}^\mathrm{solid,quasicoherent}_{\square,\Gamma^\mathrm{perfect}_{\text{Robba},X,v,I}\{t^{1/2}\}} 
\end{align}
generated by the mixed-parity almost cristalline bundles (locally finite free ones). These are defined to be the $(\infty,1)$-categories of mixed-parity cristalline complexes:
\begin{align}
\mathrm{preModule}^\mathrm{solid,quasicoherent,mixed-parityalmostcristalline}_{\square,\Gamma^\mathrm{perfect}_{\text{Robba},X,v,\infty}\{t^{1/2}\}},\\
\mathrm{preModule}^\mathrm{solid,quasicoherent,mixed-parityalmostcristalline}_{\square,\Gamma^\mathrm{perfect}_{\text{Robba},X,v,I}\{t^{1/2}\}}. 
\end{align}
\end{definition}

\indent Then the corresponding mixed-parity cristalline functors can be extended to these categories:
\begin{align}
\mathrm{preModule}^\mathrm{solid,quasicoherent,mixed-paritycristalline}_{\square,\Gamma^\mathrm{perfect}_{\text{Robba},X,v,\infty}\{t^{1/2}\}},\\
\mathrm{preModule}^\mathrm{solid,quasicoherent,mixed-paritycristalline}_{\square,\Gamma^\mathrm{perfect}_{\text{Robba},X,v,I}\{t^{1/2}\}}, 
\end{align}
and
\begin{align}
\mathrm{preModule}^\mathrm{solid,quasicoherent,mixed-parityalmostcristalline}_{\square,\Gamma^\mathrm{perfect}_{\text{Robba},X,v,\infty}\{t^{1/2}\}},\\
\mathrm{preModule}^\mathrm{solid,quasicoherent,mixed-parityalmostcristalline}_{\square,\Gamma^\mathrm{perfect}_{\text{Robba},X,v,I}\{t^{1/2}\}}. 
\end{align}

\subsubsection{Mixed-Parity Hodge Modules with Frobenius}

\noindent Now we consider the key objects in our study namely those complexes generated by certain mixed-parity Hodge modules. We start from the following definition.

\begin{remark}
All the coherent sheaves over mixed-parity Robba sheaves in this section will carry the corresponding Frobenius morphism $\varphi: F \overset{\sim}{\longrightarrow} \varphi^*F$.
\end{remark}

\begin{definition}
For any locally free coherent sheaf $F$ over
\begin{align}
\Gamma^\mathrm{perfect}_{\text{Robba},X,v,\infty}\{t^{1/2}\},\Gamma^\mathrm{perfect}_{\text{Robba},X,v,I}\{t^{1/2}\},
\end{align} 
we consider the following functor $\mathrm{dR}$ sending $F$ to the following object:
\begin{align}
f_*(F\otimes_{\Gamma^\mathrm{perfect}_{\text{Robba},X,v,\infty}\{t^{1/2}\}} \Gamma^\mathcal{O}_{\text{cristalline},X,v}\{t^{1/2}\})
\end{align}
or 
\begin{align}
f_*(F\otimes_{\Gamma^\mathrm{perfect}_{\text{Robba},X,v,I}\{t^{1/2}\}} \Gamma^\mathcal{O}_{\text{cristalline},X,v}\{t^{1/2}\}).
\end{align}
We call $F$ mixed-parity cristalline if we have the following isomorphism:
\begin{align}
f^*f_*(F\otimes_{\Gamma^\mathrm{perfect}_{\text{Robba},X,v,\infty}\{t^{1/2}\}} \Gamma^\mathcal{O}_{\text{cristalline},X,v}\{t^{1/2}\}) \otimes \Gamma^\mathcal{O}_{\text{cristalline},X,v}\{t^{1/2}\} \overset{\sim}{\longrightarrow} F \otimes \Gamma^\mathcal{O}_{\text{cristalline},X,v}\{t^{1/2}\} 
\end{align}
or 
\begin{align}
f^*f_*(F\otimes_{\Gamma^\mathrm{perfect}_{\text{Robba},X,v,\infty}\{t^{1/2}\}} \Gamma^\mathcal{O}_{\text{cristalline},X,v}\{t^{1/2}\}) \otimes \Gamma^\mathcal{O}_{\text{cristalline},X,v}\{t^{1/2}\} \overset{\sim}{\longrightarrow} F \otimes \Gamma^\mathcal{O}_{\text{cristalline},X,v}\{t^{1/2}\}. 
\end{align}
\end{definition}

\begin{definition}
For any locally free coherent sheaf $F$ over
\begin{align}
\Gamma^\mathrm{perfect}_{\text{Robba},X,v,\infty}\{t^{1/2}\},\Gamma^\mathrm{perfect}_{\text{Robba},X,v,I}\{t^{1/2}\},
\end{align} 
we consider the following functor $\mathrm{dR}^\mathrm{almost}$ sending $F$ to the following object:
\begin{align}
f_*(F\otimes_{\Gamma^\mathrm{perfect}_{\text{Robba},X,v,\infty}\{t^{1/2}\}} \Gamma^\mathcal{O}_{\text{cristalline},X,v}\{t^{1/2},\log(t)\})
\end{align}
or 
\begin{align}
f_*(F\otimes_{\Gamma^\mathrm{perfect}_{\text{Robba},X,v,I}\{t^{1/2}\}} \Gamma^\mathcal{O}_{\text{cristalline},X,v}\{t^{1/2},\log(t)\}).
\end{align}
We call $F$ mixed-parity almost cristalline if we have the following isomorphism:
\begin{align}
f^*f_*(F\otimes_{\Gamma^\mathrm{perfect}_{\text{Robba},X,v,\infty}\{t^{1/2}\}} \Gamma^\mathcal{O}_{\text{cristalline},X,v}\{t^{1/2},\log(t)\}) \otimes \Gamma^\mathcal{O}_{\text{cristalline},X,v}\{t^{1/2},\log(t)\} \\
\overset{\sim}{\longrightarrow} F \otimes \Gamma^\mathcal{O}_{\text{cristalline},X,v}\{t^{1/2},\log(t)\} 
\end{align}
or 
\begin{align}
f^*f_*(F\otimes_{\Gamma^\mathrm{perfect}_{\text{Robba},X,v,I}\{t^{1/2}\}} \Gamma^\mathcal{O}_{\text{cristalline},X,v}\{t^{1/2},\log(t)\}) \otimes \Gamma^\mathcal{O}_{\text{cristalline},X,v}\{t^{1/2},\log(t)\}\\ \overset{\sim}{\longrightarrow} F \otimes \Gamma^\mathcal{O}_{\text{cristalline},X,v}\{t^{1/2},\log(t)\}. 
\end{align}
\end{definition}

\noindent We now define the $(\infty,1)$-categories of mixed-parity cristalline modules and he corresponding mixed-parity almost cristalline modules by using the objects involved to generated these categories:

\begin{definition}
Considering all the mixed parity cristalline bundles (locally finite free) as defined above, we consider the sub-$(\infty,1)$ category of 
\begin{align}
\varphi\mathrm{preModule}^\mathrm{solid,quasicoherent}_{\square,\Gamma^\mathrm{perfect}_{\text{Robba},X,v,\infty}\{t^{1/2}\}},
\varphi\mathrm{preModule}^\mathrm{solid,quasicoherent}_{\square,\Gamma^\mathrm{perfect}_{\text{Robba},X,v,I}\{t^{1/2}\}} 
\end{align}
generated by the mixed-parity cristalline bundles (locally finite free ones). These are defined to be the $(\infty,1)$-categories of mixed-parity cristalline complexes:
\begin{align}
\varphi\mathrm{preModule}^\mathrm{solid,quasicoherent,mixed-paritycristalline}_{\square,\Gamma^\mathrm{perfect}_{\text{Robba},X,v,\infty}\{t^{1/2}\}},
\varphi\mathrm{preModule}^\mathrm{solid,quasicoherent,mixed-paritycristalline}_{\square,\Gamma^\mathrm{perfect}_{\text{Robba},X,v,I}\{t^{1/2}\}}. 
\end{align}
\end{definition}

\begin{definition}
Considering all the mixed parity almost cristalline bundles (locally finite free) as defined above, we consider the sub-$(\infty,1)$ category of 
\begin{align}
\varphi\mathrm{preModule}^\mathrm{solid,quasicoherent}_{\square,\Gamma^\mathrm{perfect}_{\text{Robba},X,v,\infty}\{t^{1/2}\}},
\varphi\mathrm{preModule}^\mathrm{solid,quasicoherent}_{\square,\Gamma^\mathrm{perfect}_{\text{Robba},X,v,I}\{t^{1/2}\}} 
\end{align}
generated by the mixed-parity almost cristalline bundles (locally finite free ones). These are defined to be the $(\infty,1)$-categories of mixed-parity cristalline complexes:
\begin{align}
\varphi\mathrm{preModule}^\mathrm{solid,quasicoherent,mixed-parityalmostcristalline}_{\square,\Gamma^\mathrm{perfect}_{\text{Robba},X,v,\infty}\{t^{1/2}\}},\\
\varphi\mathrm{preModule}^\mathrm{solid,quasicoherent,mixed-parityalmostcristalline}_{\square,\Gamma^\mathrm{perfect}_{\text{Robba},X,v,I}\{t^{1/2}\}}. 
\end{align}
\end{definition}

\indent Then the corresponding mixed-parity cristalline functors can be extended to these categories:
\begin{align}
\varphi\mathrm{preModule}^\mathrm{solid,quasicoherent,mixed-paritycristalline}_{\square,\Gamma^\mathrm{perfect}_{\text{Robba},X,v,\infty}\{t^{1/2}\}},
\varphi\mathrm{preModule}^\mathrm{solid,quasicoherent,mixed-paritycristalline}_{\square,\Gamma^\mathrm{perfect}_{\text{Robba},X,v,I}\{t^{1/2}\}}, 
\end{align}
and
\begin{align}
\varphi\mathrm{preModule}^\mathrm{solid,quasicoherent,mixed-parityalmostcristalline}_{\square,\Gamma^\mathrm{perfect}_{\text{Robba},X,v,\infty}\{t^{1/2}\}},\\
\varphi\mathrm{preModule}^\mathrm{solid,quasicoherent,mixed-parityalmostcristalline}_{\square,\Gamma^\mathrm{perfect}_{\text{Robba},X,v,I}\{t^{1/2}\}}. 
\end{align}

\subsection{Mixed-Parity Cristalline Riemann-Hilbert Correspondence}

\indent This chapter will extend the corresponding Riemann-Hilbert correspondence from \cite{Sch1}, \cite{LZ}, \cite{BL1}, \cite{BL2}, \cite{M} to the mixed-parity setting.

\begin{definition}
We define the following Riemann-Hilbert functor $\text{RH}_\text{mixed-parity}$ from the one of categories:
\begin{align}
\mathrm{preModule}^\mathrm{solid,quasicoherent,mixed-paritycristalline}_{\square,\Gamma^\mathrm{perfect}_{\text{Robba},X,v,\infty}\{t^{1/2}\}},
\mathrm{preModule}^\mathrm{solid,quasicoherent,mixed-paritycristalline}_{\square,\Gamma^\mathrm{perfect}_{\text{Robba},X,v,I}\{t^{1/2}\}}, 
\end{align}
and
\begin{align}
\mathrm{preModule}^\mathrm{solid,quasicoherent,mixed-parityalmostcristalline}_{\square,\Gamma^\mathrm{perfect}_{\text{Robba},X,v,\infty}\{t^{1/2}\}},\\
\mathrm{preModule}^\mathrm{solid,quasicoherent,mixed-parityalmostcristalline}_{\square,\Gamma^\mathrm{perfect}_{\text{Robba},X,v,I}\{t^{1/2}\}} 
\end{align}
to $(\infty,1)$-categories in image denoted by:
\begin{align}
\mathrm{preModule}_{X,\text{\'et}}
\end{align}
to be the following functors sending each $F$ in the domain to:
\begin{align}
&\text{RH}_\text{mixed-parity}(F):=f_*(F\otimes_{\Gamma^\mathrm{perfect}_{\text{Robba},X,v,\infty}\{t^{1/2}\}} \Gamma^\mathcal{O}_{\text{cristalline},X,v}\{t^{1/2}\}),\\
&\text{RH}_\text{mixed-parity}(F):=f_*(F\otimes_{\Gamma^\mathrm{perfect}_{\text{Robba},X,v,I}\{t^{1/2}\}} \Gamma^\mathcal{O}_{\text{cristalline},X,v}\{t^{1/2}\}),\\
&\text{RH}_\text{mixed-parity}(F):=f_*(F\otimes_{\Gamma^\mathrm{perfect}_{\text{Robba},X,v,\infty}\{t^{1/2}\}} \Gamma^\mathcal{O}_{\text{cristalline},X,v}\{t^{1/2},\log(t)\}),\\
&\text{RH}_\text{mixed-parity}(F):=f_*(F\otimes_{\Gamma^\mathrm{perfect}_{\text{Robba},X,v,I}\{t^{1/2}\}} \Gamma^\mathcal{O}_{\text{cristalline},X,v}\{t^{1/2},\log(t)\}),\\
\end{align}
respectively.

\end{definition}

\begin{definition}
In the situation where we have the Frobenius action we consider the follwing. We define the following Riemann-Hilbert functor $\text{RH}_\text{mixed-parity}$ from the one of categories:
\begin{align}
\varphi\mathrm{preModule}^\mathrm{solid,quasicoherent,mixed-paritycristalline}_{\square,\Gamma^\mathrm{perfect}_{\text{Robba},X,v,\infty}\{t^{1/2}\}},
\varphi\mathrm{preModule}^\mathrm{solid,quasicoherent,mixed-paritycristalline}_{\square,\Gamma^\mathrm{perfect}_{\text{Robba},X,v,I}\{t^{1/2}\}}, 
\end{align}
and
\begin{align}
\varphi\mathrm{preModule}^\mathrm{solid,quasicoherent,mixed-parityalmostcristalline}_{\square,\Gamma^\mathrm{perfect}_{\text{Robba},X,v,\infty}\{t^{1/2}\}},\\
\varphi\mathrm{preModule}^\mathrm{solid,quasicoherent,mixed-parityalmostcristalline}_{\square,\Gamma^\mathrm{perfect}_{\text{Robba},X,v,I}\{t^{1/2}\}} 
\end{align}
to $(\infty,1)$-categories in image denoted by:
\begin{align}
\mathrm{preModule}_{X,\text{\'et}}
\end{align}
to be the following functors sending each $F$ in the domain to:
\begin{align}
&\text{RH}_\text{mixed-parity}(F):=f_*(F\otimes_{\Gamma^\mathrm{perfect}_{\text{Robba},X,v,\infty}\{t^{1/2}\}} \Gamma^\mathcal{O}_{\text{cristalline},X,v}\{t^{1/2}\}),\\
&\text{RH}_\text{mixed-parity}(F):=f_*(F\otimes_{\Gamma^\mathrm{perfect}_{\text{Robba},X,v,I}\{t^{1/2}\}} \Gamma^\mathcal{O}_{\text{cristalline},X,v}\{t^{1/2}\}),\\
&\text{RH}_\text{mixed-parity}(F):=f_*(F\otimes_{\Gamma^\mathrm{perfect}_{\text{Robba},X,v,\infty}\{t^{1/2}\}} \Gamma^\mathcal{O}_{\text{cristalline},X,v}\{t^{1/2},\log(t)\}),\\
&\text{RH}_\text{mixed-parity}(F):=f_*(F\otimes_{\Gamma^\mathrm{perfect}_{\text{Robba},X,v,I}\{t^{1/2}\}} \Gamma^\mathcal{O}_{\text{cristalline},X,v}\{t^{1/2},\log(t)\}),\\
\end{align}
respectively.

\end{definition}

\newpage
\section{Geometric Family of Mixed-Parity Hodge Modules III: Semi-Stable Situations}

\noindent References: \cite{Sch1}, \cite{Sch2}, \cite{FS}, \cite{KL1}, \cite{KL2}, \cite{BL1}, \cite{BL2}, \cite{BS}, \cite{BHS}, \cite{Fon1}, \cite{CS1}, \cite{CS2}, \cite{BK}, \cite{BBK}, \cite{BBBK}, \cite{KKM}, \cite{KM}, \cite{LZ}, \cite{Shi}, \cite{M}.

\subsection{Period Rings and Sheaves}

\subsubsection{Rings}

\noindent Let $X$ be a $v$-stack over $\mathrm{Spd}\mathbb{Q}_p$, which is required to be restricted to be a diamond which is further assumed to be spacial in the local setting. We have the corresponding \'etale site and the corresponding pro-\'etale site of $X$, which we denote them by $X_{v},X_\text{\'et}$. The relationship of the two sites can be reflected by the corresponding morphism $f:X_{v}\longrightarrow X_\text{\'et}$. Then we have the corresponding semi-stable period rings and sheaves from \cite{Shi}:
\begin{align}
\Gamma_{\text{semistable},X,v}, \Gamma^\mathcal{O}_{\text{semistable},X,v}.
\end{align}
Our notations are different from \cite{Shi}, we use $\Gamma$ to mean $B$ in \cite{Shi}, while $\Gamma^\mathcal{O}$ will be the corresponding $OB$ ring in \cite{Shi}.\\

\begin{definition}
\indent Now we assume that $p>2$, following \cite{BS} we join the square root of $t$ element in $\Gamma_{\text{semistable},X,v}$ which forms the sheaves:
\begin{align}
\Gamma_{\text{semistable},X,v}\{t^{1/2}\},\Gamma^\mathcal{O}_{\text{semistable},X,v}\{t^{1/2}\}.
\end{align}
And following \cite{BL1}, \cite{BL2}, \cite{Fon1}, \cite{BHS} we further have the following sheaves of rings:
\begin{align}
\Gamma_{\text{semistable},X,v}\{t^{1/2},\log(t)\},\Gamma^\mathcal{O}_{\text{semistable},X,v}\{t^{1/2},\log(t)\}.
\end{align}
\end{definition}

\begin{definition}
We use the notations:
\begin{align}
\Gamma^\mathrm{perfect}_{\text{Robba},X,v},\Gamma^\mathrm{perfect}_{\text{Robba},X,v,\infty},\Gamma^\mathrm{perfect}_{\text{Robba},X,v,I}
\end{align}
to denote the perfect Robba rings from \cite{KL1}, \cite{KL2}, where $I\subset (0,\infty)$. Then we join $t^{1/2}$ to these sheaves we have:
\begin{align}
\Gamma^\mathrm{perfect}_{\text{Robba},X,v}\{t^{1/2}\},\Gamma^\mathrm{perfect}_{\text{Robba},X,v,\infty}\{t^{1/2}\},\Gamma^\mathrm{perfect}_{\text{Robba},X,v,I}\{t^{1/2}\}.
\end{align}
And following \cite{BL1}, \cite{BL2}, \cite{Fon1}, \cite{BHS} we have the following larger sheaves:
\begin{align}
\Gamma^\mathrm{perfect}_{\text{Robba},X,v}\{t^{1/2},\log(t)\},\Gamma^\mathrm{perfect}_{\text{Robba},X,v,\infty}\{t^{1/2},\log(t)\},\Gamma^\mathrm{perfect}_{\text{Robba},X,v,I}\{t^{1/2},\log(t)\}.
\end{align} 
\end{definition}

\begin{definition}
From now on, we use the same notation to denote the period rings involved tensored with a finite extension of $\mathbb{Q}_p$ containing square root of $p$ as in \cite{BS}.
\begin{align}
\Gamma_{\text{semistable},X,v}\{t^{1/2}\},\Gamma^\mathcal{O}_{\text{semistable},X,v}\{t^{1/2}\}.
\end{align}
\begin{align}
\Gamma_{\text{semistable},X,v}\{t^{1/2},\log(t)\},\Gamma^\mathcal{O}_{\text{semistable},X,v}\{t^{1/2},\log(t)\}.
\end{align}
\begin{align}
\Gamma^\mathrm{perfect}_{\text{Robba},X,v}\{t^{1/2}\},\Gamma^\mathrm{perfect}_{\text{Robba},X,v,\infty}\{t^{1/2}\},\Gamma^\mathrm{perfect}_{\text{Robba},X,v,I}\{t^{1/2}\}.
\end{align}
\begin{align}
\Gamma^\mathrm{perfect}_{\text{Robba},X,v}\{t^{1/2},\log(t)\},\Gamma^\mathrm{perfect}_{\text{Robba},X,v,\infty}\{t^{1/2},\log(t)\},\Gamma^\mathrm{perfect}_{\text{Robba},X,v,I}\{t^{1/2},\log(t)\}.
\end{align}
This is necessary since we to extend the action of $\varphi$ to the period rings by $\varphi(t^{1/2}\otimes 1)=\varphi(t)^{1/2}\otimes 1$.
\end{definition}

\subsubsection{Modules}

\noindent We consider quasicoherent presheaves in the following two situation:
\begin{itemize}
\item[$\square$] The solid quasicoherent modules from \cite{CS1}, \cite{CS2};
\item[$\square$] The ind-Banach quasicoherent modules from \cite{BK}, \cite{BBK}, \cite{BBBK}, \cite{KKM}, \cite{KM} with the corresponding monomorphic ind-Banach quasicoherent modules from \cite{BK}, \cite{BBK}, \cite{BBBK}, \cite{KKM}, \cite{KM}.
\end{itemize}

\begin{definition}
We use the notation:
\begin{align}
\mathrm{preModule}^\mathrm{solid,quasicoherent}_{\square,\Gamma^\mathrm{perfect}_{\text{Robba},X,v}\{t^{1/2}\}},\mathrm{preModule}^\mathrm{solid,quasicoherent}_{\square,\Gamma^\mathrm{perfect}_{\text{Robba},X,v,\infty}\{t^{1/2}\}},
\mathrm{preModule}^\mathrm{solid,quasicoherent}_{\square,\Gamma^\mathrm{perfect}_{\text{Robba},X,v,I}\{t^{1/2}\}} 
\end{align}
to denote the $(\infty,1)$-categories of solid quasicoherent presheaves over the corresonding Robba sheaves. Locally the section is defined by taking the corresponding $(\infty,1)$-categories of solid modules.
\end{definition}

\begin{definition}
We use the notation:
\begin{align}
\mathrm{preModule}^\mathrm{ind-Banach,quasicoherent}_{\Gamma^\mathrm{perfect}_{\text{Robba},X,v}\{t^{1/2}\}},\\\mathrm{preModule}^\mathrm{ind-Banach,quasicoherent}_{\Gamma^\mathrm{perfect}_{\text{Robba},X,v,\infty}\{t^{1/2}\}},\\
\mathrm{preModule}^\mathrm{ind-Banach,quasicoherent}_{\Gamma^\mathrm{perfect}_{\text{Robba},X,v,I}\{t^{1/2}\}} 
\end{align}
to denote the $(\infty,1)$-categories of solid quasicoherent presheaves over the corresonding Robba sheaves. Locally the section is defined by taking the corresponding $(\infty,1)$-categories of inductive Banach  modules. 
\end{definition}

\begin{definition}
We use the notation:
\begin{align}
\mathrm{Module}^\mathrm{solid,quasicoherent}_{\square,\Gamma^\mathrm{perfect}_{\text{Robba},X,v}\{t^{1/2}\}},\mathrm{Module}^\mathrm{solid,quasicoherent}_{\square,\Gamma^\mathrm{perfect}_{\text{Robba},X,v,\infty}\{t^{1/2}\}},
\mathrm{Module}^\mathrm{solid,quasicoherent}_{\square,\Gamma^\mathrm{perfect}_{\text{Robba},X,v,I}\{t^{1/2}\}} 
\end{align}
to denote the $(\infty,1)$-categories of solid quasicoherent sheaves over the corresonding Robba sheaves. Locally the section is defined by taking the corresponding $(\infty,1)$-categories of solid modules.
\end{definition}

\subsubsection{Mixed-Parity Hodge Modules without Frobenius}

\noindent Now we consider the key objects in our study namely those complexes generated by certain mixed-parity Hodge modules. We start from the following definition.

\begin{definition}
For any locally free coherent sheaf $F$ over
\begin{align}
\Gamma^\mathrm{perfect}_{\text{Robba},X,v,\infty}\{t^{1/2}\},\Gamma^\mathrm{perfect}_{\text{Robba},X,v,I}\{t^{1/2}\},
\end{align} 
we consider the following functor $\mathrm{dR}$ sending $F$ to the following object:
\begin{align}
f_*(F\otimes_{\Gamma^\mathrm{perfect}_{\text{Robba},X,v,\infty}\{t^{1/2}\}} \Gamma^\mathcal{O}_{\text{semistable},X,v}\{t^{1/2}\})
\end{align}
or 
\begin{align}
f_*(F\otimes_{\Gamma^\mathrm{perfect}_{\text{Robba},X,v,I}\{t^{1/2}\}} \Gamma^\mathcal{O}_{\text{semistable},X,v}\{t^{1/2}\}).
\end{align}
We call $F$ mixed-parity semi-stable if we have the following isomorphism:
\begin{align}
f^*f_*(F\otimes_{\Gamma^\mathrm{perfect}_{\text{Robba},X,v,\infty}\{t^{1/2}\}} \Gamma^\mathcal{O}_{\text{semistable},X,v}\{t^{1/2}\}) \otimes \Gamma^\mathcal{O}_{\text{semistable},X,v}\{t^{1/2}\} \overset{\sim}{\longrightarrow} F \otimes \Gamma^\mathcal{O}_{\text{semistable},X,v}\{t^{1/2}\} 
\end{align}
or 
\begin{align}
f^*f_*(F\otimes_{\Gamma^\mathrm{perfect}_{\text{Robba},X,v,I}\{t^{1/2}\}} \Gamma^\mathcal{O}_{\text{semistable},X,v}\{t^{1/2}\}) \otimes \Gamma^\mathcal{O}_{\text{semistable},X,v}\{t^{1/2}\} \overset{\sim}{\longrightarrow} F \otimes \Gamma^\mathcal{O}_{\text{semistable},X,v}\{t^{1/2}\}. 
\end{align}
\end{definition}

\begin{definition}
For any locally free coherent sheaf $F$ over
\begin{align}
\Gamma^\mathrm{perfect}_{\text{Robba},X,v,\infty}\{t^{1/2}\},\Gamma^\mathrm{perfect}_{\text{Robba},X,v,I}\{t^{1/2}\},
\end{align} 
we consider the following functor $\mathrm{dR}^\mathrm{almost}$ sending $F$ to the following object:
\begin{align}
f_*(F\otimes_{\Gamma^\mathrm{perfect}_{\text{Robba},X,v,\infty}\{t^{1/2}\}} \Gamma^\mathcal{O}_{\text{semistable},X,v}\{t^{1/2},\log(t)\})
\end{align}
or 
\begin{align}
f_*(F\otimes_{\Gamma^\mathrm{perfect}_{\text{Robba},X,v,I}\{t^{1/2}\}} \Gamma^\mathcal{O}_{\text{semistable},X,v}\{t^{1/2},\log(t)\}).
\end{align}
We call $F$ mixed-parity almost semi-stable if we have the following isomorphism:
\begin{align}
f^*f_*(F\otimes_{\Gamma^\mathrm{perfect}_{\text{Robba},X,v,\infty}\{t^{1/2}\}} \Gamma^\mathcal{O}_{\text{semistable},X,v}\{t^{1/2},\log(t)\}) \otimes \Gamma^\mathcal{O}_{\text{semistable},X,v}\{t^{1/2},\log(t)\} \\
\overset{\sim}{\longrightarrow} F \otimes \Gamma^\mathcal{O}_{\text{semistable},X,v}\{t^{1/2},\log(t)\} 
\end{align}
or 
\begin{align}
f^*f_*(F\otimes_{\Gamma^\mathrm{perfect}_{\text{Robba},X,v,I}\{t^{1/2}\}} \Gamma^\mathcal{O}_{\text{semistable},X,v}\{t^{1/2},\log(t)\}) \otimes \Gamma^\mathcal{O}_{\text{semistable},X,v}\{t^{1/2},\log(t)\}\\ \overset{\sim}{\longrightarrow} F \otimes \Gamma^\mathcal{O}_{\text{semistable},X,v}\{t^{1/2},\log(t)\}. 
\end{align}
\end{definition}

\noindent We now define the $(\infty,1)$-categories of mixed-parity semi-stable modules and he corresponding mixed-parity almost semi-stable modules by using the objects involved to generated these categories:

\begin{definition}
Considering all the mixed parity semi-stable bundles (locally finite free) as defined above, we consider the sub-$(\infty,1)$ category of 
\begin{align}
\mathrm{preModule}^\mathrm{solid,quasicoherent}_{\square,\Gamma^\mathrm{perfect}_{\text{Robba},X,v,\infty}\{t^{1/2}\}},
\mathrm{preModule}^\mathrm{solid,quasicoherent}_{\square,\Gamma^\mathrm{perfect}_{\text{Robba},X,v,I}\{t^{1/2}\}} 
\end{align}
generated by the mixed-parity semi-stable bundles (locally finite free ones). These are defined to be the $(\infty,1)$-categories of mixed-parity semi-stable complexes:
\begin{align}
\mathrm{preModule}^\mathrm{solid,quasicoherent,mixed-paritysemistable}_{\square,\Gamma^\mathrm{perfect}_{\text{Robba},X,v,\infty}\{t^{1/2}\}},
\mathrm{preModule}^\mathrm{solid,quasicoherent,mixed-paritysemistable}_{\square,\Gamma^\mathrm{perfect}_{\text{Robba},X,v,I}\{t^{1/2}\}}. 
\end{align}
\end{definition}

\begin{definition}
Considering all the mixed parity almost semi-stable bundles (locally finite free) as defined above, we consider the sub-$(\infty,1)$ category of 
\begin{align}
\mathrm{preModule}^\mathrm{solid,quasicoherent}_{\square,\Gamma^\mathrm{perfect}_{\text{Robba},X,v,\infty}\{t^{1/2}\}},
\mathrm{preModule}^\mathrm{solid,quasicoherent}_{\square,\Gamma^\mathrm{perfect}_{\text{Robba},X,v,I}\{t^{1/2}\}} 
\end{align}
generated by the mixed-parity almost semi-stable bundles (locally finite free ones). These are defined to be the $(\infty,1)$-categories of mixed-parity semi-stable complexes:
\begin{align}
\mathrm{preModule}^\mathrm{solid,quasicoherent,mixed-parityalmostsemistable}_{\square,\Gamma^\mathrm{perfect}_{\text{Robba},X,v,\infty}\{t^{1/2}\}},\\
\mathrm{preModule}^\mathrm{solid,quasicoherent,mixed-parityalmostsemistable}_{\square,\Gamma^\mathrm{perfect}_{\text{Robba},X,v,I}\{t^{1/2}\}}. 
\end{align}
\end{definition}

\indent Then the corresponding mixed-parity semi-stable functors can be extended to these categories:
\begin{align}
\mathrm{preModule}^\mathrm{solid,quasicoherent,mixed-paritysemistable}_{\square,\Gamma^\mathrm{perfect}_{\text{Robba},X,v,\infty}\{t^{1/2}\}},
\mathrm{preModule}^\mathrm{solid,quasicoherent,mixed-paritysemistable}_{\square,\Gamma^\mathrm{perfect}_{\text{Robba},X,v,I}\{t^{1/2}\}}, 
\end{align}
and
\begin{align}
\mathrm{preModule}^\mathrm{solid,quasicoherent,mixed-parityalmostsemistable}_{\square,\Gamma^\mathrm{perfect}_{\text{Robba},X,v,\infty}\{t^{1/2}\}},\\
\mathrm{preModule}^\mathrm{solid,quasicoherent,mixed-parityalmostsemistable}_{\square,\Gamma^\mathrm{perfect}_{\text{Robba},X,v,I}\{t^{1/2}\}}. 
\end{align}

\subsubsection{Mixed-Parity Hodge Modules with Frobenius}

\noindent Now we consider the key objects in our study namely those complexes generated by certain mixed-parity Hodge modules. We start from the following definition.

\begin{remark}
All the coherent sheaves over mixed-parity Robba sheaves in this section will carry the corresponding Frobenius morphism $\varphi: F \overset{\sim}{\longrightarrow} \varphi^*F$.
\end{remark}

\begin{definition}
For any locally free coherent sheaf $F$ over
\begin{align}
\Gamma^\mathrm{perfect}_{\text{Robba},X,v,\infty}\{t^{1/2}\},\Gamma^\mathrm{perfect}_{\text{Robba},X,v,I}\{t^{1/2}\},
\end{align} 
we consider the following functor $\mathrm{dR}$ sending $F$ to the following object:
\begin{align}
f_*(F\otimes_{\Gamma^\mathrm{perfect}_{\text{Robba},X,v,\infty}\{t^{1/2}\}} \Gamma^\mathcal{O}_{\text{semistable},X,v}\{t^{1/2}\})
\end{align}
or 
\begin{align}
f_*(F\otimes_{\Gamma^\mathrm{perfect}_{\text{Robba},X,v,I}\{t^{1/2}\}} \Gamma^\mathcal{O}_{\text{semistable},X,v}\{t^{1/2}\}).
\end{align}
We call $F$ mixed-parity semi-stable if we have the following isomorphism:
\begin{align}
f^*f_*(F\otimes_{\Gamma^\mathrm{perfect}_{\text{Robba},X,v,\infty}\{t^{1/2}\}} \Gamma^\mathcal{O}_{\text{semistable},X,v}\{t^{1/2}\}) \otimes \Gamma^\mathcal{O}_{\text{semistable},X,v}\{t^{1/2}\} \overset{\sim}{\longrightarrow} F \otimes \Gamma^\mathcal{O}_{\text{semistable},X,v}\{t^{1/2}\} 
\end{align}
or 
\begin{align}
f^*f_*(F\otimes_{\Gamma^\mathrm{perfect}_{\text{Robba},X,v,I}\{t^{1/2}\}} \Gamma^\mathcal{O}_{\text{semistable},X,v}\{t^{1/2}\}) \otimes \Gamma^\mathcal{O}_{\text{semistable},X,v}\{t^{1/2}\} \overset{\sim}{\longrightarrow} F \otimes \Gamma^\mathcal{O}_{\text{semistable},X,v}\{t^{1/2}\}. 
\end{align}
\end{definition}

\begin{definition}
For any locally free coherent sheaf $F$ over
\begin{align}
\Gamma^\mathrm{perfect}_{\text{Robba},X,v,\infty}\{t^{1/2}\},\Gamma^\mathrm{perfect}_{\text{Robba},X,v,I}\{t^{1/2}\},
\end{align} 
we consider the following functor $\mathrm{dR}^\mathrm{almost}$ sending $F$ to the following object:
\begin{align}
f_*(F\otimes_{\Gamma^\mathrm{perfect}_{\text{Robba},X,v,\infty}\{t^{1/2}\}} \Gamma^\mathcal{O}_{\text{semistable},X,v}\{t^{1/2},\log(t)\})
\end{align}
or 
\begin{align}
f_*(F\otimes_{\Gamma^\mathrm{perfect}_{\text{Robba},X,v,I}\{t^{1/2}\}} \Gamma^\mathcal{O}_{\text{semistable},X,v}\{t^{1/2},\log(t)\}).
\end{align}
We call $F$ mixed-parity almost semi-stable if we have the following isomorphism:
\begin{align}
f^*f_*(F\otimes_{\Gamma^\mathrm{perfect}_{\text{Robba},X,v,\infty}\{t^{1/2}\}} \Gamma^\mathcal{O}_{\text{semistable},X,v}\{t^{1/2},\log(t)\}) \otimes \Gamma^\mathcal{O}_{\text{semistable},X,v}\{t^{1/2},\log(t)\} \\
\overset{\sim}{\longrightarrow} F \otimes \Gamma^\mathcal{O}_{\text{semistable},X,v}\{t^{1/2},\log(t)\} 
\end{align}
or 
\begin{align}
f^*f_*(F\otimes_{\Gamma^\mathrm{perfect}_{\text{Robba},X,v,I}\{t^{1/2}\}} \Gamma^\mathcal{O}_{\text{semistable},X,v}\{t^{1/2},\log(t)\}) \otimes \Gamma^\mathcal{O}_{\text{semistable},X,v}\{t^{1/2},\log(t)\}\\ \overset{\sim}{\longrightarrow} F \otimes \Gamma^\mathcal{O}_{\text{semistable},X,v}\{t^{1/2},\log(t)\}. 
\end{align}
\end{definition}

\noindent We now define the $(\infty,1)$-categories of mixed-parity semi-stable modules and he corresponding mixed-parity almost semi-stable modules by using the objects involved to generated these categories:

\begin{definition}
Considering all the mixed parity semi-stable bundles (locally finite free) as defined above, we consider the sub-$(\infty,1)$ category of 
\begin{align}
\varphi\mathrm{preModule}^\mathrm{solid,quasicoherent}_{\square,\Gamma^\mathrm{perfect}_{\text{Robba},X,v,\infty}\{t^{1/2}\}},
\varphi\mathrm{preModule}^\mathrm{solid,quasicoherent}_{\square,\Gamma^\mathrm{perfect}_{\text{Robba},X,v,I}\{t^{1/2}\}} 
\end{align}
generated by the mixed-parity semi-stable bundles (locally finite free ones). These are defined to be the $(\infty,1)$-categories of mixed-parity semi-stable complexes:
\begin{align}
\varphi\mathrm{preModule}^\mathrm{solid,quasicoherent,mixed-paritysemistable}_{\square,\Gamma^\mathrm{perfect}_{\text{Robba},X,v,\infty}\{t^{1/2}\}},
\varphi\mathrm{preModule}^\mathrm{solid,quasicoherent,mixed-paritysemistable}_{\square,\Gamma^\mathrm{perfect}_{\text{Robba},X,v,I}\{t^{1/2}\}}. 
\end{align}
\end{definition}

\begin{definition}
Considering all the mixed parity almost semi-stable bundles (locally finite free) as defined above, we consider the sub-$(\infty,1)$ category of 
\begin{align}
\varphi\mathrm{preModule}^\mathrm{solid,quasicoherent}_{\square,\Gamma^\mathrm{perfect}_{\text{Robba},X,v,\infty}\{t^{1/2}\}},
\varphi\mathrm{preModule}^\mathrm{solid,quasicoherent}_{\square,\Gamma^\mathrm{perfect}_{\text{Robba},X,v,I}\{t^{1/2}\}} 
\end{align}
generated by the mixed-parity almost semi-stable bundles (locally finite free ones). These are defined to be the $(\infty,1)$-categories of mixed-parity semi-stable complexes:
\begin{align}
\varphi\mathrm{preModule}^\mathrm{solid,quasicoherent,mixed-parityalmostsemistable}_{\square,\Gamma^\mathrm{perfect}_{\text{Robba},X,v,\infty}\{t^{1/2}\}},\\
\varphi\mathrm{preModule}^\mathrm{solid,quasicoherent,mixed-parityalmostsemistable}_{\square,\Gamma^\mathrm{perfect}_{\text{Robba},X,v,I}\{t^{1/2}\}}. 
\end{align}
\end{definition}

\indent Then the corresponding mixed-parity semi-stable functors can be extended to these categories:
\begin{align}
\varphi\mathrm{preModule}^\mathrm{solid,quasicoherent,mixed-paritysemistable}_{\square,\Gamma^\mathrm{perfect}_{\text{Robba},X,v,\infty}\{t^{1/2}\}},
\varphi\mathrm{preModule}^\mathrm{solid,quasicoherent,mixed-paritysemistable}_{\square,\Gamma^\mathrm{perfect}_{\text{Robba},X,v,I}\{t^{1/2}\}}, 
\end{align}
and
\begin{align}
\varphi\mathrm{preModule}^\mathrm{solid,quasicoherent,mixed-parityalmostsemistable}_{\square,\Gamma^\mathrm{perfect}_{\text{Robba},X,v,\infty}\{t^{1/2}\}},\\
\varphi\mathrm{preModule}^\mathrm{solid,quasicoherent,mixed-parityalmostsemistable}_{\square,\Gamma^\mathrm{perfect}_{\text{Robba},X,v,I}\{t^{1/2}\}}. 
\end{align}

\subsection{Mixed-Parity semi-stable Riemann-Hilbert Correspondence}

\indent This chapter will extend the corresponding Riemann-Hilbert correspondence from \cite{Sch1}, \cite{LZ}, \cite{BL1}, \cite{BL2}, \cite{M} to the mixed-parity setting.

\begin{definition}
We define the following Riemann-Hilbert functor $\text{RH}_\text{mixed-parity}$ from the one of categories:
\begin{align}
\mathrm{preModule}^\mathrm{solid,quasicoherent,mixed-paritysemistable}_{\square,\Gamma^\mathrm{perfect}_{\text{Robba},X,v,\infty}\{t^{1/2}\}},
\mathrm{preModule}^\mathrm{solid,quasicoherent,mixed-paritysemistable}_{\square,\Gamma^\mathrm{perfect}_{\text{Robba},X,v,I}\{t^{1/2}\}}, 
\end{align}
and
\begin{align}
\mathrm{preModule}^\mathrm{solid,quasicoherent,mixed-parityalmostsemistable}_{\square,\Gamma^\mathrm{perfect}_{\text{Robba},X,v,\infty}\{t^{1/2}\}},\\
\mathrm{preModule}^\mathrm{solid,quasicoherent,mixed-parityalmostsemistable}_{\square,\Gamma^\mathrm{perfect}_{\text{Robba},X,v,I}\{t^{1/2}\}} 
\end{align}
to $(\infty,1)$-categories in image denoted by:
\begin{align}
\mathrm{preModule}_{X,\text{\'et}}
\end{align}
to be the following functors sending each $F$ in the domain to:
\begin{align}
&\text{RH}_\text{mixed-parity}(F):=f_*(F\otimes_{\Gamma^\mathrm{perfect}_{\text{Robba},X,v,\infty}\{t^{1/2}\}} \Gamma^\mathcal{O}_{\text{semistable},X,v}\{t^{1/2}\}),\\
&\text{RH}_\text{mixed-parity}(F):=f_*(F\otimes_{\Gamma^\mathrm{perfect}_{\text{Robba},X,v,I}\{t^{1/2}\}} \Gamma^\mathcal{O}_{\text{semistable},X,v}\{t^{1/2}\}),\\
&\text{RH}_\text{mixed-parity}(F):=f_*(F\otimes_{\Gamma^\mathrm{perfect}_{\text{Robba},X,v,\infty}\{t^{1/2}\}} \Gamma^\mathcal{O}_{\text{semistable},X,v}\{t^{1/2},\log(t)\}),\\
&\text{RH}_\text{mixed-parity}(F):=f_*(F\otimes_{\Gamma^\mathrm{perfect}_{\text{Robba},X,v,I}\{t^{1/2}\}} \Gamma^\mathcal{O}_{\text{semistable},X,v}\{t^{1/2},\log(t)\}),\\
\end{align}
respectively.

\end{definition}

\begin{definition}
In the situation where we have the Frobenius action we consider the follwing. We define the following Riemann-Hilbert functor $\text{RH}_\text{mixed-parity}$ from the one of categories:
\begin{align}
\varphi\mathrm{preModule}^\mathrm{solid,quasicoherent,mixed-paritysemistable}_{\square,\Gamma^\mathrm{perfect}_{\text{Robba},X,v,\infty}\{t^{1/2}\}},
\varphi\mathrm{preModule}^\mathrm{solid,quasicoherent,mixed-paritysemistable}_{\square,\Gamma^\mathrm{perfect}_{\text{Robba},X,v,I}\{t^{1/2}\}}, 
\end{align}
and
\begin{align}
\varphi\mathrm{preModule}^\mathrm{solid,quasicoherent,mixed-parityalmostsemistable}_{\square,\Gamma^\mathrm{perfect}_{\text{Robba},X,v,\infty}\{t^{1/2}\}},\\
\varphi\mathrm{preModule}^\mathrm{solid,quasicoherent,mixed-parityalmostsemistable}_{\square,\Gamma^\mathrm{perfect}_{\text{Robba},X,v,I}\{t^{1/2}\}} 
\end{align}
to $(\infty,1)$-categories in image denoted by:
\begin{align}
\mathrm{preModule}_{X,\text{\'et}}
\end{align}
to be the following functors sending each $F$ in the domain to:
\begin{align}
&\text{RH}_\text{mixed-parity}(F):=f_*(F\otimes_{\Gamma^\mathrm{perfect}_{\text{Robba},X,v,\infty}\{t^{1/2}\}} \Gamma^\mathcal{O}_{\text{semistable},X,v}\{t^{1/2}\}),\\
&\text{RH}_\text{mixed-parity}(F):=f_*(F\otimes_{\Gamma^\mathrm{perfect}_{\text{Robba},X,v,I}\{t^{1/2}\}} \Gamma^\mathcal{O}_{\text{semistable},X,v}\{t^{1/2}\}),\\
&\text{RH}_\text{mixed-parity}(F):=f_*(F\otimes_{\Gamma^\mathrm{perfect}_{\text{Robba},X,v,\infty}\{t^{1/2}\}} \Gamma^\mathcal{O}_{\text{semistable},X,v}\{t^{1/2},\log(t)\}),\\
&\text{RH}_\text{mixed-parity}(F):=f_*(F\otimes_{\Gamma^\mathrm{perfect}_{\text{Robba},X,v,I}\{t^{1/2}\}} \Gamma^\mathcal{O}_{\text{semistable},X,v}\{t^{1/2},\log(t)\}),\\
\end{align}
respectively.

\end{definition}

\begin{remark}
We now have discussed the corresponding two different morphisms:
\begin{align}
f: X_\text{pro\'et}\longrightarrow X_\text{\'et};\\
f': X_\text{v}\longrightarrow X_\text{\'et}.
\end{align}
One can consider the following relation among the sites:
\begin{align}
X_\text{v}\longrightarrow X_\text{pro\'et}\longrightarrow X_\text{\'et}
\end{align}
which produces $f'$. The map:
\begin{align}
g: X_\text{v}\longrightarrow X_\text{pro\'et}
\end{align}
can help us relate the corresponding constructions above as in \cite[Proposition 2.37]{B}. Namely we have:
\begin{align}
&\mathrm{dR}_{v}=\mathrm{dR}_{\text{pro\'et}}g_*;\\
&\mathrm{dR}_{v,\text{almost}}=\mathrm{dR}_{\text{pro\'et},\text{almost}}g_*;\\
&\mathrm{cristalline}_{v}=\mathrm{cristalline}_{\text{pro\'et}}g_*;\\
&\mathrm{cristalline}_{v,\text{almost}}=\mathrm{cristalline}_{\text{pro\'et},\text{almost}}g_*;\\
&\mathrm{semistable}_{v}=\mathrm{semistable}_{\text{pro\'et}}g_*;\\
&\mathrm{semistable}_{v,\text{almost}}=\mathrm{semistable}_{\text{pro\'et},\text{almost}}g_*.
\end{align} 
\end{remark}

\chapter{Generalized Langlands Program}

\newpage
\section{Moduli $v$-Stack}

\noindent References: \cite{FS}, \cite{FF}, \cite{Sch1},\cite{Sch2}, \cite{KL1}, \cite{KL2};

\noindent Further References:\cite{Lan1}, \cite{Drin1}, \cite{Drin2}, \cite{Zhu}, \cite{DHKM}.\\

\noindent We consider the category of all the perfectoid spaces over $\overline{{\mathbb{Q}}_p((\mu_{p^\infty}))^\wedge}^{\wedge,\flat}$ as in \cite{FS}. We use the notation Perfectoid\_v to denote the associated $v$-site after \cite{FS}, \cite{Sch2}. Let $p>2$. For any $\mathrm{Spa}(A,A^+)\in \text{perfectoid}_{v}$ we have the perfect Robba rings from \cite{KL1}, \cite{KL2}:
\begin{align}
\Gamma^\text{perfect}_{\text{Robba},\mathrm{Spa}(A,A^+),I\subset (0,\infty)}.
\end{align}
We also have the corresponding de Rham period rings:
\begin{align}
\Gamma^+_{\text{deRham},\mathrm{Spa}(A,A^+)},\Gamma_{\text{deRham},\mathrm{Spa}(A,A^+)}.
\end{align}
In the first filtration of this first de Rham period ring we have the generator $t$, we now extend the corresponding rings above by adding the square root of $t$, $t^{1/2}
$ following \cite{BS}. We then have the extended rings:
\begin{align}
\Gamma^\text{perfect}_{\text{Robba},\mathrm{Spa}(A,A^+),I\subset (0,\infty)}\{t^{1/2}\},
\end{align}
\begin{align}
\Gamma^+_{\text{deRham},\mathrm{Spa}(A,A^+)}\{t^{1/2}\},\Gamma_{\text{deRham},\mathrm{Spa}(A,A^+)}\{t^{1/2}\}.
\end{align}
Then we form the corresponding extended Fargues-Fontaine curve (after choosing a large finite extension $E$ of $\mathbb{Q}_p$ containing $\varphi(t)^{1/2}$):
\begin{align}
\mathrm{FF}_A:=\bigcup_{I\subset (0,\infty)}\mathrm{Spa}(\Gamma^\text{perfect}_{\text{Robba},\mathrm{Spa}(A,A^+),I\subset (0,\infty)}\{t^{1/2}\}\otimes_{\mathbb{Q}_p}E,\Gamma^{\text{perfect},+}_{\text{Robba},\mathrm{Spa}(A,A^+),I\subset (0,\infty)}\{t^{1/2}\}\otimes_{\mathbb{Q}_p}E)/\varphi^\mathbb{Z},
\end{align}
where the Frobenius is extended to $t^{1/2}\otimes 1$ by acting $\varphi(t)^{1/2}\otimes 1$.

\begin{definition}
Let $G$ be any $p$-adic group as in \cite{FS}\footnote{That is to say the group $G$ is defined over $\mathbb{Q}_p$. And the Robba rings are defined over $\mathbb{Q}_p$ as well, which strictly speaking are generated from Witt vectors in \cite{KL1}, but one can generalize this directly to the level of \cite{KL2} by replacing the field $\overline{{\mathbb{Q}}_p((\mu_{p^\infty}))^\wedge}^{\wedge,\flat}$ with some larger field $\overline{F((\mu_{p^\infty}))^\wedge}^{\wedge,\flat}$, where $F/\mathbb{Q}_p$ is finite extension of $\mathbb{Q}_p$.}. We now define the pre-v-stack $\text{Moduli}_G$ to be a presheaf valued in the groupoid over
\begin{align}
\text{perfectoid}_{v} 
\end{align}
sendind each $\mathrm{Spa}(A,A^+)$ perfectoid in the site to the groupoid of all the locally finite free coherent sheaves carrying $G$-bundle structure over  
\begin{align}
\mathrm{FF}_A:=\bigcup_{I\subset (0,\infty)}\mathrm{Spa}(\Gamma^\text{perfect}_{\text{Robba},\mathrm{Spa}(A,A^+),I\subset (0,\infty)}\{t^{1/2}\}\otimes_{\mathbb{Q}_p}E,\Gamma^{\text{perfect},+}_{\text{Robba},\mathrm{Spa}(A,A^+),I\subset (0,\infty)}\{t^{1/2}\}\otimes_{\mathbb{Q}_p}E)/\varphi^\mathbb{Z}.
\end{align}
\end{definition}

\begin{proposition}
This prestack is a small $v$-stack in the $v$-topology.
\end{proposition}

\begin{proof}
The proof will be the same as in \cite[Proposition III.1.3]{FS}. Our stack can also be regarded as a two components extension of the original stack in \cite{FS}.  
\end{proof}

\section{Motives over $\mathrm{Moduli}_G$}

\noindent With the notation in the previous section, we now consider the sheaves over extended Fargues-Fontain stacks:

\begin{definition}
\begin{align}
\mathrm{FF}_{\mathrm{Moduli}_G}:=\bigcup_{I\subset (0,\infty)}\mathrm{Spa}(\Gamma^\text{perfect}_{\text{Robba},{\mathrm{Moduli}_G},I\subset (0,\infty)}\{t^{1/2}\}\otimes_{\mathbb{Q}_p}E,\Gamma^{\text{perfect},+}_{\text{Robba},{\mathrm{Moduli}_G},I\subset (0,\infty)}\{t^{1/2}\}\otimes_{\mathbb{Q}_p}E)/\varphi^\mathbb{Z},
\end{align}
which has the corresonding structure map as in the following:
\[\displayindent=-0.4in
\xymatrix@R+1pc{
&\mathrm{FF}_{\mathrm{Moduli}_G} \ar[d]  \\
&\mathrm{FF}_{\mathrm{FF}_*} \ar[d]\\
&\mathrm{FF}_{\mathrm{Spd}^\diamond(\mathbb{Q}_p)} \ar[d]\\
& \mathrm{Spd}^\diamond(\mathbb{Q}_p).  
}
\]
\end{definition}

\begin{definition}
We use the notation
\begin{align}
\mathrm{Quasicoherent}^{\mathrm{solid}}_{\mathrm{FF}_{\mathrm{Moduli}_G},\mathcal{O}_{\mathrm{FF}_{\mathrm{Moduli}_G}}}
\end{align}
to denote $(\infty,1)$-category of all the solid quasicoherent sheaves over the stack $\mathrm{FF}_{\mathrm{Moduli}_G}$. For any local perfectoid $Y\in {\mathrm{Moduli}_G}_v$ we define the corresponding $(\infty,1)$-category in the local sense.\\
We use the notation
\begin{align}
\mathrm{Quasicoherent}^{\mathrm{solid,perfectcomplexes}}_{\mathrm{FF}_{\mathrm{Moduli}_G},\mathcal{O}_{\mathrm{FF}_{\mathrm{Moduli}_G}}}
\end{align}
to denote $(\infty,1)$-category of all the solid quasicoherent sheaves over the stack $\mathrm{FF}_{\mathrm{Moduli}_G}$ which are perfect complexes. For any local perfectoid $Y\in {\mathrm{Moduli}_G}_v$ we define the corresponding $(\infty,1)$-category in the local sense.
\end{definition}

\begin{definition}
We use the notation
\begin{align}
\mathrm{Quasicoherent}^{\mathrm{indBanach}}_{\mathrm{FF}_{\mathrm{Moduli}_G},\mathcal{O}_{\mathrm{FF}_{\mathrm{Moduli}_G}}}
\end{align}
to denote $(\infty,1)$-category of all the ind-Banach quasicoherent sheaves over the stack $\mathrm{FF}_{\mathrm{Moduli}_G}$. For any local perfectoid $Y\in {\mathrm{Moduli}_G}_v$ we define the corresponding $(\infty,1)$-category in the local sense.\\
We use the notation
\begin{align}
\mathrm{Quasicoherent}^{\mathrm{indBanach,perfectcomplexes}}_{\mathrm{FF}_{\mathrm{Moduli}_G},\mathcal{O}_{\mathrm{FF}_{\mathrm{Moduli}_G}}}
\end{align}
to denote $(\infty,1)$-category of all the indBanach quasicoherent sheaves over the stack $\mathrm{FF}_{\mathrm{Moduli}_G}$ which are perfect complexes. For any local perfectoid $Y\in {\mathrm{Moduli}_G}_v$ we define the corresponding $(\infty,1)$-category in the local sense.
\end{definition}

\begin{definition}
We use the notation
\begin{align}
\{\varphi\mathrm{Module}^{\mathrm{solid}}(\Gamma^\text{perfect}_{\text{Robba},{\mathrm{Moduli}_G},I}\{t^{1/2}\}\otimes_{\mathbb{Q}_p}E)\}_{I\subset (0,\infty)}
\end{align}
to denote $(\infty,1)$-category of all the solid $\varphi$-modules over the extended Robba ring. The modules satisfy the Frobenius pullback condition and glueing condition for overlapped intervals $I\subset J\subset K$. For any local perfectoid $Y\in {\mathrm{Moduli}_G}_v$ we define the corresponding $(\infty,1)$-category in the local sense.\\
We use the notation
\begin{align}
\{\varphi\mathrm{Module}^{\mathrm{solid,perfectcomplexes}}(\Gamma^\text{perfect}_{\text{Robba},{\mathrm{Moduli}_G},I}\{t^{1/2}\}\otimes_{\mathbb{Q}_p}E)\}_{I\subset (0,\infty)}
\end{align}
to denote $(\infty,1)$-category of all the solid $\varphi$-modules over the extended Robba ring which are perfect complexes. The modules satisfy the Frobenius pullback condition and glueing condition for overlapped intervals $I\subset J\subset K$. For any local perfectoid $Y\in {\mathrm{Moduli}_G}_v$ we define the corresponding $(\infty,1)$-category in the local sense. 

\end{definition}

\begin{definition}
We use the notation
\begin{align}
\{\varphi\mathrm{Module}^{\mathrm{indBanach}}(\Gamma^\text{perfect}_{\text{Robba},{\mathrm{Moduli}_G},I}\{t^{1/2}\}\otimes_{\mathbb{Q}_p}E)\}_{I\subset (0,\infty)}
\end{align}
to denote $(\infty,1)$-category of all the ind-Banach $\varphi$-modules over the extended Robba ring. The modules satisfy the Frobenius pullback condition and glueing condition for overlapped intervals $I\subset J\subset K$. For any local perfectoid $Y\in {\mathrm{Moduli}_G}_v$ we define the corresponding $(\infty,1)$-category in the local sense.\\
We use the notation
\begin{align}
\{\varphi\mathrm{Module}^{\mathrm{ind-Banach,perfectcomplexes}}(\Gamma^\text{perfect}_{\text{Robba},{\mathrm{Moduli}_G},I}\{t^{1/2}\}\otimes_{\mathbb{Q}_p}E)\}_{I\subset (0,\infty)}
\end{align}
to denote $(\infty,1)$-category of all the ind-Banach $\varphi$-modules over the extended Robba ring which are perfect complexes. The modules satisfy the Frobenius pullback condition and glueing condition for overlapped intervals $I\subset J\subset K$. For any local perfectoid $Y\in {\mathrm{Moduli}_G}_v$ we define the corresponding $(\infty,1)$-category in the local sense. 

\end{definition}

\begin{definition}
We use the notation
\begin{align}
\varphi\mathrm{Module}^{\mathrm{solid}}(\Gamma^\text{perfect}_{\text{Robba},{\mathrm{Moduli}_G},\infty}\{t^{1/2}\}\otimes_{\mathbb{Q}_p}E)
\end{align}
to denote $(\infty,1)$-category of all the solid $\varphi$-modules over the extended Robba ring. The modules satisfy the Frobenius pullback condition. For any local perfectoid $Y\in {\mathrm{Moduli}_G}_v$ we define the corresponding $(\infty,1)$-category in the local sense.\\
We use the notation
\begin{align}
\varphi\mathrm{Module}^{\mathrm{solid,perfectcomplexes}}(\Gamma^\text{perfect}_{\text{Robba},{\mathrm{Moduli}_G},\infty}\{t^{1/2}\}\otimes_{\mathbb{Q}_p}E)
\end{align}
to denote $(\infty,1)$-category of all the solid $\varphi$-modules over the extended Robba ring which are perfect complexes. The modules satisfy the Frobenius pullback condition. For any local perfectoid $Y\in {\mathrm{Moduli}_G}_v$ we define the corresponding $(\infty,1)$-category in the local sense. 

\end{definition}

\begin{definition}
We use the notation
\begin{align}
\varphi\mathrm{Module}^{\mathrm{indBanach}}(\Gamma^\text{perfect}_{\text{Robba},{\mathrm{Moduli}_G},\infty}\{t^{1/2}\}\otimes_{\mathbb{Q}_p}E)
\end{align}
to denote $(\infty,1)$-category of all the ind-Banach $\varphi$-modules over the extended Robba ring. The modules satisfy the Frobenius pullback condition. For any local perfectoid $Y\in {\mathrm{Moduli}_G}_v$ we define the corresponding $(\infty,1)$-category in the local sense.\\
We use the notation
\begin{align}
\left\{\underset{\mathrm{indBanach,perfectcomplexes},\Gamma^\text{perfect}_{\text{Robba},X,I}\{t^{1/2}\}\otimes_{\mathbb{Q}_p}E}{\varphi\mathrm{Module}}\right\}_{I\subset (0,\infty)}
\end{align}
to denote $(\infty,1)$-category of all the ind-Banach $\varphi$-modules over the extended Robba ring which are perfect complexes. The modules satisfy the Frobenius pullback condition. For any local perfectoid $Y\in {\mathrm{Moduli}_G}_v$ we define the corresponding $(\infty,1)$-category in the local sense. 

\end{definition}

\begin{proposition}
We have the following commutative diagram by taking the global section functor in the horizontal rows:\\
\[\displayindent=+0in
\xymatrix@R+7pc{
\mathrm{Quasicoherent}^{\mathrm{solid}}_{\mathrm{FF}_{\mathrm{Spd}(\mathbb{Q}_p)^\diamond},\mathcal{O}_{\mathrm{FF}_{\mathrm{Spd}(\mathbb{Q}_p)^\diamond}}} \ar[r] \ar[d] &\{\varphi\mathrm{Module}^{\mathrm{solid}}(\Gamma^\text{perfect}_{\text{Robba},{\mathrm{Spd}(\mathbb{Q}_p)^\diamond},I}\{t^{1/2}\}\otimes_{\mathbb{Q}_p}E)\}_{I\subset (0,\infty)} \ar[d]\\
\mathrm{Quasicoherent}^{\mathrm{solid}}_{\mathrm{FF}_{\mathrm{Moduli}_G},\mathcal{O}_{\mathrm{FF}_{\mathrm{Moduli}_G}}}  \ar[r] \ar[r] \ar[r] &\{\varphi\mathrm{Module}^{\mathrm{solid}}(\Gamma^\text{perfect}_{\text{Robba},{\mathrm{Moduli}_G},I}\{t^{1/2}\}\otimes_{\mathbb{Q}_p}E)\}_{I\subset (0,\infty)}.  \\  
}
\]
\end{proposition}

\begin{proposition}
We have the following commutative diagram by taking the global section functor in the horizontal rows:\\
\[\displayindent=+0in
\xymatrix@R+7pc{
\mathrm{Quasicoherent}^{\mathrm{solid,perfectcomplexes}}_{\mathrm{FF}_{\mathrm{Spd}(\mathbb{Q}_p)^\diamond},\mathcal{O}_{\mathrm{FF}_{\mathrm{Spd}(\mathbb{Q}_p)^\diamond}}} \ar[r] \ar[d] &\{\underset{\mathrm{solid,perfectcomplexes}}{\varphi\mathrm{Module}}(\Gamma^\text{perfect}_{\text{Robba},{\mathrm{Spd}(\mathbb{Q}_p)^\diamond},I}\{t^{1/2}\}\otimes_{\mathbb{Q}_p}E)\}_{I\subset (0,\infty)} \ar[d]  \\
\mathrm{Quasicoherent}^{\mathrm{solid,perfectcomplexes}}_{\mathrm{FF}_{\mathrm{Moduli}_G},\mathcal{O}_{\mathrm{FF}_{\mathrm{Moduli}_G}}}  \ar[r] \ar[r] \ar[r] &\{\underset{\mathrm{solid,perfectcomplexes}}{\varphi\mathrm{Module}}(\Gamma^\text{perfect}_{\text{Robba},{\mathrm{Moduli}_G},I}\{t^{1/2}\}\otimes_{\mathbb{Q}_p}E)\}_{I\subset (0,\infty)}.   
}
\]
\end{proposition}

\begin{proposition}
We have the following commutative diagram by taking the global section functor in the horizontal rows:\\
\[\displayindent=+0in
\xymatrix@C-0.1in@R+7pc{
\mathrm{Quasicoherent}^{\mathrm{indBanach}}_{\mathrm{FF}_{\mathrm{Spd}(\mathbb{Q}_p)^\diamond},\mathcal{O}_{\mathrm{FF}_{\mathrm{Spd}(\mathbb{Q}_p)^\diamond}}} \ar[r] \ar[d] &\{\varphi\mathrm{Module}^\text{indBanach}(\Gamma^\text{perfect}_{\text{Robba},{\mathrm{Spd}(\mathbb{Q}_p)^\diamond},I}\{t^{1/2}\}\otimes_{\mathbb{Q}_p}E)\}_{I\subset (0,\infty)}\ar[d]\\
\mathrm{Quasicoherent}^{\mathrm{indBanach}}_{\mathrm{FF}_{\mathrm{Moduli}_G},\mathcal{O}_{\mathrm{FF}_{\mathrm{Moduli}_G}}}  \ar[r] \ar[r] \ar[r] &\{\varphi\mathrm{Module}^{\mathrm{indBanach}}(\Gamma^\text{perfect}_{\text{Robba},{\mathrm{Moduli}_G},I}\{t^{1/2}\}\otimes_{\mathbb{Q}_p}E)\}_{I\subset (0,\infty)}.    
}
\]
\end{proposition}

\begin{proposition}
We have the following commutative diagram by taking the global section functor in the horizontal rows:
\[\displayindent=+0in
\xymatrix@C-0.2in@R+7pc{
\mathrm{Quasicoherent}^{\mathrm{indBanach,perfectcomplexes}}_{\mathrm{FF}_{\mathrm{Spd}(\mathbb{Q}_p)^\diamond},\mathcal{O}_{\mathrm{FF}_{\mathrm{Spd}(\mathbb{Q}_p)^\diamond}}} \ar[d]\ar[r] &
\left\{\underset{\mathrm{indBanach,perfectcomplexes},\Gamma^\text{perfect}_{\text{Robba},\mathrm{Spd}(\mathbb{Q}_p)^\diamond,I}\{t^{1/2}\}\otimes_{\mathbb{Q}_p}E}{\varphi\mathrm{Module}}\right\}_{I\subset (0,\infty)}\ar[d]\\
\mathrm{Quasicoherent}^{\mathrm{indBanach,perfectcomplexes}}_{\mathrm{FF}_{\mathrm{Moduli}_G},\mathcal{O}_{\mathrm{FF}_{\mathrm{Moduli}_G}}}  \ar[r] \ar[r] \ar[r] &
\left\{\underset{\mathrm{indBanach,perfectcomplexes},\Gamma^\text{perfect}_{\text{Robba},\mathrm{Moduli}_G,I}\{t^{1/2}\}\otimes_{\mathbb{Q}_p}E}{\varphi\mathrm{Module}}\right\}_{I\subset (0,\infty)}.    
}
\]

\end{proposition}

\indent Taking the corresponding simplicial commutative object we have the following propositions:

\begin{proposition}
We have the following commutative diagram by taking the global section functor in the horizontal rows:
\[\displayindent=+0in
\xymatrix@R+7pc{
\underset{\mathrm{Quasicoherent}^{\mathrm{solid}}_{\mathrm{FF}_{\mathrm{Spd}(\mathbb{Q}_p)^\diamond},\mathcal{O}_{\mathrm{FF}_{\mathrm{Spd}(\mathbb{Q}_p)^\diamond}}}}{\mathrm{SimplicialRings}} \ar[d] \ar[r] &\underset{\{\varphi\mathrm{Module}^{\mathrm{solid}}(\Gamma^\text{perfect}_{\text{Robba},{\mathrm{Spd}(\mathbb{Q}_p)^\diamond},I}\{t^{1/2}\}\otimes_{\mathbb{Q}_p}E)\}_{I\subset (0,\infty)}}{\mathrm{SimplicialRings}} \ar[d]\\
\underset{\mathrm{Quasicoherent}^{\mathrm{solid}}_{\mathrm{FF}_{\mathrm{Moduli}_G},\mathcal{O}_{\mathrm{FF}_{\mathrm{Moduli}_G}}}}{\mathrm{SimplicialRings}}  \ar[r] \ar[r] \ar[r] &\underset{\{\varphi\mathrm{Module}^{\mathrm{solid}}(\Gamma^\text{perfect}_{\text{Robba},{\mathrm{Moduli}_G},I}\{t^{1/2}\}\otimes_{\mathbb{Q}_p}E)\}_{I\subset (0,\infty)}}{\mathrm{SimplicialRings}}.
}
\]
\end{proposition}

\begin{proposition}
We have the following commutative diagram by taking the global section functor in the horizontal rows:
\[\displayindent=+0in
\xymatrix@C+0in@R+7pc{
\underset{\mathrm{Quasicoherent}^{\mathrm{solid,perfectcomplexes}}_{\mathrm{FF}_{\mathrm{Spd}(\mathbb{Q}_p)^\diamond},\mathcal{O}_{\mathrm{FF}_{\mathrm{Spd}(\mathbb{Q}_p)^\diamond}}}}{\mathrm{SimplicialRings}}\ar[d] \ar[r] &\underset{\{\varphi\mathrm{Module}^{\mathrm{solid,perfectcomplexes}}(\Gamma^\text{perfect}_{\text{Robba},{\mathrm{Spd}(\mathbb{Q}_p)^\diamond},I}\{t^{1/2}\}\otimes_{\mathbb{Q}_p}E)\}_{I\subset (0,\infty)}}{\mathrm{SimplicialRings}}\ar[d]\\
\underset{\mathrm{Quasicoherent}^{\mathrm{solid,perfectcomplexes}}_{\mathrm{FF}_{\mathrm{Moduli}_G},\mathcal{O}_{\mathrm{FF}_{\mathrm{Moduli}_G}}}}{\mathrm{SimplicialRings}}  \ar[r] \ar[r] \ar[r] &\underset{\{\varphi\mathrm{Module}^{\mathrm{solid,perfectcomplexes}}(\Gamma^\text{perfect}_{\text{Robba},{\mathrm{Moduli}_G},I}\{t^{1/2}\}\otimes_{\mathbb{Q}_p}E)\}_{I\subset (0,\infty)}}{\mathrm{SimplicialRings}}.  
}
\]
\end{proposition}

\begin{proposition}
We have the following commutative diagram by taking the global section functor in the horizontal rows:
\[\displayindent=+0in
\xymatrix@R+7pc{
\underset{\mathrm{Quasicoherent}^{\mathrm{indBanach}}_{\mathrm{FF}_{\mathrm{Moduli}_G},\mathcal{O}_{\mathrm{FF}_{\mathrm{Moduli}_G}}}}{\mathrm{SimplicialRings}}  \ar[r] \ar[r] \ar[r] &\underset{\{\varphi\mathrm{Module}^{\mathrm{indBanach}}(\Gamma^\text{perfect}_{\text{Robba},{\mathrm{Moduli}_G},I}\{t^{1/2}\}\otimes_{\mathbb{Q}_p}E)\}_{I\subset (0,\infty)}}{\mathrm{SimplicialRings}}   \\
\underset{\mathrm{Quasicoherent}^{\mathrm{indBanach}}_{\mathrm{FF}_{\mathrm{Spd}(\mathbb{Q}_p)^\diamond},\mathcal{O}_{\mathrm{FF}_{\mathrm{Spd}(\mathbb{Q}_p)^\diamond}}}}{\mathrm{SimplicialRings}} \ar[u]\ar[r] &\underset{\{\varphi\mathrm{Module}^\text{indBanach}(\Gamma^\text{perfect}_{\text{Robba},{\mathrm{Spd}(\mathbb{Q}_p)^\diamond},I}\{t^{1/2}\}\otimes_{\mathbb{Q}_p}E)\}_{I\subset (0,\infty)}}{\mathrm{SimplicialRings}}.\ar[u]  
}
\]
\end{proposition}

\begin{proposition}
We have the following commutative diagram by taking the global section functor in the horizontal rows:
\[\displayindent=+0in
\xymatrix@C-0.1in@R+7pc{
\underset{\mathrm{Quasicoherent}^{\mathrm{indBanach,perfectcomplexes}}_{\mathrm{FF}_{\mathrm{Moduli}_G},\mathcal{O}_{\mathrm{FF}_{\mathrm{Moduli}_G}}}}{\mathrm{SimplicialRings}}  \ar[r] \ar[r] \ar[r] &\underset{\{\varphi\mathrm{Module}^{\mathrm{indBanach,perfectcomplexes}}(\Gamma^\text{perfect}_{\text{Robba},{\mathrm{Moduli}_G},I}\{t^{1/2}\}\otimes_{\mathbb{Q}_p}E)\}_{I\subset (0,\infty)}}{\mathrm{SimplicialRings}}   \\
\underset{\mathrm{Quasicoherent}^{\mathrm{indBanach,perfectcomplexes}}_{\mathrm{FF}_{\mathrm{Spd}(\mathbb{Q}_p)^\diamond},\mathcal{O}_{\mathrm{FF}_{\mathrm{Spd}(\mathbb{Q}_p)^\diamond}}}}{\mathrm{SimplicialRings}} \ar[u]\ar[r] &\underset{\{\varphi\mathrm{Module}^{\mathrm{indBanach,perfectcomplexes}}(\Gamma^\text{perfect}_{\text{Robba},{\mathrm{Spd}(\mathbb{Q}_p)^\diamond},I}\{t^{1/2}\}\otimes_{\mathbb{Q}_p}E)\}_{I\subset (0,\infty)}}{\mathrm{SimplicialRings}}.\ar[u]  
}
\]

\end{proposition}

\newpage
\section{Moduli $v$-Stack in More General Setting}

\noindent References: \cite{FS}, \cite{FF}, \cite{Sch1},\cite{Sch2}, \cite{KL1}, \cite{KL2};

\noindent Further References: \cite{Lan1}, \cite{Drin1}, \cite{Drin2}, \cite{Zhu}, \cite{DHKM}.\\

\noindent We consider the category of all the perfectoid spaces over $\mathrm{Spd}\overline{\mathbb{F}}_p$ as in \cite{FS}. We use the notation $\mathrm{perfectoid}_v$ to denote the associated $v$-site after \cite{FS}, \cite{Sch2}. Let $p>2$. Now we fix a finite extension $K$ of $\mathbb{Q}_p$. And the Robba rings are defined over $K$ as well, namely we consider the generalized Witt vector over $\mathcal{O}_K$ as in \cite{KL2}\footnote{In \cite{FS} and \cite{KL2}, this field is denoted by $E$ where the relative $p$-adic Hodge theory in \cite{KL2} and the Langlands correspondence in \cite{FS} both happen over this field $E$. To be more precise relative p-adic Hodge theory studies $p$-adic cohomologization over analytic stacks over $E$, while the Langlands correspondences relates derived $\infty$-categories of $\mathrm{Moduli}_{G/E}$ and derived $\infty$-categories of moduli stack of representations of $W_{E,2}$ into the Langlands dual groups.}. For any $\mathrm{Spa}(A,A^+)\in \text{perfectoid}_{v}$ we have the perfect Robba rings from \cite{KL1}, \cite{KL2}:
\begin{align}
\Gamma^\text{perfect}_{\text{Robba},\mathrm{Spa}(A,A^+),I\subset (0,\infty)}.
\end{align}
We also have the corresponding de Rham period rings:
\begin{align}
\Gamma^+_{\text{deRham},\mathrm{Spa}(A,A^+)},\Gamma_{\text{deRham},\mathrm{Spa}(A,A^+)}.
\end{align}
In the first filtration of this first de Rham period ring we have the generator $t$, we now extend the corresponding rings above by adding the square root of $t$, $t^{1/2}
$ following \cite{BS}. We then have the extended rings:
\begin{align}
\Gamma^\text{perfect}_{\text{Robba},\mathrm{Spa}(A,A^+),I\subset (0,\infty)}\{t^{1/2}\},
\end{align}
\begin{align}
\Gamma^+_{\text{deRham},\mathrm{Spa}(A,A^+)}\{t^{1/2}\},\Gamma_{\text{deRham},\mathrm{Spa}(A,A^+)}\{t^{1/2}\}.
\end{align}
Then we form the corresponding extended Fargues-Fontaine curve (after choosing a large finite extension $E$ of $\mathbb{Q}_p$ containing $\varphi(t)^{1/2}$):
\begin{align}
\mathrm{FF}_A:=\bigcup_{I\subset (0,\infty)}\mathrm{Spa}(\Gamma^\text{perfect}_{\text{Robba},\mathrm{Spa}(A,A^+),I\subset (0,\infty)}\{t^{1/2}\}\otimes_{\mathbb{Q}_p}E,\Gamma^{\text{perfect},+}_{\text{Robba},\mathrm{Spa}(A,A^+),I\subset (0,\infty)}\{t^{1/2}\}\otimes_{\mathbb{Q}_p}E)/\varphi^\mathbb{Z},
\end{align}
where the Frobenius is extended to $t^{1/2}\otimes 1$ by acting $\varphi(t)^{1/2}\otimes 1$.

\begin{definition}
Let $G/K$ be any $p$-adic group as in \cite{FS}. That is to say the group $G$ is defined over $K$, where $K$ is some finite extension of $\mathbb{Q}_p$, defined as above. We now define the pre-v-stack $\text{Moduli}_G$ to be a presheaf valued in the groupoid over
\begin{align}
\text{perfectoid}_{v} 
\end{align}
sendind each $\mathrm{Spa}(A,A^+)$ perfectoid in the site to the groupoid of all the locally finite free coherent sheaves carrying $G$-bundle structure over  
\begin{align}
\mathrm{FF}_A:=\bigcup_{I\subset (0,\infty)}\mathrm{Spa}(\Gamma^\text{perfect}_{\text{Robba},\mathrm{Spa}(A,A^+),I\subset (0,\infty)}\{t^{1/2}\}\otimes_{\mathbb{Q}_p}E,\Gamma^{\text{perfect},+}_{\text{Robba},\mathrm{Spa}(A,A^+),I\subset (0,\infty)}\{t^{1/2}\}\otimes_{\mathbb{Q}_p}E)/\varphi^\mathbb{Z}.
\end{align}
\end{definition}

\begin{proposition}
This prestack is a small $v$-stack in the $v$-topology.
\end{proposition}

\begin{proof}
The proof will be the same as in \cite[Proposition III.1.3]{FS}. Our stack can also be regarded as a two components extension of the original stack in \cite{FS}.  
\end{proof}

\section{Motives over $\mathrm{Moduli}_G$ in More General Setting}

\noindent Keeping the generality in the previous section, we now consider the sheaves over extended Fargues-Fontain stacks:

\begin{definition}
\begin{align}
\mathrm{FF}_{\mathrm{Moduli}_G}:=\bigcup_{I\subset (0,\infty)}\mathrm{Spa}(\Gamma^\text{perfect}_{\text{Robba},{\mathrm{Moduli}_G},I\subset (0,\infty)}\{t^{1/2}\}\otimes_{\mathbb{Q}_p}E,\Gamma^{\text{perfect},+}_{\text{Robba},{\mathrm{Moduli}_G},I\subset (0,\infty)}\{t^{1/2}\}\otimes_{\mathbb{Q}_p}E)/\varphi^\mathbb{Z},
\end{align}
which has the corresonding structure map as in the following:
\[\displayindent=-0.4in
\xymatrix@R+1pc{
&\mathrm{FF}_{\mathrm{Moduli}_G} \ar[d]  \\
&\mathrm{FF}_{\mathrm{FF}_*} \ar[d]\\
&\mathrm{FF}_{\mathrm{Spd}^\diamond(\mathbb{Q}_p)} \ar[d]\\
& \mathrm{Spd}^\diamond(\mathbb{Q}_p).  
}
\]
\end{definition}

\begin{definition}
We use the notation
\begin{align}
\mathrm{Quasicoherent}^{\mathrm{solid}}_{\mathrm{FF}_{\mathrm{Moduli}_G},\mathcal{O}_{\mathrm{FF}_{\mathrm{Moduli}_G}}}
\end{align}
to denote $(\infty,1)$-category of all the solid quasicoherent sheaves over the stack $\mathrm{FF}_{\mathrm{Moduli}_G}$. For any local perfectoid $Y\in {\mathrm{Moduli}_G}_v$ we define the corresponding $(\infty,1)$-category in the local sense.\\
We use the notation
\begin{align}
\mathrm{Quasicoherent}^{\mathrm{solid,perfectcomplexes}}_{\mathrm{FF}_{\mathrm{Moduli}_G},\mathcal{O}_{\mathrm{FF}_{\mathrm{Moduli}_G}}}
\end{align}
to denote $(\infty,1)$-category of all the solid quasicoherent sheaves over the stack $\mathrm{FF}_{\mathrm{Moduli}_G}$ which are perfect complexes. For any local perfectoid $Y\in {\mathrm{Moduli}_G}_v$ we define the corresponding $(\infty,1)$-category in the local sense.
\end{definition}

\begin{definition}
We use the notation
\begin{align}
\mathrm{Quasicoherent}^{\mathrm{indBanach}}_{\mathrm{FF}_{\mathrm{Moduli}_G},\mathcal{O}_{\mathrm{FF}_{\mathrm{Moduli}_G}}}
\end{align}
to denote $(\infty,1)$-category of all the ind-Banach quasicoherent sheaves over the stack $\mathrm{FF}_{\mathrm{Moduli}_G}$. For any local perfectoid $Y\in {\mathrm{Moduli}_G}_v$ we define the corresponding $(\infty,1)$-category in the local sense.\\
We use the notation
\begin{align}
\mathrm{Quasicoherent}^{\mathrm{indBanach,perfectcomplexes}}_{\mathrm{FF}_{\mathrm{Moduli}_G},\mathcal{O}_{\mathrm{FF}_{\mathrm{Moduli}_G}}}
\end{align}
to denote $(\infty,1)$-category of all the indBanach quasicoherent sheaves over the stack $\mathrm{FF}_{\mathrm{Moduli}_G}$ which are perfect complexes. For any local perfectoid $Y\in {\mathrm{Moduli}_G}_v$ we define the corresponding $(\infty,1)$-category in the local sense.
\end{definition}

\begin{definition}
We use the notation
\begin{align}
\{\varphi\mathrm{Module}^{\mathrm{solid}}(\Gamma^\text{perfect}_{\text{Robba},{\mathrm{Moduli}_G},I}\{t^{1/2}\}\otimes_{\mathbb{Q}_p}E)\}_{I\subset (0,\infty)}
\end{align}
to denote $(\infty,1)$-category of all the solid $\varphi$-modules over the extended Robba ring. The modules satisfy the Frobenius pullback condition and glueing condition for overlapped intervals $I\subset J\subset K$. For any local perfectoid $Y\in {\mathrm{Moduli}_G}_v$ we define the corresponding $(\infty,1)$-category in the local sense.\\
We use the notation
\begin{align}
\{\varphi\mathrm{Module}^{\mathrm{solid,perfectcomplexes}}(\Gamma^\text{perfect}_{\text{Robba},{\mathrm{Moduli}_G},I}\{t^{1/2}\}\otimes_{\mathbb{Q}_p}E)\}_{I\subset (0,\infty)}
\end{align}
to denote $(\infty,1)$-category of all the solid $\varphi$-modules over the extended Robba ring which are perfect complexes. The modules satisfy the Frobenius pullback condition and glueing condition for overlapped intervals $I\subset J\subset K$. For any local perfectoid $Y\in {\mathrm{Moduli}_G}_v$ we define the corresponding $(\infty,1)$-category in the local sense. 

\end{definition}

\begin{definition}
We use the notation
\begin{align}
\{\varphi\mathrm{Module}^{\mathrm{indBanach}}(\Gamma^\text{perfect}_{\text{Robba},{\mathrm{Moduli}_G},I}\{t^{1/2}\}\otimes_{\mathbb{Q}_p}E)\}_{I\subset (0,\infty)}
\end{align}
to denote $(\infty,1)$-category of all the ind-Banach $\varphi$-modules over the extended Robba ring. The modules satisfy the Frobenius pullback condition and glueing condition for overlapped intervals $I\subset J\subset K$. For any local perfectoid $Y\in {\mathrm{Moduli}_G}_v$ we define the corresponding $(\infty,1)$-category in the local sense.\\
We use the notation
\begin{align}
\{\varphi\mathrm{Module}^{\mathrm{ind-Banach,perfectcomplexes}}(\Gamma^\text{perfect}_{\text{Robba},{\mathrm{Moduli}_G},I}\{t^{1/2}\}\otimes_{\mathbb{Q}_p}E)\}_{I\subset (0,\infty)}
\end{align}
to denote $(\infty,1)$-category of all the ind-Banach $\varphi$-modules over the extended Robba ring which are perfect complexes. The modules satisfy the Frobenius pullback condition and glueing condition for overlapped intervals $I\subset J\subset K$. For any local perfectoid $Y\in {\mathrm{Moduli}_G}_v$ we define the corresponding $(\infty,1)$-category in the local sense. 

\end{definition}

\begin{definition}
We use the notation
\begin{align}
\varphi\mathrm{Module}^{\mathrm{solid}}(\Gamma^\text{perfect}_{\text{Robba},{\mathrm{Moduli}_G},\infty}\{t^{1/2}\}\otimes_{\mathbb{Q}_p}E)
\end{align}
to denote $(\infty,1)$-category of all the solid $\varphi$-modules over the extended Robba ring. The modules satisfy the Frobenius pullback condition. For any local perfectoid $Y\in {\mathrm{Moduli}_G}_v$ we define the corresponding $(\infty,1)$-category in the local sense.\\
We use the notation
\begin{align}
\varphi\mathrm{Module}^{\mathrm{solid,perfectcomplexes}}(\Gamma^\text{perfect}_{\text{Robba},{\mathrm{Moduli}_G},\infty}\{t^{1/2}\}\otimes_{\mathbb{Q}_p}E)
\end{align}
to denote $(\infty,1)$-category of all the solid $\varphi$-modules over the extended Robba ring which are perfect complexes. The modules satisfy the Frobenius pullback condition. For any local perfectoid $Y\in {\mathrm{Moduli}_G}_v$ we define the corresponding $(\infty,1)$-category in the local sense. 

\end{definition}

\begin{definition}
We use the notation
\begin{align}
\varphi\mathrm{Module}^{\mathrm{indBanach}}(\Gamma^\text{perfect}_{\text{Robba},{\mathrm{Moduli}_G},\infty}\{t^{1/2}\}\otimes_{\mathbb{Q}_p}E)
\end{align}
to denote $(\infty,1)$-category of all the ind-Banach $\varphi$-modules over the extended Robba ring. The modules satisfy the Frobenius pullback condition. For any local perfectoid $Y\in {\mathrm{Moduli}_G}_v$ we define the corresponding $(\infty,1)$-category in the local sense.\\
We use the notation
\begin{align}
\left\{\underset{\mathrm{indBanach,perfectcomplexes},\Gamma^\text{perfect}_{\text{Robba},X,I}\{t^{1/2}\}\otimes_{\mathbb{Q}_p}E}{\varphi\mathrm{Module}}\right\}_{I\subset (0,\infty)}
\end{align}
to denote $(\infty,1)$-category of all the ind-Banach $\varphi$-modules over the extended Robba ring which are perfect complexes. The modules satisfy the Frobenius pullback condition. For any local perfectoid $Y\in {\mathrm{Moduli}_G}_v$ we define the corresponding $(\infty,1)$-category in the local sense. 

\end{definition}

\begin{proposition}
We have the following commutative diagram by taking the global section functor in the horizontal rows:\\
\[\displayindent=+0in
\xymatrix@R+7pc{
\mathrm{Quasicoherent}^{\mathrm{solid}}_{\mathrm{FF}_{\mathrm{Spd}(\mathbb{Q}_p)^\diamond},\mathcal{O}_{\mathrm{FF}_{\mathrm{Spd}(\mathbb{Q}_p)^\diamond}}} \ar[r] \ar[d] &\{\varphi\mathrm{Module}^{\mathrm{solid}}(\Gamma^\text{perfect}_{\text{Robba},{\mathrm{Spd}(\mathbb{Q}_p)^\diamond},I}\{t^{1/2}\}\otimes_{\mathbb{Q}_p}E)\}_{I\subset (0,\infty)} \ar[d]\\
\mathrm{Quasicoherent}^{\mathrm{solid}}_{\mathrm{FF}_{\mathrm{Moduli}_G},\mathcal{O}_{\mathrm{FF}_{\mathrm{Moduli}_G}}}  \ar[r] \ar[r] \ar[r] &\{\varphi\mathrm{Module}^{\mathrm{solid}}(\Gamma^\text{perfect}_{\text{Robba},{\mathrm{Moduli}_G},I}\{t^{1/2}\}\otimes_{\mathbb{Q}_p}E)\}_{I\subset (0,\infty)}.  \\  
}
\]
\end{proposition}

\begin{proposition}
We have the following commutative diagram by taking the global section functor in the horizontal rows:\\
\[\displayindent=+0in
\xymatrix@R+7pc{
\mathrm{Quasicoherent}^{\mathrm{solid,perfectcomplexes}}_{\mathrm{FF}_{\mathrm{Spd}(\mathbb{Q}_p)^\diamond},\mathcal{O}_{\mathrm{FF}_{\mathrm{Spd}(\mathbb{Q}_p)^\diamond}}} \ar[r] \ar[d] &\{\underset{\mathrm{solid,perfectcomplexes}}{\varphi\mathrm{Module}}(\Gamma^\text{perfect}_{\text{Robba},{\mathrm{Spd}(\mathbb{Q}_p)^\diamond},I}\{t^{1/2}\}\otimes_{\mathbb{Q}_p}E)\}_{I\subset (0,\infty)} \ar[d]  \\
\mathrm{Quasicoherent}^{\mathrm{solid,perfectcomplexes}}_{\mathrm{FF}_{\mathrm{Moduli}_G},\mathcal{O}_{\mathrm{FF}_{\mathrm{Moduli}_G}}}  \ar[r] \ar[r] \ar[r] &\{\underset{\mathrm{solid,perfectcomplexes}}{\varphi\mathrm{Module}}(\Gamma^\text{perfect}_{\text{Robba},{\mathrm{Moduli}_G},I}\{t^{1/2}\}\otimes_{\mathbb{Q}_p}E)\}_{I\subset (0,\infty)}.   
}
\]
\end{proposition}

\begin{proposition}
We have the following commutative diagram by taking the global section functor in the horizontal rows:\\
\[\displayindent=+0in
\xymatrix@C-0.1in@R+7pc{
\mathrm{Quasicoherent}^{\mathrm{indBanach}}_{\mathrm{FF}_{\mathrm{Spd}(\mathbb{Q}_p)^\diamond},\mathcal{O}_{\mathrm{FF}_{\mathrm{Spd}(\mathbb{Q}_p)^\diamond}}} \ar[r] \ar[d] &\{\varphi\mathrm{Module}^\text{indBanach}(\Gamma^\text{perfect}_{\text{Robba},{\mathrm{Spd}(\mathbb{Q}_p)^\diamond},I}\{t^{1/2}\}\otimes_{\mathbb{Q}_p}E)\}_{I\subset (0,\infty)}\ar[d]\\
\mathrm{Quasicoherent}^{\mathrm{indBanach}}_{\mathrm{FF}_{\mathrm{Moduli}_G},\mathcal{O}_{\mathrm{FF}_{\mathrm{Moduli}_G}}}  \ar[r] \ar[r] \ar[r] &\{\varphi\mathrm{Module}^{\mathrm{indBanach}}(\Gamma^\text{perfect}_{\text{Robba},{\mathrm{Moduli}_G},I}\{t^{1/2}\}\otimes_{\mathbb{Q}_p}E)\}_{I\subset (0,\infty)}.    
}
\]
\end{proposition}

\begin{proposition}
We have the following commutative diagram by taking the global section functor in the horizontal rows:
\[\displayindent=+0in
\xymatrix@C-0.2in@R+7pc{
\mathrm{Quasicoherent}^{\mathrm{indBanach,perfectcomplexes}}_{\mathrm{FF}_{\mathrm{Spd}(\mathbb{Q}_p)^\diamond},\mathcal{O}_{\mathrm{FF}_{\mathrm{Spd}(\mathbb{Q}_p)^\diamond}}} \ar[d]\ar[r] &
\left\{\underset{\mathrm{indBanach,perfectcomplexes},\Gamma^\text{perfect}_{\text{Robba},\mathrm{Spd}(\mathbb{Q}_p)^\diamond,I}\{t^{1/2}\}\otimes_{\mathbb{Q}_p}E}{\varphi\mathrm{Module}}\right\}_{I\subset (0,\infty)}\ar[d]\\
\mathrm{Quasicoherent}^{\mathrm{indBanach,perfectcomplexes}}_{\mathrm{FF}_{\mathrm{Moduli}_G},\mathcal{O}_{\mathrm{FF}_{\mathrm{Moduli}_G}}}  \ar[r] \ar[r] \ar[r] &
\left\{\underset{\mathrm{indBanach,perfectcomplexes},\Gamma^\text{perfect}_{\text{Robba},\mathrm{Moduli}_G,I}\{t^{1/2}\}\otimes_{\mathbb{Q}_p}E}{\varphi\mathrm{Module}}\right\}_{I\subset (0,\infty)}.    
}
\]

\end{proposition}

\indent Taking the corresponding simplicial commutative object we have the following propositions:

\begin{proposition}
We have the following commutative diagram by taking the global section functor in the horizontal rows:
\[\displayindent=+0in
\xymatrix@R+7pc{
\underset{\mathrm{Quasicoherent}^{\mathrm{solid}}_{\mathrm{FF}_{\mathrm{Spd}(\mathbb{Q}_p)^\diamond},\mathcal{O}_{\mathrm{FF}_{\mathrm{Spd}(\mathbb{Q}_p)^\diamond}}}}{\mathrm{SimplicialRings}} \ar[d] \ar[r] &\underset{\{\varphi\mathrm{Module}^{\mathrm{solid}}(\Gamma^\text{perfect}_{\text{Robba},{\mathrm{Spd}(\mathbb{Q}_p)^\diamond},I}\{t^{1/2}\}\otimes_{\mathbb{Q}_p}E)\}_{I\subset (0,\infty)}}{\mathrm{SimplicialRings}} \ar[d]\\
\underset{\mathrm{Quasicoherent}^{\mathrm{solid}}_{\mathrm{FF}_{\mathrm{Moduli}_G},\mathcal{O}_{\mathrm{FF}_{\mathrm{Moduli}_G}}}}{\mathrm{SimplicialRings}}  \ar[r] \ar[r] \ar[r] &\underset{\{\varphi\mathrm{Module}^{\mathrm{solid}}(\Gamma^\text{perfect}_{\text{Robba},{\mathrm{Moduli}_G},I}\{t^{1/2}\}\otimes_{\mathbb{Q}_p}E)\}_{I\subset (0,\infty)}}{\mathrm{SimplicialRings}}.
}
\]
\end{proposition}

\begin{proposition}
We have the following commutative diagram by taking the global section functor in the horizontal rows:
\[\displayindent=+0in
\xymatrix@C+0in@R+7pc{
\underset{\mathrm{Quasicoherent}^{\mathrm{solid,perfectcomplexes}}_{\mathrm{FF}_{\mathrm{Spd}(\mathbb{Q}_p)^\diamond},\mathcal{O}_{\mathrm{FF}_{\mathrm{Spd}(\mathbb{Q}_p)^\diamond}}}}{\mathrm{SimplicialRings}}\ar[d] \ar[r] &\underset{\{\varphi\mathrm{Module}^{\mathrm{solid,perfectcomplexes}}(\Gamma^\text{perfect}_{\text{Robba},{\mathrm{Spd}(\mathbb{Q}_p)^\diamond},I}\{t^{1/2}\}\otimes_{\mathbb{Q}_p}E)\}_{I\subset (0,\infty)}}{\mathrm{SimplicialRings}}\ar[d]\\
\underset{\mathrm{Quasicoherent}^{\mathrm{solid,perfectcomplexes}}_{\mathrm{FF}_{\mathrm{Moduli}_G},\mathcal{O}_{\mathrm{FF}_{\mathrm{Moduli}_G}}}}{\mathrm{SimplicialRings}}  \ar[r] \ar[r] \ar[r] &\underset{\{\varphi\mathrm{Module}^{\mathrm{solid,perfectcomplexes}}(\Gamma^\text{perfect}_{\text{Robba},{\mathrm{Moduli}_G},I}\{t^{1/2}\}\otimes_{\mathbb{Q}_p}E)\}_{I\subset (0,\infty)}}{\mathrm{SimplicialRings}}.  
}
\]
\end{proposition}

\begin{proposition}
We have the following commutative diagram by taking the global section functor in the horizontal rows:
\[\displayindent=+0in
\xymatrix@R+7pc{
\underset{\mathrm{Quasicoherent}^{\mathrm{indBanach}}_{\mathrm{FF}_{\mathrm{Moduli}_G},\mathcal{O}_{\mathrm{FF}_{\mathrm{Moduli}_G}}}}{\mathrm{SimplicialRings}}  \ar[r] \ar[r] \ar[r] &\underset{\{\varphi\mathrm{Module}^{\mathrm{indBanach}}(\Gamma^\text{perfect}_{\text{Robba},{\mathrm{Moduli}_G},I}\{t^{1/2}\}\otimes_{\mathbb{Q}_p}E)\}_{I\subset (0,\infty)}}{\mathrm{SimplicialRings}}   \\
\underset{\mathrm{Quasicoherent}^{\mathrm{indBanach}}_{\mathrm{FF}_{\mathrm{Spd}(\mathbb{Q}_p)^\diamond},\mathcal{O}_{\mathrm{FF}_{\mathrm{Spd}(\mathbb{Q}_p)^\diamond}}}}{\mathrm{SimplicialRings}} \ar[u]\ar[r] &\underset{\{\varphi\mathrm{Module}^\text{indBanach}(\Gamma^\text{perfect}_{\text{Robba},{\mathrm{Spd}(\mathbb{Q}_p)^\diamond},I}\{t^{1/2}\}\otimes_{\mathbb{Q}_p}E)\}_{I\subset (0,\infty)}}{\mathrm{SimplicialRings}}.\ar[u]  
}
\]
\end{proposition}

\begin{proposition}
We have the following commutative diagram by taking the global section functor in the horizontal rows:
\[\displayindent=+0in
\xymatrix@C-0.1in@R+7pc{
\underset{\mathrm{Quasicoherent}^{\mathrm{indBanach,perfectcomplexes}}_{\mathrm{FF}_{\mathrm{Moduli}_G},\mathcal{O}_{\mathrm{FF}_{\mathrm{Moduli}_G}}}}{\mathrm{SimplicialRings}}  \ar[r] \ar[r] \ar[r] &\underset{\{\varphi\mathrm{Module}^{\mathrm{indBanach,perfectcomplexes}}(\Gamma^\text{perfect}_{\text{Robba},{\mathrm{Moduli}_G},I}\{t^{1/2}\}\otimes_{\mathbb{Q}_p}E)\}_{I\subset (0,\infty)}}{\mathrm{SimplicialRings}}   \\
\underset{\mathrm{Quasicoherent}^{\mathrm{indBanach,perfectcomplexes}}_{\mathrm{FF}_{\mathrm{Spd}(\mathbb{Q}_p)^\diamond},\mathcal{O}_{\mathrm{FF}_{\mathrm{Spd}(\mathbb{Q}_p)^\diamond}}}}{\mathrm{SimplicialRings}} \ar[u]\ar[r] &\underset{\{\varphi\mathrm{Module}^{\mathrm{indBanach,perfectcomplexes}}(\Gamma^\text{perfect}_{\text{Robba},{\mathrm{Spd}(\mathbb{Q}_p)^\diamond},I}\{t^{1/2}\}\otimes_{\mathbb{Q}_p}E)\}_{I\subset (0,\infty)}}{\mathrm{SimplicialRings}}.\ar[u]  
}
\]

\end{proposition}

\newpage
\section{Generalized Local Langlands Correspondence after Fargues-Scholze}

\begin{reference}\mbox{}\\
$\square$ References on p-adic Hodge Theory:\\ \cite{pHodgeT}, \cite{pHodgeF}, \cite{pHodgeS1}, \cite{pHodgeS2}, \cite{pHodgeKL1}, \cite{pHodgeKL2}, \cite{pHodgeBS}, \cite{pHodgeKPX};\\
$\square$ References on foundations of p-adic analysis:\\
\cite{pToAnCS1}, \cite{pToAnCS2}, \cite{pToAnCS3}, \cite{pToAnCS4},  
\cite{pToAnBBBK};\\
$\square$ References on local Langlands program:\\ 
\cite{LPL}, \cite{LPD1}, \cite{LPLL}, \cite{LPVL}, \cite{LPC}, \cite{LPFS}, \cite{LPGL}, \cite{LPEGH}, \cite{LPEG}, \cite{LPZ}, \cite{LPDHKM}, \cite{LPD2};\\
$\square$ References on p-adic Langlands program:\\\cite{LPL}, \cite{LPD1}, \cite{LPLL}, \cite{LPVL}, \cite{LPC}, \cite{LPFS}, \cite{LPGL}, \cite{LPEGH}, \cite{LPEG}, \cite{LPZ}, \cite{LPDHKM}, \cite{LPD2}.
\end{reference}

\subsection{Perfectoid Moduli Stacks}

\noindent As above, we discussed certain mixed-parity moduli stacks over the mixed-parity Fargues-Fontaine curves where we work directly over the \textit{tilt} of some field $K:=\overline{\mathbb{Q}_p((\mu_{p^\infty}))}^\wedge$. We assume first that we still work over this field where we have the corresponding element $t$ as $\log[\varepsilon]$ as in \cite{pHodgeS1} for the de Rham sheaves over the pro-\'etale site. Then in the scenario where we have the splitting $t^{1/2}$ we then extend the underlying diamond to be a two-fold covering. To be more precise we defined the corresponding prestack 
\begin{align}
\mathrm{Bundle}_{2,G}(\mathrm{FF}_{2,\square})
\end{align}
which sends each perfectoid space $X_\square$ in the v-site $\mathrm{Perfectoid}_{K,v}$ to the corresponding groupoid of $G$-torsors over the relative extended Fargues-Fontaine curves $\mathrm{FF}_{2,X_\square}$. Here we assume $G$ is the $p$-adic group in \cite{LPFS}, namely we have a base field $K'$ such that $G$ is defined over this field and we assume other conditions from \cite{LPFS}. This is indeed a $v$-stack as in \cite[Chapter III, Proposition 1.3]{LPFS}.

\begin{proposition}
\begin{align}
\mathrm{Bundle}_{2,G}(\mathrm{FF}_{2,\square})
\end{align}
is a small $v$-stack over the site $\mathrm{Perfectoid}_{K,v}$. $\square$
\end{proposition}

\indent In the corresponding discussion as above we actually condidered the $p$-adic motives over the corresponding $v$-stack here which is also \textit{small}. Namely locally for each corresponding local perfectoid $Y$ we have the corresponding chance to attach the corresonding Robba rings\footnote{Note that we need to tensor the Robba rings with a large finite extension of $\mathbb{Q}_p$ to guarantee that $\varphi.(t^{1/2})$ makes sense in the rings.}
\begin{align}
\widetilde{\Pi}_{Y,I,2},\\
\widetilde{\Pi}_{Y,\infty,2},\\
\widetilde{\Pi}_{Y,2},
\end{align}
in the correspoding extended setting where we join the element $t^{1/2}$. We can have the corresponding Frobenius action or not actualy extended from the corresponding Robba rings. Be careful \textit{from now on} we consider the following assumption from \cite{pHodgeKL2} and \cite{LPFS}:

\begin{assumption}
We assume that the corresponding $G$ and the corresponding Robba rings are defined over a same field $K'$ which is finite extension of $\mathbb{Q}_p$, where we often implicize the corresponding field here. We assume the field is $p$-adic which is denoted by $E$ in \cite{LPFS}.
\end{assumption}

\noindent This will generalize the above discussion to the situation where we work over the field $\overline{K'((\mu_{p^\infty}))}^{\wedge,\flat}$ which contains the field $k_{K'}$, the residue field of $K'$. Then we have the corresponding definitions of the stacks as in the above after \cite{pHodgeKL2} and \cite{LPFS}.

\newpage
\subsection{Generalization of Langlands Program and the Geometrization}

\begin{assumption}
We assume that the corresponding $G$ and the corresponding Robba rings are defined over a same field $K'$ which is finite extension of $\mathbb{Q}_p$, where we often implicize the corresponding field here. We assume the field is $p$-adic which is denoted by $E$ in \cite{LPFS}. The $\mathrm{Bundle}_{2,G}$ is defined over $\mathrm{perfectoid}_{\mathrm{Spd}\overline{\mathbb{F}}_p}$. 
\end{assumption}

\indent  The main conjecture in \cite{LPFS} relates derived $\infty$-category of all the $A$-valued complexes over $\mathrm{Bundle}_G(\mathrm{FF}_{\square})$ to the corresponding derived $\infty$-category of coherent sheaves over the stack of $L$-parameter $\mathrm{Stack}_{L,{G}^\mathrm{dual},A}/{G}^\mathrm{dual}$. The former derived category in our setting can be defined directly which is also well-defined from \cite{LPFS} and \cite{pHodgeS2}:
\begin{align}
\mathrm{DerivedCat}(\mathrm{Bundle}_{2,G}(\mathrm{FF}_{2,\square}))_{\mathrm{Coeff}: A}
\end{align}
which is actually not expected to completely isomorphic to the corresponding category of the latter in our situation:
\begin{align}
\mathrm{DerivedCat}_{\mathrm{bounded},\mathrm{coherent},\mathrm{Nilpotent}}(\mathrm{Stack}_{2,L,{G}^\mathrm{dual},A}/{G}^\mathrm{dual})
\end{align}
where the stack of $L$-parameters (which origin two-fold covering of the Weil group $\mathrm{Weil}_{K',2}$) is the pull back of the stack $\mathrm{Stack}_{L,{G}^\mathrm{dual},A}/{G}^\mathrm{dual}$ along the covering map of $\mathrm{Weil}_{K',2}$ over $\mathrm{Weil}_{K'}$. Fargues-Scholze's conjecture conjectures directly that they are isomophic in the corresponding usual non-mixed-parity situation. In our situation the relationship should be clear as well but there is some tiny different due to the fact that we need to tensor the Robba rings with a large finite extension of $\mathbb{Q}_p$ to guarantee that $\varphi.(t^{1/2})$ makes sense in the rings.

\begin{conjecture} (Fargues-Scholze, \cite[Chapter I, Conjecture 10.2]{LPFS})
There is a canonical isomorphism between the two $\infty$-categories:
\begin{align}
\mathrm{DerivedCat}(\mathrm{Bundle}_{G}(\mathrm{FF}_{\square}))_{\mathrm{Coeff}: A}
\end{align}
with
\begin{align}
\mathrm{DerivedCat}_{\mathrm{bounded},\mathrm{coherent},\mathrm{Nilpotent}}(\mathrm{Stack}_{L,{G}^\mathrm{dual},A}/{G}^\mathrm{dual}).
\end{align}
\end{conjecture}

\begin{conjecture} (After Fargues-Scholze, \cite[Chapter I, Conjecture 10.2]{LPFS})
There is a canonical direct relationship between the two $\infty$-categories:
\begin{align}
\mathrm{DerivedCat}(\mathrm{Bundle}_{G}(\mathrm{FF}_{\square}))_{\mathrm{Coeff}: A}
\end{align}
with 
\begin{align}
\mathrm{DerivedCat}(\mathrm{Bundle}_{2,G}(\mathrm{FF}_{2,\square}))_{\mathrm{Coeff}: A}
\end{align}
after we consider the pull-back of the categories along the map:
\begin{align}
\mathrm{DerivedCat}_{\mathrm{bounded},\mathrm{coherent},\mathrm{Nilpotent}}&(\mathrm{Stack}_{2,L,{G}^\mathrm{dual},A}/{G}^\mathrm{dual})\\&\rightarrow \mathrm{DerivedCat}_{\mathrm{bounded},\mathrm{coherent},\mathrm{Nilpotent}}(\mathrm{Stack}_{L,{G}^\mathrm{dual},A}/{G}^\mathrm{dual}).
\end{align}
\end{conjecture}

\begin{remark}
The automorphic to Galois/Weil direction of the Langlands correspondence in \cite{LPFS} can be realized in our setting as well, but we need to consider a corresponding larger category to similarize the corresponding $\mathrm{Weil}_{K',2}$-$G$-equivariant sheaves coming from the smooth representations of $G$. For instance over $\mathrm{Bundle}_{2,G}(\mathrm{FF}_{2,\square})$ we consider such category of $\mathrm{Weil}_{K',2}$-$G$-equivariant sheaves not coming from the smooth representations of $G$, which we denote by $\mathrm{Perv}_{K',G,\mathrm{smooth},2}$, exactly as in the corresponding usual situation in \cite{LPFS}. Then we take the corresponding $\mathrm{Weil}_{K',2}$-equivariant sheaves and cohomologies obtained from this. Then all the corresponding stacks of shtukas and the cohomologies can be defined by using the pull back along:
\begin{align}
\mathrm{FF}_{2,\square} \rightarrow \mathrm{FF}_{\square}
\end{align}
after which one can extract the corresponding representations of product of $\mathrm{Weil}_{K',2}$. Then the corresponding Drinfeld's lemma in our mixed parity setting holds true as well since we just take the original Drinfeld's lemma for $K'$ over the corresponding morphism:
\begin{align}
\mathrm{Gal}_{\mathbb{Q}_p/K,2} \rightarrow \mathrm{Gal}_{\mathbb{Q}_p/K}.
\end{align}
This will gives us the corresponding mixed-parity $L$-parameter from $\mathrm{Weil}_{K',2}$ into $\widehat{G}(A)$. This is why we call our current framework \textit{generalized Langlands correspondence}. At least under the well-defined geometrization, both sides of the generalized correspondence should behave quite similar to those in \cite{LPFS}.

\end{remark}

\newpage
\subsubsection*{Acknowledgements}
The author thanks Professor Kedlaya for all those suggestions around the corresponding mixed-parity theoreticalization of the work of Kedlaya, Kedlaya-Liu, Kedlaya-Pottharst-Xiao. The author thanks Professor Sorensen for suggestions on the mixed-parity representation theoretic perspectives.

\newpage
\bibliographystyle{}

\begin{thebibliography}{}

\bibitem[Sch1]{Sch1} Scholze, Peter. "p-adic Hodge Theory for Rigid-Analytic Varieties." Forum of Mathematics. Pi, vol. 1, 2013, https://doi.org/10.1017/fmp.2013.1.

\bibitem[KL1]{KL1} Kedlaya, Kiran Sridhara, and Ruochuan Liu. "Relative p-Adic Hodge Theory: Foundations." Soci\'et\'e math\'ematique de France, 2015.

\bibitem[KL2]{KL2} Kedlaya, Kiran S., and Ruochuan Liu. "Relative p-Adic Hodge Theory, II: Imperfect Period Rings." 2016, https://doi.org/10.48550/arxiv.1602.06899.

\bibitem[BL1]{BL1} Bhatt, Bhargav and Jacob Lurie. "A p-adic Riemann-Hilbert functor I: torsion coefficients."

\bibitem[BL2]{BL2} Bhatt, Bhargav and Jacob Lurie. "A p-adic Riemann-Hilbert functor II: $\mathbb{Q}_p$-coefficients."

\bibitem[BS]{BS} Breuil, Christophe, and Peter Schneider. "First Steps Towards p-Adic Langlands Functoriality." Journal F\"ur Die Reine Und Angewandte Mathematik, vol. 2007, no. 610, 2007, pp. 149-80, https://doi.org/10.1515/CRELLE.2007.070.

\bibitem[Fon1]{Fon1} Fontaine, Jean-Marc. "Arithm\'etique des repr\'esentations galoisiennes $p$-adiques." In Cohomologie $p$-adiques et applications arithm\'etiques (III), Ast\'erisque, no. 295 (2004), pp. 1-115. http://www.numdam.org/item/AST$\_$2004$\_$295$\_$1$\_$0/.

\bibitem[BHS]{BHS} Breuil, Christophe, et al. "A Local Model for the Trianguline Variety and Applications." Publications Math\'ematiques. Institut Des Hautes \'Etudes Scientifiques, vol. 130, no. 1, 2019, pp. 299-412, https://doi.org/10.1007/s10240-019-00111-y.

\bibitem[M]{M}  Mann, Lucas. "A $p$-Adic 6-Functor Formalism in Rigid-Analytic Geometry." 2022, https://doi.org/10.48550/arxiv.2206.02022.

\bibitem[CS1]{CS1} Clausen, Dustin and Peter Scholze. "Lectures on Condensed Mathematics." https://www.math.uni-bonn.de/people/scholze/Condensed.pdf. 

\bibitem[CS2]{CS2} Clausen, Dustin and Peter Scholze. "Lectures on Analytic Geometry." https://www.math.uni-bonn.de/people/scholze/Analytic.pdf.
 
\bibitem[BK]{BK} Bambozzi, Federico, and Kobi Kremnizer. "On the Sheafyness Property of Spectra of Banach Rings." 2020, https://doi.org/10.48550/arxiv.2009.13926.

\bibitem[BBK]{BBK} Ben-Bassat, Oren, and Kobi Kremnizer. "Fr\'echet Modules and Descent." Theory and Applications of Categories, vol. 39, no. 9, 2023.

\bibitem[BBBK]{BBBK} Bambozzi, Federico, et al. "Analytic Geometry over $F_1$ and the Fargues-Fontaine Curve." Advances in Mathematics (New York. 1965), vol. 356, 2019, https://doi.org/10.1016/j.aim.2019.106815.

\bibitem[BBM]{KM} Ben-Bassat, Oren, and Devarshi Mukherjee. "Analytification, Localization and Homotopy Epimorphisms." Bulletin Des Sciences Math\'ematiques, vol. 176, 2022, https://doi.org/10.1016/j.bulsci.2022.103129.

\bibitem[KKM]{KKM} Kelly, Jack, Kobi Kremnizer, and Devarshi Mukherjee. 2021. "Analytic Hochschild-Kostant-Rosenberg Theorem." Advances in Mathematics (New York. 1965), vol. 410, 2022, https://doi.org/10.1016/j.aim.2022.108694.

\bibitem[LZ]{LZ} Liu, Ruochuan, and Xinwen Zhu. "Rigidity and a Riemann-Hilbert Correspondence for p-Adic Local Systems." Inventiones Mathematicae, vol. 207, no. 1, 2017, pp. 291-343, https://doi.org/10.1007/s00222-016-0671-7.

\bibitem[FS]{FS} Fargues, Laurent, and Peter Scholze. "Geometrization of the Local Langlands Correspondence." 2021, https://doi.org/10.48550/arxiv.2102.13459.
 
\bibitem[FF]{FF} Fargues, Laurent, Jean Marc Fontaine. "Courbes et fibr\'es vectoriels en th\'eorie de Hodge p-adique." Ast\'erisque 406 (2018): 1-382.

\bibitem[Sch2]{Sch2} Scholze, Peter. "Etale Cohomology of Diamonds." 2017, https://doi.org/10.48550/arxiv.1709.07343.

\bibitem[Sch1]{Sch1} Scholze, Peter. "$p$-adic Hodge Theory for Rigid-Analytic Varieties." Forum of Mathematics. Pi, vol. 1, 2013, https://doi.org/10.1017/fmp.2013.1.

\bibitem[KL1]{KL1} Kedlaya, Kiran Sridhara, and Ruochuan Liu. "Relative p-Adic Hodge Theory: Foundations." Soci\'et\'e math\'ematique de France, 2015.

\bibitem[KL3]{KL3} Kedlaya, Kiran, and Ruochuan Liu. "On Families of $(\varphi, \Gamma)$-Modules." Algebra and Number Theory, vol. 4, no. 7, 2010, pp. 943-967, https://doi.org/10.2140/ant.2010.4.943.

\bibitem[Ked1]{Ked1} Kedlaya, Kiran S. "A p-Adic Local Monodromy Theorem." Annals of Mathematics, vol. 160, no. 1, 2004, pp. 93-184, https://doi.org/10.4007/annals.2004.160.93.

\bibitem[KPX]{KPX} Kedlaya, Kiran S, Jonathan Pottharst, and Liang Xiao. 2012. "Cohomology of Arithmetic Families of $(\varphi,\Gamma)$-Modules." Journal of the American Mathematical Society, vol. 27, no. 4, 2014, pp. 1043-1115, https://doi.org/10.1090/S0894-0347-2014-00794-3.
 
\bibitem[KL2]{KL2} Kedlaya, Kiran S., and Ruochuan Liu. "Relative p-Adic Hodge Theory, II: Imperfect Period Rings." 2016, https://doi.org/10.48550/arxiv.1602.06899.

\bibitem[BS]{BS} Breuil, Christophe, and Peter Schneider. "First Steps Towards p-Adic Langlands Functoriality." Journal F\"ur Die Reine Und Angewandte Mathematik, vol. 2007, no. 610, 2007, pp. 149-80, https://doi.org/10.1515/CRELLE.2007.070.

\bibitem[CS1]{CS1} Clausen, Dustin and Peter Scholze. "Lectures on Condensed Mathematics." https://www.math.uni-bonn.de/people/scholze/Condensed.pdf.

\bibitem[CS2]{CS2} Clausen, Dustin and Peter Scholze. "Lectures on Analytic Geometry." https://www.math.uni-bonn.de/people/scholze/Analytic.pdf.

\bibitem[TT]{TT} Tan, Fucheng, and Jilong Tong. "Crystalline Comparison Isomorphisms in p-Adic Hodge Theory: the Absolutely Unramified Case." Algebra and Number Theory, vol. 13, no. 7, 2019, pp. 1509-81, https://doi.org/10.2140/ant.2019.13.1509.

\bibitem[LBV]{LBV} Bras, Arthur-C\'esar Le, and Alberto Vezzani. "The de Rham-Fargues-Fontaine Cohomology." 2021, https://doi.org/10.48550/arxiv.2105.13028.

\bibitem[B]{B} Bosco, Guido. "Rational p-adic Hodge Theory for Rigid Analytic Varieties." arXiv:2306.06100. 

\bibitem[Shi]{Shi} Shimizu, K. (2022). "A $p$-adic monodromy theorem for de Rham local systems." Compositio Mathematica, 158(12), 2157-2205. doi:10.1112/S0010437X2200776X.

\bibitem[AI1]{AI1} Andreatta, Fabrizio, and Adrian Iovita. "Global Applications of Relative $(\varphi,\Gamma)$-Modules I." Ast\'erisque, Volume 319, 339-420, 2008.

\bibitem[AB1]{AB1} Andreatta, Fabrizio; Brinon, Olivier. "$\mathrm {B}_{\mathrm {dR}}$-repr\'esentations dans le cas relatif." Annales scientifiques de l'\'Ecole Normale Sup\'erieure, Serie 4, Volume 43 (2010) no. 2, pp. 279-339. doi : 10.24033/asens.2121. http://www.numdam.org/articles/10.24033/asens.2121/.

\bibitem[AB2]{AB2} Andreatta, Fabrizio, and Olivier Brinon. 2013. "Acyclicit\'e G\'eom\'etrique de Bcris." Commentarii Mathematici Helvetici 88 (4): 965-1022. https://doi.org/10.4171/CMH/309.

\bibitem[AI2]{AI2} Andreatta, Fabrizio, and Adrian Iovita. 2012. "Semistable Sheaves and Comparison Isomorphisms in the Semistable Case." Rendiconti - Seminario Matematico Della Universit\`a Di Padova 128: 131-285. https://doi.org/10.4171/RSMUP/128-7.

\bibitem[Fon2]{Fon2} Fontaine, Jean Marc. "Cohomologie de De Rham, cohomologie cristalline et representations p-adiques." (1983). Algebraic Geometry, Proceedings of the Japan-France Conference, Tokyo/Kyoto, 1982.

\bibitem[Fon3]{Fon3} Fontaine, Jean-Marc. 1982. "Sur Certains Types de Representations p-Adiques Du Groupe de Galois D'un Corps Local; Construction D'un Anneau de Barsotti-Tate." Annals of Mathematics 115 (3): 529-577. https://doi.org/10.2307/2007012.

\bibitem[Fa1]{Fa1} Faltings, Gerd. 1988. "p-Adic Hodge Theory." Journal of the American Mathematical Society 1 (1): 255-299. https://doi.org/10.1090/S0894-0347-1988-0924705-1.

\bibitem[KL]{KL} Kedlaya, Kiran S, and Ruochuan Liu. 2016. "Finiteness of Cohomology of Local Systems on Rigid Analytic Spaces." https://doi.org/10.48550/arxiv.1611.06930.

\bibitem[Fa2]{Fa2} Faltings, Gerd. "Crystalline cohomology and p-adic Galois-representations." (1988). Algebraic Analysis, Geometry and Number Theory, Baltimore MD 1988, page 25-80 JHU press, Baltimore MD 1989.

\bibitem[Fa3]{Fa3} FALTINGS, Gerd. "Almost \'Etale Extensions." Ast\'erisque, Soci\'et\'e math\'ematique de France, 2002, pp. 185-270.

\bibitem[Fon4]{Fon4} Fontaine, Jean-Marc. Expos\'e II : Le corps des p\'eriodes $p$-adiques, in P\'eriodes $p$-adiques - S\'eminaire de Bures, 1988, Ast\'erisque, no. 223 (1994), Talk no. 2, 43 p. http://www.numdam.org/item/AST$\_$1994$\_$223$\_$59$\_$0/.

\bibitem[Fon5]{Fon5} Fontaine, Jean-Marc. "Repr\'esentations p-adiques semi-stables." Ast\'erisque, 223, 1994, page 113-184.

\bibitem[Fon6]{Fon6} Fontaine, Jean-Marc. "Repr\'esentations $\ell$-adiques potentiellement semi-stables." Ast\'erisque, 223, 1994, page 321-347.

\bibitem[SW]{SW} Scholze, Peter, and Jared Weinstein. Berkeley Lectures on p-Adic Geometry. Princeton University Press, 2020.

\bibitem[AI3]{AI3} Andreatta, Fabrizio, and Adrian Iovita. 2013. "Comparison Isomorphisms for Smooth Formal Schemes." Journal of the Institute of Mathematics of Jussieu 12 (1): 77-151. https://doi.org/10.1017/S1474748012000643.

\bibitem[Lan1]{Lan1} Robert Langlands. 1967. "Letter to Weil."

\bibitem[Drin1]{Drin1} Drinfel'd, Vladimir G. "Elliptic modules." Mathematics of the USSR-Sbornik 23, no. 4 (1974): 561.

\bibitem[Drin2]{Drin2} Drinfeld, Vladimir Gershonovich. "Langlands' conjecture for $\mathrm{GL}(2)$ over functional fields." In Proceedings of the International Congress of Mathematicians (Helsinki, 1978), vol. 2, pp. 565-574. 1980.

\bibitem[Zhu]{Zhu} Zhu, Xinwen. "Coherent sheaves on the stack of Langlands parameters." arXiv preprint arXiv:2008.02998 (2020).

\bibitem[DHKM]{DHKM} Dat, Jean-Fran\c{c}ois, David Helm, Robert Kurinczuk, and Gilbert Moss. "Moduli of Langlands parameters." arXiv preprint arXiv:2009.06708 (2020).


\bibitem[pHodgeT]{pHodgeT} Tate, John T. "p-Divisible groups." In Proceedings of a Conference on Local Fields: NUFFIC Summer School held at Driebergen (The Netherlands) in 1966, pp. 158-183. Berlin, Heidelberg: Springer Berlin Heidelberg, 1967.
\bibitem[pHodgeF]{pHodgeF} Fontaine, Jean-Marc. "Sur certains types de repr\'esentations p-adiques du groupe de Galois d'un corps local; construction d'un anneau de Barsotti-Tate." Annals of Mathematics 115, no. 3 (1982): 529-577.
\bibitem[pHodgeS1]{pHodgeS1} Scholze, Peter. "p-adic Hodge theory for rigid-analytic varieties." In Forum of Mathematics, Pi, vol. 1, p. e1. Cambridge University Press, 2013.
\bibitem[pHodgeS2]{pHodgeS2} Scholze, Peter. "\'Etale cohomology of diamonds." arXiv preprint arXiv:1709.07343 (2017).
\bibitem[pHodgeKL1]{pHodgeKL1} Kedlaya, K. S., and R. Liu. "Relative p-adic Hodge theory: foundations, Ast\'erisque 371." Soc. Math. France, Paris (2015).
\bibitem[pHodgeKL2]{pHodgeKL2} Kedlaya, Kiran S., and Ruochuan Liu. "Relative p-adic Hodge theory, II: Imperfect period rings." arXiv preprint arXiv:1602.06899 (2016).
\bibitem[pHodgeBS]{pHodgeBS} Breuil, Christophe and Schneider, Peter. "First steps towards p-adic Langlands functoriality" Journal f\"ur die reine und angewandte Mathematik 2007, no. 610 (2007): 149-180. https://doi.org/10.1515/CRELLE.2007.070.
\bibitem[pHodgeKPX]{pHodgeKPX} Kedlaya, Kiran, Jonathan Pottharst, and Liang Xiao. "Cohomology of arithmetic families of $(\varphi, \Gamma)$-modules." Journal of the American Mathematical Society 27, no. 4 (2014): 1043-1115.



\bibitem[pToAnCS1]{pToAnCS1} Dustin Clausen and Peter Scholze. "Lectures on Condensed Mathematics." https://www.math.uni-bonn.de/people/scholze/Condensed.pdf.
\bibitem[pToAnCS2]{pToAnCS2} Dustin Clausen and Peter Scholze. "Lectures on Analytic Geometry." https://www.math.uni-bonn.de/people/scholze/Analytic.pdf.
\bibitem[pToAnCS3]{pToAnCS3} Dustin Clausen and Peter Scholze. "Condensed Mathematics and Complex Geometry." https://people.mpim-bonn.mpg.de/scholze/Complex.pdf.
\bibitem[pToAnCS4]{pToAnCS4} Dustin Clausen and Peter Scholze. "Analytic Stacks." https://people.mpim-bonn.mpg.de/scholze/AnalyticStacks.html.
\bibitem[pToAnBBBK]{pToAnBBBK} Bambozzi, Federico, Oren Ben-Bassat, and Kobi Kremnizer. "Analytic geometry over F1 and the Fargues-Fontaine curve." Advances in Mathematics 356 (2019): 106815.



\bibitem[LPL]{LPL} Robert Langlands. "Letter to Weil." 1967. 
\bibitem[LPD1]{LPD1} Drinfeld, Vladimir Gershonovich. "Langlands' conjecture for $GL(2)$ over functional fields." In Proceedings of the International Congress of Mathematicians (Helsinki, 1978), vol. 2, pp. 565-574. 1980.
\bibitem[LPLL]{LPLL} Lafforgue, Laurent. "Chtoucas de Drinfeld et correspondance de Langlands." Inventiones mathematicae 147 (2002): 1-241.
\bibitem[LPVL]{LPVL} Lafforgue, Vincent. "Chtoucas pour les groupes r\'eductifs et param\'etrisation de Langlands globale." Journal of the American Mathematical Society 31, no. 3 (2018): 719-891. 
\bibitem[LPC]{LPC} Colmez, Pierre. "Repr\'esentations de $GL2(\mathbb{Q}_p)$ et $(\varphi, \Gamma)$-modules." Ast\'erisque 330, no. 281 (2010): 281-509.
\bibitem[LPFS]{LPFS} Fargues, Laurent, and Peter Scholze. "Geometrization of the local Langlands correspondence." arXiv preprint arXiv:2102.13459 (2021).
\bibitem[LPGL]{LPGL} Genestier, Alain, and Vincent Lafforgue. "Chtoucas restreints pour les groupes r\'eductifs et param\'etrisation de Langlands locale." arXiv preprint arXiv:1709.00978 (2017).
\bibitem[LPEGH]{LPEGH} Emerton, Matthew, Toby Gee, and Eugen Hellmann. "An introduction to the categorical p-adic Langlands program." arXiv preprint arXiv:2210.01404 (2022).
\bibitem[LPEG]{LPEG} Emerton, Matthew, and Toby Gee. Moduli Stacks of \'Etale $(\varphi, \Gamma)$-Modules and the Existence of Crystalline Lifts:(AMS-215). No. 215. Princeton University Press, 2022.
\bibitem[LPZ]{LPZ} Zhu, Xinwen. "Coherent sheaves on the stack of Langlands parameters." arXiv preprint arXiv:2008.02998 (2020).
\bibitem[LPDHKM]{LPDHKM} Dat, Jean-Fran\c{c}ois, David Helm, Robert Kurinczuk, and Gilbert Moss. "Moduli of Langlands parameters." arXiv preprint arXiv:2009.06708 (2020).

\bibitem[To2]{To2} Tong, Xin. "Arithmetic Hodge-Iwasawa Moduli Stacks." arXiv preprint arXiv:2401.09403 (2024). 
\bibitem[LPD2]{LPD2} Drinfel'd, Vladimir G. "Elliptic modules." Mathematics of the USSR-Sbornik 23, no. 4 (1974): 561.



\end{thebibliography}

\end{document}